\documentclass[10pt]{article}
\usepackage[dvipsnames]{xcolor}
\usepackage{float}
\usepackage{cancel}
\usepackage{leftidx}
\usepackage{amsmath}
\usepackage{amscd}
\usepackage{amssymb}
\usepackage{rotating}
\usepackage{amsfonts}
\usepackage{amsthm}
\usepackage{verbatim}
\usepackage{simpler-wick}
\usepackage{diagbox}
\usepackage{mathtools}
\usepackage[overload]{empheq}
\usepackage{cases}
\usepackage[top=1in, bottom=1in, left=1.25in, right=1.25in]{geometry}
\usepackage[colorlinks,linkcolor=blue, anchorcolor=blue, urlcolor=blue, citecolor=blue]{hyperref}
\usepackage{cleveref} 
 \usepackage[hyperpageref]{backref}

\input xy
\xyoption{all}

\allowdisplaybreaks[3]

\newcommand{\pqed}{\hfill\qedsymbol\\}
\newcommand{\nn}{\nonumber}
\newcommand{\sfh}{\mathsf{h}}
\newcommand{\tsfh}{\tilde{\mathsf{h}}}
\newcommand{\sfF}{\mathsf{F}}
\newcommand{\Mbar}{\overline{\mathcal{M}}}
\newcommand{\sqp}{\diamond}
\newcommand{\bqp}{\star}
\newcommand{\Spec}{\mathrm{Spec}}
\newcommand{\sfbd}{\mathbf{d}^{\mathbf{d}}}  
\newcommand{\elld}{\mathbf{d}!} 

\newtheorem{theorem}{Theorem}[section]
\newtheorem{proposition}[theorem]{Proposition}
\newtheorem{lemma}[theorem]{Lemma}
\newtheorem{corollary}[theorem]{Corollary}
\newtheorem{conjecture}[theorem]{Conjecture}
\newtheorem{theorem/definition}[theorem]{Theorem/Definition}
\newtheorem{question}[theorem]{Question}
\newtheorem{problem}[theorem]{Problem}

\theoremstyle{remark}
\newtheorem{remark}[theorem]{Remark}

\theoremstyle{definition}
 \newtheorem{example}[theorem]{Example}
\newtheorem{definition}[theorem]{Definition}
\newtheorem{construction}[theorem]{Construction}

\newtheorem{notation}[theorem]{Notation}

\title{Big quantum cohomology of Fano complete intersections}
\author{Xiaowen Hu}
\date{}

\begin{document}

\maketitle

\begin{abstract}
For smooth complete intersections in the projective spaces, we use the deformation invariance of Gromov-Witten invariants and results in classical invariant theory to study the symmetric reduction of the WDVV equation by the monodromy groups.

We propose a double root recursion conjecture for genus 0 invariants of non-exceptional Fano complete intersections other than the cubic hypersurfaces and the odd dimensional intersection of two quadrics. Based on it we develop an algorithm to compute the genus 0 invariants of any length and with any number of primitive insertions. The effectiveness of our algorithm is the main conjecture of this paper. We prove the conjecture at the first non-trivial order, which amounts to the computation of genus zero Gromov-Witten invariants with exactly 4 primitive insertions.

We show a reconstruction theorem for the odd dimensional intersection of two quadrics. There are some 4-point and 8-point genus 0 Gromov-Witten invariants for cubic hypersurfaces that cannot be computed by the monodromy group method. We compute them by studying the cohomology ring structure of the Fano variety of lines and by Zinger’s reduced genus 1 Gromov-Witten invariants. Then we get a reconstruction theorem for cubic hypersurfaces.

At the end of the paper, we give numerical examples and some conjectural closed formulae.
\end{abstract}



\tableofcontents

\setcounter{section}{0}
\section{Introduction}\label{sec:Introduction}
There have been intensive studies on the quantum cohomology of smooth complete intersections in the projective spaces, involving only ambient cohomology classes, especially for the Fano and Calabi-Yau types (see e.g., \cite{Giv96}, \cite{LLY97}).
On the other hand, we knew very little about the full quantum cohomology of $X$ with primitive cohomology classes as insertions, except for the $3$-point invariants \cite{Bea95}, \cite{CJ99}, with some restrictions on the multidegree of $X$. The  primitive cohomology classes do not come from restrictions from classes of the ambient projective spaces, and disappear on the singular fibers of the family of all complete intersections, so the Gromov-Witten invariants involving such classes cannot be computed  directly by the quantum hyperplane property, virtual torus localization, or the degeneration formula\footnote{After the 4th version of this paper, there appears a new approach \cite{ABPZ21}, which shows a new degeneration formula, and gives  an effective algorithm in all genera. However, at present, this algorithm seems not easy to give explicit computations.}.

 The  knowledge of the  big quantum cohomology involving primitive cohomology is necessary to understand the relation between quantum cohomology and derived categories such as Dubrovin's conjecture, and the numerical mirror symmetry.

The purpose of this paper is to study the big quantum cohomology of  the smooth Fano complete intersections  of dimension $\geq 3$ in the projective spaces. We use the deformation invariance of the Gromov-Witten invariants to equip the WDVV equations  with the monodromy group action, and then do  the symmetric reduction, and use some other ad hoc geometric method if necessary. The deformation invariance of Gromov-Witten invariants has been an important ingredient in e.g. \cite{CJ99} and \cite{OP06}, but its power has not yet  been fully developed. The use of big monodromy groups and their invariants was even earlier, in the study of Donaldson polynomials of algebraic surfaces \cite{FMM87}. For cubic hypersurfaces some additional geometric arguments are needed (\S \ref{sec:FanoVarietyOfLines} and \S \ref{sec:vanishing-reducedGenus1}).

 Before going to a detailed account of our main results, we mention  that in the following three cases of non-exceptional complete intersections (see Definition \ref{def-exceptional}) we have a complete reconstruction theorem for the genus 0 Gromov-Witten invariants:
\begin{itemize}
\item[(1)] Cubic hypersurfaces of dimension at least 3;
\item[(2)] Odd dimensional complete intersections of two quadrics of dimension at least 3;
\item[(3)] The cases where $\gcd (n-2, \mathsf{a}(n,\mathbf{d}))>1$, where $n$ is the dimension of the complete intersection, and $\mathsf{a}(n,\mathbf{d})$ is the Fano index.
\end{itemize}

For genus 0 Gromov-Witten invariants of non-exceptional Fano projective complete intersections other than the above cases, we find a double root recursion phenomenon. 
 \begin{conjecture}[Double root recursion]\label{conj-sqrtRecursion-unified-intro-sketchForm}
 The genus 0 Gromov-Witten invariants of a non-exceptional smooth projective complete intersection, which is not a cubic hypersurface or a complete intersection of two quadric hypersurfaces, can be reconstructed by recursively solving quadratic equations, of one variable, each of which has a unique double root.
 \end{conjecture}
The precise version of this conjecture is Conjecture \ref{conj-sqrtRecursion-unified-intro}. We prove the first nontrivial stage of this conjectural recursion. Then we give numerical verifications by examples in higher stages. 
 Based on this conjecture, we develop an algorithm to compute the genus 0 invariants of any length with any number of primitive insertions. It has the following features:
\begin{enumerate}
       \item For an even dimensional non-exceptional Fano complete intersection $X$ other than the cubic hypersurfaces, we check the double root recursion in each step, the so called $s$-order $l$, which encodes the genus 0 Gromov-Witten invariants of $X$ with exact $2l$ primitive insertions. If the double root recursion holds at this step, we can proceed on;
       \item For an odd dimensional non-exceptional Fano complete intersection $X$ other than the cubic hypersurfaces and intersections of two quadrics, we do the same as the even-dimensional ones, but for $l>1+ \frac{\mathrm{rank}\ H^n(X)}{4}$ the computation is conjectural. The reason is eventually the anti-commutativity of the Gromov-Witten with odd degree insertions.
       \item The trivial equations in the double root recursion lead to closed formulae which eventually express the generating function $F$ of genus 0 Gromov-Witten invariants in terms of the generating function $F^{(0)}$ of \emph{ambient} genus 0 Gromov-Witten invariants.
 \end{enumerate} 
 The effectiveness of our algorithm in even dimension, and its validity  in odd dimensions for invariants with more than  $\frac{\mathrm{rank}\ H^n(X)}{2}+2$ primitive insertions, is the essence of Conjecture \ref{conj-sqrtRecursion-unified-intro-sketchForm}. 

In  Appendix \ref{sec:algorithm}, we present our algorithm. The algorithm is implemented in a Macaulay2 package 
\texttt{QuantumCohomologyFanoCompleteIntersection}.
The reader can find the package in

 \url{https://github.com/huxw06/Quantum-cohomology-of-Fano-complete-intersections}

\subsection{Main results}
\label{sec:intro-1}
Let $X$ be a smooth complex projective variety. Denote the cohomology group $\oplus_{i=0}^{2\dim X}H^{i}(X,\mathbb{C})$ by $H^{*}(X)$. Let $\gamma_1,\dots,\gamma_m$ be a basis of $H^{*}(X)$. The generating function of genus zero Gromov-Witten invariants, or the quantum cohomology, of $X$ is defined to be
\begin{eqnarray}\label{eq-generating1}
\sfF(t^1,\dots,t^m,\mathsf{q})=\sum_{k=0}^{\infty}\sum_{\beta\in H_{2}(X)}\frac{1}{k!}
\langle \sum_{i=1}^{m}t^i \gamma_i, \dots,\sum_{i=1}^{m}t^i \gamma_i\rangle_{0,k,\beta}\mathsf{q}^{\beta}.
\end{eqnarray}
The multiplication of $t^1,\dots,t^m$ is graded commutative. More precisely, if $\gamma_1,\dots,\gamma_m$ are  chosen to have pure degrees 
$|\gamma_1|,\dots, |\gamma_m|$, then
\begin{eqnarray}
t^{i}t^j=(-1)^{|\gamma_i| |\gamma_j|}t^{j}t^i.
\end{eqnarray}
The Poincar\'e pairing on $H^{*}(X)$ is denoted by $(\cdot,\cdot)$. For a pure degree basis $\gamma_1,\dots,\gamma_m$, let 
$g_{ij}=(\gamma_i,\gamma_j)$, and the inverse matrix be $g^{ij}$. Then $\sfF(t^1,\dots,t^m,\mathsf{q})$ satisfies the WDVV equation
\begin{eqnarray}\label{eq-wdvv_1}
\sum_{e=1}^{m}\sum_{f=1}^{m}\frac{\partial^3 \sfF}{\partial t_a\partial t_b\partial t_e}g^{ef}\frac{\partial^3 \sfF}{\partial t_f\partial t_c\partial t_d}
=(-1)^{|\gamma_b| |\gamma_c|}\sum_{e=1}^{m}\sum_{f=1}^{m}\frac{\partial^3 \sfF}{\partial t_a\partial t_c\partial t_e}g^{ef}\frac{\partial^3 \sfF}{\partial t_f\partial t_b\partial t_d}.
\end{eqnarray} 
Thus to $\sfF(t^1,\dots,t^m,\mathsf{q})$ is associated a  \emph{formal Frobenius (super-)manifold} \cite[\S III.1.1]{Man99}, and we denote it by $\mathcal{M}_{X}$. There are also various choices of subspaces of 
$H^{*}(X)$ to obtain Frobenius manifolds. The following lemma is easily deduced from 
(\ref{eq-wdvv_1}).
\begin{lemma}\label{lem-subspace}
Let $V$ be a subspace of $H^{*}(X)$. Suppose there is another subspace $V^{\bot}$, such that
\begin{itemize}
\item[(i)] $H^{*}(X)=V\oplus V^{\bot}$, and $(v_1,v_2)=0$ for $v_1\in V, v_2\in V^{\bot}$;
\item[(ii)] $\langle \gamma_1,\dots,\gamma_k, \gamma\rangle_{0,k+1,\beta}=0$ for $\gamma_1,\dots,\gamma_k\in V$, $\gamma\in V^{\bot}$ and  $k\geq 1$, $\beta \in H_2(X)$.
\end{itemize}
Then choosing a pure degree  basis of $V$, says $ \gamma_1,\dots,\gamma_l$, the generating function
\begin{eqnarray}\label{eq-generating2}
\sfF(t^1,\dots,t^l,\mathsf{q})=\sum_{k=0}^{\infty}\sum_{\beta\in H_{2}(X)}\frac{1}{k!}
\langle \sum_{i=1}^{l}t^i \gamma_i, \dots,\sum_{i=1}^{l}t^i \gamma_i\rangle_{0,k,\beta}\mathsf{q}^{\beta}
\end{eqnarray}
satisfies the WDVV equation
\begin{eqnarray}\label{eq-wdvv_2}
\sum_{e=1}^{l}\sum_{f=1}^{l}\frac{\partial^3 \sfF}{\partial t_a\partial t_b\partial t_e}g^{ef}\frac{\partial^3 \sfF}{\partial t_f\partial t_c\partial t_d}
=(-1)^{|\gamma_b| |\gamma_c|}\sum_{e=1}^{l}\sum_{f=1}^{l}\frac{\partial^3 \sfF}{\partial t_a\partial t_c\partial t_e}g^{ef}\frac{\partial^3 \sfF}{\partial t_f\partial t_b\partial t_d}.
\end{eqnarray} 
\end{lemma}

Now let $X$ be a complete intersection of dimension $n$ in $\mathbb{P}^{n+r}$. By the Lefschetz hyperplane theorem and Poincar\'e duality, we have an injective homomorphism $j^*: H^{*}(\mathbb{P}^{n+r})\rightarrow H^{*}(X)$. We call the image of $j^*$ the \emph{ambient cohomology} of $X$ according to the embedding $j: X\rightarrow \mathbb{P}^{n+r}$, and denote it by $H^*_{\mathrm{amb}}(X)$. We denote the primitive cohomology group of $X$ by 
$ H^*_{\mathrm{prim}}(X)$. We have the decomposition
\begin{eqnarray*}
H^*(X)=H^*_{\mathrm{amb}}(X)\oplus H^n_{\mathrm{prim}}(X).
\end{eqnarray*}

The subspace $H^*_{\mathrm{amb}}(X)$ is equal to the subspace of $H^*(X)$ fixed by the monodromy group of the total family of the smooth complete intersections in $\mathbb{P}^{n+r}$ which contains $X$ as a fibre. From this and the \emph{deformation axiom} of Gromov-Witten invariants, we can deduce (see Corollary \ref{cor-mono3}, or \cite{BK05}, \cite{LP04}) that $H^*_{\mathrm{amb}}(X)$ satisfies the assumptions of Lemma \ref{lem-subspace}, thus we obtain a Frobenius manifold $\mathcal{M}_{\mathrm{amb}}$.  Note that Zinger showed in \cite{Zin14} that $\mathcal{M}_{\mathrm{amb}}$ is an analytic Frobenius manifold, not only a formal one. See also Remark \ref{rem:convergence-F(0)}  for a simple proof of this fact, as a consequence of our algorithm.  

By the \emph{quantum hyperplane property}, the quantum cohomology for ambient classes can be computed as the \emph{twisted} quantum cohomology of the ambient space $\mathbb{P}^{n+r}$. This enables us to compute  $\mathcal{M}_{\mathrm{amb}}$. For Fano and Calabi-Yau complete intersections we have the \emph{mirror theorems} \cite{Giv96}, \cite{LLY97} for the \emph{small $J$-functions}, which encode the
information of all the genus 0 invariants involving the ambient classes only.

In this paper we will see that  the WDVV equation together with the standard properties of Gromov-Witten invariants, especially the deformation invariance,  may enable us to compute the genus 0 Gromov-Witten invariants involving the primitive classes. In general, if $\mathcal{X}\rightarrow S$ is a smooth family of projective varieties, and $\mathcal{X}_0=X$ is a special fibre, 
then $G=\pi_{1}(S,0)$ acts on $H^{*}(X)$ as the monodromy action. Denote the \emph{Novikov ring} of $X$ by
$\Lambda$, i.e., $\Lambda$ is spanned over $\mathbb{C}$ by $\{\mathsf{q}^{\alpha}: \alpha\in H_{2}(X)\}$ with the relation $\mathsf{q}^{\alpha+\beta}=\mathsf{q}^{\alpha}\cdot \mathsf{q}^{\beta} $. 
The following lemma is a consequence of the deformation axiom and the $S_n$ invariance axiom (see Corollary \ref{cor-mono2}).
\begin{lemma}\label{lem-inv1}
For any genus $g$, the generating function 
$\mathcal{F}_{g}(t^1,\dots,t^m,\mathsf{q})$ of genus $g$ primary Gromov-Witten invariants of $X$
lies in the $G$-invariant part $$\big(\mathrm{Sym}(H^*(X)^{\vee})[[\Lambda]]\big)^{G}, $$
where  $\mathrm{Sym}(H^*(X)^{\vee})$ is the $\mathbb{Z}/2 \mathbb{Z}$-\emph{graded symmetric product} of $H^*(X)^{\vee}$. 
\end{lemma}

 When $X$ is a smooth complete intersection in $\mathbb{P}^{n}$, we have good knowledge about the Zariski closure of $G$, and can find its invariants.
 The crucial observation is that  Lemma \ref{lem-inv1} implies some symmetric reduction of the tautological relations, and thus  gives us new equations
  for $\mathcal{F}_g$. In this paper we focus on $\sfF=\mathcal{F}_0$ and study the symmetric reduction of  
  (\ref{eq-wdvv_1}). \\

In the following of this paper, by \emph{complete intersections} we always mean smooth complete intersections in the projective spaces, and without loss of generality we assume that every component $d_i$ of the multi-degree $\mathbf{d}=(d_1,\dots, d_r)$ is at least 2. For complete intersections of dimension 2, the primary Gromov-Witten invariants involving primitive classes can be reduced to those without primitive classes, by the divisor equation. For complete intersections of Calabi-Yau type or of general type, a simple degree counting shows that the Gromov-Witten invariants involving primitive classes are trivial. So we consider only the \emph{Fano complete intersections of dimension at least 3}, and thus the Picard number is 1. Except for a small class of $\mathbf{d}$ 
(for such $\mathbf{d}$, $X$ is called \emph{exceptional}, see Section \ref{sec:monodromyGroup} for the definition), the Zariski closure of $G$ is the orthogonal group or the symplectic group, whose ring of invariants is very simple. Let $\gamma_{0},\dots, \gamma_{n}$ be a basis of $H_{\mathrm{amb}}^{*}(X)$, 
 $\gamma_{n+1},\dots, \gamma_{n+m}$ a basis of  $H_{\mathrm{prim}}^{n}(X)$, and $t^0,\dots,t^{n+m}$  the dual basis. Let $g_{ab}=g(\gamma_a,\gamma_b)$ be
 the Poincar\'e pairing. We introduce a new variable
 \begin{eqnarray}
s=\frac{(-1)^n}{2}\sum_{i=n+1}^{n+m}\sum_{j=n+1}^{n+m}g_{ij}t^{i}t^{j}.
 \end{eqnarray}
 Note that when the dimension of $X$ is odd, $m$ is even and $s^{\frac{m}{2}+1}=0$ because $t^i$ and $t^j$ are skew-commutative.
\begin{theorem}\label{thm-mainthm1}(= Theorem \ref{thm-wdvveventhm} + Theorem \ref{thm-wdvvoddthm})
Suppose $X$ is a non-exceptional complete intersection of dimension $n \geq 3$. Let $m=\mathrm{rank}\ H^{n}_{\mathrm{prim}}(X)$.
\begin{itemize}
\item[(i)] $\sfF$ is a series in $t^0, \dots, t^n$ and $s$.
\item[(ii)] If $\dim X$ is even, the  WDVV equation for $\sfF$ is equivalent to the WDVV equation for $\sfF^{(0)}$ together with
\begin{eqnarray}\label{eq-wdvv230}
&&\sum_{e=0}^{n}\sum_{f=0}^{n}\frac{\partial^3 \sfF}{\partial t^a \partial t^b\partial t^e}g^{ef}\frac{\partial^2 \sfF}{\partial s\partial t^f}+2s\frac{\partial^3 \sfF}{\partial s\partial t^a\partial t^b}\frac{\partial^2 \sfF}{\partial s^2}\nn\\
&=& \frac{\partial \sfF}{\partial s\partial t^a}\frac{\partial \sfF}{\partial s\partial t^b}, \quad 0\leq a,b\leq n,
\end{eqnarray}
and
\begin{eqnarray}\label{eq-wdvv240}
\sum_{e=0}^{n}\sum_{f=0}^{n}\frac{\partial \sfF}{\partial s\partial t^e}g^{ef}\frac{\partial \sfF}{\partial s\partial t^f}+2s\Big(\frac{\partial^2 \sfF}{\partial s^2}\Big)^2= 0 .
\end{eqnarray}
\item[(iii)] If $\dim X$ is odd, 
the  WDVV equation for $\sfF$ is equivalent to the  WDVV equation for $\sfF^{(0)}$ together with
\begin{eqnarray}\label{eq-wdvv230odd}
&&\sum_{e=0}^{n}\sum_{f=0}^{n}\frac{\partial^3 \sfF}{\partial t^a \partial t^b\partial t^e}g^{ef}\frac{\partial^2 \sfF}{\partial s\partial t^f}+2s\frac{\partial^3 \sfF}{\partial s\partial t^a\partial t^b}\frac{\partial^2 \sfF}{\partial s^2}\nn\\
&=& \frac{\partial \sfF}{\partial s\partial t^a}\frac{\partial \sfF}{\partial s\partial t^b}\mod s^{\frac{m}{2}}, \quad 0\leq a,b\leq n,
\end{eqnarray}
and
\begin{eqnarray}\label{eq-wdvv240odd}
\sum_{e=0}^{n}\sum_{f=0}^{n}\frac{\partial \sfF}{\partial s\partial t^e}g^{ef}\frac{\partial \sfF}{\partial s\partial t^a}+2s\Big(\frac{\partial^2 \sfF}{\partial s^2}\Big)^2= 0 \mod s^{\frac{m}{2}}.
\end{eqnarray}
\end{itemize}
\end{theorem}

We also study the consequence of the deformation axiom and the $S_n$-invariance for the descendant invariants. For the result see Appendix \ref{sec:symmetricReduction-JFunction}. In particular, it turns out that the generating function of one-point descendants of primitive classes is related to $\sfF$ in a simple way (Theorem \ref{thm-Jreconstruction1} (ii)).\\

Let
\begin{eqnarray*}
\sfF^{(k)}(t^0, t^1,\dots,t^{n})=\Big(\frac{\partial ^{k}}{\partial s^{k}}\sfF\Big) \Big |_{s=0},
\end{eqnarray*}
so we have the expansion
\begin{eqnarray*}
\sfF=\sfF^{(0)}+s\sfF^{(1)}+\frac{s^{2}}{2}\sfF^{(2)}+\dots
\end{eqnarray*}
where $\sfF^{(0)}$ is the generating function of the genus 0 primary Gromov-Witten invariants of the ambient classes.  

We need also the Euler vector field for $\sfF$. 
Suppose $\gamma_{i}=\sfh_i$ is the $i$-th power (in the ordinary cohomology ring) of the hyperplane class, 
$0\leq i\leq n$. Let $\mathsf{a}(n,\mathbf{d})=n+r+1-\sum_{i=1}^{r}d_i$ be the Fano index of $X$.
Recall that the Euler vector field is
\begin{gather*}
E:=\sum_{i=0}^{n}(1-i)t^{i}\frac{\partial}{\partial t^i}+\sum_{j=n+1}^{n+m}(1-\frac{n}{2})t^j\frac{\partial }{\partial t^j}+\mathsf{a}(n,\mathbf{d})\frac{\partial}{\partial t^1}.
\end{gather*}
We introduce another differential operator
\begin{gather}\label{eq-effective-EulerVectorField}
\sum_{i=0}^{n}(1-i)t^{i}\frac{\partial}{\partial t^i}+(2-n)s\frac{\partial }{\partial s}+\mathsf{a}(n,\mathbf{d})\frac{\partial}{\partial t^1}.
\end{gather}
This operator has the same effect as $E$ when they act on $\sfF$. In this paper, we abusively refer to (\ref{eq-effective-EulerVectorField}) as  the Euler vector field as well and denote it also by $E$.
Let $c(t_0, \dots, t^{n+m})$ be the classical triple intersection form, i.e,
\begin{eqnarray}\label{eq-tripleIntersectionForm}
c(t_0, \dots, t^{n+m})=\sum_a\sum_b\sum_c\frac{t^{a}t^{b}t^{c}}{6}\int_{X}\gamma_a \gamma_b \gamma_c.
\end{eqnarray}
Then $\sfF$ satisfies the following system
\begin{subequations}\label{eq-intro-system-even}
\begin{align}[left ={\empheqlbrace}]
&\sfF_{abe}g^{ef}\sfF_{sf}+2s\sfF_{sab}\sfF_{ss}=\sfF_{sa}\sfF_{sb},\ (\mbox{mod}\ s^{\frac{m}{2}}\ \mbox{when $n$ is odd}) \quad \mbox{for}\ 0\leq a,b\leq n,\label{eq-intro-system-even-1}\\
&\sfF_{se}g^{ef}\sfF_{sf}+2s\sfF_{ss}\sfF_{ss}=0,\ (\mbox{mod}\ s^{\frac{m}{2}}\ \mbox{when $n$ is odd}) \label{eq-intro-system-even-2}\\
&E\sfF=(3-n)\sfF+\mathsf{a}(n,\mathbf{d})\frac{\partial}{\partial t^1}c. \label{eq-intro-system-even-3}
\end{align}
\end{subequations}
Here we have written (\ref{eq-wdvv230}) and (\ref{eq-wdvv240}) in compact forms, where  the Einstein summation is taken over $e,f=0,1,\dots,n$. 
The goal of this paper can be presented as:\\

\textbf{Goal}: Reconstruct $\sfF$ from $\sfF^{(0)}$ using (\ref{eq-intro-system-even}).\\

Before proceeding to discuss the solution of (\ref{eq-intro-system-even}), we first make some comments 
on the non-semisimplicity of the quantum cohomology of $X$. A smooth K\"{a}hler manifold $X$ is called of Hodge-Tate type if it has the hodge numbers
$h^{p,q}=0$ for  $p\neq q$. Let $\mathcal{M}_{\mathrm{even}}$ be the Frobenius manifold produced by the quantum cohomology of $X$ with even degree insertions. By \cite{BaM04}, \cite{HMT09},  $\mathcal{M}_{\mathrm{even}}$ is generically semisimple implies that $X$ is of Hodge-Tate type. A complete intersection is of Hodge-Tate type if and only if it is exceptional. Thus for non-exceptional complete intersections,
 $\mathcal{M}_{\mathrm{even}}$ is not generically semisimple. Here we concern the ambient locus $\mathcal{M}_{\mathrm{amb}}$ in $\mathcal{M}_{\mathrm{even}}$. When we solve $\sfF^{(1)}$ we obtain the following  corollary of 
 Theorem \ref{thm-mainthm1} in passing. (This corollary is contrary to  the result of \cite{TX97}, where it was stated that $\mathcal{M}_{\mathrm{amb}}$ is generically semisimple if the degree of $X$ is relatively lower to $n$. 
There is an error in the proof of \cite{TX97}, page 487, line 15, where they incorrectly quoted from \cite{Bea95} that $c_1=\frac{1}{2}l_0^2$ and deduced that some determinant is nonzero modulo $(t^2)$. The correct computation is $c_1=l_0^2$ and thus the determinant is zero modulo $(t^2)$. I am indebted to Nicolas Perrin for pointing out this.) For a general account of the semisimplicity of Fano varieties I refer the reader to \cite{Per14}.

\begin{proposition}\label{cor-nonsemisimple}(= Corollary \ref{cor-ss4})
For every non-exceptional complete intersection $X$ with dimension $\geq 3$, $\mathcal{M}_{\mathrm{amb}}$ is nowhere semisimple.
\end{proposition}

In the proof of Proposition  \ref{cor-nonsemisimple}, we introduce a vector field
\begin{eqnarray}\label{eq-intro-ss2.5}
\Theta:=\sum_{e=0}^{n}\sum_{f=0}^{n}\sfF_{e}^{(1)}g^{ef}\gamma_{f}.
\end{eqnarray}
This vector field $\Theta$ has nice properties (see Proposition \ref{prop-nonsemisimple}); in particular, it is a common eigenvector of the quantum multiplication by all ambient cohomology classes. This enables us to determine $\sfF^{(1)}$ (see Proposition \ref{prop-Theta-inSmallQuantumCohomology} and Theorem \ref{thm-reconstruction-II-F(1)}).

Now we discuss the solution of (\ref{eq-intro-system-even}). Recall  the dimension constraint:  
$$\langle \gamma_{a},\gamma_{b},\dots, \gamma_{c}\rangle_{0,k,\beta}\neq 0$$  implies 
\begin{equation}\label{eq-intro-Dim}
	\deg \gamma_a+ \deg \gamma_b+\dots+ \deg \gamma_c=n-3+k+\beta\cdot \mathsf{a}(n,\mathbf{d}).
\end{equation}
Here $\deg$ means the complex degree, see Notations and Conventions \ref{sec:nota-conv}.
As we have observed, when $\mathsf{a}(n,\mathbf{d})\leq 0$ the reconstruction from $\sfF^{(0)}$ to $\sfF$ is trivial. So  we assume 
$\mathsf{a}(n,\mathbf{d})>0$, i.e. $X$ is a Fano complete intersection. Then by (\ref{eq-intro-Dim}), 
\begin{equation*}
	F:=\sfF|_{\mathsf{q}=1}
\end{equation*}
makes sense as a formal series and determines $\sfF$. The previous differential equations for $\sfF$ hold also for $F$, and it is convenient to state the following results in terms of $F$. 

Expanding both sides of the differential equations (\ref{eq-intro-system-even}) into series of $s$, we obtain many  equations of genus 0 Gromov-Witten invariants. Then by studying the structure of the Frobenius algebra of the small quantum cohomology of ambient classes of $X$ we obtain the following reconstruction theorem.

\begin{theorem}\label{thm-intro-reconstruction-II-F(1)}(= Theorem \ref{thm-reconstruction-II-F(1)})
Let $X$ be a non-exceptional Fano complete intersections in a projective space. Then
$F^{(1)}$ can be reconstructed by (\ref{eq-system1-even}) when $n$ is even (resp. (\ref{eq-system1-odd}) when $n$ is odd) and (\ref{eq-Dim}) from the generating function $F^{(0)}$ of the ambient quantum cohomology.
\end{theorem}
In particular,  Proposition \ref{prop-Theta-inSmallQuantumCohomology} gives all the 3-point invariants, which gives a new proof of, and  generalizes, the results in \cite{Bea95} and \cite{CJ99}. Then we show:
\begin{theorem}\label{thm-induction1}(= Theorem \ref{thm-reconstruction-I})
For non-exceptional Fano complete intersections with dimension $\geq 3$, $F^{(k)}$ can be reconstructed by the equations (\ref{eq-intro-system-even-1}), (\ref{eq-intro-system-even-3}), and (\ref{eq-Dim}) from the following data:
\begin{itemize}
\item[(i)] The generating function $F^{(0)}$ of ambient quantum cohomology;
\item[(ii)] The generating function $F^{(1)}$ of  quantum cohomology with exactly two primitive insertions;
\item[(iii)] The constant term $F^{(i)}(0)$ of $F^{(i)}$ for $2\leq i\leq k$.
\end{itemize}
\end{theorem}

By the dimension constraint, if $\gcd (n-2, \mathsf{a}(n,\mathbf{d}))>1$, then $F^{(k)}(0)=0$ for $k\geq 1$. So an immediate corollary is
\begin{corollary}\label{cor-reconstruction-intro}
For non-exceptional Fano complete intersections with  $\gcd (n-2, \mathsf{a}(n,\mathbf{d}))>1$, $F$ can be reconstructed by (\ref{eq-intro-system-even-1}), (\ref{eq-intro-system-even-3}) and (\ref{eq-Dim}) from $F^{(0)}$.
\end{corollary}
The classical  way to use WDVV (\ref{eq-wdvv_1}) to get recursions is to use the leading terms.
 Namely, selecting a monomial $t^I$, where $I$ is a multi-index, and extracting the coefficients of $t^I$ on both sides, we get an equation of the form 
 \begin{eqnarray}
      &&\mathrm{Coeff}_{t^I}(\partial_{t^a}\partial_{t^b}\partial_{t^e}F)g^{ef}(\partial_{t^f}\partial_{t^c}\partial_{t^d}F)(0)
      +     (\partial_{t^a}\partial_{t^b}\partial_{t^e}F)(0)g^{ef}\mathrm{Coeff}_{t^I}(\partial_{t^f}\partial_{t^c}\partial_{t^d}F)\nn\\
      &&-\mathrm{Coeff}_{t^I}(\partial_{t^a}\partial_{t^c}\partial_{t^e}F)g^{ef}(\partial_{t^f}\partial_{t^b}\partial_{t^d}F)(0)
      -(\partial_{t^a}\partial_{t^c}\partial_{t^e}F)(0)g^{ef}\mathrm{Coeff}_{t^I}(\partial_{t^f}\partial_{t^b}\partial_{t^d}F)\nn\\
      &=& \mbox{combinations of coefficients of lower order terms}.
 \end{eqnarray}
 Here we have omitted the signs in (\ref{eq-wdvv_1}), and  have adopted Einstein's summation convention, i.e. omitting the summation notations of the repeated indices $e$ and $f$. More generally, for a fixed length $l$, one can use the knowledge of the invariants of length $\leq l$ to get recursions. We call the resulted recursions \emph{essentially linear recursions}. The application of (\ref{eq-intro-system-even-1}) and (\ref{eq-intro-system-even-3}) in the proof of both Theorem \ref{thm-intro-reconstruction-II-F(1)} and Corollary \ref{cor-reconstruction-intro} is  similar to this and we call it  essentially linear recursions as well.

When $\gcd (n-2, \mathsf{a}(n,\mathbf{d}))=1$, there may be nonzero $F^{(k)}(0)$ for $k\geq 2$. In Theorem \ref{thm-intro-reconstruction-II-F(1)} and Theorem \ref{thm-induction1} the equation (\ref{eq-intro-system-even-2}) is only used to determine $F^{(1)}(0)$. We hope that (\ref{eq-intro-system-even-2}) will also help to determine  $F^{(k)}(0)$ for 
$k\geq 2$.  This seems to be the most difficult aspect of the system (\ref{eq-intro-system-even}). In fact I regard (\ref{eq-intro-system-even-2}) as the essential consequence of the symmetric reduction by the monodromy group; see Remark \ref{rmk-reconstruction-II}. 

So our central task is to compute $F^{(k)}(0)$ from the system (\ref{eq-intro-system-even}) or  some additional geometric tools if necessary. The first nontrivial one is $F^{(2)}(0)$. With  the above results obtained by (\ref{eq-intro-system-even-1}) and (\ref{eq-intro-system-even-3}), the equation (\ref{eq-intro-system-even-2}) will give a quadratic equation for $F^{(2)}(0)$.
It is quite a miracle that the involved quadratic equations for $F^{(2)}(0)$ have two equal roots, except for the case of  cubic hypersurfaces. 
This is the first occurrence of the double root recursion phenomenon.

To state the result,
we need to introduce some notations. For $\mathbf{d}=(d_1,\dots,d_r)$, let
\begin{eqnarray}\label{eq-alb-intro}
|\mathbf{d}|=\sum_{i=1}^{r}d_i,\
\elld:=\prod_{i=1}^{r}d_i !, & \sfbd:=d_1^{d_1}\cdots d_{r}^{d_r},
\end{eqnarray}
For a smooth Fano complete intersection $X$ of multidegree $\mathbf{d}=(d_1,\dots,d_r)\in \mathbb{P}^{n+r}$, let
$\sfh$ be the hyperplane class of $X$. We define
\begin{equation*}
	\sfh_i=\underbrace{\sfh\cup\cdots\cup\sfh}_{\mbox{$i$ factors}},
\end{equation*}
\begin{eqnarray}\label{eq-ss5-intro}
\tilde{\sfh}=\left\{
\begin{array}{cc}
\sfh, & \mathsf{a}(n,\mathbf{d})\geq 2,\\ 
\sfh+\elld, & \mathsf{a}(n,\mathbf{d})=1,
\end{array}\right.
\end{eqnarray}
and
\begin{equation*}
	\tsfh_i=\underbrace{\tsfh\sqp\cdots\sqp\tsfh}_{\mbox{$i$ factors}}
\end{equation*}
where $\sqp$ stands for the small quantum product. Let $M$ and $W$ be the transition matrices between $\sfh_i$ and $\tsfh_i$:
\begin{equation*}
	\sfh_i=\sum_{j=0}^{n}M_i^j \tsfh_j,\ \tsfh_i=\sum_{j=0}^{n}W_i^j \sfh_j.
\end{equation*}

\begin{theorem}\label{thm-induction2}(= Theorem \ref{thm-higher10.1})
Let $X_n(\mathbf{d})$ be a non-exceptional complete intersection.
Then
\begin{eqnarray}\label{eq-intro-higher22}
F^{(2)}(0)=
\begin{dcases}
1, & \mathrm{if}\  \mathbf{d}=(2,2);\\
1\ \mathrm{or}\ 4, & \mathrm{if}\  \mathbf{d}= (3);\\
\frac{-\sum_{j=0}^n j M_{j}^{1}W_{n}^{j}
+\sfbd\sum_{j=0}^n 
j M_{j}^{1}W_{n- \mathsf{a}(n,\mathbf{d})}^{j}}{\mathsf{a}(n,\mathbf{d})\prod_{i=1}^r d_i}
	\
, & \mathrm{if}\ l=\frac{n-1}{\mathsf{a}(n,\mathbf{d})}\in \mathbb{Z}_{\geq 2};\\
0, & \mathrm{otherwise}.
\end{dcases}
\end{eqnarray}
\end{theorem}

The formula (\ref{eq-intro-higher22}) for $F^{(2)}(0)$, when $l=\frac{n-1}{\mathsf{a}(n,\mathbf{d})}\in \mathbb{Z}_{\geq 2}$, may be not satisfactory. In the  case $l=2$, we have a closed formula.  
\begin{theorem}\label{thm-intro-F20-a(n,d)=(n-1)/2}(= Theorem \ref{thm-F20-a(n,d)=(n-1)/2})
Suppose $\frac{n-1}{\mathsf{a}(n,\mathbf{d})}=2$. Then
\begin{equation*}
      F^{(2)}(0)=\frac{\prod_{i=1}^r d_i!(d_i-1)!}{2}.
\end{equation*}
\end{theorem}

For cubic hypersurfaces, to determine $F^{(2)}(0)$ we need more geometric inputs. Since the involved invariants are of degree 1, by an easy vanishing result of the genus 1 reduced invariants defined by Zinger \cite{Zin09},  we can use the  genus 1 standard versus reduced formula  of \cite{Zin08} to compute $F^{(2)}(0)$. We can also study the structure of the cohomology ring of the Fano variety of lines $\Omega_X$ on cubic hypersurfaces $X$, with the help from the result of \cite{GS14} on the Betti numbers of $\Omega_X$, and by the way obtain 
$F^{(2)}(0)$.

\begin{theorem}\label{thm-intro6}(= Theorem \ref{thm-cubic5} (iv) or Theorem \ref{thm-4points-fanoIndex-(n-1)})
For the cubic hypersurfaces of dimension $\geq 3$, $F^{(2)}(0)=1$. 
\end{theorem}

As a byproduct, we obtain a complete description of the ring structure of $H^*(\Omega_X)$, see Theorem \ref{thm-cubic} for details.  

 For $\mathbf{d}= (3)$ or $(2,2)$, we have a complete reconstruction theorem, by an essentially linear recursion on $F^{(k)}(0)$ for $k\geq 3$.

\begin{theorem}\label{thm-intro5}(= Theorem \ref{thm-reconstructcubicandquadric})
\begin{itemize}
\item[(i)] For the cubic threefold $X$, $F$ can be reconstructed by (\ref{eq-intro-system-even}) and the dimension constraint (\ref{eq-Dim}) from $F^{(0)}$ and  $F^{(2)}(0)$, $F^{(4)}(0)$.
\item[(ii)] For cubic hypersurfaces $X$ with $\dim X\geq 4$, $F$ can be reconstructed by (\ref{eq-intro-system-even}) and the dimension constraint (\ref{eq-Dim}) from $F^{(0)}$ and  $F^{(2)}(0)$.
\item[(iii)] For odd-dimensional intersections of two quadrics with $\dim X\geq 2$, $F$ can be reconstructed by (\ref{eq-intro-system-even}) and  the dimension constraint (\ref{eq-Dim}) from  $F^{(0)}$ and  $F^{(2)}(0)$.
\end{itemize}
\end{theorem}

After the first version of this paper, this trick to compute degree 1 invariants was  applied by Hua-Zhong Ke   in his proof \cite{Ke18} of the conjecture $\mathcal{O}$ for Fano complete intersections.\\

As we see in Theorem \ref{thm-intro5} (ii), the cubic threefold is special. To compute $F^{(4)}(0)$, we show a vanishing theorem for certain degree 2 reduced genus 1 invariants on cubic hypersurfaces:
\begin{theorem}\label{thm-intro-vanishing-deg2-cubics}(= Theorem \ref{thm-vanishing-deg2-cubics})
Let $X$ be a cubic hypersurface in $\mathbb{P}^N$. Let $\alpha_1,\dots,\alpha_k\in H^*(X)$. Then 
\begin{equation}\label{eq-intro-vanishing-deg2-cubics}
\langle \alpha_1,\dots,\alpha_k\rangle_{1,2}^0=0=\langle \psi \alpha_1,\alpha_2,\dots,\alpha_k\rangle_{1,2}^0.
\end{equation}
\end{theorem}
The idea of the proof is to show first that the evaluation map restricted to $\Mbar_{1,k}^{0}(X,2)$ factors through $\Mbar_{0,k}(X,1)$, and then count the virtual dimensions.
Consequently  we can use Zinger's standard versus reduced formula to obtain:
\begin{theorem}\label{thm-8points-cubic3fold-intro}(= Theorem \ref{thm-8points-cubic3fold-final})
For cubic 3-folds $X$, $F^{(4)}(0)=0$. Equivalently, for any $\gamma_i\in H^3(X;\mathbb{Q})$, $1\leq i\leq 8$,
\begin{equation*}
 	\langle \gamma_1,\dots,\gamma_8\rangle_{0,8,2}^X=0.
 \end{equation*} 
\end{theorem}

Theorem  \ref{thm-intro6}, \ref{thm-intro5}, and \ref{thm-8points-cubic3fold-intro} together  give a complete reconstruction procedure  for all smooth cubic hypersurfaces of dimension $\geq 3$, and odd dimensional complete intersections of two quadrics. \\

Now we begin to discuss the computation of $F^{(k)}(0)$ for $k\geq 3$, for Fano complete intersections other than the cases in Theorem \ref{thm-intro5}. 
We give an algorithm in Appendix \ref{sec:algorithm} to explicitly compute $F^{(0)}$, $F^{(1)}$ and also $F^{(k)}$ where we regard $F^{(k)}(0)$ as unknowns. We implement this algorithm as a Macaulay2 package 
\begin{equation*}
      \mbox{\texttt{QuantumCohomologyFanoCompleteIntersection}. }
\end{equation*}
In this algorithm we use the coordinates $\tau^0,\dots,\tau^n$ dual to the small quantum cohomology basis $1,\tsfh,\dots,\tsfh_n$. It has the advantage that the linear recursion of the highest order terms in the WDVV equations is quite simple. The cost is that the expression of the Euler field becomes complicated.

The use of $\tau$-coordinates is also essential for the proof of Theorem \ref{thm-induction2}.  In fact the expression of $F^{(2)}(0)$ in (\ref{eq-intro-higher22}) (see also (\ref{eq-functionOfWM})) essentially comes from the Euler field in terms of $\tau^i$. But note that  $F^{(k)}(0)$ is independent of the choice of coordinates.  So the occurrence of such an expression  justifies the naturality of our use of the coordinates $\tau^i$.\\

With this package in hand, we are able to investigate $F^{(k)}(0)$ for $k>2$.
From  (\ref{eq-intro-system-even-2}) written in $\tau$-coordinates, we get  for $l\geq 2$, ($2\leq l\leq \frac{m}{2}$ when $n$ is odd)
\begin{eqnarray}\label{eq-wdvv24expand-sqrt-intro}
\sum_{k=1}^{l}\sum_{a=0}^n\sum_{b=0}^n
\frac{\partial_{\tau^a}F^{(k)}\eta^{ab}\partial_{\tau^b}F^{(l+1-k)}}{(k-1)!(l-k)!}
+2\sum_{k=2}^{l}\frac{F^{(k)}F^{(l+2-k)}}{(k-2)!(l-k)!}=0.
\end{eqnarray}

 Let $I=(i_0,\dots,i_n)\in \mathbb{Z}_{\geq 0}^{n+1}$. We define a function  $\mathrm{Eqc}(n,\mathbf{d},l,I)$ to extract the coefficient of $(\tau^0)^{i_0}\dots (\tau^n)^{i_n}$ in (\ref{eq-wdvv24expand-sqrt-intro}) and uses Theorem \ref{thm-induction1} to obtain an equation on 
 \begin{equation*}
	F^{(2)}(0),F^{(3)}(0),\dots,F^{(l)}(0).
\end{equation*}
In our package this function is implemented by 
\begin{equation*}
       \mbox{\texttt{equationOfConstTerm}}
\end{equation*}
with the input $\{n,\mathbf{d},l,I\}$. We denote $F^{(k)}(0)$ by $z_k$ in the package.
We take a quintic 4-fold as an example. Running
\begin{equation*}
      \mbox{\texttt{equationOfConstTerm}}\  \{4,\{5\},2,\{0,0,0,0,0\}\}
\end{equation*}
returns
\begin{equation*}
      2\,{z}_{2}^{2}-8352000\,{z}_{2}+8719488000000,
\end{equation*}
which factors as
\begin{equation*}
     2 \left({{z}_{2}-2088000}\right)^{2}.
\end{equation*}
So $F^{(2)}(0)=2088000$. 
The matches the formula in Theorem  \ref{thm-induction2}. Then running
\begin{equation*}
      \mbox{\texttt{equationOfConstTerm}}\  \{4,\{5\},3,\{0,0,0,0,0\}\}
\end{equation*}
returns
\begin{equation}\label{eq-sqrtRecursion-3}
      46080\,{z}_{2}^{2}+8\,{z}_{2}{z}_{3}+3119454720000\,{z}_{2}-16704000\,{z}_{3}-6714318458880000000.
\end{equation}
Substituting $z_2=2088000$ into (\ref{eq-sqrtRecursion-3}) we get 0, i.e. a trivial equation. Running
\begin{equation*}
      \mbox{\texttt{equationOfConstTerm}}\  \{4,\{5\},4,\{0,0,0,0,0\}\}
\end{equation*}
returns
\begin{eqnarray}\label{eq-sqrtRecursion-4}
      &-586224\,{z}_{2}^{3}+3190863801600\,{z}_{2}^{2}+1644480\,{z}_{2}{z}_{3}+12\,{z}_{3}^{2}+12\,{z}_{2}{z}_{4}-7369983201945600000\,{z}_{2}\nn\\
      &+6501980160000\,{z}_{3}-25056000\,{z}_{4}+8870266887085670400000000.
\end{eqnarray}
Substituting $z_2=2088000$ into (\ref{eq-sqrtRecursion-4}) we get
\begin{equation*}
     12 \left({{z}_{3}+413985600000}\right)^{2},
\end{equation*}
again a quadratic equation with two equal roots! So we get $F^{(3)}(0)=-413985600000$. Proceeding in this way, we find, in all examples that we checked, that when the dimension $n$ is even and when $\mathbf{d}\neq (3)$, we can always compute $F^{(k)}(0)$ by recursively solving a quadratic equation with two equal roots. We call this phenomenon, and the resulted method to compute the leading terms $F^{(k)}(0)$, the \emph{double root recursion}.

In odd dimensions, one can see from (\ref{eq-wdvv24expand-sqrt-intro}) that there are no quadratic equations to compute all $F^{(k)}(0)$, but only $2\leq k\leq \lfloor \frac{m}{4}\rfloor +1$. For a conjectural remedy we refer the reader to Section \ref{sec:sqrtRecursion-oddDim}.

We state the main conjecture in this paper, both the cases $n$ is even or odd, in a uniform way.
\begin{conjecture}[Double root recursion]\label{conj-sqrtRecursion-unified-intro}(= Conjecture \ref{conj-sqrtRecursion}+\ref{conj-sqrtRecursion-oddDim})
Let $X=X_n(\mathbf{d})$ be a non-exceptional smooth complete intersection, with $n\geq 3$, and multidegree $\mathbf{d}$. Suppose $\mathbf{d}\neq (3)$ or $(2,2)$. Let $m=\mathrm{rank}\ H_{\mathrm{prim}}^n(X)$. Let $\gamma_0,\dots,\gamma_n$ be a basis of $H^*_{\mathrm{amb}}(X)$, and $t^0,\dots,t^n$ be the dual basis. Let $g_{e,f}=(\gamma_e,\gamma_f)$ be the Poincaré pairing, and $(g^{e,f})_{0\leq e,f\leq n}$ the dual matrix of $(g_{e,f})$. Let $E$ be the Euler vector field. Let $G$ be a series in $t^0,\dots,t^n$ and $s$. Let 
\begin{equation*}
      G^{(i)}(t^0,\dots,t^n):= \Big(\big(\frac{\partial}{\partial s}\big)^i G\Big)|_{s=0}.
\end{equation*}
Then the system 
\begin{subequations}\label{eq-system1-G-intro}
\begin{align}[left ={\empheqlbrace}]
& G^{(0)}=F^{(0)},\\
&\sum_{e=0}^n \sum_{f=0}^n \frac{\partial^3 G}{\partial t^a\partial t^b \partial t^e} g^{ef}\frac{\partial^2 G}{\partial s \partial t^f}+2s \frac{\partial^3 G}{\partial s\partial t^a \partial t^b} \frac{\partial^2 G}{\partial s\partial s} =\frac{\partial^2 G}{\partial s\partial t^a} \frac{\partial^2 G}{\partial s\partial t^b},\quad \mbox{for}\ 0\leq a,b\leq n,\\
&\sum_{e=0}^n \sum_{f=0}^n \frac{\partial^2 G}{\partial s\partial t^e}g^{ef}\frac{\partial^2 G}{\partial s\partial t^f}+2s \big(\frac{\partial^2 G}{\partial s\partial s})^2=0, \\
&EG=(3-n)G+\mathsf{a}(n,\mathbf{d})\frac{\partial}{\partial t^1}c, \\
&G^{(k)}(0)=0\ \mbox{unless}\ \frac{(n-2)k-n+3}{\mathsf{a}(n,\mathbf{d})}\in \mathbb{Z},
\end{align}
\end{subequations}
has a unique solution, and the solution can be given by the double root recursion as Conjecture \ref{conj-sqrtRecursion}. More precisely,
\begin{enumerate}
      \item[(i)] For even $l\geq 2$, the polynomial $\mathrm{Eqc}(n,\mathbf{d},l,0)$ of $z_2,\dots,z_l$, after substituting the values of $z_i$ for  $2\leq i\leq \frac{l}{2}$, becomes a complete square in $z_{\frac{l+2}{2}}$;
      \item[(ii)] For odd $l\geq 2$, the polynomial $\mathrm{Eqc}(n,\mathbf{d},l,0)$ of $z_2,\dots,z_l$, after substituting the values of $z_i$ for  $2\leq i\leq \frac{l+1}{2}$, becomes 0;
      \item[(iii)] For $l\geq 2$ and $I\neq 0$, the polynomial $\mathrm{Eqc}(n,\mathbf{d},l,I)$ of $z_2,\dots,z_l$, after substituting the values of $z_i$ for  $2\leq i\leq \lfloor\frac{l+2}{2}\rfloor$, becomes 0.
\end{enumerate}
Finally, denote by $\widetilde{F}$ the unique solution to the above system, and let $F$ be the generating function of genus zero primary Gromov-Witten invariants of $X$. Then
\begin{equation}\label{eq-sqrtRecursion-unified}
F=\begin{cases}
	\widetilde{F},& \mbox{if $n$ is even};\\
     \sum_{k=0}^{\frac{m}{2}}\frac{s^k}{k!}\widetilde{F}^{(k)},& \mbox{if $n$ is odd}.
   \end{cases}
\end{equation}
\end{conjecture}

For a detailed account of this conjecture in even and odd dimensions separately, we refer the reader to Section \ref{sec:sqrtRecursion} and \ref{sec:sqrtRecursion-oddDim}. Here we have the following remarks.
      \begin{enumerate}
            \item In even dimensions, the generating function $F$ satisfies the system (\ref{eq-system1-G-intro})  automatically, as a consequence of Theorem \ref{thm-mainthm1} (ii). So Conjecture  \ref{conj-sqrtRecursion-unified-intro}  is purely a formal statement on $F$.
            \item In odd dimensions, the generating function $F$ has at most $s$-degree $\frac{m}{2}$, and the system  (\ref{eq-system1-G-intro}) holds for $F$ only in the  "mod $s^{\frac{m}{2}}$" sense. But the solution $\widetilde{F}$ can have terms of arbitrarily high $s$-degrees, whose geometric meaning is absent at present.
            \item The parts (ii) and (iii) mean that   the essentially linear recursions will never give nontrivial equations, for the cases other than those in Theorem \ref{thm-intro5}.
            \item For $n$ odd and $\mathbf{d}=(2,2)$, i.e. odd dimensional complete intersections of two quadrics, there is only one possible nonzero $F^{(k)}(0)$ due to the dimension reason. This is $F^{(2)}(0)$, which is computed by the double root recursion at $s$-order 2 as we have seen in Theorem \ref{thm-induction2}. On the other hand  the system (\ref{eq-system1-G-intro}) has a unique formal solution by an essentially linear recursion on the leading terms $G^{(k)}(0)$. In other words, in this case the (genuine) genus 0 Gromov-Witten invariants is also determined by the double root recursion, but we must exclude this case in Conjecture \ref{conj-sqrtRecursion-unified-intro}.
      \end{enumerate}

We sketch in the following table our knowledge and tools on the leading terms $F^{(k)}(0)$ of non-exceptional smooth complete intersections of dimension $\geq 3$. The term \emph{drr} stands for double root recursion.
\vspace{-0.8cm}
\begin{table}[H]
\centering
\caption{$F^{(k)}(0)$}
\begin{tabular}{*{5}{|c}|}
\hline
\diagbox{$(n,\mathbf{d})$}{$F^{(k)}(0)$} & $k=1$ & $k=2$ & $3\leq k\leq \lfloor \frac{m}{4}\rfloor +1$ & $ k>  \frac{m}{4} +1$ \\
\hline
$\mathbf{d}=(3)$, $n=3$ & eigenvector & geometric method & \multicolumn{2}{c|}{geometric method}  \\
\hline
$\mathbf{d}=(3)$, $n\geq 4$ & eigenvector & geometric method &\multicolumn{2}{c|}{essentially linear recursion} \\
\hline
$\mathbf{d}\neq (3)$, even $n$ & eigenvector & drr & \multicolumn{2}{c|}{\textcolor{green}{drr}} \\
\hline
$\mathbf{d}\neq (3)$, odd $n$ & eigenvector & drr & \textcolor{green}{drr}  & \textcolor{red}{drr}\\ 
\hline
\end{tabular}
\label{table:F(k)(0)}
\end{table}
In this table, an item in black color means this term has been computed or  shown to be computable by the indicated method. The items double root recursion in green mean that for given $(n,\mathbf{d})$ one can check whether the double root recursion holds for $F^{(i)}(0)$ inductively from $i$ to $k$, and if the double root recursion holds, then $F^{(k)}(0)$ is obtained. Red color means that for $k$ in this range, even if the  double root recursion holds, the value obtained is \emph{hypothetical}.

We would like to remind the reader that to verify, by examples, that the computed values in the red range match the genuine values of the Gromov-Witten invariants, is quite hard. The first (i.e. the involved genus 0 Gromov-Witten invariant has the smallest length) nontrivial verification is to show, for a complete intersection $X$ of dimension $3$ and  multidegree $(2,2,2)$, 
\begin{equation*}
F^{(9)}(0)=-4251528=-2^{3}3^{12}.
\end{equation*}
This amounts to computing genus zero Gromov-Witten invariants of $X$ with exactly 18 primitive insertions and of (stable map) degree 9.

In Section \ref{sec:sqrtRecursion-examples} we display some $F^{(k)}(0)$ that we compute by our package. From the results at there we make the following conjecture.
\begin{conjecture}\label{conj-integrality-positivity-intro}(= Conjecture \ref{conj-integrality-positivity})
When $\mathbf{d}= (3)$, $F^{(n+1)}(0)=0$. When $\mathbf{d}\neq (3)$, or $\mathbf{d}=(3)$ and $l\neq n+1$, $F^{(l)}(0)$ is a positive integer if $l$ is even, and is a negative integer if $l$ is odd.
\end{conjecture}
The integrality is related to the integrality of genus 0 primary Gromov-Witten invariants of semipositive symplectic manifolds (\cite[Theorem A]{Ruan96}; see also \cite[Theorem 7.1.1]{MS12}). We show the integrality in the odd dimension case. The positivity seems mysterious and we can say nothing.

Our final goal is not only the reconstruction theorems. We wish to find closed formulae. Even the generating function $F^{(0)}$ of ambient quantum cohomology seems hopeless to have a closed formula. But the parts (ii) and (iii) in Conjecture \ref{conj-sqrtRecursion-unified-intro} imply that there are closed formulae of $F^{(k)}$ for $k\geq 2$, in terms of $F^{(i)}$ for $i<k$. We illustrate this by  giving a conjectural formula of $F^{(2)}$ in Section \ref{sec:closedFormula}. 

\vspace{0.3cm}

In the above we are concerned with only the non-exceptional complete intersections. For completeness let us summarize the results for exceptional complete intersections.  The monodromy groups of the exceptional complete intersections are finite groups.

If $X$ is an odd dimensional quadric, $H^*_{\mathrm{prim}}(X)=0$, so we have nothing to do.

If $X$ is an even dimensional quadric, $H^*_{\mathrm{prim}}(X)=1$, and the monodromy group is $\mathbb{Z}/2 \mathbb{Z}\cong \mathrm{O}(1)$. So the above results, especially Theorem \ref{thm-induction1}, remain valid. By the dimension constraint one easily finds that $F^{(k)}(0)=0$. So $F$ is reconstructible. 

If $X$ is an even dimensional complete intersection of two quadrics, the monodromy group is the Weyl group $D_{n+3}$, where $n=\dim X$. In this case we obtain partial results for the 4 point invariants in Section \ref{sec:4points-Invariants}, which is used in  \cite{Hu21} to obtain all 4 point invariants. Then based on this result, we show in \cite{Hu21} the computation of $F$ is reduced to a single unknown invariant of length $n+3$, which we call the \emph{special correlator}.\\

Besides the above-mentioned conjectures and questions, we propose several problems.

\begin{problem}
Use Conjecture \ref{conj-sqrtRecursion-unified-intro} (ii) and (iii) to find closed formulae of $F^{(k)}$ in terms of $F^{(i)}$ for $i<k$, for $k\geq 3$. Prove such formulae for $k\geq 2$.
\end{problem}

\begin{problem}
Find a formula of $F^{(k)}(0)$ for $k>2$, in terms of the matrices $W$ and $M$, as (\ref{eq-intro-higher22}). 
\end{problem}
These two problems together are almost equivalent to proving Conjecture \ref{conj-sqrtRecursion-unified-intro}.

A closely related problem is to find a closed formula for 
\begin{equation}\label{eq-functionOfWM}
	-\sum_{j=0}^n j M_{j}^{1}W_{n}^{j}
+\sfbd\sum_{j=0}^n 
j M_{j}^{1}W_{n- \mathsf{a}(n,\mathbf{d})}^{j}.
\end{equation}
This quantity, which appears in Theorem \ref{thm-induction2} for $F^{(2)}(0)$, looks bizarre at first sight. But it and the related quantities $\mathsf{c}(n,l,\mathbf{d})$ (see (\ref{eq-def-functionC})) seem ubiquitous in the computations in Section \ref{sec:4pointsInvariantsAmbient} and Section \ref{sec:ReconstructionII-F(2)}. In principle for each $(n,\mathbf{d})$ it can be computed from the mirror formula. It is desirable to find a simpler expression, but I am not sure whether this is possible in general.
In the case $\mathsf{a}(n,\mathbf{d})=\frac{n-1}{2}$, we obtain a closed formula in Theorem \ref{thm-F20-WM-a(n,d)=(n-1)/2}, from which Theorem \ref{thm-intro-F20-a(n,d)=(n-1)/2} follows.

\begin{problem}\label{ques-ZarClosure-monodromyGroup}
For a given smooth projective variety $X$, find the Zariski closure of all possible monodromies in $\mathrm{Aut}(H^*(X))$. When this (closure) group is the orthogonal group or symplectic in each degree, check whether the double root recursion is true.
\end{problem}
Here all possible monodromies mean the monodromies induced by all  smooth proper families over connected bases where $X$ sits in as a fibre. Even a rigid variety may have nontrivial monodromies, e.g. the quadric hypersurfaces. For smooth complete intersections in $\mathbb{P}^n$, by \cite{Del73} this closure coincides with the closure of the subgroup preserving the Poincar\'{e} pairing, the integral structure and the hyperplane class.  In any case, the monodromy actions on $H^*(X)$ lie in  the actions by the \emph{mapping class group} of $X$. What can we say about their Zariski closures? In general they are not equal, e.g. $X=\mathbb{P}^2$. Nevertheless, I wish that for the frequently appearing varieties $X$ one can determine the former as a subgroup of the latter defined by some conditions.


\begin{problem}\label{prob-mirrorSymmetry}
Establish the full (numerical) mirror symmetry for a Fano complete intersection $X$ in $\mathbb{P}^n$. Namely, find a B-model construction of the Frobenius manifold induced by the big quantum cohomology of $X$. 
\end{problem}

The numerical mirror symmetry for Fano complete intersections, from the view of A-model, was restricted to, at most, the Frobenius manifold arising from quantum cohomology with only  ambient cohomology classes. Our work provides a possibility to compute the full Frobenius manifold in the A-model. Conversely, 
I wish that mirror symmetry sheds light on Conjecture \ref{conj-sqrtRecursion-unified-intro}, for example gives evidence for the values of $F^{(k)}(0)$ in the red range in Table \ref{table:F(k)(0)}, or gives hints for the unknown special correlator of an even dimensional complete intersection of two quadrics.

\begin{problem}
The monodromy invariance of Gromov-Witten invariants in all genera, or the cohomological field theory, of Fano complete intersections enables us to make symmetric reduction of the induced tautological relations. What information can we extract from this? Can we do this in the Gromov-Witten theory valued in the stable cohomology of $\Mbar_{g,n}$?
\end{problem}
We made progress in genus 1 in this direction (\cite{Hu22}).\\

This paper is organized as follows.

In Section \ref{sec:GW-invariants} we recall the definition and  axiomatic properties of Gromov-Witten invariants. Then we recall some classical results on the small quantum cohomology of complete intersections, including the mirror theorem. 

 In Section \ref{sec:monodromyGroup} we recall some results on the monodromy group action on the primitive cohomology groups of smooth complete intersections in the projective spaces, and apply the invariant theory of orthogonal groups and symplectic groups to obtain a theorem on the form of the generating function of Gromov-Witten invariants.

 In Section \ref{sec:red}
 we compute the symmetric reduction of the WDVV equation by the Zariski closure of the monodromy group. We need an estimate of the dimension of primitive cohomology of smooth complete intersections, which is given in Appendix \ref{sec:estimate-dimPrim}.
 The  symmetric reduction of  the quantum differential equation for the big $J$-function is given in  Appendix \ref{sec:symmetricReduction-JFunction}. 

In Section \ref{sec:structureOfFrobAlg-reconstructionTheorems} we use the equations after symmetric reduction to study the vector field $\Theta$, and use its property to show the non-semisimplicity of $\mathcal{M}_{\mathrm{amb}}$ and compute all the 3-point invariants involving two primitive classes. Then we prove the  reconstruction  I of the invariants involving both primitive and ambient classes from the invariants with only primitive insertions.

Section \ref{sec:4pointsInvariantsAmbient} is a preparation for the computation of 4 point invariant with primitive insertions. We compute certain sums of 4 point invariants with ambient insertions, then the 4 point invariants with exactly two primitive insertions. The use of $\tau$-coordinates is essential in this and the next section.

 In Section \ref{sec:ReconstructionII-F(2)}, we show Theorem \ref{thm-induction2} and Theorem \ref{thm-intro5}.

 In Section \ref{sec:FanoVarietyOfLines} we compute the cohomology ring of the Fano variety of lines on cubic hypersurfaces and by the way we obtain the 4-point invariants with 4 primitive insertions for the cubic hypersurfaces. This provides the initial data for the reconstruction of genus zero GW invariants of all cubic hypersurfaces of dimension $\geq 4$. 

 In Section \ref{sec:genusReduction} We show some vanishing results on Zinger's reduced genus 1 invariants, and then by Zinger's standard versus reduced formula we compute the invariants with only primitive insertions for $n$-dimensional complete intersections of Fano index $n-1$, and the 8 point invariants with only primitive insertions for cubic 3-folds. 
  This provides the initial data for the reconstruction of genus zero GW invariants of all cubic hypersurfaces of dimension $\geq 3$.

In Section \ref{sec:HigherOrderLeadingTerms} we present the double root recursion Conjecture \ref{conj-sqrtRecursion-unified-intro} in even and odd dimensions. Then we present examples and make some observations. Finally we deduce a conjectural closed formula of $F^{(2)}$ from Conjecture \ref{conj-sqrtRecursion-unified-intro} (ii) and (iii). Note that in Conjecture \ref{conj-sqrtRecursion-unified-intro} we exclude the case of cubic hypersurfaces. For the latters, we present also a conjectural closed formula of $F^{(2)}$ by an analogy and checking examples. 

The algorithms in our Macaulay2 package \texttt{QuantumCohomologyFanoCompleteIntersection} is described in Appendix \ref{sec:algorithm}. For the reader  interested in the algorithm to compute the ambient generating function $F^{(0)}$,   Sections \ref{sec:algorithm-summary} to \ref{sec:algorithm-recursionOfAmbientGeneratingFunction} can be read independently, except for notations referred to other sections on only a few occasions. 
     Intermediate formulae are computed explicitly so that they directly match the corresponding functions in our package.

In Appendix \ref{sec:proof-F2(0)-closedFormula} we show Theorem \ref{thm-intro-F20-a(n,d)=(n-1)/2}. We first reduce the statement to the evaluation of a two point invariant with ambient classes and a $\psi$-class as insertions, then we compute this invariant by torus localization.\\


\subsection{Notations and Conventions}\label{sec:nota-conv}
\begin{enumerate}
      \item The ground field will always be $\mathbb{C}$ in this paper.
      \item We denote $H^i(X)=H^i(X;\mathbb{Q})$ for a topological space $X$. For a Kähler manifold we denote by $H^*_{\mathrm{prim}}(X)\subset H^*(X,\mathbb{C})$ the primitive cohomology group.
      \item For $\gamma\in H^i(X;\mathbb{Q})$ or $H^i(X;\mathbb{C})$, we denote $\deg \gamma=\frac{i}{2}$, i.e. the \emph{complex degree} of $\gamma$, even when $i$ is odd. The \emph{real degree} $i$ is denoted by $\deg_{\mathbb{R}} \gamma$, which is used explicitly only in Section \ref{sec:red-odd}.
      \item A basis, or a subset of the basis, of the cohomology group $H^i(X;\mathbb{C})$ is always indexed by  subscripts, while the superscript on one of the base means a power of it, either in the ordinary product or the quantum product, which will be indicated in the context. In some cases,  a symbol with a subscript is equal to the same symbol with the same superscript. For example, if $X$ is a Fano complete intersection of dimension $n$ in a projective space,  and $0\leq i\leq n$, then $\tilde{\sfh}_i=\tilde{\sfh}^i$ (see \S \ref{sec:mirror-theorem}) for the definition of these terms). In this case, we use $\tilde{\sfh}_i$ ($0\leq i\leq n$) when we do computations involving a particular basis, and use $\tilde{\sfh}^i$ ($i\geq 0$) when we deal with the ring structure. I hope that this choice of notations will not confuse the reader.
      \item In this paper, for $\mathbf{d}=(d_1,\dots,d_r)\in \mathbb{Z}_{>0}^r$, we use $X_n(\mathbf{d})$ to denote a smooth complete intersection of multi-degree $\mathbf{d}$ in $\mathbb{P}^{n+r}$, and $m=\mathrm{rank}\ H^{n}_{\mathrm{prim}}(X)$.
      When we say \emph{the} complete intersection of multi-degree $\mathbf{d}$, we mean that we have chosen an arbitrary smooth complete intersection of multi-degree $d$, and what we are doing is independent of the choice. When we are studying a family of complete intersections, we will make the terminology precise. 
      \item Throughout this paper, $n$ denotes the dimension of the complete intersection in consideration, and $\mathbf{d}$ the multi-degree, $\mathsf{a}(n,\mathbf{d})$ the Fano index.
      \item For a smooth projective variety $X$, $\langle \alpha_1\psi^{k_1},\dots,\alpha_l\psi^{k_l}\rangle_{g,l,\beta}^{X}$ means the Gromov-Witten invariant of $X$ of genus $g$, $l$ marked points and degree $\beta$. If the target space is clear from the context, we drop the superscript $X$. The genus $g$ in the subscript will never be omitted in this paper.  If the number of insertions is obvious we usually drop the subscript  $l$. Since the target variety in this paper are Fano complete intersections, the degree of a nonzero invariant is determined by the genus $g$, the dimension of the target variety, and the number of its insertions, and thus when the insertions are  given, the degree can be dropped in principle. However we will still present the degree $\beta$ when it is a definite number; this turns out convenient, e.g. for the application of the divisor equation. So when there are only two numbers in the subscript, the second one is the degree $\beta$.
      \item The generating function with Novikov variables $\mathsf{q}$ of  genus 0 primary Gromov-Witten invariants is denoted by $\sfF$. The specialization at $\mathsf{q}=1$, which is legitimate when $X$ is Fano, is denoted by $F$. 
      \item For the brevity of expressions, we often use the Einstein summation convention. The range of the indices in the summation will be specified at the beginning of the corresponding section.

      \item For a real number $a$, $\lfloor a\rfloor$ denotes the greatest integer not larger than $a$, and $\lceil a\rceil$ denotes the smallest integer not less than $a$.
      
      \item  We denote the small quantum multiplication by $\sqp$, and the big quantum multiplication by $\bqp$. The symbol $\circ$ is preserved to stand for  the composition of operators.
\end{enumerate}

\section{Quick recap of Gromov-Witten invariants}\label{sec:GW-invariants}
We recall the definition of the Gromov-Witten invariants, and their properties that we need to use in this paper. Our main references are
\cite[Chapter VI]{Man99} and
\cite[\S 26.3-26.4]{HKKPTVVZ03}.

Let $X$ be a smooth projective scheme over $\mathbb{C}$ of dimension $n$. Let $k\in \mathbb{Z}_{\geq 0}$, and $\beta\in H_2(X;\mathbb{Z})/\mathrm{tor}$. The stack $\Mbar_{g,k}(X,\beta)$ of stable maps of degree $\beta$ from  genus $g$  semistable curves with $k$ marked points to $X$  is a proper Deligne-Mumford stack and carries a virtual fundamental class (\cite{BF97}, \cite{LT98}) $[\Mbar_{g,k}(X,\beta)]^{\mathrm{vir}}$ of dimension $(1-g)(n-3)+k+c_1(T_X)\cdot \beta$. For each $1\leq i\leq k$, the section $\sigma_i$ pulls back the relative cotangent line bundle of the universal curve to form a line bundle on $\Mbar_{g,k}(X,\beta)$, whose first Chern class is denoted by $\psi_i$; moreover there is an associated \emph{evaluation map} $\mathrm{ev}_i=f\circ \sigma_i$, where $f$ is the universal stable map. For $\gamma_1,\dots,\gamma_k\in H^*(X;\mathbb{Q})$ and $a_1,\dots,a_k\in \mathbb{Z}_{\geq 0}$, there is an associated \emph{Gromov-Witten invariant}
\begin{equation*}
	\langle \psi_1^{a_1}\gamma_1,\dots,\psi_k^{a_k}\gamma_k\rangle_{g,k,\beta}^X:=\int_{[\Mbar_{g,k}(X,\beta)]^{\mathrm{vir}}}\prod_{i=1}^{k}\psi_i^{a_i}\mathrm{ev}_i^*\gamma_i\in \mathbb{Q}.
\end{equation*}
A term like $\psi_i^{a_i}\gamma_i$ in $\langle \psi_1^{a_1}\gamma_1,\dots,\psi_k^{a_k}\gamma_k\rangle_{g,k,\beta}$ is called an \emph{insertion} of this invariant. We say $(g,k,\beta)$ is in the \emph{stable range} if either $2g-2+k>0$ or $\beta$ is a nonzero effective curve class.
It is convenient to use simplified notations in the following occasions:
\begin{enumerate}
	\item[(i)] The superscript $X$ will be omitted when it is obvious;
	\item[(ii)] the subscript $i$ of $\psi_i$ in $\langle \psi_1^{a_1}\gamma_1,\dots,\psi_k^{a_k}\gamma_k\rangle_{g,k,\beta}$ might be dropped when there is no confusion; 
	\item[(iii)] the subscript $k$  might be dropped when $k$ is obvious from the expression;
	\item[(iv)] the subscript  $\beta$ might be dropped when it can be uniquely determined, when it is nonzero, by the insertions and the following condition (\ref{eq-Dim}), which is always the case for Fano complete intersections in projective spaces.
\end{enumerate}

The GW invariants $\langle \psi_1^{a_1}\gamma_1,\dots,\psi_k^{a_k}\gamma_k\rangle_{g,k,\beta}$ with $a_1=\dots=a_k=0$ are called \emph{primary}. 
For brevity we will call a genus 0 primary Gromov-Witten invariant with $k$ marked points a \emph{correlator of length $k$}.

For  two cohomology classes $\gamma_1$ and $\gamma_2$, we denote the Poincaré pairing by
\begin{equation*}
	(\gamma_1,\gamma_2):=\int_X \gamma_1\cup \gamma_2.
\end{equation*}
The GW invariants satisfy the following properties. 
\begin{gather}\label{eq-Deg0}\tag{Deg0}
\langle \gamma_1,\dots,\gamma_k\rangle_{g,k,0}=
\begin{cases}
\int_{X}\gamma_1\cup \gamma_2\cup \gamma_3,& \mbox{if}\ g=0, k=3;\nn\\
-\frac{1}{24}\int_X \gamma_1\cup c_{n-1}(T_X), & \mbox{if}\ g=1, k=1,\\
0, & \mbox{if}\ 2g-2+k\geq 2.
\end{cases}
\end{gather}
The dimension constraint:
\begin{multline}\label{eq-Dim}\tag{Dim}
	\langle \psi_1^{a_1}\gamma_{1},\dots,\psi_k^{a_k} \gamma_{k}\rangle_{g,k,\beta}=0\ \mbox{unless}\\
	\sum_{i=1}^k(a_i+\deg \gamma_{b_k})=(1-g)(n-3)+k+c_1(T_X)\cap \beta.
\end{multline}
The $S_n$-equivariance:
\begin{multline}\label{eq-Sym}\tag{Sym}
	\langle \psi_1^{a_1}\gamma_{1},\dots, \psi_{i-1}^{a_{i-1}}\gamma_{i-1}, \psi_{i}^{a_{i}}\gamma_{i},\dots,\psi_k^{a_k} \gamma_{b_k}\rangle_{g,k,\beta}\\
	=(-1)^{|\gamma_{i-1}|\cdot |\gamma_{i}|}
	\langle \psi_1^{a_1}\gamma_{1},\dots, \psi_{i-1}^{a_{i}}\gamma_{i}, \psi_{i}^{a_{i-1}}\gamma_{i-1},\dots,\psi_k^{a_k} \gamma_{b_k}\rangle_{g,k,\beta}.
\end{multline}
The divisor equation: for $\gamma\in H^2(X)$,
\begin{gather}\label{eq-Div}\tag{Div}
	\langle \psi^{a_1}_1\gamma_1,\dots,\psi^{a_k}_k\gamma_k,\gamma\rangle_{g,k+1,\beta}=(\gamma\cap \beta)\langle \psi^{a_1}_1\gamma_1,\dots,\psi^{a_k}_k\gamma_k\rangle_{g,k,\beta}\\
	+\sum_{i=1}^{k}\langle \psi^{a_1}_1\gamma_1,\dots,\psi^{a_{i-1}}_{i-1}\gamma_{i-1},\psi^{a_i-1}_i\gamma_{a_i}\cup \gamma,\psi^{a_{i+1}}_{i+1}\gamma_{i+1},\dots, \psi^{a_k}_k\gamma_k\rangle_{g,k,\beta}.
\end{gather}
In particular, when $\gamma=c_1(T_X)$, we call this equation the \emph{divisor equation with the first Chern class}.\\
\noindent The string equation: for $(g,k,\beta)$ in the stable range,
\begin{gather}\label{eq-String}\tag{Str}
	\langle \psi^{a_1}_1\gamma_1,\dots,\psi^{a_k}_k\gamma_k,1\rangle_{g,k+1,\beta}\\
	=\sum_{i=1}^{k}\langle \psi^{a_1}_1\gamma_1,\dots,\psi^{a_{i-1}}_{i-1}\gamma_{i-1},\psi^{a_i-1}_i\gamma_{a_i},\psi^{a_{i+1}}_{i+1}\gamma_{i+1},\dots, \psi^{a_k}_k\gamma_k\rangle_{g,k,\beta},
\end{gather}
and in particular, the \emph{fundamental class axiom}:
\begin{gather}\label{eq-FCA}\tag{FCA}
	\langle 1, \gamma_{1},\dots,\gamma_{k-1}\rangle_{g,k,\beta}=
	\begin{cases}
	(\gamma_1,\gamma_2),& \mbox{if}\ g=0, k=3, \beta=0;\\
	0, & \mbox{if}\ 3g-3+k\geq 1\ \mbox{or}\ \beta\neq 0.
	\end{cases}
\end{gather}
The dilaton equation:
\begin{gather}\label{eq-Dilaton}\tag{Dil}
		\langle \psi_1^{a_1}\gamma_1,\dots,\psi_k^{a_k}\gamma_k,\psi_{k+1}\rangle_{g,k+1,\beta}=(2g-2+k)\langle\psi_1^{a_1}\gamma_1,\dots,\psi_k^{a_k}\gamma_k\rangle_{g,k,\beta}.
\end{gather}

Now suppose that  $\gamma_0,\dots,\gamma_N$ is a basis of $H^*(X;\mathbb{Q})$. Let $g_{i,j}=(\gamma_i,\gamma_j)$,and $(g^{ij})$ the inverse matrix. Then there is the genus 0 topological recursion relation
\begin{align}
	&\langle \psi^{a_1+1}\gamma_{b_1},\psi^{a_2}\gamma_{b_2},\psi^{a_3}\gamma_{b_3},
	\prod_{i\in S}\psi_i^{a_i}\mathrm{ev}_i^*\gamma_{b_i}\rangle_{0,\beta}\nn\\
	=&\sum_{\begin{subarray}{c}\beta_1+\beta_2=\beta\\ S_1\sqcup S_2=S\end{subarray}} \sum_{e,f}\Big((\pm)\langle \psi^{a_1}\gamma_{b_1},\prod_{i\in S_1}\psi_i^{a_i}\mathrm{ev}_i^*\gamma_{b_i}, \gamma_e\rangle_{0,\beta_1}\nn\\
      &g^{ef}
	\langle \gamma_f, \psi^{a_2}\gamma_{b_2},\psi_i^{a_3}\gamma_{b_3},
	\prod_{i\in S_2}\psi^{a_i}\mathrm{ev}_i^*\gamma_{b_i}\rangle_{0,\beta_2}\Big),\label{eq-TRR0}\tag{TRR0}
\end{align}
and the genus 1 topological recursion relation (see e.g., \cite{Get98}):
\begin{align}
	&\langle \psi^{a+1}\gamma_{b_1},
	\prod_{i\in S}\psi_i^{a_i}\mathrm{ev}_i^*\gamma_{b_i}\rangle_{1,\beta}\nn\\
	=&\sum_{\begin{subarray}{c}\beta_1+\beta_2=\beta\\ S_1\sqcup S_2=S\end{subarray}} \sum_{e,f}\Big((\pm)\langle \psi^{a}\gamma_{b_1},\prod_{i\in S_1}\psi_i^{a_i}\mathrm{ev}_i^*\gamma_{b_i}, \gamma_e\rangle_{0,\beta_1}g^{ef}
	\langle \gamma_f, \prod_{i\in S_2}\psi^{a_i}\mathrm{ev}_i^*\gamma_{b_i}\rangle_{1,\beta_2}\nn\\
	&+ \frac{1}{24}\sum_{e,f}\langle \prod_{i\in S}\psi_i^{a_i}\mathrm{ev}_i^*\gamma_{b_i}, \gamma_e,g^{ef}\gamma_f
	\rangle_{0,\beta}
	\Big),\label{eq-TRR1}\tag{TRR1}
\end{align}
where $(\pm)$ means the sign arising from the transpositions of odd degree classes. For example, if $\gamma_b$ and $\gamma_c$ are classes of pure degrees, then $(\pm)=(-1)^{|\gamma_b|\cdot |\gamma_c|}$, where $|\cdot|$ stands for the degree of cohomology classes.

Let $T^0,\dots,T^N$ be the dual basis with respect to $\gamma_0,\dots,\gamma_N$, then the genus $g$ generating function is defined as
 \begin{equation}\label{eq-def-generatingFunction}
 	\mathcal{F}_g(T^0,\dots,T^{N},\mathsf{q})=\sum_{k\geq 0} \sum_{\beta} \frac{1}{k!}\big\langle \sum_{i=0}^N \gamma_i T^i,\dots,\sum_{i=0}^N \gamma_i T^i\big\rangle_{g,k,\beta} \mathsf{q}^{\beta},
 \end{equation}
where the invariants outside of the stable range are defined to be zero, by convention. 

\subsection{Genus zero generating functions}
The reference for this subsection is \cite[Chapter I to III]{Man99}. 
 In this paper, only the genus zero generating function will be manifestly used. With the notations as above, we denote
\begin{equation}\label{eq-generatingFunction-with-q}
	\mathsf{F}(T^0,\dots,T^{N},\mathsf{q})=\mathcal{F}_0(T^0,\dots,T^{N},\mathsf{q}),
\end{equation}
and
\begin{equation}\label{eq-generatingFunction-without-q}
	F(T^0,\dots,T^{N})=\mathcal{F}_0(T^0,\dots,T^{N},\mathsf{q}=1).
\end{equation}
Note that $\mathsf{F}(T^0,\dots,T^{N},\mathsf{q})$ always makes sense as a formal series, while for $F(T^0,\dots,T^{N})$ to make sense, one needs to address the convergence issue.
For Fano complete intersections of dimension at least 3 in projective spaces, $H_2(X;\mathbb{Z})$ is free of rank 1, and for a nonzero invariant the degree $\beta$ is determined by its insertions, as can be seen from (\ref{eq-Dim}). So in these cases $F(T^0,\dots,T^{N})$ also makes sense as a formal series, and it determines $\mathsf{F}(T^0,\dots,T^{N},\mathsf{q})$. However, in any case, $\mathsf{F}(T^0,\dots,T^{N},\mathsf{q})$ has the advantage that the divisor equation can be written as a differential equation. In fact, let $\beta_1,\dots,\beta_r$ be a  basis of $H_2(X;\mathbb{Z})/\mathrm{tor}$, and $\mathsf{q}^i=\mathsf{q}^{\beta_i}$. Suppose that $\gamma_0,\dots, \gamma_N\in H^*(X)$ is a basis of pure degrees.  Then the divisor equation (\ref{eq-Div}) for primary genus zero GW invariants is equivalent to the system
\begin{equation*}
	\frac{\partial \mathsf{F}}{\partial T^{j}}= \sum_{i=1}^r (\gamma_j\cap \beta_i) \mathsf{q}^i\frac{\partial \mathsf{F}}{\partial \mathsf{q}^i},
\end{equation*}
where $T^j$ runs over the coordinates in $T^0,\dots, T^N$ that $\gamma_j\in H^2(X)$.

 Both types of generating functions satisfy the \emph{WDVV} equation
\begin{gather}\label{eq-WDVV}\tag{WDVV}
\sum_{e=0}^N \sum_{f=0}^N \frac{\partial^3 \mathsf{F}}{\partial T^a \partial T^b\partial T^e}g^{ef}\frac{\partial^3 \mathsf{F}}{\partial T^f \partial T^c\partial T^d}
=\sum_{e=0}^N \sum_{f=0}^N (\pm)\frac{\partial^3 \mathsf{F}}{\partial T^a \partial T^c\partial T^e}g^{ef}\frac{\partial^3 \mathsf{F}}{\partial T^f \partial T^b\partial T^d}.
\end{gather}

The \emph{big quantum product} is defined as
\begin{equation*}
 	\gamma_a\bqp \gamma_b=\sum_{e}\sum_f\frac{\partial^3 \mathsf{F}}{\partial T^a \partial T^b\partial T^e}g^{ef}\gamma_f,
 \end{equation*} 
 and the \emph{small quantum product} is defined as
\begin{equation*}
 	\gamma_a\sqp \gamma_b=\gamma_a\sqp \gamma_b|_{T^0=\dots=T^N=0}.
 \end{equation*} 

Now suppose $\gamma_0,\dots,\gamma_N$ have pure degrees. Let
\[
c_1(T_X)=\sum_{i=0}^N a_i \gamma_i.
\]
Of course $a_i=0$ unless $|\gamma_i|=2$.
Then (\ref{eq-Dim}) and the divisor equation (\ref{eq-Div}) with first Chern class implies that 
\begin{gather}\label{eq-EulerVectorField}\tag{EV}
	EF=(3-n)F+  \sum_{i=0}^N a_i \frac{\partial }{\partial T^i} c,
\end{gather}
where the \emph{Euler vector field}
\begin{equation}\label{eq-EV-0}
	E=\sum_{i=0}^N(1-\frac{|\gamma_i|}{2})\frac{\partial }{\partial T^i}+\sum_{i=0}^N a_i \frac{\partial }{\partial T^i}
\end{equation}
and $c$ is the classical triple intersection form (\ref{eq-tripleIntersectionForm}). 
We make the following observation: for primary genus zero GW invariants, in the three properties (\ref{eq-Dim}), (\ref{eq-Div}) with first Chern class, and (\ref{eq-EulerVectorField}), any two of the three imply the third one.

The WDVV equation and the Euler vector field for $\mathsf{F}$ give rise to a \emph{formal Frobenius (super-)manifold}. When $F$ is convergent, we get a Frobenius manifold defined on the convergence domain.

\subsection{Quantum cohomology of complete intersections with ambient insertions}
\label{sec:mirror-theorem}
We recall some results of Beauville and Givental, which will be used in this paper.

For $n\in \mathbb{Z}_{\geq 1}$ and $\mathbf{d}=(d_1,\dots,d_r)\in \mathbb{Z}_{\geq 2}^r$, we denote by $X_n(\mathbf{d})$ a smooth complete intersection in $\mathbb{P}^{n+r+1}$ of multi-degree $\mathbf{d}$. Throughout this paper, we denote by $\sfh$ the hyperplane cohomology class on the relevant projective space and its restriction on the complete intersections.
Recall the notation (\ref{eq-alb-intro}).
Let $\mathsf{a}(n,\mathbf{d}):=n+r+1-|\mathbf{d}|$ be the Fano index of $X_n(\mathbf{d})$, and
\begin{eqnarray}\label{eq-ss5}
\tilde{\sfh}=\left\{
\begin{array}{cc}
\sfh, & \mathsf{a}(n,\mathbf{d})\geq 2,\\ 
\sfh+\elld\mathsf{q}, & \mathsf{a}(n,\mathbf{d})=1.
\end{array}\right.
\end{eqnarray}
By Lemma \ref{lem-subspace} and the following Corollary \ref{cor-mono3}, the ambient subgroup  $H^*_{\mathrm{amb}}(X;\mathbb{C})$ is closed under (small and big) quantum products.
The structure of the small quantum products of ambient classes is  computed in \cite[Corollary 9.3 and 10.9]{Giv96}.
\begin{theorem}[Givental]\label{thm-Givental-smallQuantumCohomology}
The ambient small quantum cohomology ring of $X_n(\mathbf{d})$ is 
\begin{equation}\label{eq-Givental-smallQuantumCohomology}
	\mathbb{C}[\tilde{\sfh}]/(\tilde{\sfh}^{n+1}-\sfbd\tilde{\sfh}^{n+1-\mathsf{a}(n,\mathbf{d})}\mathsf{q}).
\end{equation}
\end{theorem}
 Note that here $\tilde{\sfh}^i$ is the small quantum product 
\begin{equation}\label{eq-quantumPower}
\underbrace{\tilde{\sfh}\sqp \cdots \sqp \tilde{\sfh}}_{i\ \mbox{factors}}.
\end{equation}
Let $\sfh_i$ be the ordinary $i$-th power of $\sfh$ in the singular cohomology. For the consistency of our use of subscript for a basis of the cohomology groups (see \S \ref{sec:nota-conv}), we denote by $\tsfh_i$ the $i$-th quantum power (\ref{eq-quantumPower}).
To the best knowledge of the author, for general $(n,\mathbf{d})$ there is  no explicit formula of the transformation matrix between  $\tsfh_i$ and $\sfh_i$. But when the Fano index is  large enough, this is possible for the degree reason (\cite[Main Theorem, (1.8) and (2.1)-(2.3)]{Bea95}):
\begin{theorem}[Beauville]\label{thm-Beuville}
Suppose $2 \mathsf{a}(n,\mathbf{d})\geq n+1$. For  $0\leq p\leq n+1-\mathsf{a}(n,\mathbf{d})$, let 
\begin{equation}\label{eq-def-Beauville-ell}
	\ell_p=\frac{1}{\prod_{i=1}^{r}d_i}\langle \mathsf{h},\mathsf{h}_{n-p},\mathsf{h}_{\mathsf{a}(n,\mathbf{d})-1+p}\rangle_{0,3,1}^{X_{n}(\mathbf{d})}. 
\end{equation}
Then $\ell_p$ are integers that  are explicitly given by \cite[(2.1)]{Bea95}. In particular, 
\begin{equation}\label{eq-Beauville-ell-l0-l1}
	\ell_0=\elld,\
	\ell_1=\elld\sum_{\begin{subarray}{c}1\leq i\leq r\\
	1\leq j<d_i\end{subarray}}\frac{d_i-j}{j},
\end{equation}
They satisfy
\begin{equation}\label{eq-Beauville-dual}
	\ell_p=\ell_{n+1- \mathsf{a}(n,\mathbf{d})-p},
\end{equation}
and 
\begin{equation}\label{eq-Beauville-ell-sum}
	\sum_{p=0}^{n+1-\mathsf{a}(n,\mathbf{d})}\ell_p=\sfbd.
\end{equation}
Moreover,  
\begin{equation}\label{eq-Beauville-tsfhTosfh}
	\tsfh_p=\begin{cases}
	\sfh_p,& \mbox{if}\ 0\leq p< \mathsf{a}(n,\mathbf{d});\\
	\sfh_p+(\sum_{i=0}^{p-\mathsf{a}(n,\mathbf{d})}\ell_i)\sfh_{p- \mathsf{a}(n,\mathbf{d})},& \mbox{if}\ \mathsf{a}(n,\mathbf{d})\leq p\leq n,
	\end{cases}
\end{equation}
and
\begin{equation}\label{eq-Beauville-sfhTotsfh}
	\sfh_p=\begin{cases}
	\tsfh_p,& \mbox{if}\ 0\leq p< \mathsf{a}(n,\mathbf{d});\\
	\tsfh_p-(\sum_{i=0}^{p-\mathsf{a}(n,\mathbf{d})}\ell_i)\tsfh_{p- \mathsf{a}(n,\mathbf{d})},& \mbox{if}\ \mathsf{a}(n,\mathbf{d})\leq p\leq n.
	\end{cases}
\end{equation}
\end{theorem}

Both Theorem \ref{thm-Givental-smallQuantumCohomology} and Theorem \ref{thm-Beuville}  can be deduced  from the mirror theorem. 
For brevity we recall in the following the mirror theorem in the case that the Fano index $\mathsf{a}(n,\mathbf{d})> 1$. This suffices for its only manifest use in the main body of this paper, in Section \ref{sec:application-of-SvR-4points}.
For a complete treatment we refer the reader to Appendix \ref{sec:algorithm}.

Denote by $j$ the embedding $X_{n}(\mathbf{d})\hookrightarrow \mathbb{P}^{n+r}$.  
The small $J$-function of degree $\beta$ is defined as the generating function of the 1-point descendant invariants of $ \mathbb{P}^{n+r}$ twisted by 
$\mathcal{O}_{\mathbb{P}^{n+r}}(\mathbf{d})$, or equivalently by the quantum hyperplane theorem, 
\begin{eqnarray*}
J_{\beta}(z):=
j_*\mathrm{ev}_{*}\frac{[\Mbar_{0,1}(X_{n}(\mathbf{d}),\beta)]^{\mathrm{vir}}}{z(z-\psi)}
=\frac{1}{\prod_{i=1}^{r}d_i}\sum_{i=0}^{n}
\langle \sfh_{i}\psi^{n-2+\mathsf{a}(n,\mathbf{d})\beta-i}\rangle_{0,1,\beta}^{X_{n}(\mathbf{d})}\frac{\sfh_{n-i}}{z^{n+\mathsf{a}(n,\mathbf{d})\beta-i}}.
\end{eqnarray*}
The mirror theorem (\cite{Giv96}, \cite{LLY97}) says
\begin{eqnarray}\label{eq-mirrorFormula-a(n,d)>1}
J_{\beta}(z)=\frac{\prod_{i=1}^{r}\prod_{k=1}^{d_i \beta}(d_i \sfh+kz)}{\prod_{k=1}^{\beta}(\sfh+kz)^{n+r+1}},
\end{eqnarray}
where the power is taken in the ordinary cohomology ring.

\section{The monodromy group and its invariants}\label{sec:monodromyGroup}
For $n\geq 1$, $r\geq 1$, $\mathbf{d}=(d_{1},\dots,d_{r})\in \mathbb{Z}^{r}$, where $d_{i}\geq 2$ for $1\leq i\leq r$, denote the
family of all nonsingular complete intersections in $\mathbb{P}^{n+r}$ of multi-degree $(d_{1},\dots,d_{r})$ by $\pi_{\mathbf{d}}:\mathcal{X}_{n}(\mathbf{d})\rightarrow S_{n}(\mathbf{d})$. The cohomology groups $\coprod_{s\in S}H^{*}(\mathcal{X}_{n}(\mathbf{d})_{s},\mathbb{Q})$ form a local system over $S_{n}(\mathbf{d})$. In this context, the \emph{deformation invariance} (e.g. \cite[Theorem 4.2]{LT98}) states that
\begin{theorem}\label{thm-mono1}
Let $C$ be an irreducible smooth curve, and $f:C\rightarrow S_n(\mathbf{d})$ a morphism. Let $f^*\mathcal{X}_{n}(\mathbf{d}) $ be the pulled back family over $C$.
Let $\gamma_{1},\dots, \gamma_{k}$ be  sections of the local system 
$\coprod_{c\in C}H^{*}(f^*\mathcal{X}_{n}(\mathbf{d})_{c},\mathbb{Q})$. Then the genus $g$ Gromov-Witten invariant
\begin{eqnarray*}
\langle \gamma_{1,c},\dots, \gamma_{k,c}\rangle^{f^*\mathcal{X}_{n}(\mathbf{d})_{c}}_{g,k,d}
\end{eqnarray*}
is independent of $c\in C$.
\end{theorem}
The proof of  \cite[Theorem 4.2]{LT98} works for any  coefficient  group $\mathbb{Z}/N \mathbb{Z}$ where $N\in \mathbb{N}$, over which the virtual fundamental class can be defined (a priori,  it is defined over $\mathbb{Q}$). In particular, for fixed $(g,k,d)$ one can replace $\mathbb{Q}$ by $\mathbb{Q}_{\ell}$ or even $\mathbb{Z}/\ell^N \mathbb{Z}$ for sufficiently large primes $\ell$.

For  a  fixed fibre $ \mathcal{X}_{n}(\mathbf{d})_{s}$, the fundamental group  $\pi_{1}(S_n(\mathbf{d}),s)$ acts on $H^{*}(\mathcal{X}_{n}(\mathbf{d})_{s})$.  The following is a consequence of Theorem \ref{thm-mono1}. 
\begin{corollary}\label{cor-mono2}
For any closed point $s\in S$, and $\gamma_1,\dots,\gamma_k\in H^{*}(\mathcal{X}_{n}(\mathbf{d})_{s})$, and $h\in \pi_{1}(S_n(\mathbf{d}),s)$ we have
\begin{equation}\label{eq-mono2}
\langle \gamma_{1},\dots, \gamma_{k}\rangle_{g,k,d}=\langle h.\gamma_{1},\dots, h.\gamma_{k}\rangle_{g,k,d}, 
\end{equation}
where $h.\gamma$ denotes the monodromy action of $h$ on $\gamma$.
\end{corollary}
To Theorem \ref{thm-mono1}, we need to find a cover of the base scheme such that  the selected classes extend to global sections of the local system. Such covers in general may not be algebraic. So we need to use 
$\ell$-adic cohomology.
\begin{proof}
Let $\mathcal{X}=\mathcal{X}_{n}(\mathbf{d})$,  $X=\mathcal{X}_{n}(\mathbf{d})_{s}$, $S=S_n(\mathbf{d})$, and $\pi=\pi_{\mathbf{d}}$. By the classical Zariski-Van Kampen theorem, there exists a pencil $L$ containing $s$ such that $\pi_1(L\cap S,s)\rightarrow \pi_1(S,s) $ is surjective. So we replace $S$ by $L\cap S$, and show the statement for the family restricted to it. 

The monodromy action of $\pi_{1}(S,s)$ factors through the étale fundamental group, which is the profinite completion of $\pi_{1}(S,s)$. More precisely, there is a natural commutative diagram
\[
\xymatrix{
	\pi_{1}(S,s) \ar[r] \ar[d] & \mathrm{Aut}\big(H^*(X;\mathbb{Z})\big) \ar[d] \\
\pi_{1,\mathrm{\acute{e}t}}(S,s) \ar[r] & \mathrm{Aut}\big(H_{\mathrm{\acute{e}t}}^*(X;\mathbb{Z}_{\ell})\big).
}
\]
In view of the remark after Theorem \ref{thm-mono1}, we then need only arbitrarily choose a sufficiently large prime, and show the statement for all $\gamma_i\in H^*(X; \mathbb{Z}/\ell^N)$ for all $N>0$.

Let $\eta=\Spec(K)$ be the generic point of $S$, and $\overline{K}$ an algebraic closure of $K$, and $\bar{\eta}=\Spec(\overline{K})$. Choose an imbedding $K\hookrightarrow \overline{K}$, and thus a morphism $\bar{\eta}\rightarrow \eta$. By the smooth and proper base change theorems, there is a canonical specialization isomorphism
\begin{equation}\label{eq-specializationIsomorphism}
 	 H_{\mathrm{\acute{e}t}}^i(X; \mathbb{Z}/\ell^N)\cong (R^i\pi_* \mathbb{Z}/\ell^N)_{s}
       \xrightarrow{\sim} (R^i\pi_* \mathbb{Z}/\ell^N)_{\bar{\eta}}\cong H_{\mathrm{\acute{e}t}}^*(\mathcal{X}_{\bar{\eta}}; \mathbb{Z}/\ell^N).
 \end{equation} 
 There is a surjective homomorphism (recall that $S$ is now a curve)
 \begin{equation}\label{eq-surjectiveHomomorphism-fundGroups}
 	\mathrm{Gal}(\overline{K}/K)\twoheadrightarrow \pi_{1,\mathrm{\acute{e}t}}(S,s).
 \end{equation}
The lifting of the action of  $\pi_{1,\mathrm{\acute{e}t}}(C,s)$ to $\mathrm{Gal}(\overline{K}/K)$, on $H_{\mathrm{\acute{e}t}}^*(\mathcal{X}_{\bar{\eta}}; \mathbb{Z}/\ell^N)$ via (\ref{eq-specializationIsomorphism}), coincides with the action induced by the automorphisms of $\overline{K}$, as  $\mathcal{X}_{\bar{\eta}}=\mathcal{X}\times_K \overline{K}$. 
We denote still by $\gamma_i$ the element of $H_{\mathrm{\acute{e}t}}^*(\mathcal{X}_{\bar{\eta}}; \mathbb{Z}/\ell^N)$ via (\ref{eq-specializationIsomorphism}). Then we are left to show (\ref{eq-mono2}) for $h\in \mathrm{Gal}(\overline{K}/K)$. By definition of stalks in étale cohomology, there is an open subset $U$ of $S$ and an étale cover $\widetilde{U}\rightarrow U$ such that $\gamma_i$ are induced by sections in $H^*(\mathcal{X}\times_S \tilde{U},\mathbb{Z}/\ell^N)$. Applying Theorem \ref{thm-mono1} to $\widetilde{U}\rightarrow S\rightarrow S_n(\mathbf{d})$ we are done.
\end{proof}

\begin{remark}\label{rem:symplectic-mono}
The above algebraic proof is somewhat clumsy. One can also use the symplectic definition of Gromov-Witten invariants, and then Corollary \ref{cor-mono2} is straightforward. In fact, for the later use in this paper, we need only the statement for Fano complete intersections. So  the definition of Gromov-Witten invariants for \emph{semipositive} symplectic manifolds (e.g. \cite{Ruan96}, \cite{MS12}) suffices for us.
\end{remark}

We focus on  the middle dimensional cohomology groups $H^{n}(\mathcal{X}_{n}(\mathbf{d})_{s})$. The Poincar\'{e} pairing induces a bilinear form $Q$ on $H^{n}(\mathcal{X}_{n}(\mathbf{d})_{s},\mathbb{C})$, which is symmetric when $n$ is even, and skew-symmetric when $n$ is odd. By the invariant cycle theorem, the subspace of the invariant classes under the monodromy action of $\pi_{1}(S_{n}(\mathbf{d}),s)$ is zero when 
$n$ is odd, and is spanned by $\sfh_{\frac{n}{2}}$ when $n$ is even.
In the former case, denote $H^{n}_{\mathrm{var}}(\mathcal{X}_{n}(\mathbf{d})_{s})=H^{n}(\mathcal{X}_{n}(\mathbf{d})_{s})$.
 In the latter case, denote $H^{n}_{\mathrm{var}}(\mathcal{X}_{n}(\mathbf{d})_{s})=(\mathbb{C}\sfh^{n/2})^{\perp}$, i.e., the orthogonal complement of $H^{n/2}$ with respect to $Q$. We have
 \begin{eqnarray*}
 H^{n}_{\mathrm{var}}(\mathcal{X}_{n}(\mathbf{d})_{s})=
 H^{n}_{\mathrm{prim}}(\mathcal{X}_{n}(\mathbf{d})_{s}).
 \end{eqnarray*}
Then since $H^{*}_{\mathrm{amb}}(\mathcal{X}_{n}(\mathbf{d})_{s})$ is fixed by the monodromy action of $\pi_{1}(S_{n}(\mathbf{d}),s)$,  Corollary \ref{cor-mono2} implies
\begin{corollary}\label{cor-mono3}
For $\gamma\in  H^{n}_{\mathrm{prim}}(\mathcal{X}_{n}(\mathbf{d})_{s})$, 
and $\gamma_{i}\in H^{*}_{\mathrm{amb}}(\mathcal{X}_{n}(\mathbf{d})_{s})$, $1\leq i\leq k$,
\begin{eqnarray*}
\langle \gamma_{1},\dots, \gamma_{k},\gamma\rangle_{g,k+1,d}=0.
\end{eqnarray*}
\end{corollary}

Now we fix a fibre over $0\in S$ and denote it by $X_n(\mathbf{d})$. Then 
 $V_n(\mathbf{d}):=H^{n}_{\mathrm{prim}}(X_n(\mathbf{d}))$ forms a representation of the fundamental group $\pi_{1}(S_n(\mathbf{d}),0)$,
 and we denote it by $\rho_n(\mathbf{d})$. Let $M_{n}(\mathbf{d})$ be the Zariski closure of
   $\rho_{n}(\mathbf{d})\big(\pi_{1}(S_n(\mathbf{d}),0)\big)$ in $\mathrm{GL}(V_n(\mathbf{d}))$. Then  $M_{n}(\mathbf{d})$ lies in $\mathrm{O}(V_n(\mathbf{d}))$ or $\mathrm{Sp}(V_n(\mathbf{d}))$ when $n$ is even or odd, respectively, where 
   $\mathrm{O}(V_n(\mathbf{d}))$ and $\mathrm{Sp}(V_n(\mathbf{d}))$  are the orthogonal group or the symplectic group with respect to the nondegenerate bilinear form $Q$. 
   
   For (abstract) subgroups of complex linear algebraic groups, taking images by homomorphisms of algebraic groups and taking Zariski closures are commuting. So to find the decomposition of the representations of $\pi_{1}(S_{\mathbf{d}},0)$ on the tensor products, symmetric powers and exterior powers of $V_n(\mathbf{d})$, we only need to consider the corresponding representations of the closed algebraic group $M_n(\mathbf{d})$.\\

\begin{definition}\label{def-exceptional}
$X_n(\mathbf{d})$ is called an \emph{exceptional complete intersection} (in the projective spaces), if it is one of the following cases:
(i) $X_{n}(2)$;
(ii) $X_{n}(2,2)$, $n$ is even;
(iii) $X_2(3)$, i.e., a cubic surface.
Otherwise $X_n(\mathbf{d})$ is called \emph{non-exceptional}.
\end{definition}
Then we recall:
\begin{theorem}[{\cite{Del73}, \cite{PS03}}]\label{thm-DPS}
 If $X_n(\mathbf{d})$ is a non-exceptional complete intersection, we have $M_{n}(\mathbf{d})=\mathrm{O}(V_{n}(\mathbf{d}))$ when $n$ is even,
  and $M_{n}(\mathbf{d})=\mathrm{Sp}(V_{n}(\mathbf{d}))$ when $n$ is odd. For the  exceptional cases, the ranks of the primitive cohomology are:
\begin{align*}
\dim H_{\mathrm{prim}}^{n}(X_{n}(2))&=\begin{cases}
0, & 2\nmid n\\
1, & 2| n
\end{cases},\\
\dim H_{\mathrm{prim}}^{n}(X_n(2,2))&=n+3,\ \mbox{where}\ 2| n,\\
\dim H_{\mathrm{prim}}^{2}(X_2(3))&=6,
\end{align*}
and the corresponding $M_{n}(\mathbf{d})$ and $\rho_{n}(\mathbf{d})$ are:
\begin{itemize}
\item[(i)] $X_{n}(2)$, $n$ is even, $M_{n}(2)\cong \mathbb{Z}/2\mathbb{Z}$, $\rho_{n}(2)$ is  the unique nontrivial 1-dim representation;
\item[(ii)] $X_{n}(2,2)$, $n$ is even, $M_{n}(2,2)$ is the Weyl group of the root system D$_{n+3}$, $\rho_{n+3}(2,2)$ is the standard representation of this Weyl group ;
\item[(iii)] $X_2(3)$, $M_{2}(3)$ is the Weyl group of the root system E$_{6}$, $\rho_{2}(3)$ is the standard representation of this Weyl group.
\end{itemize}
\end{theorem}

\begin{theorem}\phantomsection\label{thm-monodromythm}
\begin{itemize}
\item[(i)] Let $X$ be an even-dimensional non-exceptional complete intersection of dimension $n\geq 4$. Suppose $\gamma_0,\dots,\gamma_{n+m}$ be a basis of $H^{*}(X,\mathbb{C})$, where $\gamma_{0},\dots,\gamma_{n}$ is a basis of $H_{\mathrm{amb}}^{*}(X,\mathbb{C})$, and 
$\gamma_{n+1},\dots,\gamma_{n+m}$ is an orthonormal basis of $H^{n}_{\mathrm{prim}}(X)$, in the sense $Q(\gamma_{i},\gamma_{j})=\delta_{ij}$. Let 
$t^{0},\dots, t^{n+m}$ be the dual basis of $\gamma_{0},\dots,\gamma_{n+m}$. Then for  $g\geq 0$, the genus $g$ generating function $\mathcal{F}_{g}$ can be written in a unique way as a series in $t^{0},\dots,t^{n}$ and $s$, where
\begin{eqnarray*}
s=\sum_{\mu=n+1}^{n+m}\frac{(t^{\mu})^2}{2}.
\end{eqnarray*}
\item[(ii)] Let $X$ be an odd-dimensional complete intersection of dimension $n\geq 3$. Suppose $\gamma_0,\dots,\gamma_{n+m}$ be a basis of $H^{*}(X,\mathbb{C})$, 
where $\gamma_{0},\dots,\gamma_{n}$ is a basis of $H_{\mathrm{amb}}^{*}(X,\mathbb{C})$, and $\gamma_{n+1},\dots,\gamma_{n+m}$ is a symplectic basis of $H^{n}_{\mathrm{prim}}(X)$. Let 
$t^{0},\dots, t^{n+m}$ be the dual basis of $\gamma_{0},\dots,\gamma_{n+m}$. Then for any $g\geq 0$, the genus $g$ generating function $\mathcal{F}_{g}$ can be written in a unique way as a series in $t^{0},\dots,t^{n}$ and $s$ with the degree of $s$ not greater than $\frac{m}{2}$, where
\begin{eqnarray*}
s=-\sum_{\mu=n+1}^{n+\frac{m}{2}}t^{\mu}t^{\mu+\frac{m}{2}}.
\end{eqnarray*}
\end{itemize}
\end{theorem}
\begin{proof} When $\dim X\geq 3$, the Picard group is of rank one. Then by Lemma \ref{lem-inv1}, our theorem is an immediate consequence of Theorem \ref{thm-DPS} and the first fundamental theorems of polynomial invariants of the orthogonal groups and symplectic groups. When $\dim X$ is even, see e.g., \cite[Theorem 5.2.2]{GW09}. When $\dim X$ is odd, what we need is the  first fundamental theorem of skew-polynomial invariants of symplectic groups, which has been developed in \cite[\S 3, the totally real case]{Rib83} (see also \cite[Theorem 3.4]{Tho07}).
\end{proof}

\begin{remark}
The exceptional complete intersections are exactly the complete intersections of Hodge-Tate type. By the non-semisimplicity theorem of \cite{BaM04} 
(see also \cite{HMT09}), the big quantum cohomology of non-exceptional complete intersections are not generically semisimple. Theorem \ref{thm-monodromythm} gives another view to this result. Let $E$ be the Euler vector field for the genus zero generating function, and 
suppose $t^i$ and $t^j$ are as in Theorem \ref{thm-monodromythm} with $n+1\leq i,j\leq n+m$.
Then it is easy to check that 
\begin{eqnarray*}
\tilde{E}=E+t^j\frac{\partial}{\partial t^i}-t^i\frac{\partial}{\partial t^j}
\end{eqnarray*}
is another Euler vector field and commuting with $E$, simply because the dependence of $F$ on $t^i$ and $t^j$ is via the dependence on a single variable $s$. By the argument of \cite[Theorem 1.3]{HMT09} this implies the non-semisimplicity.

For the exceptional complete intersections, on the contrary, the monodromy groups are finite groups, so the rings of invariants have the same dimensions as the ranks of the primitive cohomology groups. For example, for the even dimensional complete intersection $X_n(2,2)$, let $t^1,\dots,t^{n+3}$ be a dual orthonormal basis, then the symmetric invariants of the Weyl group $W_{n+3}$ are generated by $\prod_{i=1}^{n+3}t^{i}$ and symmetric polynomials of $(t^1)^2,\dots,(t^{n+3})^2$. Therefore to show that the big quantum cohomology of $X_n(2,2)$ is semisimple (if it were), one has to compute the invariants to length \emph{at least $l$} such that there is no continuous family of automorphisms preserving the degree $2k$ forms
\[
\sum_{i=1}^{n+3}(t^i)^{2k}
\]
for $2\leq 2k\leq l$. In \cite{Hu21} we elaborate this approach to show that the quantum cohomology of $X_n(2,2)$ is generically semisimple.
\end{remark}


\section{The reduction of the WDVV equation by the monodromy group action}
\label{sec:red}
In this section we study the WDVV equation for the complete intersections. In the notations of Theorem \ref{thm-monodromythm}, the genus zero generating function  can be written as $F=F(t^0,\dots,t^n, s)$. The WDVV equation for the original genus 0 generating function 
can be written as an equivalent system of differential equations for  $F(t^{0},\dots,t^{n},s)$. 
We study the even and odd dimensional cases separately. For later use we define
\begin{eqnarray*}
F^{(k)}(t^0,\dots,t^{n})=\Big(\frac{\partial ^{k}}{\partial s^{k}}F\Big) \Big |_{s=0},
\end{eqnarray*}
so $F$ can be expanded as
\begin{eqnarray*}
F=F^{(0)}+sF^{(1)}+\frac{s^{2}}{2}F^{(2)}+\dots
\end{eqnarray*}
Note that $F^{(0)}$ is the  generating function of the \emph{ambient quantum cohomology} of $X$. By Lemma \ref{lem-subspace} and Corollary \ref{cor-mono3}, we know that $F^{(0)}$ satisfies the following WDVV equation, for $0\leq a,b,c,d\leq n$,
\begin{eqnarray}\label{eq-wdvv0}
\sum_{e=0}^{n}\sum_{f=0}^{n}\frac{\partial^3 F^{(0)}}{\partial t^a\partial t^b\partial t^e}g^{ef}\frac{\partial^3 F^{(0)}}{\partial t^f\partial t^c\partial t^d}
=\sum_{e=0}^{n}\sum_{f=0}^{n}\frac{\partial^3 F^{(0)}}{\partial t^a\partial t^c\partial t^e}g^{ef}\frac{\partial^3 F^{(0)}}{\partial t^f\partial t^b\partial t^d}.
\end{eqnarray} 

\subsection{Even dimensions}
\label{sec:red-even}
Suppose $n=\dim X$ is even.
Let $\gamma_{0}=1,\gamma_{1},\dots,\gamma_{n}$ be a basis of $H_{\mathrm{amb}}^{*}(X)$ , and $\gamma_{n+1},\dots, \gamma_{n+m}$
be an orthonormal basis of $H_{\mathrm{prim}}^{n}(X)$. The corresponding dual basis is denoted by $t^{0},\dots, t^{n+m}$.
Let
\begin{eqnarray*}
s=\sum_{\mu=n+1}^{n+m}\frac{(t^{\mu})^2}{2},
\end{eqnarray*}
then by Theorem \ref{thm-monodromythm} the genus 0 generating function   can be written as a power series of $t^0,\dots,t^n$ and $s$. We denote this series by $F(t^{0},\dots,t^{n},s)$. 
The WDVV equation for the original genus 0 generating function  is
\begin{eqnarray}\label{eq-wdvv-even}
\sum_{e=0}^{n+m}\sum_{f=0}^{n+m}\frac{\partial^3 F}{\partial t^a\partial t^b\partial t^e}g^{ef}\frac{\partial^3 F}{\partial t^f\partial t^c\partial t^d}
=\sum_{e=0}^{n+m}\sum_{f=0}^{n+m}\frac{\partial^3 F}{\partial t^a\partial t^c\partial t^e}g^{ef}\frac{\partial^3 F}{\partial t^f\partial t^b\partial t^d}
\end{eqnarray} 
for $0\leq a,b,c,d\leq n+m$. 
For $0\leq a,b,c\leq n$, we use $F_{abc}$ to denote $$\frac{\partial^3 F}{\partial t^a\partial t^b\partial t^c},$$
and for brevity use $F_{abs}$ or $F_{sab}$ to denote $$\frac{\partial^3 F}{\partial t^a\partial t^b\partial s}.$$
Similarly the meanings of $F_{ssa}$ and $F_{ss}$, $F_{sss}$ are obvious. We will use  the Einstein convention, where \emph{the summation indices run from $0$ to $n$}.\\

Some choices of the $4$-tuples $(a,b,c,d)$ in (\ref{eq-wdvv-even}) give trivial equations, e.g., $b=c=d$.
We consider all the possible choices of $4$-tuples $(a,b,c,d)$ that may give nontrivial equations. They are listed as follows.
\begin{enumerate}
\item[(1)] $0\leq a,b,c,d\leq n$. In this case the LHS of  (\ref{eq-wdvv-even}) is 
\begin{eqnarray}\label{eq-detail1}
&&F_{abe}g^{ef}F_{fcd}+\sum_{e=n+1}^{n+m}\sum_{f=n+1}^{n+m}\frac{\partial^3 F}{\partial t^a\partial t^b\partial t^e}g^{ef}\frac{\partial^3 F}{\partial t^f\partial t^c\partial t^d}\nn\\
&=&F_{abe}g^{ef}F_{fcd}+\sum_{e=n+1}^{n+m}\sum_{f=n+1}^{n+m}F_{abs}t_{e}g^{ef}t_{f}F_{scd}\nn\\
&=&F_{abe}g^{ef}F_{fcd}+2sF_{abs}F_{scd},
\end{eqnarray}
and the RHS of  (\ref{eq-wdvv-even}) is 
\begin{eqnarray*}
F_{ace}g^{ef}F_{fbd}+2sF_{acs}F_{sbd}.
\end{eqnarray*}
So we obtain
\begin{eqnarray}\label{eq-wdvv1}\nn\\
F_{abe}g^{ef}F_{fcd}+2sF_{abs}F_{scd}=F_{ace}g^{ef}F_{fbd}+2sF_{acs}F_{sbd}.
\end{eqnarray}
\item[(2)] $0\leq a,b,c\leq n$, $n+1\leq d\leq n+m$. In this case the  the LHS of  (\ref{eq-wdvv-even}) is 
\begin{eqnarray}\label{eq-detail2}
&&F_{abe}g^{ef}F_{fcd}+\sum_{e=n+1}^{n+m}\sum_{f=n+1}^{n+m}\frac{\partial^3 F}{\partial t^a\partial t^b\partial t^e}g^{ef}\frac{\partial^3 F}{\partial t^f\partial t^c\partial t^d}\nn\\
&=&F_{abe}g^{ef}F_{fcd}+\sum_{e=n+1}^{n+m}\sum_{f=n+1}^{n+m}F_{abs}t^{e}g^{ef}t^{f}t^{d}F_{ssc}
+\sum_{e=n+1}^{n+m}F_{abs}t^{e}g^{ed}F_{sc}\nn\\
&=&t^{d}F_{abe}g^{ef}F_{sfc}
+2t^{d}sF_{sab}F_{ssc}
+t^{d}F_{abs}F_{sc}.
\end{eqnarray}
Interchanging $b$ and $c$ we obtain the RHS of  (\ref{eq-wdvv-even}), so
\begin{eqnarray}\label{eq-wdvv2}
&&t^{d}F_{abe}g^{ef}F_{sfc}
+2t^{d}sF_{sab}F_{ssc}
+t^{d}F_{abs}F_{sc}\nn\\
&=&t^{d}F_{ace}g^{ef}F_{sfb}
+2t^{d}sF_{sac}F_{ssb}
+t^{d}F_{acs}F_{sb},
\end{eqnarray}

i.e.,
\begin{eqnarray}\label{eq-wdvv3}
&&F_{abe}g^{ef}F_{sfc}
+2sF_{sab}F_{ssc}
+F_{abs}F_{sc}\nn\\
&=&F_{ace}g^{ef}F_{sfb}
+2sF_{sac}F_{ssb}
+F_{acs}F_{sb}.
\end{eqnarray}
From now on we omit the details such as (\ref{eq-detail2}). 
\item[(3)] $0\leq a,b\leq n$, $n+1\leq c, d\leq n+m$, and $c\neq d$, we obtain
\begin{eqnarray}\label{eq-wdvv5}
&&F_{abe}g^{ef}F_{ssf}
+2sF_{sab}F_{sss}
+2F_{sab}F_{ss}\nn\\
&=&F_{sae}g^{ef}F_{sfb}
+2sF_{ssa}F_{ssb}+F_{sa}F_{ssb}
+F_{ssa}F_{sb}.
\end{eqnarray}

\item[(4)] $0\leq a,b\leq n$, $n+1\leq c= d\leq n+m$, we obtain
\begin{eqnarray*}
&&F_{abe}g^{ef}F_{sf}+2sF_{sab}F_{ss}
+(t^{c})^{2}(F_{abe}g^{ef}F_{ssf}
+2sF_{sab}F_{sss}+2F_{sab}F_{ss})\nn\\
&=&(t^{c})^{2}(F_{sae}g^{ef}F_{sfb}
+2sF_{ssa}F_{ssb}+F_{sa}F_{ssb}+F_{ssa}F_{sb})
+F_{sa}F_{sb}.
\end{eqnarray*}
When $m\geq 2$ this is equivalent to  (\ref{eq-wdvv5}) together with
\begin{eqnarray}\label{eq-wdvv6}
F_{abe}g^{ef}F_{sf}
+2sF_{sab}F_{ss}
=F_{sa}F_{sb}.
\end{eqnarray}

\item[(5)](When $m\geq 3$)  $0\leq a\leq n$, $n+1\leq b,c, d\leq n+m$, and $b\neq c$, $c\neq d$, $b\neq d$, we obtain
\begin{eqnarray*}
&&t^{b}t^{c}t^{d}(F_{sae}g^{ef}F_{ssf}
+2sF_{ssa}F_{sss}
+F_{sa}F_{sss}
+2F_{ssa}F_{ss})\nn\\
&=&t^{b}t^{c}t^{d}(F_{sae}g^{ef}F_{ssf}
+2sF_{ssa}F_{sss}
+F_{sa}F_{sss}
+2F_{ssa}F_{ss}),
\end{eqnarray*}
which is trivial.

\item[(6)] (When $m\geq 2$) $0\leq a\leq n$, $n+1\leq b,c\leq n+m$, $c=d$ and $b\neq c$, we obtain
\begin{eqnarray*}
&&t^{b}(t^{c})^{2}(F_{sae}g^{ef}F_{ssf}+2sF_{ssa}F_{sss}+2F_{ssa}F_{ss}+F_{sa}F_{sss})\nn\\
&&+t^{b}(F_{sae}g^{ef}F_{sf}
+2sF_{ssa}F_{ss}+F_{sa}F_{ss})\nn\\
&=&t^{b}(t^{c})^{2}(F_{sae}g^{ef}F_{ssf}+2sF_{ssa}F_{sss}
+2F_{ssa}F_{ss}+F_{sa}F_{sss})
+t^{b}F_{sa}F_{ss},
\end{eqnarray*}
which is equivalent to
\begin{eqnarray}\label{eq-wdvv8}
F_{sae}g^{ef}F_{sf}
+2sF_{ssa}F_{ss}=0.
\end{eqnarray}

\item[(7)] (When $m\geq 4$) $n+1\leq a,b,c, d\leq n+m$, and $a,b,c,d$ are pairwise distinct,
\begin{eqnarray*}
&&t^{a}t^{b}t^{c}t^{d}(F_{sse}g^{ef}F_{ssf}
+2sF_{sss}F_{sss}
+4F_{ss}F_{sss})\nn\\
&=&t^{a}t^{b}t^{c}t^{d}(F_{sse}g^{ef}F_{ssf}
+2sF_{sss}F_{sss}
+4F_{ss}F_{sss}),
\end{eqnarray*}
which is trivial.
\item[(8)] (When $m\geq 3$) $n+1\leq a,b,c\leq n+m$, $c=d$, and $a,b,c$ are pairwise distinct, we obtain
\begin{eqnarray*}
&&t^{a}t^{b}(t^{c})^2(F_{sse}g^{ef}F_{ssf}+2sF_{sss}F_{sss}+2F_{sss}F_{ss}
+F_{ss}F_{sss}+F_{ss}F_{sss})\nn\\
&&+t^{a}t^{b}(F_{sse}g^{ef}F_{sf}
+2sF_{sss}F_{ss}+F_{ss}F_{ss}+F_{ss}F_{ss})\nn\\
&=&t^{a}t^{b}(t^{c})^2(F_{sse}g^{ef}F_{ssf}
+2sF_{sss}F_{sss}+2F_{sss}F_{ss}+F_{ss}F_{sss}+F_{ss}F_{sss})\nn\\
&&+t^{a}t^{b}F_{ss}F_{ss},
\end{eqnarray*}
which is equivalent to
\begin{eqnarray}\label{eq-wdvv9}
F_{sse}g^{ef}F_{sf}+2sF_{sss}F_{ss}+F_{ss}F_{ss}
=0.
\end{eqnarray}

\item[(9)] (When $m\geq 2$) $n+1\leq a,b\leq n+m$, $a\neq b$, $c=a$, $d=b$, we obtain
\begin{eqnarray}
&&((t^{a})^{2}+(t^{b})^{2})\Big(F_{sse}g^{ef}F_{sf}
+2sF_{sss}F_{ss}
+F_{ss}F_{ss}\Big)\nn\\
&&+F_{se}g^{ef}F_{sf}
+2sF_{ss}F_{ss}
=0,
\end{eqnarray}
which when $m\geq 3$ is equivalent to (\ref{eq-wdvv9}) together with
\begin{eqnarray}\label{eq-wdvv10}
F_{se}g^{ef}F_{sf}
+2sF_{ss}F_{ss}=0.
\end{eqnarray}
\end{enumerate}
It is not hard to see that the other WDVV equations do not give  new equations.\\

By Corollary \ref{cor-estimate-rank-primCoh},
when $X$ is an even dimensional non-exceptional complete intersection of dimension $\geq 3$, we have $m=\mathrm{rk}\ H^{n}_{\mathrm{prim}}\big(X_{n}(\mathbf{d})\big)\geq 3$. So the WDVV equation (\ref{eq-wdvv-even}) is equivalent to the collection of 
(\ref{eq-wdvv1}), (\ref{eq-wdvv3}), (\ref{eq-wdvv5}), (\ref{eq-wdvv6}), (\ref{eq-wdvv8}), (\ref{eq-wdvv9}), (\ref{eq-wdvv10}). These equations are not independent.
Differentiating (\ref{eq-wdvv10}) by $t^{a}$, $0\leq a\leq n$,  we obtain  (\ref{eq-wdvv8}), and differentiating (\ref{eq-wdvv10}) by $s$ we obtain  (\ref{eq-wdvv9}).

Differentiating (\ref{eq-wdvv6}) by $s$ yields
\begin{eqnarray}\label{eq-wdvv21}
&&F_{sabe}g^{ef}F_{sf}+F_{abe}g^{ef}F_{ssf}+2F_{sab}F_{ss}+2sF_{ssab}F_{ss}+2sF_{sab}F_{sss}\nn\\
&=&F_{sa}F_{ssb}+F_{ssa}F_{sb}.
\end{eqnarray}
For $0\leq a,b\leq n$, differentiating (\ref{eq-wdvv10}) by $t^{a}$ and $t^{b}$ yields
\begin{eqnarray}\label{eq-wdvv22}
F_{sabe}g^{ef}F_{sf}+F_{sae}g^{ef}F_{sbf}+2sF_{ssab}F_{ss}+2sF_{ssa}F_{ssb}=0.
\end{eqnarray}
Substituting (\ref{eq-wdvv21}) and (\ref{eq-wdvv22}) into both sides of (\ref{eq-wdvv5}), we see that (\ref{eq-wdvv5}) is equivalent to
\begin{eqnarray*}
F_{sa}F_{ssb}+F_{ssa}F_{sb}-F_{sabe}g^{ef}F_{sf}-2sF_{ssab}F_{ss}\\
=-F_{sabe}g^{ef}F_{sf}-2sF_{ssab}F_{ss}+F_{sa}F_{ssb}+F_{ssa}F_{sb},
\end{eqnarray*}
which is trivial.\\

For $0\leq c\leq n$, differentiating (\ref{eq-wdvv6}) by $t^{c}$,
\begin{eqnarray}\label{eq-wdvv22.1}
&&F_{abce}g^{ef}F_{sf}+F_{abe}g^{ef}F_{sfc}+2sF_{sabc}F_{ss}+2sF_{sab}F_{ssc}\nn\\
&=&F_{sac}F_{sb}+F_{sa}F_{sbc}.
\end{eqnarray}
Similarly for $0\leq b\leq n$, we have
\begin{eqnarray}\label{eq-wdvv22.2}
&&F_{abce}g^{ef}F_{sf}+F_{ace}g^{ef}F_{sfb}+2sF_{sabc}F_{ss}+2sF_{sac}F_{ssb}\nn\\
&=&F_{sab}F_{sc}+F_{sa}F_{sbc}.
\end{eqnarray}
Then taking the difference of (\ref{eq-wdvv22.1}) and (\ref{eq-wdvv22.2}) we obtain (\ref{eq-wdvv3}). Finally, for $0\leq a,b,c,d\leq n$, computing
\begin{eqnarray*}
0=\frac{\partial^{2}}{\partial t^{c}\partial t^{d}}\Big(F_{abe}g^{ef}F_{sf}+2sF_{sab}F_{ss}-F_{sa}F_{sb}\Big)\\
+\frac{\partial^{2}}{\partial t^{a}\partial t^{b}}\Big(F_{cde}g^{ef}F_{sf}+2sF_{scd}F_{ss}-F_{sc}F_{sd}\Big)
\end{eqnarray*}
and
\begin{eqnarray*}
0=\frac{\partial^{2}}{\partial t^{b}\partial t^{d}}\Big(F_{ace}g^{ef}F_{sf}+2sF_{sac}F_{ss}-F_{sa}F_{sc}\Big)\\
+\frac{\partial^{2}}{\partial t^{a}\partial t^{c}}\Big(F_{bde}g^{ef}F_{sf}+2sF_{sbd}F_{ss}-F_{sb}F_{sd}\Big)
\end{eqnarray*}
then taking the difference, we obtain
\begin{eqnarray*}
\frac{\partial}{\partial s}\Big(F_{abe}g^{ef}F_{fcd}+2sF_{sab}F_{scd}\Big)
= \frac{\partial}{\partial s}\Big(F_{ace}g^{ef}F_{fbd}+2sF_{sac}F_{sbd}\Big).
\end{eqnarray*}
This means that, expanding both sides of  (\ref{eq-wdvv1}) as series of $s$, the part of  positive powers are equal.
On the other hand,  $0$-th part of (\ref{eq-wdvv1}) is (\ref{eq-wdvv0}), the WDVV equation for $F^{(0)}$. So (\ref{eq-wdvv1}) can be derived from (\ref{eq-wdvv6}) and the  WDVV equation for $F^{(0)}$. 

Summarizing, we obtain the following:
\begin{theorem}\label{thm-wdvveventhm}
For even dimensional non-exceptional complete intersections of $\dim$ at least 4, the collection of WDVV equations (\ref{eq-wdvv-even}) for $F$ is equivalent to the WDVV (\ref{eq-wdvv0}) for $F^{(0)}$ together with
\begin{eqnarray}\label{eq-wdvv23}
F_{abe}g^{ef}F_{sf}+2sF_{sab}F_{ss}=F_{sa}F_{sb}, & 0\leq a,b\leq n,
\end{eqnarray}
and
\begin{eqnarray}\label{eq-wdvv24}
F_{se}g^{ef}F_{sf}+2sF_{ss}F_{ss}=0.
\end{eqnarray}
\end{theorem}

\subsection{Odd dimensions}
\label{sec:red-odd}
Suppose $\dim X$ is odd.
Let $\gamma_{0}=1,\gamma_{1},\dots,\gamma_{n}$ be a basis of $H_{\mathrm{amb}}^{*}(X)$ , and $\gamma_{n+1},\dots, \gamma_{n+m}$
be a symplectic basis of $H_{\mathrm{prim}}^{*}(X)$, i.e,
\begin{eqnarray*}
g_{i,i+\frac{m}{2}}=1, & g_{i+\frac{m}{2},i}=-1, & n+1\leq i\leq n+\frac{m}{2},
\end{eqnarray*}
thus
\begin{eqnarray*}
g^{i,i+\frac{m}{2}}=-1, & g^{i+\frac{m}{2},i}=1, & n+1\leq i\leq n+\frac{m}{2}.
\end{eqnarray*}

 The corresponding dual basis is denoted by $t^{0},\dots, t^{n+m}$.
Let
\begin{eqnarray*}
s=-\sum_{i=n+1}^{n+\frac{m}{2}}t^i t^{i+\frac{m}{2}}.
\end{eqnarray*}
then by Theorem \ref{thm-monodromythm} the genus 0 generating function   can be written as a series $F(t^{0},\dots,t^{n},s)$. 
The WDVV equation for the original genus 0 generating function
\begin{eqnarray}\label{eq-wdvv-odd}
\sum_{e=0}^{n+m}\sum_{f=0}^{n+m}\frac{\partial^3 F}{\partial t^a\partial t^b\partial t^e}g^{ef}\frac{\partial^3 F}{\partial t^f\partial t^c\partial t^d}
=(-1)^{bc}\sum_{e=0}^{n+m}\sum_{f=0}^{n+m}\frac{\partial^3 F}{\partial t^a\partial t^c\partial t^e}g^{ef}\frac{\partial^3 F}{\partial t^f\partial t^b\partial t^d}.
\end{eqnarray} 
for $0\leq a,b,c,d\leq n+m$, where $(-1)^{bc}:=(-1)^{\deg_{\mathbb{R}} \gamma_b\cdot \deg_{\mathbb{R}} \gamma_c}$.  Here 
$\deg_{\mathbb{R}} \gamma_b$ means the real degree of the cohomology class $\gamma_b$ (see \S \ref{sec:nota-conv}).
As in the previous subsection we use the compact notations  $F_{abc}$, $F_{sab}$ etc., and  the Einstein convention summation is from $0$ to $n$.\\

The \emph{graded commutativity} makes the final differential equations and their derivation slightly different from those in even dimensions. We record several facts resulted from the graded commutativity.
\begin{enumerate}
      \item[(i)] Although the variable $s$ is of even degree, its $(\frac{m}{2}+1)$-th power is 0.  Thus a power series of $s$ is in fact a polynomial of $s$ with degree $\leq \frac{m}{2}$.
      \item[(ii)]  Suppose $g(s)$ is a polynomial in $s$, then for some $n+1\leq c\leq n+\frac{m}{2}$, 
\begin{eqnarray*}
t^c g(s)=0\quad 
(\mbox{or}\ t^ct^{c+\frac{m}{2}} g(s)=0)
\end{eqnarray*}
only implies that the coefficients of $s^k$ of $g(s)$ is 0 for $0\leq k\leq \frac{m}{2}-1$.
      \item[(iii)] Similarly, 
\begin{eqnarray*}
t^c t^{d} g(s)=0
\end{eqnarray*}
for some $n+1\leq c,d\leq n+m$ and $c\not\equiv d \mod \frac{m}{2}$ implies that the coefficients of $s^k$ of $g(s)$ is 0 for $0\leq k\leq \frac{m}{2}-2$. 
\end{enumerate}
These facts will be used in the following of this subsection implicitly. Moreover, we use the notation
$$f(s)= g(s)\mod s^k $$ to denote that the coefficients of $s^{i}$ in $f$ and $g$ are equal for $0\leq i\leq k-1$.\\

In odd dimensions, there are more choices of  4-tuples $(a,b,c,d)$ in (\ref{eq-wdvv-odd}) than the even dimensions. 
We list all the possible choices of 4-tuples $(a,b,c,d)$ that may lead to nontrivial equations. For each choice we omit the details which are similar to the even dimensions, and give only the final equations.
\begin{enumerate}
\item[(1)] $0\leq a,b,c,d\leq n$.
\begin{eqnarray}\label{eq-wdvv1odd}
F_{abe}g^{ef}F_{fcd}+2sF_{abs}F_{scd}=F_{ace}g^{ef}F_{fbd}+2sF_{acs}F_{sbd}.
\end{eqnarray}
Note that $s^{\frac{m}{2}+1}=0$. So if we regard $F$ as a series of $s$, (\ref{eq-wdvv1odd}) means an equality mod $s^{\frac{m}{2}+1}$.

\item[(2)] $0\leq a,b,c\leq n$, $n+1\leq d\leq n+\frac{m}{2}$. We get
\begin{eqnarray}\label{eq-wdvv3odd-0}
&&t^{d+\frac{m}{2}}\big(F_{abe}g^{ef}F_{sfc}
+2sF_{sab}F_{ssc}
+F_{abs}F_{sc}\big)\nn\\
&=&t^{d+\frac{m}{2}}\big(F_{ace}g^{ef}F_{sfb}
+2sF_{sac}F_{ssb}
+F_{acs}F_{sb}\big).
\end{eqnarray}
Since $t^i$ are skew-commutative for $n+1\leq i\leq n+m$, (\ref{eq-wdvv3odd-0}) is equivalent to
\begin{eqnarray}\label{eq-wdvv3odd}
&&F_{abe}g^{ef}F_{sfc}
+2sF_{sab}F_{ssc}
+F_{abs}F_{sc}\nn\\
&=&F_{ace}g^{ef}F_{sfb}
+2sF_{sac}F_{ssb}
+F_{acs}F_{sb} \mod s^{\frac{m}{2}}.
\end{eqnarray}
\item[(3)] (When $m\geq 4$) $0\leq a,b\leq n$, $n+1\leq c,d\leq n+\frac{m}{2}$, and $c\neq d$, we get
\begin{eqnarray}
&&t^{c+\frac{m}{2}}t^{d+\frac{m}{2}}\big(F_{abe}g^{ef}F_{ssf}
+2sF_{sab}F_{sss}
+2F_{sab}F_{ss}\big)\nn\\
&=&t^{c+\frac{m}{2}}t^{d+\frac{m}{2}}\big(F_{sae}g^{ef}F_{sfb}
+2sF_{ssa}F_{ssb}+F_{sa}F_{ssb}
+F_{ssa}F_{sb}\big).
\end{eqnarray}
which is equivalent to 
\begin{eqnarray}\label{eq-wdvv4odd}
&&F_{abe}g^{ef}F_{ssf}
+2sF_{sab}F_{sss}
+2F_{sab}F_{ss}\nn\\
&=&F_{sae}g^{ef}F_{sfb}
+2sF_{ssa}F_{ssb}+F_{sa}F_{ssb}
+F_{ssa}F_{sb} \mod s^{\frac{m}{2}-1}.
\end{eqnarray}
Such reasoning will be omitted in the following.
\item[(4)]  (When $m\geq 4$)  $0\leq a,b\leq n$, $n+1\leq c\leq n+\frac{m}{2}$, $n+\frac{m}{2}+1\leq d\leq n+m$, $d\neq c+\frac{m}{2}$.
\begin{eqnarray}\label{eq-wdvv5odd}
&&F_{abe}g^{ef}F_{ssf}
+2sF_{sab}F_{sss}
+2F_{sab}F_{ss}\nn\\
&=&F_{sae}g^{ef}F_{sfb}
+2sF_{ssa}F_{ssb}+F_{sa}F_{ssb}
+F_{ssa}F_{sb} \mod s^{\frac{m}{2}-1}.
\end{eqnarray}

\item[(5)] $0\leq a,b\leq n$, $n+1\leq c\leq n+\frac{m}{2}$, $d= c+\frac{m}{2}$.

\begin{eqnarray}\label{eq-wdvv6odd}
&&F_{abe}g^{ef}F_{sf}
+2sF_{sab}F_{ss}+t^{c}t^{c+\frac{m}{2}}\big(F_{abe}g^{ef}F_{ssf}
+2sF_{sab}F_{sss}+2F_{sab}F_{ss}\big)\nn\\
&=&F_{sa}F_{sb}+t^{c}t^{c+\frac{m}{2}}\big( F_{ase}g^{ef}F_{fbs}
+2sF_{ssa}F_{ssb}
+F_{sa}F_{ssb}
+F_{ssa}F_{sb}\big).
\end{eqnarray}

\item[(6)]  (When $m\geq 4$)  $0\leq a\leq n$, $n+1\leq b,c\leq n+\frac{m}{2}$, and $b\neq c$, $d=c+\frac{m}{2}$. 
\begin{eqnarray}\label{eq-wdvv8odd}
F_{sae}g^{ef}F_{sf}
+2sF_{ssa}F_{ss}=0 \mod s^\frac{m}{2}.
\end{eqnarray}

\item[(7)] $0\leq a\leq n$, $n+1\leq b=c\leq n+\frac{m}{2}$, $d=c+\frac{m}{2}$,
\begin{eqnarray}
F_{aes}g^{ef}F_{fs}+2sF_{ass}F_{ss}=0 \mod s^{\frac{m}{2}}.
\end{eqnarray}

\item[(8)]  (When $m\geq 6$)  $n+1\leq a,b,c\leq n+\frac{m}{2}$, $d=c+\frac{m}{2}$, and $a,b,c$ are pairwise distinct,
\begin{eqnarray}\label{eq-wdvv9odd}
F_{sse}g^{ef}F_{sf}+2sF_{sss}F_{ss}+F_{ss}F_{ss}
=0 \mod s^{\frac{m}{2}-1}.
\end{eqnarray}

\item[(9)]  (When $m\geq 4$)  $n+1\leq a, b\leq n+\frac{m}{2}$, $a\neq b$, $c=a+\frac{m}{2}$, $d=b+\frac{m}{2}$,
\begin{eqnarray}\label{eq-wdvv10odd}
&&(t^{a}t^{a+\frac{m}{2}}+t^{b}t^{b+\frac{m}{2}})\Big(F_{sse}g^{ef}F_{sf}
+2sF_{sss}F_{ss}
+F_{ss}F_{ss}\Big)\nn\\
&&+F_{se}g^{ef}F_{sf}
+2sF_{ss}F_{ss}
=0.
\end{eqnarray}

\item[(10)] $n+1\leq a\leq n+\frac{m}{2}$, $b=a$, $c=d=a+\frac{m}{2}$,
\begin{eqnarray}\label{eq-wdvv10.5odd}
F_{es}g^{ef}F_{fs}+2sF_{ss}F_{ss}+2t^{a}t^{a+\frac{m}{2}}(F_{es}g^{ef}F_{fss}
+F_{ss}F_{ss}+2sF_{ss}F_{sss})=0.
\end{eqnarray}
Note that (\ref{eq-wdvv10.5odd}) implies (\ref{eq-wdvv10odd}).
\end{enumerate}

For the other choices of $a,b,c,d$ the resulted equations are trivial.

\begin{theorem}\label{thm-wdvvoddthm}
For odd dimensional non-exceptional complete intersections of $\dim$ at least 3, the collection of  WDVV equations (\ref{eq-wdvv-odd})  for $F$ is equivalent to WDVV (\ref{eq-wdvv0}) for $F^{(0)}$ together with
\begin{eqnarray}\label{eq-wdvv23odd}
F_{abe}g^{ef}F_{sf}+2sF_{sab}F_{ss}\equiv F_{sa}F_{sb} \mod s^{\frac{m}{2}}, & 0\leq a,b\leq n,
\end{eqnarray}
and
\begin{eqnarray}\label{eq-wdvv24odd}
F_{se}g^{ef}F_{sf}+2sF_{ss}F_{ss}\equiv 0 \mod s^{\frac{m}{2}}.
\end{eqnarray}
\end{theorem}
\begin{proof}
By Corollary \ref{cor-estimate-rank-primCoh}, for odd dimensional complete intersections of dimension $\geq 3$ we have $m=\mathrm{rk}\ H_{\mathrm{prim}}^{n}(X)\geq 4$. 
First we suppose $m\geq 6$.
Taking $s=0$ in (\ref{eq-wdvv1odd}) yields (\ref{eq-wdvv0}). 
(\ref{eq-wdvv23odd}) follows from (\ref{eq-wdvv5odd}) and (\ref{eq-wdvv6odd}). (\ref{eq-wdvv24odd}) follows from (\ref{eq-wdvv10.5odd}) and (\ref{eq-wdvv9odd}). 

When $m=4$ the equation (\ref{eq-wdvv9odd}) is missing, so we need some ad hoc reasoning for (\ref{eq-wdvv24odd}). Suppose
\begin{eqnarray*}
F_{sse}g^{ef}F_{sf}+2sF_{sss}F_{ss}+F_{ss}F_{ss}
=x_0+x_1s+x_2 s^2,
\end{eqnarray*}
then (\ref{eq-wdvv10odd}) implies
\begin{eqnarray*}
F_{se}g^{ef}F_{sf}
+2sF_{ss}F_{ss}
=x_0 s+x_1 s^2.
\end{eqnarray*}
But 
\begin{eqnarray*}
2(F_{sse}g^{ef}F_{sf}+2sF_{sss}F_{ss}+F_{ss}F_{ss})
=\frac{\partial }{\partial s}\big(F_{se}g^{ef}F_{sf}
+2sF_{ss}F_{ss}\big),
\end{eqnarray*}
so 
\begin{eqnarray*}
2(x_0+x_1 s+2x_2 s^2)=x_0+2x_1s,
\end{eqnarray*}
thus
\begin{eqnarray*}
x_0=x_2=0.
\end{eqnarray*}
Therefore in this case (\ref{eq-wdvv9odd}) and (\ref{eq-wdvv24odd}) are still valid.

Conversely,  the equations (\ref{eq-wdvv1odd})-(\ref{eq-wdvv10.5odd}) can be deduced from (\ref{eq-wdvv-even}), (\ref{eq-wdvv23odd}) and (\ref{eq-wdvv24odd}) almost verbatim as in Section \ref{sec:red-even}; one needs only additionally record the power of $s$ that is mod out in each equation and their derivatives.
\end{proof}


\subsection{Incorporating the Euler vector field}
If $\gamma_0,\dots,\gamma_n$ are taken to be $1,\sfh,\dots,\sfh_n$, then the Euler vector field (\ref{eq-EulerVectorField}) can be written (abusively, as our convention after \eqref{eq-effective-EulerVectorField}) as
\begin{eqnarray}\label{eq-EV-symmetricReducted}
E=\sum_{i=0}^{n}(1-i)t^{i}\frac{\partial}{\partial t^i}+(2-n)s\frac{\partial }{\partial s}+\mathsf{a}(n,\mathbf{d})\frac{\partial}{\partial t^1}.
\end{eqnarray}
We obtain, for even $n$,
\begin{subequations}\label{eq-system1-even}
\begin{align}[left ={\empheqlbrace}]
&F_{abe}g^{ef}F_{sf}+2sF_{sab}F_{ss}=F_{sa}F_{sb},\quad \mbox{for}\ 0\leq a,b\leq n,\label{eq-system1-even-1}\\
&F_{se}g^{ef}F_{sf}+2sF_{ss}F_{ss}=0,\label{eq-system1-even-2}\\
&EF=(3-n)F+\mathsf{a}(n,\mathbf{d})\frac{\partial}{\partial t^1}c,\label{eq-system1-even-3}
\end{align}
\end{subequations}
and for odd $n$,
\begin{subequations}\label{eq-system1-odd}
\begin{align}[left ={\empheqlbrace}]
&F_{abe}g^{ef}F_{sf}+2sF_{sab}F_{ss}=F_{sa}F_{sb} \mod s^{\frac{m}{2}},\quad\mbox{for}\ 0\leq a,b\leq n,\label{eq-system1-odd-1}\\
&F_{se}g^{ef}F_{sf}+2sF_{ss}F_{ss}=0 \mod s^{\frac{m}{2}}, \label{eq-system1-odd-2}\\
&EF=(3-n)F+\mathsf{a}(n,\mathbf{d})\frac{\partial}{\partial t^1}c, \label{eq-system1-odd-3}
\end{align}
\end{subequations}
where $c$ is defined in (\ref{eq-tripleIntersectionForm}). Our aim in this paper is to study to what extent the system (\ref{eq-system1-even}) or (\ref{eq-system1-odd}) with $F|_{s=0}=F^{(0)}$ as initial given data, can be solved,  possibly with the help of some properties listed in Section \ref{sec:GW-invariants} or some other geometric reason.

\subsection{ \texorpdfstring{$F^{(l)}(0)$}{F(l)(0)} as ratios}\label{sec:meaningOfF(l)(0)}
It will turn out to be crucial to compute the constant leading term $F^{(l)}(0)$ of $F^{(l)}$. The following lemma gives a basis-independent description of $F^{(l)}(0)$. 
Let $A_{2l}$ be the set 
\begin{multline}\label{eq-def-A2l}
A_{2l}=\big\{\big((i_1,j_1),(i_2,j_2),\dots,(i_l,j_l)\big)|\{i_1,j_1,i_2,j_2,\dots,i_l,j_l\}=\{1,\dots,2l\},\\
i_k<j_k\ \mbox{for } 1\leq k\leq l, i_1<i_2<\dots<i_l\big\}.
\end{multline}
In other words, the elements of $A_{2l}$ parametrize the unordered pairings in a set of cardinality $2l$. For example, the elements of $A_{4}$ can be depicted as
\begin{equation*}
	\wick{\c1 1 \c2 2 \c1 3 \c2 4}\quad
	\wick{\c1 1 \c1 2 \c2 3 \c2 4}\quad
	\wick{\c1 1 \c2 2 \c2 3 \c1 4}\quad.
\end{equation*}

For $\sigma=\big((i_1,j_1),(i_2,j_2),\dots,(i_l,j_l)\big)\in A_{2l}$, and $G=(g_{i,j})_{1\leq i,j\leq 2l}$ a $2l\times 2l$ symmetric matrix (resp. a $2l\times 2l$ skew-symmetric matrix), we define
\begin{eqnarray}
  \mathrm{P}_{\sigma}(G)= \prod_{k=1}^{l}g_{i_k,j_k}.\\
  \big(\mbox{resp. }  \mathrm{Pf}_{\sigma}(G)=\mathrm{sgn}(\sigma) \prod_{k=1}^{l}g_{i_k,j_k}.\big)
\end{eqnarray}
Then define
\begin{eqnarray}
\mathrm{P}(G)=\sum_{\sigma\in A_{2l}}\mathrm{P}_{\sigma}(G).\label{eq-pfaffian-def-even}\\
\big(\mbox{resp. } \mathrm{Pf}(G)=\sum_{\sigma\in A_{2l}}\mathrm{Pf}_{\sigma}(G).\big)\label{eq-pfaffian-def-odd}
\end{eqnarray}
For skew-symmetric $G$, $\mathrm{Pf}(G)$ is no other than the Pfaffian of $G$. For symmetric $G$, we call $\mathrm{P}(G)$ the \emph{permanent Pfaffian} of $G$.

Now for  $\alpha_1,\dots,\alpha_{2l}\in H^*_{\mathrm{prim}}(X)$, we define $G(\alpha_1,\dots,\alpha_{2l})$ to be the matrix $G=(g_{i,j})_{1\leq i,j\leq 2l}$ with $g_{i,j}=(\alpha_i,\alpha_j)$. Thus $G$ is symmetric when $n$ is even,  and is skew-symmetric when $n$ is odd.
\begin{proposition}\label{prop-initialValues-meaning}\label{prop-F(l)-ratio}
Let $\gamma_0,\dots,\gamma_n$ be a basis of $H^{*}_{\mathrm{amb}}(X)$, and $t^{0},\dots,t^{n}$ the dual basis.  Then for any $k\geq 0$ and $i_1,\dots,i_k\in [0,n]$ and $\alpha_1,\dots,\alpha_{2l}\in H^*_{\mathrm{prim}}(X)$,
\begin{enumerate}
      \item[(i)] When $n$ is even,
            \begin{equation}\label{eq-initialValues-meaning-even}
                  \langle \gamma_{i_1},\dots,\gamma_{i_k}, \alpha_1,\dots,\alpha_{2l}\rangle_{0,k+2l}=\partial_{t^{i_1}}\circ\dots\circ\partial_{t^{i_k}}F^{(l)}(0)\cdot \mathrm{P}\big(G(\alpha_1,\dots,\alpha_{2l})\big);
            \end{equation}
      \item[(ii)]      When $n$ is odd,
                  \begin{equation}\label{eq-initialValues-meaning-odd}
                        \langle \gamma_{i_1},\dots,\gamma_{i_k}, \alpha_1,\dots,\alpha_{2l}\rangle_{0,k+2l}=\partial_{t^{i_1}}\circ\dots\circ\partial_{t^{i_k}}F^{(l)}(0)\cdot \mathrm{Pf}\big(G(\alpha_1,\dots,\alpha_{2l})\big).
                  \end{equation}
\end{enumerate}
\end{proposition}
\begin{proof}
For brevity we assume $k=0$; the argument carries over verbatim to $k\geq 0$. 
First we consider the even dimensional case. As at the beginning of Section \ref{sec:red-even}, let $\gamma_{n+1},\dots,\gamma_{n+m}$ be an orthonormal basis of $H^{*}_{\mathrm{prim}}(X)$, with dual coordinates $t^{n+1},\dots,t^{n+m}$. Since both sides of (\ref{eq-initialValues-meaning-even}) are invariant under permutations, we can assume that, without loss of generality, each $\alpha_i$ equals some  $\gamma_j$, where $n+1\leq j\leq n+m$. If some $\gamma_j$ appears in $\alpha_1,\dots, \alpha_{2l}$ exactly odd times, then both sides of  (\ref{eq-initialValues-meaning-even}) vanish. Assume that $\gamma_j$ appears exactly $2l_j$ times, for $n+1\leq j\leq n+m$. Then
\begin{equation}\label{eq-lem-initialValues-meaning-1}
	\langle \alpha_1,\dots,\alpha_{2l}\rangle=\frac{\partial^{2l}F}{(\partial t^{n+1})^{2l_{n+1}}\dots (\partial t^{n+m})^{2l_{n+m}}}(0),
\end{equation}
and
\begin{equation}\label{eq-lem-initialValues-meaning-2}
	\mathrm{P}\big(G(\alpha_1,\dots,\alpha_{2l})\big)=\prod_{j=n+1}^{n+m}|A_{2l_j}|,
\end{equation}
where $|S|$ denote the cardinality of a set $S$. Since 
\[
s=\frac{1}{2}\sum_{j=n+1}^{n+m}(t^j)^2,
\]
by the description of $A_{2l}$ as the  set of unordered pairings as above, we have
\begin{equation*}
	\frac{\partial^{2l} F}{(\partial t^j)^{2l}}(0)=|A_{2l}| \frac{\partial^l F}{\partial s^l}(0),
\end{equation*}
and more generally
\begin{equation}\label{eq-lem-initialValues-meaning-3}
	\frac{\partial^{2l}F}{(\partial t^{n+1})^{2l_{n+1}}\dots (\partial t^{n+m})^{2l_{n+m}}}(0)=\prod_{j=n+1}^{n+m}|A_{2l_{j}}|\cdot \frac{\partial^l F}{\partial s^l}(0).
	\end{equation}
So (\ref{eq-initialValues-meaning-even}) follows from 	(\ref{eq-lem-initialValues-meaning-1})-(\ref{eq-lem-initialValues-meaning-3}).

In the odd dimension case we take a symplectic basis as at the beginning of Section \ref{sec:red-odd}. Then (\ref{eq-initialValues-meaning-odd}) is obvious because in a nonzero invariant every pair $\gamma_j,\gamma_{j+\frac{m}{2}}$ appears at most one time.
\end{proof}

\begin{corollary}
$F^{(l)}(0)\in \mathbb{Q}$.
\end{corollary}
\begin{proof}
When $n$ is even (resp. $n$ is odd), by choosing $\alpha_i$'s in an orthogonal basis (resp. a symplectic basis) of $H^n_{\mathrm{prim}}(X)$, there exists $\alpha_1,\dots,\alpha_{2l}$ such that 
\begin{equation}\label{eq-nonzero-pfaffian}
    \mathrm{P}\big(G(\alpha_1,\dots,\alpha_{2l})\big)\neq 0.\ (\mbox{resp. } \mathrm{Pf}\big(G(\alpha_1,\dots,\alpha_{2l})\big)\neq 0   \mbox{ for } 1\leq l\leq \frac{m}{2}.)
\end{equation}
 So there exists also $\alpha_1,\dots,\alpha_{2l}\in H^n_{\mathrm{prim}}(X)\cap H^n(X;\mathbb{Q})$ such that (\ref{eq-nonzero-pfaffian}) holds. Since for such $\alpha_i$'s $\langle \alpha_1,\dots,\alpha_{2l}\rangle\in \mathbb{Q}$, the conclusion follows from (\ref{eq-initialValues-meaning-even}) (resp. \ref{eq-initialValues-meaning-odd}).
\end{proof}

Numerical experiments suggest that $F^{(l)}(0)\in \mathbb{Z}$. For a discussion see Remark \ref{rem:integrality}.

\subsection{An extension problem of Frobenius manifolds}
In this subsection we try to give an interpretation for the reconstruction problem from $F^{(0)}$ to $F$. The result will not be used in the remaining sections, while we think that the discussions in the subsection will be helpful to understand the nature of the reconstruction. For brevity we restrict to the even dimensional cases. 
The idea is to interpret (\ref{eq-wdvv23}) and (\ref{eq-wdvv24}) as some classical PDEs. If we make the change of variables
\begin{eqnarray}\label{eq-changeofvariable1}
 s=\frac{r^2}{2},
 \end{eqnarray}
 then (\ref{eq-wdvv24}) is transformed into
 \begin{eqnarray}\label{eq-changeofvariable2}
\frac{\partial F_s}{\partial t^e}g^{ef}\frac{\partial F_s}{\partial t^f}+\frac{\partial F_s}{\partial r}\frac{\partial F_s}{\partial r}=0.
 \end{eqnarray}
 If we define the pairing
 \begin{eqnarray}\label{eq-changeofvariable3}
g(\partial_r,\partial_r)=1, & g(\partial_r,\partial_{t^a})=0, & 0\leq a\leq n,
 \end{eqnarray}
and view $F_s$ as a series in $t^0,\dots,t^n, r$,
 then (\ref{eq-changeofvariable2}) is an \emph{eikonal equation} for $F_s$. We will see in Theorem \ref{thm-reconstruction-II-F(1)} that $F^{(1)}=F_s|_{s=0}=F_s|_{r=0}$ can be regarded as also input data, so we need only to solve  the Cauchy problem of the eikonal equation
\begin{eqnarray}\label{eq-changeofvariable4}
\left\{
\begin{array}{c}
\frac{\partial F_s}{\partial t^e}g^{ef}\frac{\partial F_s}{\partial t^f}+\frac{\partial F_s}{\partial r}\frac{\partial F_s}{\partial r}=0,\\
F_s|_{r=0}=\mathrm{Known\ data}.
\end{array}
\right.
 \end{eqnarray} 
 However as $\frac{\partial F_s}{\partial r}|_{r=0}=0$, (\ref{eq-changeofvariable4}) is in fact a \emph{characteristic Cauchy problem}, for which
  the solution is in general not unique. Now we interpret the system (\ref{eq-system1-even}) as an extension problem of Frobenius manifolds. For this purpose we need to introduce the notion of \emph{Frobenius submanifolds}. 
  
\begin{definition}
Let $(M,g,\bqp,e,E)$ be a Frobenius manifold, $N$ a submanifold. $N$ is called a \emph{Frobenius submanifold} if the following are satisfied:
\begin{itemize}
\item[(i)] $g|_N$ is nondegenerate and flat;
\item[(ii)] $e|_{N}, E|_{N}\in TN$;
\item[(iii)] $TN\bqp TN \subset TN$;
\item[(iv)] $(N,g|_N,\bqp |_{TN},e|_N,E|_N)$ is a Frobenius manifold, and an affine flat structure on $N$ can be locally extended to be an affine flat structure on $M$.
\end{itemize}
\end{definition}  

\begin{remark}
The Frobenius submanifold in our definition is called a \emph{natural Frobenius submanifold} in \cite{Str01}.
\end{remark}
  
As we have seen in the introduction, there is a Frobenius manifold $\mathcal{M}_{\mathrm{amb}}$ corresponding to the ambient quantum cohomology, with generating function  $F^{(0)}$, flat coordinates $t^0, \dots,t^n$, and Euler vector field
\begin{eqnarray}\label{eq-Euleramb}
E_{\mathrm{amb}}=\sum_{i=0}^{n}(1-i)t^{i}\frac{\partial}{\partial t^i}
+\mathsf{a}(n,\mathbf{d})\frac{\partial}{\partial t^1}.
\end{eqnarray}
We introduce a new variable $r$ as \eqref{eq-changeofvariable1} and extend the pairing as (\ref{eq-changeofvariable3}), and define the generating function  to be
 $F$ for the whole quantum cohomology with the change of variables (\ref{eq-changeofvariable1}).
\begin{lemma}
 The system of  WDVV equations for $F$ in the flat coordinates $t^0,\dots,t^n,r$ is equivalent to
\begin{numcases}{}
F_{abe}g^{ef}F_{fcd}+2sF_{abs}F_{scd}=F_{ace}g^{ef}F_{fbd}+2sF_{acs}F_{sbd}, \label{eq-wdvvred1}\\
F_{abe}g^{ef}F_{sfc}
+2sF_{sab}F_{ssc}
+F_{abs}F_{sc}
=F_{ace}g^{ef}F_{sfb}
+2sF_{sac}F_{ssb}
+F_{acs}F_{sb}, \label{eq-wdvvred2}\\
F_{abe}g^{ef}F_{sf}
+2sF_{abe}g^{ef}F_{ssf}
+4s^2F_{sab}F_{sss}+6sF_{sab}F_{ss}\nn\\
=2s F_{sae}g^{ef}F_{sfb}
+4s^2 F_{ssa}F_{ssb}
+2s F_{sa}F_{ssb}
+2sF_{ssa}F_{sb}
+F_{sa}F_{sb},\label{eq-wdvvred3}
\end{numcases}
where the indices $a,b$ runs over $0,\dots,n$, and the Einstein summation convention runs also over $0,\dots,n$.
\end{lemma}
\begin{proof}
The pairing $g(\partial_r,\partial_r)$ is equivalent to $g(\partial_s,\partial_s)=\frac{1}{r^2}=\frac{1}{2s}$. Since
\[
F_{rr}=\frac{\partial}{\partial r}\big(r\frac{\partial F}{\partial s}\big)=F_s+2s F_{ss},
\]
the WDVV equation
\[
F_{abe}g^{ef}F_{fcr}+F_{abr}F_{rcr}=F_{ace}g^{ef}F_{fbr}+F_{acr}F_{rbr}
\]
is equivalent to (\ref{eq-wdvvred2}). The other equations are derived in a similar way.
\end{proof}

We define the degree of $\partial_r$ to be $\frac{n}{2}$, then the Euler field is
\begin{eqnarray}\label{eq-Eulerred}
E=\sum_{i=0}^{n}(1-i)t^{i}\frac{\partial}{\partial t^i}+(1-\frac{n}{2})r\frac{\partial }{\partial r}
+\mathsf{a}(n,\mathbf{d})\frac{\partial}{\partial t^1}.
\end{eqnarray}
It is easily seen that the system (\ref{eq-system1-even}) implies (\ref{eq-wdvvred1})-(\ref{eq-Eulerred}); in fact (\ref{eq-wdvvred1}) is (\ref{eq-wdvv1}), (\ref{eq-wdvvred2}) is (\ref{eq-wdvv3}), and for (\ref{eq-wdvvred3}) see the proof of the following Lemma \ref{lem-equiv-WDVV-submanifolds}. So we have constructed a Frobenius manifold 
$\mathcal{M}$, which has $\mathcal{M}_{\mathrm{amb}}$ as a codimension 1 Frobenius submanifold.
 Conversely, we have:

\begin{lemma}\label{lem-equiv-WDVV-submanifolds}
The system of equations (\ref{eq-wdvvred1})-(\ref{eq-Eulerred}) together with (\ref{eq-changeofvariable2}) is equivalent to the system (\ref{eq-system1-even}) together with the WDVV equations for $F|_{s=0}$.
\end{lemma}
\begin{proof} It suffices to show the direction ``$\Longrightarrow$". The equations (\ref{eq-wdvvred1}) restricted to $s=0$ imply the WDVV for $F|_{s=0}$.
We write (\ref{eq-wdvvred2}) as
\begin{eqnarray*}
\frac{\partial}{\partial t^c}\big(F_{abe}g^{ef}F_{sf}+2sF_{sab}F_{ss}-F_{sa}F_{sb}\big)
=\frac{\partial}{\partial t^b}\big(F_{ace}g^{ef}F_{sf}+2sF_{sac}F_{ss}-F_{sa}F_{sc}\big).
\end{eqnarray*}
Thus there exists $\Phi_a(\mathbf{t},s)$ for $a=0,\dots,n$ such that
\begin{eqnarray}\label{eq-wdvvred4-0}
F_{abe}g^{ef}F_{sf}+2sF_{sab}F_{ss}-F_{sa}F_{sb}=\frac{\partial\Phi_a}{\partial t^b},\ \mbox{for } 0\leq b\leq n.
\end{eqnarray}
Since the LHS of (\ref{eq-wdvvred4-0}) is symmetric in $a$ and $b$, there exits $\Phi(\mathbf{t},s)$ such that
\begin{equation*}
	\Phi_a=\frac{\partial\Phi}{\partial t^a},
\end{equation*}
and thus
\begin{eqnarray}\label{eq-wdvvred4}
F_{abe}g^{ef}F_{sf}+2sF_{sab}F_{ss}-F_{sa}F_{sb}=\frac{\partial^2 \Phi}{\partial t^a\partial t^b}.
\end{eqnarray}
We write (\ref{eq-wdvvred3}) as
\begin{eqnarray*}
&&F_{abe}g^{ef}F_{sf}+2sF_{sab}F_{ss}-F_{sa}F_{sb}\\
&=&2s(F_{sae}g^{ef}F_{sfb}+2s F_{ssa}F_{ssb}+F_{sa}F_{ssb}+F_{ssa}F_{sb}\\
&&-F_{abe}g^{ef}F_{ssf}-2sF_{sab}F_{sss}-2F_{sab}F_{ss}),
\end{eqnarray*}
i.e.,
\begin{eqnarray*}
&&F_{abe}g^{ef}F_{sf}+2sF_{sab}F_{ss}-F_{sa}F_{sb}\\
&=&s\frac{\partial^2}{\partial t^a\partial t^b}(F_{se}g^{ef}F_{sf}+2s F_{ss}F_{ss})
-2s\frac{\partial}{\partial s}\big(F_{abe}g^{ef}F_{sf}+2sF_{sab}F_{ss}-F_{sa}F_{sb}\big).
\end{eqnarray*}
Thus there exist $\phi(s)$ and $\phi_{a}(s)$ such that
\begin{eqnarray*}
s(F_{se}g^{ef}F_{sf}+2s F_{ss}F_{ss})=\Phi+2s\frac{\partial \Phi}{\partial s}+\sum_{a}\phi_a(s) t^a+\phi(s).
\end{eqnarray*}
Solving $f_{a}(s)+2sf'_{a}(s)=\phi_{a}(s)$ and $f(s)+2sf'(s)=\phi(s)$, modifying $\Phi$ by $\Phi+\sum_{a}f_a(s) t^a+f(s)$, we have
\begin{eqnarray}\label{eq-wdvvred5}
s(F_{se}g^{ef}F_{sf}+2s F_{ss}F_{ss})=\Phi+2s\frac{\partial \Phi}{\partial s},
\end{eqnarray}
and (\ref{eq-wdvvred4}) still holds. By (\ref{eq-changeofvariable2}), the LHS of (\ref{eq-wdvv24}) equals 0, so 
\[
\Phi+2s\frac{\partial \Phi}{\partial s}=0.
\]
Hence $\Phi=0$, and we obtain  (\ref{eq-wdvv23}). 
\end{proof}

We call  (\ref{eq-changeofvariable2}) the  \emph{eikonal equation in the normal direction} for $F$. One can of course write it in the coordinates
$t^0, \dots,t^n,r$. In summary, we can state the problem of finding a solution for (\ref{eq-system1-even}) with the initial condition $F|_{s=0}=F^{(0)}$ as follows.
\begin{proposition}
Solving the system (\ref{eq-system1-even}) with the initial condition at $s=0$ is equivalent to finding a codimension 1  embedding of $\mathcal{M}_{amb}$ into a Frobenius manifold with the Euler vector field (\ref{eq-Eulerred}) satisfying the eikonal equation in the normal direction.
\end{proposition}

\section{Structure of Frobenius algebras and reconstruction theorems}\label{sec:structureOfFrobAlg-reconstructionTheorems}
Expanding both sides of (\ref{eq-wdvv23}) and (\ref{eq-wdvv24}) (resp. (\ref{eq-wdvv23odd}) and (\ref{eq-wdvv24odd}) when $n$ is odd) as power series of $s$,  we obtain  equations for $F^{(i)}$ and their derivatives:
\begin{eqnarray}\label{eq-wdvv23expand}
\sum_{j=0}^{k}\frac{F_{ abe}^{(j)}g^{ef}F_{ f}^{(k-j+1)}}{j!(k-j)!}
+\sum_{j=1}^{k}\frac{2F_{ ab}^{(j)}F^{(k-j+2)}}{(j-1)!(k-j)!}
=\sum_{j=1}^{k+1}\frac{F_{ a}^{(j)}F_{b}^{(k-j+2)}}{(j-1)!(k-j+1)!},\\
  (\mbox{resp. for $k\leq \frac{m}{2}-1$ when $n$ is odd})\nn
\end{eqnarray}
\begin{eqnarray}\label{eq-wdvv24expand}
\sum_{j=1}^{k+1}\frac{F_{e}^{(j)}g^{ef}F_{f}^{(k+2-j)}}{(j-1)!(k+1-j)!}
+2\sum_{j=2}^{k+1}\frac{F^{(j)}F^{(k+3-j)}}{(j-2)!(k+1-j)!}=0,\\
  (\mbox{resp. for $k\leq \frac{m}{2}-1$ when $n$ is odd})\nn
\end{eqnarray}
where $0\leq a,b\leq n$.
The Euler vector field gives, for $k\geq 1$,
\begin{equation*}
E_{\mathrm{amb}}F^{(k)}+(2-n)k F^{(k-1)}=(3-n)F^{(k)}.
\end{equation*}

For later use, we rewrite (\ref{eq-wdvv23expand}) as, for $k\geq 1$,
\begin{eqnarray}\label{eq-higher24}
&&F_{ abe}^{(0)}g^{ef}F_{ f}^{(k+1)}+2k F_{ ab}^{(1)}F^{(k+1)}-F_{ a}^{(k+1)}F_{ b}^{(1)}-F_{ a}^{(1)}F_{ b}^{(k+1)}\nn\\
&=&\sum_{j=2}^{k}\binom{k}{j-1}F_{ a}^{(j)}F_{b}^{(k-j+2)}-\sum_{j=1}^{k}\binom{k}{j}F_{ abe}^{(j)}g^{ef}F_{ f}^{(k-j+1)}\nn\\
&&-2k\sum_{j=2}^{k}\binom{k-1}{j-1}F_{ ab}^{(j)}F^{(k-j+2)}.
\end{eqnarray}
We rewrite (\ref{eq-wdvv24expand}) as
\begin{eqnarray}\label{eq-higher25-F2}
F_{ e}^{(1)}g^{ef}F_{ f}^{(2)}+F^{(2)}F^{(2)}=0,
\end{eqnarray}
and for $k\geq 2$,
\begin{eqnarray}\label{eq-higher25}
&&F_{ e}^{(1)}g^{ef}F_{ f}^{(k+1)}+2kF^{(2)}F^{(k+1)}\nn\\
&=&-\frac{1}{2}\sum_{j=2}^{k}\binom{k}{j-1}F_{ e}^{(j)}g^{ef}F_{ f}^{(k+2-j)}
-k\sum_{j=3}^{k}\binom{k-1}{j-2}F^{(j)}F^{(k+3-j)}.
\end{eqnarray}

\subsection{\texorpdfstring{$F^{(1)}$}{F(1)} and non-semisimplicity of \texorpdfstring{$\mathcal{M}_{\mathrm{amb}}$}{Mamb}}
Taking $k=0$ in (\ref{eq-wdvv23expand}) and (\ref{eq-wdvv24expand}) we get
\begin{eqnarray}\label{eq-ss1}
\sum_{e=0}^{n}\sum_{f=0}^{n}F_{abe}^{(0)}g^{ef}F_{f}^{(1)}=F_{a}^{(1)}F_{b}^{(1)},
\end{eqnarray}
and
\begin{eqnarray}\label{eq-ss2}
\sum_{e=0}^{n}\sum_{f=0}^{n} F_{ e}^{(1)}g^{ef}F_{ f}^{(1)}=0.
\end{eqnarray}
Define 
\begin{eqnarray}\label{eq-ss2.5}
\Theta:=\sum_{e=0}^{n}\sum_{f=0}^{n}F_{e}^{(1)}g^{ef}\gamma_{f}.
\end{eqnarray}
Then $\Theta$ is a vector field on $\mathcal{M}_{\mathrm{amb}}$, and is independent of the choice of  flat coordinates.
\begin{proposition}\phantomsection\label{prop-nonsemisimple}
\begin{itemize}
\item[(i)] As a  vector field on the ambient Frobenius manifold $\mathcal{M}_{\mathrm{amb}}$ for $X$, 
 $\Theta$ is a common eigenvector for the multiplications by 
$\gamma_{a}$ with eigenvalue $F_{ a}^{(1)}$, for $0\leq a\leq n$. 
\item[(ii)] Denote the ambient big quantum multiplication by $\bqp_{\mathrm{amb}}$, then $\Theta\bqp_{\mathrm{amb}}\Theta=0$.
\item[(iii)] $(\Theta,1)=1$.
\end{itemize}
\end{proposition}
\begin{proof} (i) follows from (\ref{eq-ss1}):
\begin{eqnarray*}
&& \gamma_a\bqp \Theta=\sum_{e=0}^{n}\sum_{f=0}^{n}\sum_{c=0}^{n}\sum_{d=0}^{n}F_{e}^{(1)}g^{ef}F_{afc}^{(0)}g^{cd}\gamma_d
=\sum_{c=0}^{n}\sum_{d=0}^{n}F_{a}^{(1)}F_{c}^{(1)}g^{cd}\gamma_d.
\end{eqnarray*}
Then by (i) and (\ref{eq-ss2}),
\begin{eqnarray*}
\Theta\bqp \Theta= \sum_{a=0}^{n}\sum_{b=0}^{n}F_{a}^{(1)}g^{ab}\gamma_{b}\bqp \sum_{e=0}^{n}\sum_{f=0}^{n}F_{e}^{(1)}g^{ef}\gamma_{f}
=\sum_{a=0}^{n}\sum_{b=0}^{n}F_{a}^{(1)}g^{ab} F^{(1)}_b \cdot \Theta=0.
\end{eqnarray*}
Since 
\[
(\Theta,1)=\sum_{e=0}^{n}\sum_{f=0}^{n}F_{e}^{(1)}g^{ef}(\gamma_{f},1)=\frac{\partial F^{(1)}}{\partial t^0},
\]
(iii) is equivalent to 
$\frac{\partial F^{(1)}}{\partial t^0}=1$. By  (\ref{eq-String}) and the choice of the coordinate $s$, $\frac{\partial F^{(1)}}{\partial t^0}$ is a constant such that 
\[
\langle 1,\gamma_i,\gamma_j\rangle_{0,3,0}=\frac{\partial F^{(1)}}{\partial t^0}\cdot(\gamma_i,\gamma_j),
\]
hence $\frac{\partial F^{(1)}}{\partial t^0}=1$.
\end{proof}

Therefore we obtain a square-zero element $\Theta$ which is itself nonzero everywhere. So we obtain the following non-semisimplicity result. Note that usually one talks about semisimplicity either in the formal sense, or assuming a convergent region of the generating function. This corollary holds in both senses.

\begin{corollary}\label{cor-ss4}
Let $X$ be a non-exceptional complete intersection of dimension at least 3. Then the Frobenius manifold associated to the ambient big quantum cohomology of $X$ is not semisimple. 
\end{corollary}

\subsection{Quasi-canonical bases and reconstruction theorems}
\begin{lemma}\label{lem-artinAlgebra}
Suppose $ 0\neq b\in \mathbb{C}$, $1\leq k\leq n$, $n\geq 2$. Then we have an isomorphism of Artin algebras 
\begin{equation}\label{eq-artinAlgebra-isomorphism-0}
\varphi:  \mathbb{C}[\epsilon]/(\epsilon^{k})\oplus \mathbb{C}^{n-k+1}\cong \mathbb{C}[w]/(w^{n+1}-bw^{k}),
\end{equation}
such that $\varphi(\epsilon)=w^{n-k+2}-bw$.
\end{lemma}
\begin{proof} First we  check
 $\varphi(\epsilon)^k=0$:
\begin{eqnarray*}
\varphi(\epsilon)^k= w^k \sum_{i=0}^{k}(-b)^{k-i}\binom{k}{i}w^{(n-k+1)i}
= \sum_{i=0}^{k}(-b)^{k-i}\binom{k}{i}b^{i}w^{k}=0.
\end{eqnarray*}
Since $w^{n-k+2}-bw$ has $n-k+2$ distinct roots, we have an isomorphism of $\mathbb{C}$-algebras
$\mathbb{C}[w]/\big(\varphi(\epsilon)\big)\cong \mathbb{C}^{n-k+2}$. 
\end{proof}

 Let $e_0=1$, $e_i=\epsilon^i$ for $1\leq i\leq k-1$.  For $1\leq i\leq n-k+1$, let $e_{i+k-1}$ be the identity element of the $i$-th copy of $\mathbb{C}$ on the LHS of (\ref{eq-artinAlgebra-isomorphism-0}). Then $\{e_0,\dots,e_{n}\}$ is  a basis of the LHS of (\ref{eq-artinAlgebra-isomorphism-0}). 
Let $\zeta=\exp\big(\frac{2\pi \sqrt{-1}}{n-k+1}\big)$.  By the Chinese remainder theorem,  one can explicitly define $\varphi$ by
\begin{equation}\label{eq-artinAlgebra-isomorphism-1}
	\varphi(e_{i+k-1})=\frac{1}{(n-k+1)b^{\frac{n}{n-k+1}}\zeta^{ni}}\frac{w^{n+1}-bw^k}{w-b^{\frac{1}{n-k+1}}\zeta^i},\ \mbox{for } 1\leq i\leq n-k+1.
\end{equation}

\begin{lemma}\label{lem-artinAlg-eigenvector}
Up to a scalar factor, $\epsilon^{k-1}$ is the unique non-zero element in $\mathbb{C}[\epsilon]/(\epsilon^{k})\oplus \mathbb{C}^{n-k+1}$ satisfying both the following conditions: 
\begin{enumerate}
      \item[(i)] it is a common eigenvector in for multiplications by any element in $\mathbb{C}[\epsilon]/(\epsilon^{k})\oplus \mathbb{C}^{n-k+1}$;
      \item[(ii)] it is nilpotent.
\end{enumerate}
Moreover, we have
\begin{equation}\label{eq-artinAlg-eigenvector}
\varphi(\epsilon^{k-1})= (-1)^{k}b^{k-2}(w^n-bw^{k-1}).
\end{equation}
\end{lemma}
\begin{proof}
The first statement is obvious from the algebra structure of $\mathbb{C}[\epsilon]/(\epsilon^{k})\oplus \mathbb{C}^{n-k+1}$. For the second, we compute as
\begin{eqnarray*}
\varphi(\epsilon^{k-1})&=&(-b)^{k-1} w^{k-1}+w^{k-1} \sum_{i=1}^{k}(-b)^{k-1-i}\binom{k-1}{i}w^{(n-k+1)i}\\
&=&(-b)^{k-1} w^{k-1}+b^{k-2}w^{n} \sum_{i=1}^{k}(-1)^{k-1-i}\binom{k-1}{i}\\
&=& (-1)^{k}b^{k-2}(w^n-bw^{k-1}).
\end{eqnarray*}
\end{proof}

Replacing $w$ by $\tsfh$, $k$ by $n+1-\mathsf{a}(n,\mathbf{d})$, and $b$ by $\sfbd$, we define an isomorphism 
\begin{equation*}
	\varphi: \mathbb{C}[\epsilon]/(\epsilon^{n+1-\mathsf{a}(n,\mathbf{d})})\oplus \mathbb{C}^{\mathsf{a}(n,\mathbf{d})} \xrightarrow{\sim} \mathbb{C}[\tilde{\sfh}]/(\tilde{\sfh}^{n+1}-\sfbd\tilde{\sfh}^{n+1-\mathsf{a}(n,\mathbf{d})})
\end{equation*}
as (\ref{eq-artinAlgebra-isomorphism-0}) and (\ref{eq-artinAlgebra-isomorphism-1}).
\begin{definition}\label{lem-smallCanonicalBasis}
 Via $\varphi$, we identify $e_0,\dots,e_n$ with their images in 
\[
\mathbb{C}[\tilde{\sfh}]/(\tilde{\sfh}^{n+1}-\sfbd\tilde{\sfh}^{n+1-\mathsf{a}(n,\mathbf{d})})\cong H^*_{\mathrm{amb}}(X;\mathbb{C})\ \mbox{as } \mbox{complex vector spaces}.
\]
  We call $e_0,\dots,e_n$ the \emph{quasi-canonical basis} of $H^*_{\mathrm{amb}}(X;\mathbb{C})$.
  We denote by $u^0,\dots,u^n$  the dual basis of $e_0,\dots,e_n$. 
\end{definition}
\begin{definition}\label{lem-BGbasis}
We introduce the notations (see also \S \ref{sec:nota-conv} and \S \ref{sec:mirror-theorem})
\[
\tsfh_0=1,\ 
\tsfh_i=\underbrace{\tilde{\sfh}\sqp \cdots \sqp \tilde{\sfh}}_{i\ \mbox{factors}},\ 
\mbox{for } 1\leq i\leq n.
\]
We call $\tsfh_0,\dots,\tsfh_n$ the \emph{Beauville-Givental basis} of $H^*_{\mathrm{amb}}(X;\mathbb{C})$. 
We denote by $\tau^0,\dots,\tau^n$ the dual basis of $\tsfh^0,\dots,\tsfh^n$.
\end{definition} 
Let $M,L\in \mathrm{GL}_{n+1}(\mathbb{Q})$ such that
\begin{equation*}
	\sfh_i=\sum_{j=0}^{n}M_i^j \tsfh_j,\ \tsfh_i=\sum_{j=0}^{n} L_i^j e_j,\ 
\mbox{for } 0\leq i\leq n.
\end{equation*}
Then
\begin{equation*}
	u^j=\sum_{i=0}^n (ML)_i^j t^i,\
	t^i=\sum_{j=0}^{n} \big((ML)^{-1}\big)_j^i u^j.
\end{equation*}

\begin{lemma}\label{lem-EulerVectorField-coordinate-u}
In the coordinates $u^0,\dots,u^n$, the Euler vector field is
\begin{eqnarray}\label{eq-EulerVectorField-coordinate-u}
E&=& \sum_{i=0}^{n}\sum_{j=0}^{n}\sum_{l=0}^n
(1-i)\big((ML)^{-1}\big)_j^i (ML)_i^l u^j \frac{\partial}{\partial u^l}\nn\\
&&+(2-n)s\frac{\partial }{\partial s}
+\mathsf{a}(n,\mathbf{d})\sum_{i=0}^{n}L_1^i\frac{\partial}{\partial u^i}
-\delta_{\mathsf{a}(n,\mathbf{d}),1}\elld \frac{\partial}{\partial u^0}.
\end{eqnarray}
\end{lemma}
\begin{proof}
We have
\begin{eqnarray*}
E&=& \sum_{i=0}^{n}(1-i)t^{i}\frac{\partial}{\partial t^i}+(2-n)s\frac{\partial }{\partial s}
+\mathsf{a}(n,\mathbf{d})\frac{\partial}{\partial t^1}\\
&=& \sum_{i=0}^{n}\sum_{j=0}^{n}\sum_{l=0}^n
(1-i)\big((ML)^{-1}\big)_j^i (ML)_i^l u^j \frac{\partial}{\partial u^l}\\
&&+(2-n)s\frac{\partial }{\partial s}
+\mathsf{a}(n,\mathbf{d})\sum_{i=0}^{n}(ML)_1^i\frac{\partial}{\partial u^i}.
\end{eqnarray*}
Since
\begin{eqnarray*}
&&(ML)_1^i=\sum_{j=0}^{n} M_1^j L_j^i\\
&=& \begin{cases}
L_1^i,& \mbox{if } \mathsf{a}(n,\mathbf{d})\geq 2,\\
L_1^i-\elld L_0^i,& \mbox{if } \mathsf{a}(n,\mathbf{d})=1
\end{cases}\\
&=& \begin{cases}
L_1^i,& \mbox{if } \mathsf{a}(n,\mathbf{d})\geq 2,\\
L_1^i-\elld\delta_{i,0},& \mbox{if } \mathsf{a}(n,\mathbf{d})=1.
\end{cases}
\end{eqnarray*}
we get (\ref{eq-EulerVectorField-coordinate-u}).
\end{proof}

\begin{lemma}\label{lem-nonvanishing-coefficient-L11}
$L_1^1\neq 0$.
\end{lemma}
\begin{proof}
$\tsfh$ generates the $\mathbb{C}$-algebra $\mathbb{C}[\tilde{\sfh}]/(\tilde{\sfh}^{n+1}-\sfbd\tilde{\sfh}^{n+1-\mathsf{a}(n,\mathbf{d})})$, but $e_0$ and $e_2,\dots,e_n$ cannot. So by definition $L_1^1\neq 0$.
\end{proof}

\begin{proposition}\label{prop-Theta-inSmallQuantumCohomology}
\begin{eqnarray}\label{eq-Theta-inSmallQuantumCohomology}
\Theta|_{t^0=\dots=t^n=0}=\frac{1}{\prod_{i=1}^{r}d_i}
\big(\tilde{\sfh}^{n}-\sfbd\tilde{\sfh}^{n-\mathsf{a}(n,\mathbf{d})}\big).
\end{eqnarray}
\end{proposition}
\begin{proof}
By Proposition \ref{prop-nonsemisimple} (i), 
$\Theta|_{t^0=\dots=t^n=0}$ is a common eigenvector. By Proposition \ref{prop-nonsemisimple} (ii), $\Theta|_{t^0=\dots=t^n=0}$ is square zero (One can also deduce this from $\tilde{\sfh}^{n+1}=\sfbd\tilde{\sfh}^{n+1-\mathsf{a}(n,\mathbf{d})}$). So by Lemma \ref{lem-artinAlg-eigenvector}, $\Theta|_{t^0=\dots=t^n=0}$ is a multiple of $\tilde{\sfh}^{n}-\sfbd\tilde{\sfh}^{n-\mathsf{a}(n,\mathbf{d})}$. By  Proposition \ref{prop-nonsemisimple} (iii), $(\Theta|_{t^0=\dots=t^n=0},1)=1$.  
Then the coefficient in (\ref{eq-Theta-inSmallQuantumCohomology}) follows from a computation of the Poincaré pairing of the basis $\tilde{h}^i$, which we postpone to (\ref{eq-pairing1}) in \S \ref{sec:4point-only-ambient-insertions}.
\end{proof}

\begin{remark}\label{rmk-F1}
Proposition \ref{prop-Theta-inSmallQuantumCohomology} gives a simple proof of the results of \cite{Bea95} (for Fano complete intersections with $2 \mathsf{a}(n,\mathbf{d})\geq n+1$) and of \cite{CJ99} (for Fano hypersurfaces), and generalizes these results to all Fano complete intersections.
\end{remark}

Let $\lambda_i$ be the eigenvalue of $e_i$ on $\Theta$. 
\begin{lemma}\label{lem-eigenvalue-ei}
\[
\lambda_0=1,\ \lambda_i=0\ \mbox{for } 1\leq i\leq n.
\]
\end{lemma}
\begin{proof}
By Lemma \ref{lem-artinAlg-eigenvector}, $\lambda_i$ is equal to the eigenvalue of $e_i$ on $\epsilon^{n-\mathsf{a}(n,\mathbf{d})}$. Then the conclusion follows.
\end{proof}

\begin{definition}\label{def-order-gwInvariant}
Let $v^0,\dots,v^n$ be an arbitrary system of linear coordinates on $H^*_{\mathrm{amb}}(X;\mathbb{C})$. Then we can expand $F$ as a series of $v^0,\dots,v^n$ and $s$. 
We define a preorder\footnote{A preorder is weaker than a partial order in that the former does not demand \emph{anti-symmetry}, i.e. $a\preceq b$ and $b\preceq a$ does not imply $a=b$.} on the coefficients of the series:
the coefficient (as a non-evaluated symbol) of $v_0^{k_0}\cdots v_n^{k_n}s^k$ $\prec$ the coefficient of $v_0^{l_0}\cdots v_n^{l_n}s^l$, if and only if (i) $k<l$, or (ii) $k=l$ and $\sum_{i=0}^n k_i<\sum_{i=0}^n l_i$.
\end{definition}

\begin{theorem}[Reconstruction II of $F^{(1)}$]\label{thm-reconstruction-II-F(1)}
Let $X$ be a non-exceptional Fano complete intersections in a projective space. Then
$F^{(1)}$ can be reconstructed by (\ref{eq-system1-even}) when $n$ is even (resp. (\ref{eq-system1-odd}) when $n$ is odd) and (\ref{eq-Dim}) from the generating function $F^{(0)}$ of the ambient quantum cohomology.
\end{theorem}
\begin{proof}
We are going to show that the coefficients of $F^{(1)}$ expanded as a  series of $u^i$'s can be inductively determined by the given data.

\textsc{Step 1}: 
From Proposition \ref{prop-Theta-inSmallQuantumCohomology} one gets
\begin{equation*}
	\frac{\partial F^{(1)}}{\partial u^i}(0)
\end{equation*}
for $0\leq i\leq n$.

\textsc{Step 2}: 
By (\ref{eq-Dim}), one finds that $F^{(1)}(0)=0$ when $\mathsf{a}(n,\mathbf{d})>1$. When $\mathsf{a}(n,\mathbf{d})=1$, $F^{(1)}(0)$ is determined by
\begin{equation*}
\langle \gamma_i,\gamma_j\rangle_{0,2,1}=	F^{(1)}(0)\cdot (\gamma_i,\gamma_j)
\end{equation*}
for $\gamma_i,\gamma_j\in H^*_{\mathrm{prim}}(X_n(\mathbf{d}))$. As we observed below (\ref{eq-EV-0}), the Euler vector field (\ref{eq-system1-even-3}) and (\ref{eq-Dim}) implies the divisor equation with first Chern class. By the divisor equation
\begin{equation*}
	\langle \gamma_i,\gamma_j,\sfh\rangle_{0,2,1}=\langle \gamma_i,\gamma_j\rangle_{0,2,1}.
\end{equation*}
So
\[
F^{(1)}(0)=\frac{\partial F^{(1)}}{\partial t^1}(0),
\]
which has been reconstructed. 

\textsc{Step 3}:
 To compute a coefficient of $F^{(1)}$ of higher orders we differentiate (\ref{eq-ss1}). 
 Then we always are left  to solve linear systems of the form
\begin{eqnarray*}
\sum_{f=0}^n C_{ab}^{f}x_{f}-\lambda_{a}x_b-\lambda_{b}x_a=\mathrm{Lower\ order\ terms},
\end{eqnarray*}
where 
\[
C_{ab}^f=\frac{\partial^3 F}{\partial u^a \partial u^b \partial u^e}(0)g^{ef}
\]
are the structure constants of the Frobenius algebra at $0$, with respect to the basis $e_0,\dots,e_n$, and 
\[
\lambda_a=\frac{\partial F^{(1)}}{\partial u^a}(0)
\]
is equal to the eigenvalue of $e_a$ on $\Theta|_{u^0=\dots=u^n=0}$, and $x_j$, for $0\leq j\leq n$, are
\begin{gather*}
\frac{\partial^{|I|+1}F^{(1)}}{\partial u^a \partial u^{I}}(0)
\end{gather*}
for some common multi-index $I$.
For $n-\mathsf{a}(n,\mathbf{d})+1\leq j\leq n$, since 
\begin{eqnarray*}
\sum_{l=0}^n
C_{e_j e_j }^{l}x_l-\lambda_{e_j}x_{e_j}-\lambda_{e_j}x_{e_j}=x_{e_j},
\end{eqnarray*}
we can also solve $x_{e_j}$. For $2\leq j\leq n-\mathsf{a}(n,\mathbf{d})$, since 
\begin{eqnarray*}
\sum_{l=0}^n
C_{e_1, e_{j-1} }^{l}x_l-\lambda_{e_1}x_{e_j}-\lambda_{e_{j-1}}x_{e_j}=
\sum_{l=0}^n
C_{\epsilon, \epsilon^{j-1}}^{l}x_l-\lambda_{\epsilon}x_{ \epsilon^{j-1}}-\lambda_{ \epsilon^{j-1}}x_{ \epsilon}=x_{ \epsilon^{j}}=x_{e_j},
\end{eqnarray*}
we can solve $x_{e_{j}}$.   Since 
\begin{eqnarray*}
\sum_{l=0}^n
C_{e_0 e_0}^{l}x_l-\lambda_{e_0}x_{e_0}-\lambda_{e_0}x_{e_0}=-x_{e_0},
\end{eqnarray*}
we can solve $x_{e_0}$. 
Finally, by Lemma \ref{lem-EulerVectorField-coordinate-u} and Lemma \ref{lem-nonvanishing-coefficient-L11} we can use the Euler vector field $E$ to solve $x_{e_1}$.
\end{proof}

\begin{theorem}[Reconstruction I]\label{thm-reconstruction-I}
Let $X$ be a non-exceptional Fano complete intersections in a projective space. Then
for $k\geq 2$, $F^{(k)}$ can be reconstructed by (\ref{eq-system1-even-1}), (\ref{eq-system1-even-3}) when $n$ is even (resp. (\ref{eq-system1-odd-1}), (\ref{eq-system1-odd-3}) when $n$ is odd) from $F^{(0)}$, $F^{(1)}$ and the constant term $F^{(i)}(0)$ of $F^{(i)}$ for $2\leq i\leq k$.
\end{theorem}
\begin{proof} 
By induction on $k$ it suffices to show that $F^{(k+1)}$ can be reconstructed from (\ref{eq-system1-even-1}), (\ref{eq-system1-even-3}) and $F^{(i)}$ for $0\leq i\leq k$. Then we use (\ref{eq-wdvv23expand}) and the induction is the same as \textsc{Step 3} of the proof of Theorem \ref{thm-reconstruction-II-F(1)}.
 \end{proof}

\begin{remark}\label{rmk-reconstruction-II}
The labels of the above theorems may seem peculiar, as Reconstruction II comes ahead of Reconstruction I. My reason is that the equation (\ref{eq-system1-even-2}) is the essential one that  gives information that cannot be obtained without using the \emph{monodromy group}. Reconstruction I does not need (\ref{eq-system1-even-2}), thus one can deduce this type of theorem without knowing that the generating function $F$ has the form in Theorem \ref{thm-monodromythm}, and our proof above can be regarded merely as simplifying the proof that one does  without using the monodromy group.

Another feature of Reconstruction II is the inevitable use of (\ref{eq-Dim}), contrary to Reconstruction I, which only needs to use the Euler vector field as a whole (recall the observation below (\ref{eq-EV-0})). As we emphasized in the proof of Theorem \ref{thm-reconstruction-II-F(1)}, we use (\ref{eq-Dim})  only for the determination of $F^{(1)}(0)$. This is typical; see the proof of Theorem \ref{thm-reconstructcubicandquadric}. 
\end{remark}

\begin{corollary}
For non-exceptional Fano complete intersections with  $\gcd (n-2, \mathsf{a}(n,\mathbf{d}))>1$, $F$ can be reconstructed by (\ref{eq-system1-even-1}), (\ref{eq-system1-even-3}) when $n$ is even (resp. (\ref{eq-system1-odd-1}), (\ref{eq-system1-odd-3}) when $n$ is odd) and the dimension constraint (\ref{eq-Dim}) from $F^{(0)}$.
\end{corollary}
\begin{proof} By the dimension constraint, a necessary condition for $F^{(k)}(0)\neq 0$ is
\begin{eqnarray*}
\beta(k):=\frac{k(n-2)-(n-3)}{\mathsf{a}(n,\mathbf{d})}\in \mathbb{Z}^{>0}.
\end{eqnarray*}
There is no such $\beta$ if $\gcd (n-2, \mathsf{a}(n,\mathbf{d}))>1$.
\end{proof}

\section{Correlators of length 4 with ambient insertions}\label{sec:4pointsInvariantsAmbient}
This section serves as a preparation for Section \ref{sec:ReconstructionII-F(2)}, where we will compute $F^{(2)}(0)$ of all non-exceptional Fano complete intersections. For this goal, we need to compute $F^{(1)}$ up to second order, and for this in turn we need first  compute certain (sums of) invariants with only ambient insertions.

In this section we consider  Fano complete intersections $X=X_n(\mathbf{d})$ of dimension $n\geq 3$ and multi-degree $\mathbf{d}$ in projective spaces. Note that in this section, we allow $X$ to be exceptional. 
 Recall that $\mathsf{a}(n,\mathbf{d})$ denotes the Fano index of $X_n(\mathbf{d})$. 

.

\subsection{Correlators of length 4 with only ambient insertions}
\label{sec:4point-only-ambient-insertions}
Denote the $i$-th power of the hyperplane class in the ordinary cohomology ring by $\sfh_i$, and the dual basis by $t^0,\dots, t^n$. 
Recall the class $\tsfh$ in (\ref{eq-ss5}):
\begin{eqnarray}\label{ss5-with-q}
\tilde{\sfh}=\left\{
\begin{array}{cc}
\sfh, & \mathsf{a}(n,\mathbf{d})\geq 2,\\ 
\sfh+\elld\mathsf{q}, & \mathsf{a}(n,\mathbf{d})=1.
\end{array}\right.
\end{eqnarray}
Then Theorem \ref{thm-Givental-smallQuantumCohomology} of Givental says that 
the ambient small quantum cohomology ring of $X$ is
\begin{equation}\label{eq-smallQuantumCohomologyRing}
 	\mathbb{C}[\tilde{\sfh}]/(\tilde{\sfh}^{n+1}-\sfbd\tilde{\sfh}^{n+1-\mathsf{a}(n,\mathbf{d})}\mathsf{q}).
 \end{equation} 
As in \S \ref{sec:mirror-theorem}, we denote the $i$-th power of $\tsfh$ in the small quantum cohomology ring by $\tsfh_i$, and the dual basis by $\tau^0,\dots,\tau^n$.

In this section, we use the version (\ref{eq-generatingFunction-with-q}) of the generating function $\sfF$ with factors $\mathsf{q}^{\beta}$, so that the divisor equation can be written as
\begin{equation}\label{eq-divisorEquation-1}
	\frac{\partial}{\partial t^1}\sfF=\mathsf{q}\frac{\partial}{\partial \mathsf{q}}.
\end{equation}

\begin{lemma}\label{lem-transform-basis-1}
There is a matrix $M=(M_i^j)_{0\leq i,j\leq n} \in \mathrm{GL}_{n+1}(\mathbb{Q})$, satisfying 
\begin{subequations}\label{eq-matrix-M}
\begin{eqnarray}
	&&M_i^i=1\ \mbox{for}\ 0\leq i\leq n,\label{eq-matrix-M-1}
	\\
	&&M_i^j=0\ \mbox{if}\ \frac{i-j}{\mathsf{a}(n,\mathbf{d})}\not\in \mathbb{Z}_{\geq 0},
	\label{eq-matrix-M-2}
\end{eqnarray}
\end{subequations}
such that the transformations between the two bases are  of the following forms
\begin{subequations}
\begin{eqnarray}\label{eq-transform1}
\sfh_i=\sum_{j=0}^n M_{i}^{j}\tsfh_{j} \mathsf{q}^{\frac{i-j}{\mathsf{a}(n,\mathbf{d})}},\ \mbox{for } 0\leq i\leq n,
\label{eq-transform1-1}\\
\tsfh_i=\sum_{j=0}^n W_{i}^{j}\sfh_{j} \mathsf{q}^{\frac{i-j}{\mathsf{a}(n,\mathbf{d})}},\ \mbox{for } 0\leq i\leq n,
\label{eq-transform1-2}
\end{eqnarray}
\end{subequations}
where $W=M^{-1}$ also satisfies
\begin{subequations}\label{eq-matrix-W}
\begin{eqnarray}
	&&W_i^i=1\ \mbox{for}\ 0\leq i\leq n,\label{eq-matrix-W-1}
	\\
	&&W_i^j=0\ \mbox{if}\ \frac{i-j}{\mathsf{a}(n,\mathbf{d})}\not\in \mathbb{Z}_{\geq 0},
	\label{eq-matrix-W-2}
\end{eqnarray}
\end{subequations}
\end{lemma}
\begin{proof}
It suffices to show that the matrix $W$ defined by (\ref{eq-transform1-2}) satisfies (\ref{eq-matrix-W}). 
By definition $\sfh_0=\tsfh_0=1$, and  $\tsfh_1=\tsfh$ is defined by (\ref{eq-ss5}), so (\ref{eq-matrix-W}) is true for $i\leq 1$. Then note that
\begin{equation*}
	\sfh_i\sqp \sfh_j=\sum_{l=0}^n\langle \sfh_i,\sfh_j,\sfh_l\rangle_{0,3,\beta} \mathsf{q}^{ \beta} g^{0,n-l}\sfh_{n-l}.
\end{equation*}
By (\ref{eq-Dim}), $\langle \sfh_i,\sfh_j,\sfh_l\rangle_{0,3,\beta}\neq 0$ needs
\[
i+j+l=n+\beta\cdot \mathsf{a}(n,\mathbf{d}).
\]
So $\sfh_i\sqp \sfh_j$ has the form
\[
\sfh_i\sqp \sfh_j=\sum_{k\in \mathbb{Z}_{\geq 0}}
b_{k}\sfh_{i+j-k\mathsf{a}(n,\mathbf{d})}
\]
for some rational numbers $b_k$. So by induction on $i$, $W_i^j=0$ if $i<j$ or $\mathsf{a}(n,\mathbf{d})\nmid (i-j)$. Since
\[
\langle \sfh_i,\sfh_j,\sfh_{n-i-j}\rangle_{0,3,0}g^{n-i-j,i+j}=1,
\]
the diagonal entries are 1.
\end{proof}

The transformations between the corresponding coordinates are then
\begin{subequations}
\begin{eqnarray}\label{eq-transform2}
\tau^i=\sum_{j=0}^n M_j^i t^j=
\sum_{k\geq 0}^{k\mathsf{a}(n,\mathbf{d})\leq n-i}M_{i+k\mathsf{a}(n,\mathbf{d})}^{i}\mathsf{q}^k t^{i+k\mathsf{a}(n,\mathbf{d})},\label{eq-transform2-1}\\
t^i=\sum_{j=0}^n W_j^i \tau^j=
\sum_{k\geq 0}^{k\mathsf{a}(n,\mathbf{d})\leq n-i}W_{i+k\mathsf{a}(n,\mathbf{d})}^{i}\mathsf{q}^k \tau^{i+k\mathsf{a}(n,\mathbf{d})}.\label{eq-transform2-2}
\end{eqnarray}
\end{subequations}
Note that  $\{\tau^i\}$ are still flat coordinates, and thus the WDVV equation holds with coordinates $\{\tau^i\}$.
Denote the Poincar\'{e} pairing of $\tsfh_e$ and $\tsfh_f$ by $\eta_{ef}$ for $0\leq e,f\leq n$.
\begin{lemma}\label{lem-pairings-tsfh}
\begin{eqnarray}\label{eq-pairing1}
\eta_{e,f}=(\tsfh_e,\tsfh_f)=
\left\{
\begin{array}{cc}
(\sfbd)^{\frac{e+f-n}{\mathsf{a}(n,\mathbf{d})}} \mathsf{q}^{\frac{e+f-n}{\mathsf{a}(n,\mathbf{d})}}\prod_{i=1}^{r}d_i & 
\mathrm{if}\ \frac{e+f-n}{\mathsf{a}(n,\mathbf{d})}\in \mathbb{Z}_{\geq 0};\\
0, & \mathrm{otherwise}.
\end{array}
\right.
\end{eqnarray}
\begin{eqnarray}\label{eq-pairing2}
\eta^{e,f}=\frac{1}{\prod_{i=1}^{r}d_i}\cdot
\left\{
\begin{array}{cc}
-\sfbd\mathsf{q}, & \mathrm{if}\ e+f=n-\mathsf{a}(n,\mathbf{d});\\
1, & \mathrm{if}\ e+f=n;\\
0, & \mathrm{otherwise}.
\end{array}
\right.
\end{eqnarray}
\end{lemma}
\begin{proof} Recall that (\cite[Corollary 9.3 and 10.9]{Giv96})
\begin{eqnarray}\label{eq-higher26}
\tsfh_{n+1}=\sfbd\tsfh_{n+1- \mathsf{a}(n,\mathbf{d})}\mathsf{q}.
\end{eqnarray}
Note that $\eta_{ef}=(\tsfh_e,\tsfh_f)=(1,\tsfh_{e+f})$.  By the lower-triangularity of the linear transform (\ref{eq-transform1-2}),
\begin{equation}\label{eq-lem-pairings-sfh-2}
	(1,\tsfh_{i})=\begin{cases}
	0,& \mbox{if}\ i<n;\\
	(1,\sfh_{n})=\prod_{i=1}^r d_i, & \mbox{if}\ i=n.
	\end{cases}
\end{equation}
Then (\ref{eq-pairing1}) follows from (\ref{eq-higher26}) and (\ref{eq-lem-pairings-sfh-2}) by induction on $e+f$. Finally one checks directly that  (\ref{eq-pairing1}) is inverse to (\ref{eq-pairing2}).
\end{proof}

For $0\leq a,b\leq n$, we have the small quantum products
\begin{eqnarray}\label{eq-qp1}
\tsfh_{a}\sqp\tsfh_b=
\left\{
\begin{array}{ccc}
\tsfh_{a+b}, &  & \mbox{if } a+b\leq n;\\
(\sfbd)^k \mathsf{q}^k \tsfh_c,& \begin{array}{c} \text{where } n-\mathsf{a}(n,\mathbf{d})+1\leq c\leq n,\\ \text{and}\  k=\frac{a+b-c}{\mathsf{a}(n,\mathbf{d})}\in \mathbb{Z}_{\geq 0}, \end{array}  & \mbox{if } a+b>n.
\end{array}\right.
\end{eqnarray}
In performing the calculations, it proves to be convenient to introduce the following notations.
\begin{notation} The subscript in $\sfF_{\tau^a \tau^b\dots}$ means taking derivatives with respect to $\tau^{a}, \tau^b,\dots$.
\end{notation}
\begin{notation}
\begin{equation*}
	\big((\sfbd)^{k}\big)^{\vee}:=\begin{cases}
	(\sfbd)^{k}, & \mathrm{if}\ k\in \mathbb{Z}_{\geq 0};\\
	0,& \mathrm{otherwise}.
	\end{cases}
\end{equation*}
\end{notation}

We define a function of $(n,\mathbf{d})$ and $l\in \mathbb{Z}_{\geq 0}$:
\begin{eqnarray}\label{eq-def-functionC}
\mathsf{c}(n,l,\mathbf{d})&:=&
1+\sum_{i=n-l\mathsf{a}(n,\mathbf{d})}^{n}\sum_{j=0}^n 
\frac{j-i}{\mathsf{a}(n,\mathbf{d})} M_{j}^{i}W_{n}^{j}(\sfbd)^{\frac{i-n}{\mathsf{a}(n,\mathbf{d})}}  \nn\\
&&-\sum_{i=n-l\mathsf{a}(n,\mathbf{d})}^{n}\sum_{j=0}^n 
\frac{j-i}{\mathsf{a}(n,\mathbf{d})} M_{j}^{i}W_{n- \mathsf{a}(n,\mathbf{d})}^{j}(\sfbd)^{1+\frac{i-n}{\mathsf{a}(n,\mathbf{d})}},
\end{eqnarray}
where we set $M_{j}^{i}=0$ if $i<0$.

\begin{lemma}\label{lem-3pointInvariants-ambient}
\begin{eqnarray}\label{eq-qp1.5}
\sfF_{\tau^a \tau^b \tau^c}^{(0)}(0)=\left\{
\begin{array}{cc}
(\sfbd)^{\frac{a+b+c-n}{\mathsf{a}(n,\mathbf{d})}} \mathsf{q}^{\frac{a+b+c-n}{\mathsf{a}(n,\mathbf{d})}}\prod_{i=1}^{r}d_i, & \mathrm{if}\ \frac{a+b+c-n}{\mathsf{a}(n,\mathbf{d})}\in \mathbb{Z}_{\geq 0}; \\
0, & \mathrm{otherwise}.
\end{array}
\right.
\end{eqnarray}
\end{lemma}
\begin{proof}
It suffices to verify by (\ref{eq-pairing1}) that (\ref{eq-qp1.5}) implies (\ref{eq-qp1}):
\begin{eqnarray*}
&&\sum_{e=0}^n \sfF_{\tau^a \tau^b \tau^e}^{(0)}(0)\eta^{ec}\\
&=&\begin{cases}
\frac{1}{\prod_{i=1}^r d_i}\big(\sfF_{\tau^a \tau^b \tau^{n-c}}^{(0)}(0)-\sfF_{\tau^a \tau^b \tau^{n-c- \mathsf{a}(n,\mathbf{d})}}^{(0)}(0) \sfbd \mathsf{q}
\big),& \mbox{if}\ c\leq n-\mathsf{a}(n,\mathbf{d});\\
\frac{1}{\prod_{i=1}^r d_i}\sfF_{\tau^a \tau^b \tau^{n-c}}^{(0)}(0),& \mbox{if}\ c> n-\mathsf{a}(n,\mathbf{d})
\end{cases}\\
&=&\begin{cases}
\big((\sfbd)^{\frac{a+b-c}{\mathsf{a}(n,\mathbf{d})}}\big)^{\vee} \mathsf{q}^{\frac{a+b-c}{\mathsf{a}(n,\mathbf{d})}}
-\big((\sfbd)^{\frac{a+b-c}{\mathsf{a}(n,\mathbf{d})}-1}\big)^{\vee} \mathsf{q}^{\frac{a+b-c}{\mathsf{a}(n,\mathbf{d})}-1}\cdot\sfbd \mathsf{q}
,& \mbox{if}\ c\leq n-\mathsf{a}(n,\mathbf{d});\\
\big((\sfbd)^{\frac{a+b-c}{\mathsf{a}(n,\mathbf{d})}}\big)^{\vee} \mathsf{q}^{\frac{a+b-c}{\mathsf{a}(n,\mathbf{d})}},& \mbox{if}\ c> n-\mathsf{a}(n,\mathbf{d})
\end{cases}\\
&=&\begin{cases}
\delta_{a+b,c}
,& \mbox{if}\ c\leq n-\mathsf{a}(n,\mathbf{d});\\
\big((\sfbd)^{\frac{a+b-c}{\mathsf{a}(n,\mathbf{d})}}\big)^{\vee} \mathsf{q}^{\frac{a+b-c}{\mathsf{a}(n,\mathbf{d})}},& \mbox{if}\ c> n-\mathsf{a}(n,\mathbf{d}),
\end{cases}
\end{eqnarray*}
which is equivalent to (\ref{eq-qp1}).
\end{proof}

\begin{lemma}\label{lem-divisorEquation-tauCoordinate}
\begin{equation}\label{eq-divisorEquation-tauCoordinate}
	\frac{\partial}{\partial \tau_1}\sfF^{(0)}
=\Big(\mathsf{q}\frac{\partial}{\partial \mathsf{q}}+\mathsf{q}\sum_{i=0}^{n}\frac{\partial \tau^i}{\partial \mathsf{q}}\frac{\partial}{\partial \tau^i}\Big)\sfF^{(0)}
+\delta_{\mathsf{a}(n,\mathbf{d}),1}\cdot
\frac{\elld\mathsf{q}}{2}\sum_{e=0}^n \sum_{f=0}^n \eta_{ef}\tau^e\tau^f,
\end{equation}
where $\frac{\partial \tau^i}{\partial \mathsf{q}}$ is computed according to (\ref{eq-transform2-1}).
\end{lemma}
\begin{proof}
 In the coordinates $(t^0,t^1,\dots, t^n, \mathsf{q})$, the divisor equation for $\sfF^{(0)}$ is
\begin{eqnarray*}
\frac{\partial}{\partial t^1}\sfF^{(0)}=\mathsf{q}\frac{\partial}{\partial \mathsf{q}}\sfF^{(0)}.
\end{eqnarray*}
Thus in the coordinates $(\tau^0,\tau^1,\dots, \tau^n, \mathsf{q})$,  the divisor equation for $\sfF^{(0)}$ is
\begin{eqnarray*}
\sum_{i=0}^{n}\frac{\partial \tau_i}{\partial t^1}\frac{\partial}{\partial \tau^i}\sfF^{(0)}
=\Big(\mathsf{q}\frac{\partial}{\partial \mathsf{q}}+\mathsf{q}\sum_{i=0}^{n}\frac{\partial \tau^i}{\partial \mathsf{q}}\frac{\partial}{\partial \tau^i}\Big)\sfF^{(0)}.
\end{eqnarray*}
We treat the cases $\mathsf{a}(n,\mathbf{d})\geq 2$ and $\mathsf{a}(n,\mathbf{d})=1$ separately. First the case $\mathsf{a}(n,\mathbf{d})\geq 2$.
By  (\ref{eq-transform2}) 
\begin{eqnarray*}
\frac{\partial \tau^0}{\partial t^1}=0,\ \frac{\partial \tau^1}{\partial t^1}=1,\ \mbox{and } \frac{\partial \tau^b}{\partial t^1}=0\ \mbox{for } 2\leq b\leq n.
\end{eqnarray*}
So we get (\ref{eq-divisorEquation-tauCoordinate}). In the case $\mathsf{a}(n,\mathbf{d})=1$,
\begin{eqnarray*}
\frac{\partial \tau^0}{\partial t^1}=-\elld\mathsf{q},\ \frac{\partial \tau^1}{\partial t^1}=1,\ \mbox{and}\ \frac{\partial \tau^b}{\partial t^1}=0\ \mbox{for}\ 2\leq b\leq n,
\end{eqnarray*}
so
\begin{eqnarray*}
\frac{\partial}{\partial \tau^1}\sfF^{(0)}-\elld\mathsf{q}\frac{\partial}{\partial \tau^0}\sfF^{(0)}
=\Big(\mathsf{q}\frac{\partial}{\partial \mathsf{q}}+\mathsf{q}\sum_{i=0}^{n}\frac{\partial \tau^i}{\partial \mathsf{q}}\frac{\partial}{\partial \tau^i}\Big)\sfF^{(0)}.
\end{eqnarray*}
By  (\ref{eq-String}) we have
\begin{eqnarray*}
\frac{\partial}{\partial \tau^0}\sfF^{(0)}=\frac{1}{2}\sum_{e=0}^n \sum_{f=0}^n \eta_{ef}\tau^e\tau^f.
\end{eqnarray*}
So we get (\ref{eq-divisorEquation-tauCoordinate}).
\end{proof}

\begin{theorem}\label{thm-qp3}
\begin{eqnarray}\label{eq-qp4}
\sum_{e=0}^{n}\sfF_{\tau^a \tau^b \tau^c \tau^e}^{(0)}(0)\eta^{e0}=\left\{
\begin{array}{cc}
\mathsf{c}(n,l,\mathbf{d})(\sfbd)^{l} \mathsf{q}^l, & \mathrm{if}\ a,b,c\geq 1\ \mathrm{and}\  l=\frac{a+b+c-1}{\mathsf{a}(n,\mathbf{d})}\in \mathbb{Z}_{\geq 0};\\
0, &  \mathrm{otherwise}.
\end{array}
\right.
\end{eqnarray}
\end{theorem}
\begin{proof} 
By  (\ref{eq-String}) if one of $a,b,c$ equals $0$, say $a=0$, then
\[
\frac{\partial}{\partial t^0}
F^{(0)}_{\tau^b \tau^c \tau^e}=0.
\]
By (\ref{eq-transform2-2}) this yields
\[
\sfF_{\tau^0 \tau^b \tau^c \tau^e}^{(0)}=0.
\]
In the following of the proof we assume that $a,b,c\geq 1$. 
First we prove (\ref{eq-qp4}) in the case $a=1$. 
By (\ref{eq-divisorEquation-tauCoordinate}) and (\ref{eq-transform2}),
\begin{eqnarray*}
\frac{\partial}{\partial \tau^1}\sfF^{(0)}
&=&\Big(\mathsf{q}\frac{\partial}{\partial \mathsf{q}}+\mathsf{q}\sum_{i=0}^{n}\frac{\partial \tau_i}{\partial \mathsf{q}}\frac{\partial}{\partial \tau^i}\Big)F^{(0)}
+\delta_{\mathsf{a}(n,\mathbf{d}),1}\cdot
\frac{\elld\mathsf{q}}{2}\sum_{e=0}^n \sum_{f=0}^n g_{ef}\tau^e\tau^f\\
&=&\Big(\mathsf{q}\frac{\partial}{\partial \mathsf{q}}
+\sum_{i=0}^{n}\sum_{j=0}^n
\frac{j-i}{\mathsf{a}(n,\mathbf{d})} M_{j}^{i}\mathsf{q}^{\frac{j-i}{\mathsf{a}(n,\mathbf{d})}} t^{j}\frac{\partial}{\partial \tau^i}\Big)\sfF^{(0)}
+\delta_{\mathsf{a}(n,\mathbf{d}),1}\cdot
\frac{\elld\mathsf{q}}{2}\sum_{e=0}^n \sum_{f=0}^n g_{ef}\tau^e\tau^f\\
&=&\Big(\mathsf{q}\frac{\partial}{\partial \mathsf{q}}
+\sum_{i=0}^{n}\sum_{j=0}^n \sum_{k=0}^n
\frac{j-i}{\mathsf{a}(n,\mathbf{d})} M_{j}^{i}W_{k}^{j}\mathsf{q}^{\frac{k-i}{\mathsf{a}(n,\mathbf{d})}} 
\tau^{k}\frac{\partial}{\partial \tau^i}\Big)\sfF^{(0)}\\
&&+\delta_{\mathsf{a}(n,\mathbf{d}),1}\cdot
\frac{\elld\mathsf{q}}{2}\sum_{e=0}^n \sum_{f=0}^n g_{ef}\tau^e\tau^f.
\end{eqnarray*}
So
\begin{eqnarray*}
&&\frac{\partial}{\partial \tau^a}\frac{\partial}{\partial \tau^b}\frac{\partial}{\partial \tau^c}\frac{\partial}{\partial \tau^1}\sfF^{(0)}\\
&=&\frac{\partial}{\partial \tau^a}\frac{\partial}{\partial \tau^b}\frac{\partial}{\partial \tau^c}\Big(\mathsf{q}\frac{\partial}{\partial \mathsf{q}}
+\sum_{i=0}^{n}\sum_{j=0}^n \sum_{k=0}^n
\frac{j-i}{\mathsf{a}(n,\mathbf{d})} M_{j}^{i}W_{k}^{j}\mathsf{q}^{\frac{k-i}{\mathsf{a}(n,\mathbf{d})}} 
\tau^{k}\frac{\partial}{\partial \tau^i}\Big)\sfF^{(0)}\\
&=&\Big(\mathsf{q}\frac{\partial}{\partial \mathsf{q}}\frac{\partial}{\partial \tau^a}\frac{\partial}{\partial \tau^b}\frac{\partial}{\partial \tau^c}
+\frac{\partial}{\partial \tau^a}\frac{\partial}{\partial \tau^b}\frac{\partial}{\partial \tau^c}\sum_{i=0}^{n}\sum_{j=0}^n \sum_{k=0}^n
\frac{j-i}{\mathsf{a}(n,\mathbf{d})} M_{j}^{i}W_{k}^{j}\mathsf{q}^{\frac{k-i}{\mathsf{a}(n,\mathbf{d})}} 
\tau^{k}\frac{\partial}{\partial \tau^i}\Big)\sfF^{(0)}\\
&=&\Big(\mathsf{q}\frac{\partial}{\partial \mathsf{q}}\frac{\partial}{\partial \tau^a}\frac{\partial}{\partial \tau^b}\frac{\partial}{\partial \tau^c}
+\frac{\partial}{\partial \tau^a}\frac{\partial}{\partial \tau^b}
\sum_{i=0}^{n}\sum_{j=0}^n 
\frac{j-i}{\mathsf{a}(n,\mathbf{d})} M_{j}^{i}W_{c}^{j}\mathsf{q}^{\frac{c-i}{\mathsf{a}(n,\mathbf{d})}}\frac{\partial}{\partial \tau^i}\\
&&+\frac{\partial}{\partial \tau^a}\frac{\partial}{\partial \tau^c}
\sum_{i=0}^{n}\sum_{j=0}^n 
\frac{j-i}{\mathsf{a}(n,\mathbf{d})} M_{j}^{i}W_{b}^{j}\mathsf{q}^{\frac{b-i}{\mathsf{a}(n,\mathbf{d})}} \frac{\partial}{\partial \tau^i}
+\frac{\partial}{\partial \tau^b}\frac{\partial}{\partial \tau^c}
\sum_{i=0}^{n}\sum_{j=0}^n 
\frac{j-i}{\mathsf{a}(n,\mathbf{d})} M_{j}^{i}W_{a}^{j}\mathsf{q}^{\frac{a-i}{\mathsf{a}(n,\mathbf{d})}} \frac{\partial}{\partial \tau^i}\\
&&+\sum_{i=0}^{n}\sum_{j=0}^n \sum_{k=0}^n
\frac{j-i}{\mathsf{a}(n,\mathbf{d})} M_{j}^{i}W_{k}^{j}\mathsf{q}^{\frac{k-i}{\mathsf{a}(n,\mathbf{d})}} 
\tau^{k}\frac{\partial}{\partial \tau^i}\frac{\partial}{\partial \tau^a}\frac{\partial}{\partial \tau^b}\frac{\partial}{\partial \tau^c}\Big)\sfF^{(0)}.
\end{eqnarray*}
By (\ref{eq-qp1.5}), if
\[
\frac{a+b+c-n}{\mathsf{a}(n,\mathbf{d})}\not \in \mathbb{Z}_{\geq 0},
\]
then
\begin{equation}\label{eq-4pointFunction-1}
	\frac{\partial}{\partial \tau^a}\frac{\partial}{\partial \tau^b}\frac{\partial}{\partial \tau^c}\frac{\partial}{\partial \tau^1}\sfF^{(0)}|_{\tau=0}=0.
\end{equation}
If $a+b+c=n+l\cdot\mathsf{a}(n,\mathbf{d})$ with $l\in \mathbb{Z}_{\geq 0}$, from Lemma \ref{lem-3pointInvariants-ambient}  we get
\begin{eqnarray}\label{eq-4pointFunction-2}
&&\frac{\partial}{\partial \tau^a}\frac{\partial}{\partial \tau^b}\frac{\partial}{\partial \tau^c}\frac{\partial}{\partial \tau^1}\sfF^{(0)}|_{\tau=0}\nonumber\\
&=& \Big(l (\sfbd)^l
+\sum_{i=0}^{n}\sum_{j=0}^n 
\frac{j-i}{\mathsf{a}(n,\mathbf{d})} M_{j}^{i}W_{a}^{j}\big((\sfbd)^{\frac{b+c+i-n}{\mathsf{a}(n,\mathbf{d})}}\big)^{\vee}
+\sum_{i=0}^{n}\sum_{j=0}^n 
\frac{j-i}{\mathsf{a}(n,\mathbf{d})} M_{j}^{i}W_{b}^{j}\big((\sfbd)^{\frac{a+c+i-n}{\mathsf{a}(n,\mathbf{d})}}\big)^{\vee}\nonumber\\
&&+\sum_{i=0}^{n}\sum_{j=0}^n 
\frac{j-i}{\mathsf{a}(n,\mathbf{d})} M_{j}^{i}W_{c}^{j}\big((\sfbd)^{\frac{a+b+i-n}{\mathsf{a}(n,\mathbf{d})}}\big)^{\vee}\Big)
\mathsf{q}^{l}\prod_{i=1}^r d_i\nonumber\\
&=& \Big(l (\sfbd)^l
+\sum_{i=0}^{n}\sum_{j=0}^n 
\frac{j-i}{\mathsf{a}(n,\mathbf{d})} M_{j}^{i}W_{a}^{j}\big((\sfbd)^{l-\frac{a-i}{\mathsf{a}(n,\mathbf{d})}}\big)^{\vee}
+\sum_{i=0}^{n}\sum_{j=0}^n 
\frac{j-i}{\mathsf{a}(n,\mathbf{d})} M_{j}^{i}W_{b}^{j}\big((\sfbd)^{l-\frac{b-i}{\mathsf{a}(n,\mathbf{d})}}\big)^{\vee}\nonumber\\
&&+\sum_{i=0}^{n}\sum_{j=0}^n 
\frac{j-i}{\mathsf{a}(n,\mathbf{d})} M_{j}^{i}W_{c}^{j}\big((\sfbd)^{l-\frac{c-i}{\mathsf{a}(n,\mathbf{d})}}\big)^{\vee}\Big)
\mathsf{q}^{l}\prod_{i=1}^r d_i.
\end{eqnarray}
Combining the two cases, we can unify (\ref{eq-4pointFunction-1}) and (\ref{eq-4pointFunction-2}) into a formula that holds for all $0\leq a,b,c\leq n$:
 \begin{eqnarray}\label{eq-4pointFunction-3}
&&\frac{\partial}{\partial \tau^a}\frac{\partial}{\partial \tau^b}\frac{\partial}{\partial \tau^c}\frac{\partial}{\partial \tau^1}\sfF^{(0)}|_{\tau=0}\nonumber\\
&=& \Big(l \big((\sfbd)^l\big)^{\vee}
+\sum_{i=0}^{n}\sum_{j=0}^n 
\frac{j-i}{\mathsf{a}(n,\mathbf{d})} M_{j}^{i}W_{a}^{j}\big((\sfbd)^{l-\frac{a-i}{\mathsf{a}(n,\mathbf{d})}}\big)^{\vee}
+\sum_{i=0}^{n}\sum_{j=0}^n 
\frac{j-i}{\mathsf{a}(n,\mathbf{d})} M_{j}^{i}W_{b}^{j}\big((\sfbd)^{l-\frac{b-i}{\mathsf{a}(n,\mathbf{d})}}\big)^{\vee}\nonumber\\
&&+\sum_{i=0}^{n}\sum_{j=0}^n 
\frac{j-i}{\mathsf{a}(n,\mathbf{d})} M_{j}^{i}W_{c}^{j}\big((\sfbd)^{l-\frac{c-i}{\mathsf{a}(n,\mathbf{d})}}\big)^{\vee}\Big)
\mathsf{q}^{l}\prod_{i=1}^r d_i
\end{eqnarray}
where $l=\frac{a+b+c-n}{\mathsf{a}(n,\mathbf{d})}$. 
Now let
\[
l=\frac{a+b}{\mathsf{a}(n,\mathbf{d})}.
\]
Then (\ref{eq-pairing2}) and (\ref{eq-4pointFunction-3}) yields
\begin{eqnarray}
&&\sum_{e=0}^{n}\sfF_{\tau^1 \tau^a \tau^b \tau^e}^{(0)}(0)\eta^{e0}\nn\\
&=&\frac{1}{\prod_{i=1}^r d_i}\big(\sfF_{\tau^1 \tau^a \tau^b \tau^{n}}^{(0)}(0)
-\sfF_{\tau^1 \tau^a \tau^b \tau^{n- \mathsf{a}(n,\mathbf{d})}}^{(0)}(0) \sfbd \mathsf{q}\big)\nn\\
&=&\Big(l \big((\sfbd)^l\big)^{\vee}
+\sum_{i=0}^{n}\sum_{j=0}^n 
\frac{j-i}{\mathsf{a}(n,\mathbf{d})} M_{j}^{i}W_{a}^{j}\big((\sfbd)^{l-\frac{a-i}{\mathsf{a}(n,\mathbf{d})}}\big)^{\vee}
+\sum_{i=0}^{n}\sum_{j=0}^n 
\frac{j-i}{\mathsf{a}(n,\mathbf{d})} M_{j}^{i}W_{b}^{j}\big((\sfbd)^{l-\frac{b-i}{\mathsf{a}(n,\mathbf{d})}}\big)^{\vee}\nonumber\\
&&+\sum_{i=0}^{n}\sum_{j=0}^n 
\frac{j-i}{\mathsf{a}(n,\mathbf{d})} M_{j}^{i}W_{n}^{j}\big((\sfbd)^{l-\frac{n-i}{\mathsf{a}(n,\mathbf{d})}}\big)^{\vee}\Big)
\mathsf{q}^{l}\label{eq-4pointFunction-4}\\
&&-\Big((l-1) \big((\sfbd)^{l-1}\big)^{\vee}
+\sum_{i=0}^{n}\sum_{j=0}^n 
\frac{j-i}{\mathsf{a}(n,\mathbf{d})} M_{j}^{i}W_{a}^{j}\big((\sfbd)^{l-1-\frac{a-i}{\mathsf{a}(n,\mathbf{d})}}\big)^{\vee}\nn\\
&&+\sum_{i=0}^{n}\sum_{j=0}^n 
\frac{j-i}{\mathsf{a}(n,\mathbf{d})} M_{j}^{i}W_{b}^{j}\big((\sfbd)^{l-1-\frac{b-i}{\mathsf{a}(n,\mathbf{d})}}\big)^{\vee}\nonumber\\
&&+\sum_{i=0}^{n}\sum_{j=0}^n 
\frac{j-i}{\mathsf{a}(n,\mathbf{d})} M_{j}^{i}W_{n- \mathsf{a}(n,\mathbf{d})}^{j}\big((\sfbd)^{l-1-\frac{n- \mathsf{a}(n,\mathbf{d})-i}{\mathsf{a}(n,\mathbf{d})}}\big)^{\vee}\Big)
\mathsf{q}^{l}\cdot \sfbd.\label{eq-4pointFunction-5}
\end{eqnarray}
We compute the differences of the terms in (\ref{eq-4pointFunction-4}) and (\ref{eq-4pointFunction-5}) separately as follows.
\begin{equation}\label{eq-4pointFunction-4&5-terms-1}
	l \big((\sfbd)^l\big)^{\vee}-(l-1) \big((\sfbd)^{l-1}\big)^{\vee}\cdot \sfbd
	=\begin{cases}
	(\sfbd)^l, & \mbox{if}\ l\in \mathbb{Z}_{>0};\\
	0,& \mbox{if}\ l\leq 0.
	\end{cases}
\end{equation}
Since
\begin{equation*}
	\big((\sfbd)^k\big)^{\vee}-\big((\sfbd)^{k-1}\big)^{\vee}\cdot \sfbd=
	\begin{cases}
	1,& \mbox{if}\ k=0;\\
	0, & \mbox{if}\ k\in \mathbb{R}\setminus\{0\},
	\end{cases}
\end{equation*}
we have
\begin{eqnarray}\label{eq-4pointFunction-4&5-terms-2}
&&\sum_{i=0}^{n}\sum_{j=0}^n 
\frac{j-i}{\mathsf{a}(n,\mathbf{d})} M_{j}^{i}W_{a}^{j}\big((\sfbd)^{l-\frac{a-i}{\mathsf{a}(n,\mathbf{d})}}\big)^{\vee}
-\sum_{i=0}^{n}\sum_{j=0}^n 
\frac{j-i}{\mathsf{a}(n,\mathbf{d})} M_{j}^{i}W_{a}^{j}\big((\sfbd)^{l-1-\frac{a-i}{\mathsf{a}(n,\mathbf{d})}}\big)^{\vee}\cdot \sfbd\nn\\
&=&\begin{cases}
\sum_{j=0}^n 
\frac{j}{\mathsf{a}(n,\mathbf{d})} M_{j}^{0}W_{a}^{j}, & \mbox{if } l-\frac{a}{\mathsf{a}(n,\mathbf{d})}=0;\\
0,& \mbox{if } l-\frac{a}{\mathsf{a}(n,\mathbf{d})}\neq 0
\end{cases}\nn\\
&=&\begin{cases}
\sum_{j=0}^n 
\frac{j}{\mathsf{a}(n,\mathbf{d})} M_{j}^{0}W_{a}^{j}, & \mbox{if } b=0;\\
0,& \mbox{if } b>0.
\end{cases}
\end{eqnarray}
In the same way we have
\begin{eqnarray}\label{eq-4pointFunction-4&5-terms-3}
&&\sum_{i=0}^{n}\sum_{j=0}^n 
\frac{j-i}{\mathsf{a}(n,\mathbf{d})} M_{j}^{i}W_{b}^{j}\big((\sfbd)^{l-\frac{b-i}{\mathsf{a}(n,\mathbf{d})}}\big)^{\vee}
-\sum_{i=0}^{n}\sum_{j=0}^n 
\frac{j-i}{\mathsf{a}(n,\mathbf{d})} M_{j}^{i}W_{b}^{j}\big((\sfbd)^{l-1-\frac{b-i}{\mathsf{a}(n,\mathbf{d})}}\big)^{\vee}\cdot \sfbd\nn\\
&=&\begin{cases}
\sum_{j=0}^n 
\frac{j}{\mathsf{a}(n,\mathbf{d})} M_{j}^{0}W_{b}^{j}, & \mbox{if } a=0;\\
0,& \mbox{if } a>0.
\end{cases}
\end{eqnarray}
Moreover by (\ref{eq-matrix-M}) one has
\begin{eqnarray}\label{eq-4pointFunction-4&5-terms-4}
&&\sum_{i=0}^{n}\sum_{j=0}^n 
\frac{j-i}{\mathsf{a}(n,\mathbf{d})} M_{j}^{i}W_{n}^{j}\big((\sfbd)^{l-\frac{n-i}{\mathsf{a}(n,\mathbf{d})}}\big)^{\vee}\nn\\
&&-\sum_{i=0}^{n}\sum_{j=0}^n 
\frac{j-i}{\mathsf{a}(n,\mathbf{d})} M_{j}^{i}W_{n- \mathsf{a}(n,\mathbf{d})}^{j}\big((\sfbd)^{l-1-\frac{n- \mathsf{a}(n,\mathbf{d})-i}{\mathsf{a}(n,\mathbf{d})}}\big)^{\vee}\cdot \sfbd\nn\\
&=&\sum_{i=n-a-b}^{n}\sum_{j=0}^n 
\frac{j-i}{\mathsf{a}(n,\mathbf{d})} M_{j}^{i}W_{n}^{j}(\sfbd)^{l-\frac{n-i}{\mathsf{a}(n,\mathbf{d})}}\nn\\
&&-\sum_{i=n-a-b}^{n}\sum_{j=0}^n 
\frac{j-i}{\mathsf{a}(n,\mathbf{d})} M_{j}^{i}W_{n- \mathsf{a}(n,\mathbf{d})}^{j}(\sfbd)^{l+1-\frac{n-i}{\mathsf{a}(n,\mathbf{d})}}.
\end{eqnarray}
So by (\ref{eq-4pointFunction-4&5-terms-1}), (\ref{eq-4pointFunction-4&5-terms-2}), (\ref{eq-4pointFunction-4&5-terms-3}) and (\ref{eq-4pointFunction-4&5-terms-4}), when $1\leq a,b\leq n$ and $l\in \mathbb{Z}_{>0}$, we obtain
\begin{eqnarray*}
&&\sum_{e=0}^{n}\sfF_{\tau^1 \tau^a \tau^b \tau^e}^{(0)}(0)\eta^{e0}\\
&=& \mathsf{q}^{l}\Big(
(\sfbd)^l 
+\sum_{i=n-a-b}^{n}\sum_{j=0}^n 
\frac{j-i}{\mathsf{a}(n,\mathbf{d})} M_{j}^{i}W_{n}^{j}(\sfbd)^{\frac{a+b+i-n}{\mathsf{a}(n,\mathbf{d})}}  \\
&&-\sum_{i=n-a-b}^{n}\sum_{j=0}^n 
\frac{j-i}{\mathsf{a}(n,\mathbf{d})} M_{j}^{i}W_{n- \mathsf{a}(n,\mathbf{d})}^{j}(\sfbd)^{1+\frac{a+b+i-n}{\mathsf{a}(n,\mathbf{d})}}\Big)\\
&=&
\mathsf{q}^{l}(\sfbd)^l\Big(1+\sum_{i=n-a-b}^{n}\sum_{j=0}^n 
\frac{j-i}{\mathsf{a}(n,\mathbf{d})} M_{j}^{i}W_{n}^{j}(\sfbd)^{\frac{i-n}{\mathsf{a}(n,\mathbf{d})}} \\
&&-\sum_{i=n-a-b}^{n}\sum_{j=0}^n 
\frac{j-i}{\mathsf{a}(n,\mathbf{d})} M_{j}^{i}W_{n- \mathsf{a}(n,\mathbf{d})}^{j}(\sfbd)^{1+\frac{i-n}{\mathsf{a}(n,\mathbf{d})}}\Big).
\end{eqnarray*}
This proves (\ref{eq-qp4}) when at least one of $a,b,c$ is equal to 1. 
Now we prove  (\ref{eq-qp4}) for general $a,b,c$.  Differentiating the WDVV equation for $\sfF^{(0)}$ once, and taking values at $\tau=0$, one gets
\begin{eqnarray}\label{eq-qp5}
&& \sfF_{\tau^a \tau^b \tau^p\tau^e }^{(0)}(0)\eta^{ef}\sfF_{\tau^f \tau^c \tau^d}^{(0)}(0)+\sfF_{\tau^a \tau^b \tau^e}^{(0)}(0)\eta^{ef}\sfF_{\tau^f \tau^c \tau^d \tau^p}^{(0)}(0)\nn\\
&=& \sfF_{\tau^a \tau^c \tau^p\tau^e }^{(0)}(0)\eta^{ef}\sfF_{\tau^f \tau^b \tau^d}^{(0)}(0)+\sfF_{\tau^a \tau^c \tau^e}^{(0)}(0)\eta^{ef}\sfF_{\tau^f \tau^b \tau^d \tau^p}^{(0)}(0).
\end{eqnarray}
To show (\ref{eq-qp4}) for all $0\leq a,b,c\leq n$, for convenience, we extend the range $0\leq a,b,c,d\leq n$ in the subscript of $F_{\tau^a \tau^b \tau^c \tau^d}$ to the set of nonnegative integers. Note that 
\[
\mathsf{a}(n,\mathbf{d})\leq n.
\]
So for and integer $a>n$, there is a unique way to write $a$ as
\begin{equation*}
 a=a_1+k\mathsf{a}(n,\mathbf{d}),
\end{equation*}
  where $1\leq a_1\leq n$, and $k\in \mathbb{Z}_{\geq 0}$. Then we define 
\begin{equation*}
	\sfF^{(0)}_{\tau^a \tau^b \tau^c \tau^d}(0):=(\sfbd)^k \mathsf{q}^k \sfF^{(0)}_{\tau^{a_1} \tau^b \tau^c \tau^d}(0).
\end{equation*}
 We do the same  to the other superscripts 
$b,c,d$ in $\sfF^{(0)}_{\tau^a \tau^b \tau^c \tau^d}(0)$. With this convention, we have shown that (\ref{eq-qp4}) is valid for $a=1$ and all integers $b,c\geq 1$. By the definition of our basis $\tsfh_i$, (\ref{eq-qp1}) and (\ref{eq-qp5}) imply, for $a,b,c,d,p\geq 0$,
\begin{eqnarray}\label{eq-qp6}
\sfF_{\tau^a \tau^b \tau^p \tau^{c+d}}^{(0)}(0)+\sfF_{\tau^{a+b} \tau^{c} \tau^{d} \tau^p}^{(0)}(0)=\sfF_{\tau^a \tau^c \tau^p \tau^{b+d}}^{(0)}(0)+\sfF_{\tau^{a+c}\tau^b \tau^d \tau^p}^{(0)}(0).
\end{eqnarray}
Taking $b=1$, we get
\begin{eqnarray}\label{eq-qp7}
\sfF_{\tau^{a+1} \tau^c \tau^d \tau^p}^{(0)}(0)=\sfF_{\tau^a \tau^c \tau^p \tau^{1+d}}^{(0)}(0)+\sfF_{\tau^{a+c} \tau^1 \tau^d \tau^p}^{(0)}(0)-\sfF_{\tau^a \tau^1 \tau^p \tau^{c+d}}^{(0)}(0),
\end{eqnarray}
and thus
\begin{eqnarray}\label{eq-qp8}
&&\sum_{p=0}^{n}\sfF_{\tau^{a+1} \tau^c \tau^d \tau^p}^{(0)}(0)\eta^{p0}\nn\\
&=&\sum_{p=0}^{n}\sfF_{\tau^a \tau^c\tau^{1+d} \tau^{p}}^{(0)}(0)\eta^{p0}+\sum_{p=0}^{n}\sfF_{\tau^{a+c} \tau^1 \tau^d \tau^p}^{(0)}(0)\eta^{p0}-\sum_{p=0}^{n}\sfF_{\tau^a \tau^1 \tau^{c+d}\tau^{p}}^{(0)}(0)\eta^{p0}.
\end{eqnarray}
Hence (\ref{eq-qp4}) follows by induction on $a$.
\end{proof}

\subsection{Leading terms of \texorpdfstring{$F^{(1)}$}{F(1)}}

As a corollary of Theorem \ref{thm-qp3}, we can compute $F^{(1)}$ to degree 2. We treat the cases $\mathsf{a}(n,\mathbf{d})\geq 2$ and $\mathsf{a}(n,\mathbf{d})=1$ separately.
\begin{theorem}\label{thm-F(1)}
Suppose $\mathsf{a}(n,\mathbf{d})\geq 2$. Then
\begin{equation}\label{eq-thm-F(1)}
\sfF^{(1)}(\tau)
=\tau^0-\frac{1}{2}\sum_{k\geq 1}\Big(\mathsf{c}(n,k,\mathbf{d})(\sfbd)^{k}\mathsf{q}^{k}\sum_{\begin{subarray}{c}1\leq i,j\leq n\\
i+j=1+k\mathsf{a}(n,\mathbf{d})\end{subarray}
}
\tau^i \tau^{j}\Big)
+O\big((\tau)^3\big).
\end{equation}
\end{theorem}
\begin{proof}
First we assume that $X$ is a non-exceptional complete intersection.
 By (\ref{eq-ss2.5}) and Proposition \ref{prop-Theta-inSmallQuantumCohomology}, 
\begin{equation}\label{eq-cor-F(1)-1}
	\sum_{e=0}^{n}\sum_{f=0}^n \sfF_{\tau^e}^{(1)}(0)\eta^{ef}\tsfh_{f}=
	\frac{1}{\prod_{i=1}^{r}d_i}
\big(\tilde{\sfh}^{n}-\sfbd\tilde{\sfh}^{n-\mathsf{a}(n,\mathbf{d})}\mathsf{q}\big),
\end{equation}
and by Proposition \ref{prop-nonsemisimple}, $\sfF_{\tau^e}^{(1)}(0)$ is the eigenvalue of $\tsfh_{e}$ on the RHS of (\ref{eq-cor-F(1)-1}). So 
\begin{equation}\label{eq-cor-F(1)-2}
	\sfF_{\tau^e}^{(1)}(0)=\begin{cases}
	1,& \mbox{if}\ e=0;\\
	0,& \mbox{otherwise}.
	\end{cases}
\end{equation}
This gives the first term of (\ref{eq-thm-F(1)}). 
 Differentiating (\ref{eq-ss1}) once, we get
\begin{eqnarray*}
\sfF_{\tau^a \tau^b \tau^c \tau^e}^{(0)}(0)
\eta^{ef}\sfF_{\tau^f}^{(1)}(0)+\sfF_{\tau^a \tau^b \tau^e}^{(0)}(0)\eta^{ef}\sfF_{\tau^f \tau^c}^{(1)}(0)=\sfF_{\tau^a \tau^c}^{(1)}(0)\sfF_{\tau^b}^{(1)}(0)+\sfF_{\tau^a}^{(1)}(0)\sfF_{\tau^b \tau^c}^{(1)}(0).
\end{eqnarray*}
Using  (\ref{eq-cor-F(1)-2}) we get
\begin{equation}\label{eq-qp2}
\sfF_{\tau^a \tau^b \tau^c \tau^e}^{(0)}(0)\eta^{e0}+\sfF_{\tau^a \tau^b \tau^e}^{(0)}(0)\eta^{ef}\sfF_{\tau^f \tau^c}^{(1)}(0)=0,\ \mbox{for}\ 1\leq a,b\leq n.
\end{equation}
Then by (\ref{eq-qp1}) and Theorem \ref{thm-qp3}, when $a,b\geq 1$, $a+b\leq n$ and $c\geq 1$, (\ref{eq-qp2}) yields 
\begin{eqnarray}\label{eq-cor-F(1)-3}
	&&\sfF_{\tau^{a+b} \tau^c}^{(1)}(0)=
	-\sum_{e=0}^n \sfF_{\tau^a \tau^b \tau^c \tau^e}^{(0)}(0)\eta^{e0}\nonumber\\
	&=&-\mathsf{c}(n,\frac{a+b+c-1}{\mathsf{a}(n,\mathbf{d})},\mathbf{d})\big((\sfbd)^{\frac{a+b+c-1}{\mathsf{a}(n,\mathbf{d})}}\big)^{\vee} \mathsf{q}^{\frac{a+b+c-1}{\mathsf{a}(n,\mathbf{d})}}.
\end{eqnarray}
So (\ref{eq-cor-F(1)-3}) gives all $\sfF_{\tau^{i} \tau^{j}}^{(1)}(0)$ if $i\geq 2$ or $j\geq 2$. It remains to compute $\sfF_{\tau^{1} \tau^{1}}^{(1)}(0)$. But  $\sfF_{\tau^{i} \tau^{j}}^{(1)}(0)=0$ unless $\mathsf{a}(n,\mathbf{d})$ divides $i+j-1$; since $\mathsf{a}(n,\mathbf{d})\geq 2$ by assumption, this cannot happen if $i=j=1$. Finally by (\ref{eq-transform2})  $\partial_{t^0}=\partial_{\tau^0}$, so if $i=0$ or $j=0$, say $i$,  (\ref{eq-String}) yields 
\begin{equation*}
 	\frac{\partial^2 F^{(1)}}{\partial \tau^0 \partial \tau^j}=0.
 \end{equation*} 
So we complete the proof of (\ref{eq-thm-F(1)}) for non-exceptional complete intersections. 

Suppose now that $X$ is exceptional, i.e. $X$ is a quadric, or an even dimensional complete intersection of two quadrics. Then by Remark \ref{rmk-F1}, Proposition \ref{prop-Theta-inSmallQuantumCohomology} is known to be valid by \cite{Bea95}. So similar computations yield (\ref{eq-thm-F(1)}).
\end{proof}

\begin{lemma}\label{lem-specialValue-c(n,1,d)}
Suppose $\mathsf{a}(n,\mathbf{d})=1$, then
\begin{equation}\label{eq-lem-specialValue-c(n,1,d)}
	\mathsf{c}(n,1,\mathbf{d})=\frac{\elld}{\sfbd}.
\end{equation}
\end{lemma}
\begin{proof}
By the definition (\ref{eq-def-functionC}), 
\begin{equation}\label{eq-specialValue-c(n,1,d)-1}
	\mathsf{c}(n,1,\mathbf{d})=1+M_n^{n-1} (\sfbd)^{-1}. 
\end{equation}
By (\ref{eq-matrix-M}), 
\begin{equation}\label{eq-specialValue-c(n,1,d)-2}
	M_{n}^{n-1}=-W_n^{n-1}.
\end{equation}
As in the proof of Lemma \ref{lem-transform-basis-1}, by (\ref{eq-Deg0}) and (\ref{eq-FCA}), the small quantum product in the ordinary basis satisfies, for $i+j\leq n$,
\begin{equation}\label{eq-specialValue-c(n,1,d)-3}
	\sfh_i\sqp \sfh_j=\sfh_{i+j}+\sum_{k<i+j}c_k\sfh_k
\end{equation}
for some $c_k\in \mathbb{Q}$ depending on $i,j$. Similarly,
\begin{equation}\label{eq-specialValue-c(n,1,d)-4}
	\sfh\sqp \sfh_n=\sum_{k<n}b_k\sfh_k
\end{equation}
for some $b_k\in \mathbb{Q}$. So 
\begin{equation}\label{eq-specialValue-c(n,1,d)-5}
	\tsfh\sqp \sfh_n=\sfh\sqp \sfh_{n}+\elld\sfh_n
	=\elld\sfh_n+\sum_{i<n}b_i \sfh_i,
\end{equation}
and
\begin{eqnarray}\label{eq-specialValue-c(n,1,d)-6}
&& \tsfh\sqp \tsfh_n=\tsfh\sqp (\sfh_n+W_n^{n-1}\sfh_{n-1}+\sum_{i=0}^{n-2}W_n^{i}\sfh_{i})\nn\\
&&=\big(\elld+W_n^{n-1}\big)\sfh_n+\sum_{i<n}a_i \sfh_i.
\end{eqnarray}
for some $a_i\in \mathbb{Q}$. 
On the other hand by (\ref{eq-smallQuantumCohomologyRing}) we have
\begin{equation}\label{eq-specialValue-c(n,1,d)-7}
	\tsfh\sqp \tsfh_n=\tsfh_{n+1}=\sfbd \tsfh_{n}
	=\sfbd \sfh_n+ \sfbd\sum_{i=0}^{n-1}W_n^{i}\sfh_{i}.
\end{equation}
Comparing (\ref{eq-specialValue-c(n,1,d)-6}) and (\ref{eq-specialValue-c(n,1,d)-6}) we get
\begin{equation}\label{eq-specialValue-c(n,1,d)-8}
	W_n^{n-1}=\sfbd-\elld.
\end{equation}
Then by (\ref{eq-specialValue-c(n,1,d)-1}) and (\ref{eq-specialValue-c(n,1,d)-2}) we obtain (\ref{eq-lem-specialValue-c(n,1,d)}).
\end{proof}

\begin{theorem}\label{thm-F(1)-fanoIndex=1}
Suppose $\mathsf{a}(n,\mathbf{d})=1$. Then
\begin{equation}\label{eq-thm-F(1)-fanoIndex=1}
\sfF^{(1)}(\tau)
=-\elld \mathsf{q}+
\tau^0-\frac{1}{2}\sum_{k\geq 1}\Big(\mathsf{c}(n,k,\mathbf{d})(\sfbd)^{k}\mathsf{q}^{k}\sum_{\begin{subarray}{c}1\leq i,j\leq n\\
i+j=1+k\mathsf{a}(n,\mathbf{d})\end{subarray}
}
\tau^i \tau^{j}\Big)+O\big((\tau)^3\big).
\end{equation}
\end{theorem}
\begin{proof}
We need only in addition to compute the coefficient of 1 and $(\tau^1)^2$.
\begin{equation*}
	\frac{\partial \sfF^{(1)}}{\partial t^1}=\frac{\partial \sfF^{(1)}}{\partial \tau^1}-\elld\mathsf{q} \frac{\partial \sfF^{(1)}}{\partial \tau^0},
\end{equation*}
so by (\ref{eq-Div}),
\[
F^{(1)}(0)=-\elld\mathsf{q}.
\]
Since 
\begin{equation*}
	\frac{\partial}{\partial \tau^1}= \mathsf{q}\frac{\partial t^0}{\partial \tau^1}\frac{\partial}{\partial t^0}+\frac{\partial}{\partial t^1},
\end{equation*}
by (\ref{eq-String}) and (\ref{eq-Div}) we get
\[
\frac{\partial^2 \sfF^{(1)}}{\partial \tau^1 \partial \tau^1}(0)=F^{(1)}(0)=-\elld\mathsf{q}.
\]
So 
\begin{eqnarray}\label{eq-cor-F(1)-fanoIndex=1}
\sfF^{(1)}(\tau)
&=&-\elld \mathsf{q}+
\tau^0-\frac{1}{2}\sum_{k\geq 1}\Big(\mathsf{c}(n,k,\mathbf{d})(\sfbd)^{k}\mathsf{q}^{k}\sum_{\begin{subarray}{c}1\leq i,j\leq n\\
i+j=1+k\mathsf{a}(n,\mathbf{d})\\
(i,j)\neq (1,1)\end{subarray}
}
\tau^i \tau^{j}\Big)\nn\\
&&-\frac{\elld\mathsf{q}}{2} (\tau^1)^2
+O\big((\tau)^3\big).
\end{eqnarray}
Then by Lemma \ref{lem-specialValue-c(n,1,d)} we obtain (\ref{eq-thm-F(1)-fanoIndex=1}).
\end{proof}
\begin{example}\label{example-F1-leadingTerm-FanoIndex=n-1}
Suppose furthermore $\mathsf{a}(n,\mathbf{d})=n-1$.  Then (\ref{eq-thm-F(1)}) reads
\begin{eqnarray*}
\sfF^{(1)}(\tau)
&=&\tau^0-\frac{1}{2}\mathsf{c}(n,1,\mathbf{d})\sfbd\mathsf{q}\sum_{\begin{subarray}{c}1\leq i,j\leq n\\
i+j=n\end{subarray}
}
\tau^i \tau^{j}\\
&&-\mathsf{c}(n,2,\mathbf{d})(\sfbd)^2\mathsf{q}^2 \tau^{n-1}\tau^{n}
+O\big((\tau)^3\big).
\end{eqnarray*}
By (\ref{eq-def-functionC}) and (\ref{eq-matrix-M}), 
\begin{equation*}
\mathsf{c}(n,2,\mathbf{d})=\mathsf{c}(n,1,\mathbf{d})=
1+M_{n}^{n-\mathsf{a}(n,\mathbf{d})}(\sfbd)^{-1}.
\end{equation*}
Recall the notation in Theorem \ref{thm-Beuville}. 
By (\ref{eq-Beauville-sfhTotsfh}),
\begin{equation*}
	M_{n}^{n-\mathsf{a}(n,\mathbf{d})}=
	-\sum_{i=0}^{n-\mathsf{a}(n,\mathbf{d})}\ell_i. 
\end{equation*}
But by (\ref{eq-Beauville-ell-sum}), (\ref{eq-Beauville-dual}) and (\ref{eq-Beauville-ell-l0-l1}),
\begin{equation*}
	\sum_{i=0}^{n-\mathsf{a}(n,\mathbf{d})}\ell_i=\sfbd-\ell_{n+1-\mathsf{a}(n,\mathbf{d})}=\sfbd-\ell_{0}=\sfbd-\elld.
\end{equation*}
So
\begin{equation}\label{eq-c(n,1,d)-largeFanoIndex}
\mathsf{c}(n,2,\mathbf{d})=\mathsf{c}(n,1,\mathbf{d})=\frac{\elld}{\sfbd},
\end{equation}
and thus
\begin{equation}\label{eq-F1-leadingTerms-largeFanoIndex}
\sfF^{(1)}(\tau)
=\tau^0-\frac{1}{2}\elld\mathsf{q}\sum_{\begin{subarray}{c}1\leq i,j\leq n\\
i+j=n\end{subarray}
}
\tau^i \tau^{j}
-\elld\sfbd\mathsf{q}^2 \tau^{n-1}\tau^{n}
+O\big((\tau)^3\big).
\end{equation}
By (\ref{eq-Beauville-sfhTotsfh}), (\ref{eq-Beauville-ell-sum}), (\ref{eq-Beauville-dual}) and (\ref{eq-Beauville-ell-l0-l1}), one finds
\begin{equation}\label{eq-M-FanoIndex=n-1}
\begin{cases}
	M_{n-1}^0=-\ell_0=-\elld,\\
	M_{n}^1=-\ell_0-\ell_1=\ell_2- \sfbd=\ell_0- \sfbd
	=\elld- \sfbd.
\end{cases}	
\end{equation}
So by (\ref{eq-transform2-1}), 
\begin{equation}\label{eq-tauTot-FanoIndex=n-1}
\begin{cases}
\tau^0=t^0-\elld \mathsf{q}t^{n-1},\\
\tau^1=t^1+\big(\elld- \sfbd\big) \mathsf{q}t^{n},\\
\tau^i=t^i\ \mbox{for}\ i\geq 2.
\end{cases}
\end{equation}
Then we can write $F^{(1)}$ in the $t$-coordinates:
\begin{eqnarray}\label{eq-F1-leadingTerm-FanoIndex=n-1}
\sfF^{(1)}(\tau^0,\dots,\tau^n)
=t^0-\elld t^{n-1}-\frac{\elld}{2}\mathsf{q}\sum_{i=1}^{n-1}t^{i}t^{n-i}-\elld^2\mathsf{q}^{2}t^{n-1}t^{n}+O((t)^3).
\end{eqnarray}
In particular, for $X_n(2,2)$ where $n$ is odd and $n\geq 3$,
\begin{eqnarray}\label{eq-F^122}
\sfF^{(1)}(\tau^0,\dots,\tau^n)
=t^0-4t^{n-1}-2\mathsf{q}\sum_{i=1}^{n-1}t^{i}t^{n-i}-16\mathsf{q}^{2}t^{n-1}t^{n}+O((t)^3).
\end{eqnarray}
In the derivation of $\sfF^{(1)}$ we only use the first order  truncation of the WDVV with respect to $s$. By the invariant theory of $D_{n+3}$ (see the beginning of Section \ref{sec:4points-Invariants}), 
(\ref{eq-F^122}) is also valid when $n$ is even.
\end{example}

\section{Reconstruction II of \texorpdfstring{$F^{(2)}$}{F(2)}}\label{sec:ReconstructionII-F(2)}

In this section, we assume that $X=X_n(\mathbf{d})$ is a non-exceptional complete intersection in projective spaces of multi-degree $\mathbf{d}$. We use version (\ref{eq-generatingFunction-without-q}) of the generating function of genus 0 GW invariants of $X$. For the definitions of the coordinates $\tau^i$, the pairings $\eta_{ij}$, and the matrix $(M_{i}^j)$ and $(W_{i}^j)$, we refer the reader to Section \ref{sec:4point-only-ambient-insertions}; in this section we specify the parameter $\mathsf{q}$ in them to $\mathsf{q}=1$.

\subsection{Leading terms of \texorpdfstring{$F^{(2)}$}{F(2)}}

We are ready to compute $F^{(2)}(0)$ and $F_{\tau^b}^{(2)}(0)$.
\begin{lemma}
\item[(i)]
\begin{equation}\label{eq-higher5}
\sum_{e=0}^n \sum_{f=0}^nF_{\tau^e}^{(1)}\eta^{ef}F_{\tau^f}^{(2)}+F^{(2)}F^{(2)}=0.
\end{equation}
\begin{eqnarray}\label{eq-higher11}
&&-\sum_{e=0}^n \sum_{f=0}^nF_{\tau^a \tau^e}^{(1)}\eta^{ef}F_{\tau^f \tau^b}^{(1)}+\sum_{e=0}^n \sum_{f=0}^nF_{\tau^a \tau^b \tau^e}^{(0)}\eta^{ef}F_{\tau^f}^{(2)}+2F_{\tau^a \tau^b}^{(1)}F^{(2)}\nn\\
&=&F_{\tau^a}^{(2)}F_{\tau^b}^{(1)}+F_{\tau^a}^{(1)}F_{\tau^b}^{(2)}.
\end{eqnarray}
\item[(ii)]  For $2\leq b\leq n$,
\begin{eqnarray}\label{eq-higher21}
F_{\tau^b}^{(2)}(0)+2F_{\tau^1 \tau^{b-1}}^{(1)}(0)F^{(2)}(0)
=\sum_{e=0}^n \sum_{f=0}^n F_{\tau^1 \tau^e}^{(1)}(0)\eta^{ef}F_{\tau^f \tau^{b-1}}^{(1)}(0).
\end{eqnarray}
\end{lemma}
\begin{proof} (i) Taking $k=1$ in (\ref{eq-wdvv24expand})  we obtain (\ref{eq-higher5}). Differentiating (\ref{eq-ss2}) by $\partial_{\tau^a} \partial_{\tau^b}$, we have
\begin{eqnarray}\label{eq-higher10}
F_{\tau^a \tau^e}^{(1)}\eta^{ef}F_{\tau^f \tau^b}^{(1)}+F_{\tau^e}^{(1)}\eta^{ef}F_{\tau^f \tau^a \tau^b}^{(1)}=0.
\end{eqnarray}
Substituting (\ref{eq-higher10}) into (\ref{eq-wdvv23expand}) for $k=1$ , we obtain (\ref{eq-higher11}).  \\
(ii) By (\ref{eq-higher11}),
\begin{eqnarray*}
-F_{\tau^a \tau^e}^{(1)}(0)\eta^{ef}F_{\tau^f \tau^b}^{(1)}(0)+F_{\tau^a \tau^b \tau^e}^{(0)}(0)\eta^{ef}F_{\tau^f}^{(2)}(0)+2F_{\tau^a \tau^b}^{(1)}(0)F^{(2)}(0)
=0, & 1\leq a,b\leq n.
\end{eqnarray*}
Taking $a=1$, noting that $F_{\tau^a \tau^b \tau^e}^{(0)}(0)\eta^{ef}$ is the structure constant of the small quantum cohomology under the basis $\tsfh_i$, from  (\ref{eq-qp1}) we obtain (\ref{eq-higher21}).
\end{proof} 

By (\ref{eq-Dim}), $F^{(2)}(0)\neq 0$ forces
\[
4\cdot \frac{n}{2}=n-3+4+\beta\cdot \mathsf{a}(n,\mathbf{d}),
\]
for some $\beta\in \mathbb{Z}_{\geq 0}$. So if $\frac{n-1}{\mathsf{a}(n,\mathbf{d})}\not\in \mathbb{Z}$,  $F^{(2)}(0)$=0. 
\begin{lemma}
Suppose $\frac{n-1}{\mathsf{a}(n,\mathbf{d})}\in \mathbb{Z}$. Then
\begin{eqnarray}\label{eq-higher20}
&&F^{(2)}(0)F^{(2)}(0)+\Big(\frac{n-1}{\mathsf{a}(n,\mathbf{d})}\eta^{01}
-2\sum_{b=2}^{n}\eta^{0b}F_{\tau^1 \tau^{b-1}}^{(1)}(0)\Big)F^{(2)}(0)\nn\\
&&+\sum_{e=0}^n \sum_{f=0}^n\sum_{b=2}^{n}\eta^{0b}F_{\tau^1 \tau^e}^{(1)}(0)\eta^{ef}F_{\tau^f \tau^{b-1}}^{(1)}(0)=0.
\end{eqnarray}
\end{lemma}
\begin{proof}
By (\ref{eq-Div}) we have 
\begin{eqnarray}\label{eq-higher15}
\frac{\partial F^{(2)}}{\partial t^1}(0)=\frac{n-1}{\mathsf{a}(n,\mathbf{d})}F^{(2)}(0).
\end{eqnarray}
So $F_{\tau^1}^{(2)}(0)$ can be deduced from (\ref{eq-higher15}) by the change of coordinates (\ref{eq-transform2-1}) and using (\ref{eq-matrix-M}):
\begin{equation}\label{eq-higher17}
F_{\tau^1}^{(2)}(0)=\frac{n-1}{\mathsf{a}(n,\mathbf{d})}F^{(2)}(0)-\sum_{b=2}^{n}\frac{\partial \tau^b}{\partial t^1}F_{\tau^b}^{(2)}(0)
=\frac{n-1}{\mathsf{a}(n,\mathbf{d})}F^{(2)}(0).
\end{equation}

By (\ref{eq-higher5}),
\begin{eqnarray}\label{eq-higher18}
F_{\tau^e}^{(1)}(0)\eta^{ef}F_{\tau^f}^{(2)}(0)+F^{(2)}(0)F^{(2)}(0)=0,
\end{eqnarray}
which yields from(\ref{eq-cor-F(1)-2})
\begin{eqnarray}\label{eq-higher19}
\eta^{0f}F_{\tau^f}^{(2)}(0)+F^{(2)}(0)F^{(2)}(0)=0.
\end{eqnarray}
By (\ref{eq-higher21}), for $2\leq b\leq n$,
\begin{equation}\label{eq-higher21-1}
F_{\tau^b}^{(2)}(0)
=-2F_{\tau^1 \tau^{b-1}}^{(1)}(0)F^{(2)}(0)+\sum_{e=0}^n \sum_{f=0}^n F_{\tau^1 \tau^e}^{(1)}(0)\eta^{ef}F_{\tau^f \tau^{b-1}}^{(1)}(0).
\end{equation}
Substituting (\ref{eq-higher17}) and (\ref{eq-higher21-1}) into (\ref{eq-higher19}), we get (\ref{eq-higher20}).
\end{proof}
In the following we suppose $\frac{n-1}{\mathsf{a}(n,\mathbf{d})}\in \mathbb{Z}$, and  treat the cases $\frac{n-1}{\mathsf{a}(n,\mathbf{d})}=1$ and $\frac{n-1}{\mathsf{a}(n,\mathbf{d})}>1$ separately.
First we recall, from (\ref{eq-pairing2}),
\begin{equation}\label{eq-eta0i}
\eta^{0i}=
\begin{cases}
\frac{1}{\prod_{i=1}^{r}d_i}, & \mbox{if}\ i=n;\\
-\frac{\sfbd}{\prod_{i=1}^{r}d_i},  & \mbox{if}\ i=n-\mathsf{a}(n,\mathbf{d});\\
0,& \mbox{otherwise},
\end{cases}
\end{equation}
and
\begin{eqnarray}\label{eq-eta01}
\frac{n-1}{\mathsf{a}(n,\mathbf{d})}\eta^{01}=\left\{
\begin{array}{cc}
-\frac{\sfbd}{\prod_{i=1}^{r}d_i}, & \frac{n-1}{\mathsf{a}(n,\mathbf{d})}=1,\\
0, & \frac{n-1}{\mathsf{a}(n,\mathbf{d})}>1.
\end{array}\right.
\end{eqnarray}

\subsubsection{The case \texorpdfstring{$\frac{n-1}{\mathsf{a}(n,\mathbf{d})}=1$}{(n-1)/a(n,d)=1}}
In this case, by (\ref{eq-eta0i}) and  (\ref{eq-thm-F(1)}) we have
\begin{eqnarray}\label{eq-F2-case1-1}
\sum_{b=2}^{n}\eta^{0b}F_{\tau^1 \tau^{b-1}}^{(1)}(0)
=\frac{1}{\prod_{i=1}^{r}d_i}\cdot F_{\tau^1 \tau^{n-1}}^{(1)}(0)
=-\frac{c(n,1,\mathbf{d})b(\mathbf{d})}{\prod_{i=1}^{r}d_i},
\end{eqnarray}
and
\begin{eqnarray}\label{eq-F2-case1-2}
&&\sum_{b=2}^{n}\eta^{0b}F_{\tau^1 \tau^e}^{(1)}(0)\eta^{ef}F_{\tau^f \tau^{b-1}}^{(1)}(0)
= \frac{1}{\prod_{i=1}^{r}d_i}\cdot F_{\tau^1 \tau^e}^{(1)}(0)\eta^{ef}F_{\tau^f \tau^{n-1}}^{(1)}(0)\nn\\
&=&  \frac{1}{(\prod_{i=1}^{r}d_i)^2}\Big(
F_{\tau^1 \tau^{n-1}}^{(1)}(0)\Big)^2
=  \frac{c(n,1,\mathbf{d})^2b(\mathbf{d})^2}{(\prod_{i=1}^{r}d_i)^2}.
\end{eqnarray}
Putting (\ref{eq-eta01}), (\ref{eq-F2-case1-1}) and (\ref{eq-F2-case1-2}) into (\ref{eq-higher20}) we get
\begin{equation}\label{eq-F2-case1}
F^{(2)}(0)F^{(2)}(0)+\frac{(2c(n,1,\mathbf{d})-1)\sfbd}{\prod_{i=1}^{r}d_i}F^{(2)}(0)+ \frac{c (n,1,\mathbf{d})^2b(\mathbf{d})^2}{(\prod_{i=1}^{r}d_i)^2}=0.
\end{equation}
The case $\frac{n-1}{\mathsf{a}(n,\mathbf{d})}=1$ happens if and only if $\mathbf{d}=(2,2)$ or $(3)$.
For $X_n(3)$,
\[
M_{n}^1=-21,\ \sfbd=27,
\]
By (\ref{eq-matrix-M}), (\ref{eq-matrix-W}) and the definition (\ref{eq-def-functionC}), we compute
\begin{eqnarray*}
\mathsf{c}(n,1,3)&=&1+\sum_{i=1}^{n}\sum_{j=0}^n 
\frac{j-i}{n-1} M_{j}^{i}W_{n}^{j}27^{\frac{i-n}{n-1}}
-\sum_{i=1}^{n}\sum_{j=0}^n 
\frac{j-i}{n-1} M_{j}^{i}W_{1}^{j}27^{1+\frac{i-n}{n-1}}\\
&=&1+\frac{M_n^1}{27}=\frac{2}{9}.
\end{eqnarray*}
Hence (\ref{eq-F2-case1}) reads
\begin{eqnarray}\label{eq-F2-case1-cubic}
F^{(2)}(0)F^{(2)}(0)-5F^{(2)}(0)+4=0,
\end{eqnarray}
so $F^{(2)}(0)=1$ or $4$. For $X_n(2,2)$, where $n$ is odd,
\[
M_{n}^1=-12,\ \sfbd=16,
\]
By (\ref{eq-matrix-M}), (\ref{eq-matrix-W}) and the definition (\ref{eq-def-functionC}), we compute
\begin{eqnarray*}
&&\mathsf{c}(n,1,(2,2))\\
&=& 1+\sum_{i=1}^{n}\sum_{j=0}^n 
\frac{j-i}{n-1} M_{j}^{i}W_{n}^{j}16^{\frac{i-n}{n-1}}
-\sum_{i=1}^{n}\sum_{j=0}^n 
\frac{j-i}{n-1} M_{j}^{i}W_{1}^{j}16^{1+\frac{i-n}{n-1}}\\
&=&1+\frac{M_n^1}{16}=\frac{1}{4}.
\end{eqnarray*}
Hence (\ref{eq-F2-case1}) reads
\begin{eqnarray}\label{eq-F2-case1-(2,2)}
F^{(2)}(0)F^{(2)}(0)-2F^{(2)}(0)+1=0,
\end{eqnarray}
and thus $F^{(2)}(0)=1$.\\

\subsubsection{The case  \texorpdfstring{$\frac{n-1}{\mathsf{a}(n,\mathbf{d})}>1$}{(n-1)/a(n,d)>1}}
In this case by (\ref{eq-eta0i}) and  (\ref{eq-thm-F(1)}) (resp. by (\ref{eq-thm-F(1)-fanoIndex=1}) when $\mathsf{a}(n,\mathbf{d})=1$) we compute
\begin{eqnarray}\label{eq-F2-case2-1}
&&\sum_{b=2}^{n}\eta^{0b}F_{\tau^1 \tau^{b-1}}^{(1)}(0)\nn\\
&=&\frac{1}{\prod_{i=1}^{r}d_i}\Big(-\sfbd F_{\tau^1 \tau^{n-\mathsf{a}(n,\mathbf{d})-1}}^{(1)}(0)
+F_{\tau^1 \tau^{n-1}}^{(1)}(0)\Big)\nn\\
&=&\frac{1}{\prod_{i=1}^{r}d_i}\Big(c(n,\frac{n-1}{\mathsf{a}(n,\mathbf{d})}-1,\mathbf{d})\sfbd\cdot (\sfbd)^{\frac{n-1}{\mathsf{a}(n,\mathbf{d})}-1}
-c(n,\frac{n-1}{\mathsf{a}(n,\mathbf{d})},\mathbf{d})(\sfbd)^{\frac{n-1}{\mathsf{a}(n,\mathbf{d})}}\Big)\nn\\
&=&\frac{(\sfbd)^{\frac{n-1}{\mathsf{a}(n,\mathbf{d})}}}{\prod_{i=1}^{r}d_i}\Big(c(n,\frac{n-1}{\mathsf{a}(n,\mathbf{d})}-1,\mathbf{d})
-c(n,\frac{n-1}{\mathsf{a}(n,\mathbf{d})},\mathbf{d})\Big).
\end{eqnarray}

Then we compute the last group of terms in  (\ref{eq-higher20}).
Note that $\frac{n-1}{\mathsf{a}(n,\mathbf{d})}>1$  implies $n\geq 3$.
Using (\ref{eq-pairing2}) and (\ref{eq-thm-F(1)}) (resp. by (\ref{eq-thm-F(1)-fanoIndex=1}) when $\mathsf{a}(n,\mathbf{d})=1$) repeatedly, we compute
\begin{eqnarray}\label{eq-F2-case2-4}
&&\sum_{b=2}^{n}\eta^{0b}F_{\tau^1 \tau^e}^{(1)}(0)\eta^{ef}F_{\tau^f \tau^{b-1}}^{(1)}(0)\nn\\
&=& \frac{1}{\prod_{i=1}^{r}d_i}\Big(-\sfbd F_{\tau^1 \tau^e}^{(1)}(0)\eta^{ef}F_{\tau^f \tau^{n-\mathsf{a}(n,\mathbf{d})-1}}^{(1)}(0)
+F_{\tau^1 \tau^e}^{(1)}(0)\eta^{ef}F_{\tau^f \tau^{n-1}}^{(1)}(0)\Big)\nn\\
&=& \frac{1}{\prod_{i=1}^{r}d_i}\Big(-\sfbd\sum_{k=1}^{\frac{n-1}{\mathsf{a}(n,\mathbf{d})}}F_{\tau^1 \tau^{k\mathsf{a}(n,\mathbf{d})}}^{(1)}(0)\eta^{k\mathsf{a}(n,\mathbf{d}),f}F_{\tau^f \tau^{n-\mathsf{a}(n,\mathbf{d})-1}}^{(1)}(0)\nn\\
&&+\sum_{k=1}^{\frac{n-1}{\mathsf{a}(n,\mathbf{d})}}F_{\tau^1 \tau^{k\mathsf{a}(n,\mathbf{d})}}^{(1)}(0)\eta^{k\mathsf{a}(n,\mathbf{d}),f}F_{\tau^f \tau^{n-1}}^{(1)}(0)\Big)\nn\\
&=&  \frac{1}{(\prod_{i=1}^{r}d_i)^2}
\Big(\sfbd
\sum_{k=1}^{\frac{n-1}{\mathsf{a}(n,\mathbf{d})}}c (n,k,\mathbf{d})(\sfbd)^{k}F_{\tau^{n-k\mathsf{a}(n,\mathbf{d})} \tau^{n-\mathsf{a}(n,\mathbf{d})-1}}^{(1)}(0)\nn\\
&&-(\sfbd)^2
\sum_{k=1}^{\frac{n-1}{\mathsf{a}(n,\mathbf{d})}-1}c (n,k,\mathbf{d})(\sfbd)^{k}F_{\tau^{n-(k+1)\mathsf{a}(n,\mathbf{d})} \tau^{n-\mathsf{a}(n,\mathbf{d})-1}}^{(1)}(0)\nn\\
&&-\sum_{k=1}^{\frac{n-1}{\mathsf{a}(n,\mathbf{d})}}c (n,k,\mathbf{d})(\sfbd)^{k}F_{\tau^{n-k\mathsf{a}(n,\mathbf{d})} \tau^{n-1}}^{(1)}(0)\nn\\
&&+\sfbd
\sum_{k=1}^{\frac{n-1}{\mathsf{a}(n,\mathbf{d})}-1}c (n,k,\mathbf{d})(\sfbd)^{k}F_{\tau^{n-(k+1)\mathsf{a}(n,\mathbf{d})} \tau^{n-1}}^{(1)}(0)\Big)\nn\\
&=&  \frac{1}{(\prod_{i=1}^{r}d_i)^2}\Big(-\sfbd
\sum_{k=1}^{\frac{n-1}{\mathsf{a}(n,\mathbf{d})}}c(n,k,\mathbf{d})c(n,\frac{2(n-1)}{\mathsf{a}(n,\mathbf{d})}-k-1,\mathbf{d})
(\sfbd)^{k}(\sfbd)^{\frac{2(n-1)}{\mathsf{a}(n,\mathbf{d})}-k-1}\nn\\
&&+(\sfbd)^2
\sum_{k=1}^{\frac{n-1}{\mathsf{a}(n,\mathbf{d})}-1}c(n,k,\mathbf{d})c(n,\frac{2(n-1)}{\mathsf{a}(n,\mathbf{d})}-k-2,\mathbf{d})(\sfbd)^{k}(\sfbd)^{\frac{2(n-1)}{\mathsf{a}(n,\mathbf{d})}-k-2}\nn\\
&&+\sum_{k=1}^{\frac{n-1}{\mathsf{a}(n,\mathbf{d})}}c(n,k,\mathbf{d})c(n,\frac{2(n-1)}{\mathsf{a}(n,\mathbf{d})}-k,\mathbf{d})(\sfbd)^{k}(\sfbd)^{\frac{2(n-1)}{\mathsf{a}(n,\mathbf{d})}-k}\nn\\
&&-\sfbd
\sum_{k=1}^{\frac{n-1}{\mathsf{a}(n,\mathbf{d})}-1}c(n,k,\mathbf{d})c(n,\frac{2(n-1)}{\mathsf{a}(n,\mathbf{d})}-k-1,\mathbf{d})(\sfbd)^{k}(\sfbd)^{\frac{2(n-1)}{\mathsf{a}(n,\mathbf{d})}-k-1}\Big)\nn\\
&=& \frac{(\sfbd)^{\frac{2(n-1)}{\mathsf{a}(n,\mathbf{d})}}}{(\prod_{i=1}^{r}d_i)^2}\Big(-\sum_{k=1}^{\frac{n-1}{\mathsf{a}(n,\mathbf{d})}}c(n,k,\mathbf{d})c(n,\frac{2(n-1)}{\mathsf{a}(n,\mathbf{d})}-k-1,\mathbf{d})\nn\\
&&+
\sum_{k=1}^{\frac{n-1}{\mathsf{a}(n,\mathbf{d})}-1}c(n,k,\mathbf{d})c(n,\frac{2(n-1)}{\mathsf{a}(n,\mathbf{d})}-k-2,\mathbf{d})\nn\\
&&+\sum_{k=1}^{\frac{n-1}{\mathsf{a}(n,\mathbf{d})}}c(n,k,\mathbf{d})c(n,\frac{2(n-1)}{\mathsf{a}(n,\mathbf{d})}-k,\mathbf{d})
-\sum_{k=1}^{\frac{n-1}{\mathsf{a}(n,\mathbf{d})}-1}c(n,k,\mathbf{d})c(n,\frac{2(n-1)}{\mathsf{a}(n,\mathbf{d})}-k-1,\mathbf{d})\Big).\nn\\
\end{eqnarray}
Let $l=\frac{n-1}{\mathsf{a}(n,\mathbf{d})}$. By the definition (\ref{eq-def-functionC}) we have
\begin{equation}\label{eq-c(n,k,d)-stablization}
	c(n,k,\mathbf{d})=c(n,l,\mathbf{d})\ \mbox{for}\ k\geq l.
\end{equation}
 So we get, when $l\geq 3$, 
\begin{eqnarray*}
&&-\sum_{k=1}^{\frac{n-1}{\mathsf{a}(n,\mathbf{d})}}c(n,k,\mathbf{d})c(n,\frac{2(n-1)}{\mathsf{a}(n,\mathbf{d})}-k-1,\mathbf{d})
+
\sum_{k=1}^{\frac{n-1}{\mathsf{a}(n,\mathbf{d})}-1}c(n,k,\mathbf{d})c(n,\frac{2(n-1)}{\mathsf{a}(n,\mathbf{d})}-k-2,\mathbf{d})\nn\\
&&+\sum_{k=1}^{\frac{n-1}{\mathsf{a}(n,\mathbf{d})}}c(n,k,\mathbf{d})c(n,\frac{2(n-1)}{\mathsf{a}(n,\mathbf{d})}-k,\mathbf{d})
-\sum_{k=1}^{\frac{n-1}{\mathsf{a}(n,\mathbf{d})}-1}c(n,k,\mathbf{d})c(n,\frac{2(n-1)}{\mathsf{a}(n,\mathbf{d})}-k-1,\mathbf{d})\\
&=&-\sum_{k=1}^{l-2}c(n,k,\mathbf{d})c(n,2l-k-1,\mathbf{d})
-c(n,l-1,\mathbf{d})c(n,l,\mathbf{d})-c(n,l,\mathbf{d})c(n,l-1,\mathbf{d})\\
&&+\sum_{k=1}^{l-3}c(n,k,\mathbf{d})c(n,2l-k-2,\mathbf{d})
+c(n,l-2,\mathbf{d})c(n,l,\mathbf{d})+c(n,l-1,\mathbf{d})c(n,l-1,\mathbf{d})\\
&&+\sum_{k=1}^{l-1}c(n,k,\mathbf{d})c(n,2l-k,\mathbf{d})
+c(n,l,\mathbf{d})c(n,l,\mathbf{d})\\
&&-\sum_{k=1}^{l-2}c(n,k,\mathbf{d})c(n,2l-k-1,\mathbf{d})
-c(n,l-1,\mathbf{d})c(n,l,\mathbf{d})\\
&\stackrel{\mbox{by (\ref{eq-c(n,k,d)-stablization})}}{=}&-\sum_{k=1}^{l-2}c(n,k,\mathbf{d})c(n,l,\mathbf{d})
-c(n,l-1,\mathbf{d})c(n,l,\mathbf{d})-c(n,l,\mathbf{d})c(n,l-1,\mathbf{d})\\
&&+\sum_{k=1}^{l-3}c(n,k,\mathbf{d})c(n,l,\mathbf{d})
+c(n,l-2,\mathbf{d})c(n,l,\mathbf{d})+c(n,l-1,\mathbf{d})c(n,l-1,\mathbf{d})\\
&&+\sum_{k=1}^{l-1}c(n,k,\mathbf{d})c(n,l,\mathbf{d})
+c(n,l,\mathbf{d})c(n,l,\mathbf{d})\\
&&-\sum_{k=1}^{l-2}c(n,k,\mathbf{d})c(n,l,\mathbf{d})
-c(n,l-1,\mathbf{d})c(n,l,\mathbf{d})\\
&=&-c(n,l-2,\mathbf{d})c(n,l,\mathbf{d})
-c(n,l-1,\mathbf{d})c(n,l,\mathbf{d})-c(n,l,\mathbf{d})c(n,l-1,\mathbf{d})\\
&&
+c(n,l-2,\mathbf{d})c(n,l,\mathbf{d})+c(n,l-1,\mathbf{d})c(n,l-1,\mathbf{d})\\
&&+c(n,l-1,\mathbf{d})c(n,l,\mathbf{d})
+c(n,l,\mathbf{d})c(n,l,\mathbf{d})\\
&&-c(n,l-1,\mathbf{d})c(n,l,\mathbf{d})\\
&=& \big(c(n,l-1,\mathbf{d})-c(n,l,\mathbf{d})\big)^2.
\end{eqnarray*}
When $l=2$, we have similarly
\begin{eqnarray*}
&&-\sum_{k=1}^{\frac{n-1}{\mathsf{a}(n,\mathbf{d})}}c(n,k,\mathbf{d})c(n,\frac{2(n-1)}{\mathsf{a}(n,\mathbf{d})}-k-1,\mathbf{d})
+
\sum_{k=1}^{\frac{n-1}{\mathsf{a}(n,\mathbf{d})}-1}c(n,k,\mathbf{d})c(n,\frac{2(n-1)}{\mathsf{a}(n,\mathbf{d})}-k-2,\mathbf{d})\nn\\
&&+\sum_{k=1}^{\frac{n-1}{\mathsf{a}(n,\mathbf{d})}}c(n,k,\mathbf{d})c(n,\frac{2(n-1)}{\mathsf{a}(n,\mathbf{d})}-k,\mathbf{d})
-\sum_{k=1}^{\frac{n-1}{\mathsf{a}(n,\mathbf{d})}-1}c(n,k,\mathbf{d})c(n,\frac{2(n-1)}{\mathsf{a}(n,\mathbf{d})}-k-1,\mathbf{d})\\
&=&-c(n,1,\mathbf{d})c(n,2,\mathbf{d})-c(n,2,\mathbf{d})c(n,1,\mathbf{d})
+c(n,1,\mathbf{d})c(n,1,\mathbf{d})\\
&&+c(n,1,\mathbf{d})c(n,3,\mathbf{d})
+c(n,2,\mathbf{d})c(n,2,\mathbf{d})-c(n,1,\mathbf{d})c(n,2,\mathbf{d})\\
&\stackrel{\mbox{by (\ref{eq-c(n,k,d)-stablization})}}{=}&-3c(n,1,\mathbf{d})c(n,2,\mathbf{d})+c(n,1,\mathbf{d})^2+c(n,1,\mathbf{d})c(n,2,\mathbf{d})
+c(n,2,\mathbf{d})^2\\
&=& \big(c(n,1,\mathbf{d})-c(n,2,\mathbf{d})\big)^2.
\end{eqnarray*}
Therefore
\begin{equation}\label{eq-F2-case2-5}
\sum_{b=2}^{n}\eta^{0b}F_{\tau^1 \tau^e}^{(1)}(0)\eta^{ef}F_{\tau^f \tau^{b-1}}^{(1)}(0)
=\frac{(\sfbd)^{\frac{2(n-1)}{\mathsf{a}(n,\mathbf{d})}}}{(\prod_{i=1}^{r}d_i)^2}\big(c(n,l-1,\mathbf{d})-c(n,l,\mathbf{d})\big)^2.
\end{equation}

Putting (\ref{eq-eta01}), (\ref{eq-F2-case2-1}) and (\ref{eq-F2-case2-5}) together into
 (\ref{eq-higher20}), we get
\begin{equation}\label{eq-F2-case2-final}
\Big(F^{(2)}(0)-\frac{(\sfbd)^{\frac{n-1}{\mathsf{a}(n,\mathbf{d})}}}{\prod_{i=1}^{r}d_i}\big(c(n,l-1,\mathbf{d})-c(n,l,\mathbf{d})\big)\Big)^2=0.
\end{equation}

\subsubsection{The formula of \texorpdfstring{$F^{(2)}(0)$}{F(2)(0)}}
By the definition (\ref{eq-def-functionC}),
\begin{eqnarray}\label{eq-F2-case2-2}
&&c(n,\frac{n-1}{\mathsf{a}(n,\mathbf{d})}-1,\mathbf{d})
-c(n,\frac{n-1}{\mathsf{a}(n,\mathbf{d})},\mathbf{d})\nn\\
&=&
-\sum_{j=0}^n 
\frac{j-\big(n-\frac{n-1}{\mathsf{a}(n,\mathbf{d})}\cdot\mathsf{a}(n,\mathbf{d})\big)}{\mathsf{a}(n,\mathbf{d})} M_{j}^{n-\frac{n-1}{\mathsf{a}(n,\mathbf{d})}\cdot\mathsf{a}(n,\mathbf{d})}W_{n}^{j}(\sfbd)^{\frac{n-\frac{n-1}{\mathsf{a}(n,\mathbf{d})}\cdot\mathsf{a}(n,\mathbf{d})-n}{\mathsf{a}(n,\mathbf{d})}}  \nn\\
&&+\sum_{j=0}^n 
\frac{j-\big(n-\frac{n-1}{\mathsf{a}(n,\mathbf{d})}\cdot\mathsf{a}(n,\mathbf{d})\big)}{\mathsf{a}(n,\mathbf{d})} M_{j}^{n-\frac{n-1}{\mathsf{a}(n,\mathbf{d})}\cdot\mathsf{a}(n,\mathbf{d})}W_{n- \mathsf{a}(n,\mathbf{d})}^{j}(\sfbd)^{1+\frac{n-\frac{n-1}{\mathsf{a}(n,\mathbf{d})}\cdot\mathsf{a}(n,\mathbf{d})-n}{\mathsf{a}(n,\mathbf{d})}}\nn\\
&=&-\sum_{j=0}^n 
\frac{j-1}{\mathsf{a}(n,\mathbf{d})} M_{j}^{1}W_{n}^{j}(\sfbd)^{\frac{1-n}{\mathsf{a}(n,\mathbf{d})}}
+\sum_{j=0}^n 
\frac{j-1}{\mathsf{a}(n,\mathbf{d})} M_{j}^{1}W_{n- \mathsf{a}(n,\mathbf{d})}^{j}(\sfbd)^{1+\frac{1-n}{\mathsf{a}(n,\mathbf{d})}}\nn\\
&=&-\frac{(\sfbd)^{\frac{1-n}{\mathsf{a}(n,\mathbf{d})}}}{\mathsf{a}(n,\mathbf{d})}\sum_{j=0}^n j M_{j}^{1}W_{n}^{j}
+\frac{(\sfbd)^{1+\frac{1-n}{\mathsf{a}(n,\mathbf{d})}}}{\mathsf{a}(n,\mathbf{d})}\sum_{j=0}^n 
j M_{j}^{1}W_{n- \mathsf{a}(n,\mathbf{d})}^{j}\nn\\
&&+\frac{(\sfbd)^{\frac{1-n}{\mathsf{a}(n,\mathbf{d})}}}{\mathsf{a}(n,\mathbf{d})}\sum_{j=0}^n M_{j}^{1}W_{n}^{j}
-\frac{(\sfbd)^{1+\frac{1-n}{\mathsf{a}(n,\mathbf{d})}}}{\mathsf{a}(n,\mathbf{d})}\sum_{j=0}^n M_{j}^{1}W_{n- \mathsf{a}(n,\mathbf{d})}^{j}\nn\\
&=&-\frac{(\sfbd)^{\frac{1-n}{\mathsf{a}(n,\mathbf{d})}}}{\mathsf{a}(n,\mathbf{d})}\sum_{j=0}^n j M_{j}^{1}W_{n}^{j}
+\frac{(\sfbd)^{1+\frac{1-n}{\mathsf{a}(n,\mathbf{d})}}}{\mathsf{a}(n,\mathbf{d})}\sum_{j=0}^n 
j M_{j}^{1}W_{n- \mathsf{a}(n,\mathbf{d})}^{j},
\end{eqnarray}
where in the last step we use that $M\cdot W=\mathrm{Id}$ and $n-1-\mathsf{a}(n,\mathbf{d})>0$. So we obtain
\begin{eqnarray}\label{eq-F2-case2-3}
&&\frac{(\sfbd)^{\frac{n-1}{\mathsf{a}(n,\mathbf{d})}}}{\prod_{i=1}^{r}d_i}\big(c(n,l-1,\mathbf{d})-c(n,l,\mathbf{d})\big)\nn\\
&=&\frac{1}{\mathsf{a}(n,\mathbf{d})\prod_{i=1}^r d_i}
	\Big(-\sum_{j=0}^n j M_{j}^{1}W_{n}^{j}
+\sfbd\sum_{j=0}^n 
j M_{j}^{1}W_{n- \mathsf{a}(n,\mathbf{d})}^{j}\Big).
\end{eqnarray}

Summarizing (\ref{eq-F2-case1-cubic}), (\ref{eq-F2-case1-(2,2)}) and (\ref{eq-F2-case2-final}) we obtain the following:
\begin{theorem}\label{thm-higher10.1}
Let $X_n(\mathbf{d})$ be a non-exceptional complete intersection.
Then
\begin{eqnarray}\label{eq-higher22}
F^{(2)}(0)
\begin{dcases}
1, & \mathrm{if}\  \mathbf{d}=(2,2);\\
1\ \mathrm{or}\ 4, & \mathrm{if}\  \mathbf{d}= (3);\\
\frac{-\sum_{j=0}^n j M_{j}^{1}W_{n}^{j}
+\sfbd\sum_{j=0}^n 
j M_{j}^{1}W_{n- \mathsf{a}(n,\mathbf{d})}^{j}}{\mathsf{a}(n,\mathbf{d})\prod_{i=1}^r d_i}
	\
, & \mathrm{if}\ l=\frac{n-1}{\mathsf{a}(n,\mathbf{d})}\in \mathbb{Z}_{\geq 2};\\
0, & \mathrm{otherwise}.
\end{dcases}
\end{eqnarray}
\end{theorem}

We need to take an ad hoc way to determine $F^{(2)}(0)$ of the cubic hypersurfaces. 
\begin{theorem}\label{thm-F(2)(0)-cubicHypersurfaces}
For smooth cubic hypersurfaces, $F^{(2)}(0)=1$. 
\end{theorem}
We will give two proofs, respectively in Section \ref{sec:FanoVarietyOfLines} and Section \ref{sec:4points-Invariants}. 
For later use in Section \ref{sec:8points-Invariants-Cubic3Fold}, we record the following consequence of Theorem \ref{thm-F(2)(0)-cubicHypersurfaces}.
\begin{corollary}\label{cor-lowerTerms-Cubic-F(2)}
Let $X=X_n(3)$, a cubic hypersurface of dimension $n$. Then
\begin{equation}\label{eq-lowerTerms-Cubic-F(2)}
F^{(2)}=1+t^1+3t^n+O((t)^2).
\end{equation}
\end{corollary}
\begin{proof}
By (\ref{eq-higher17}), 
\begin{equation}\label{eq-lowerTerms-Cubic-F(2)-1}
	F_{\tau^1}^{(2)}(0)=F^{(2)}(0)=1.
\end{equation}
Then by (\ref{eq-higher21-1}) and (\ref{eq-thm-F(1)}), 
\begin{eqnarray*}
F_{\tau^n}^{(2)}(0)
&=&-2F_{\tau^1 \tau^{n-1}}^{(1)}(0)F^{(2)}(0)+ F_{\tau^1 \tau^{n-1}}^{(1)}(0)\eta^{n-1,1}F_{\tau^1 \tau^{n-1}}^{(1)}(0)\\
&=& 2 \mathsf{c}(n,1,3) \times 27+\frac{\mathsf{c}(n,1,3)^2\times 27^2}{3}.
\end{eqnarray*}
By (\ref{eq-c(n,1,d)-largeFanoIndex}), $\mathsf{c}(n,1,3)=\frac{2}{9}$. So 
\begin{equation}\label{eq-lowerTerms-Cubic-F(2)-2}
	F_{\tau^n}^{(2)}(0)=24.
\end{equation}
 By (\ref{eq-transform2-1}) and (\ref{eq-Beauville-ell-l0-l1}), 
\[
\tau^1=t^1-21 t^n,\ \tau^n=t^n.
\]
So from (\ref{eq-lowerTerms-Cubic-F(2)-1}) and (\ref{eq-lowerTerms-Cubic-F(2)-2}) we get (\ref{eq-lowerTerms-Cubic-F(2)}).
\end{proof}

In the case $\frac{n-1}{\mathsf{a}(n,\mathbf{d})}\in \mathbb{Z}_{\geq 2}$, one can use the mirror formula (\ref{eq-mirrorFormula-a(n,d)>1}) for the small $J$-functions to compute the matrices $W$ and $M$, and thus evaluate $F^{(2)}(0)$ by (\ref{eq-higher22}). We illustrate this in Appendix \ref{sec:computation-twoPointInvariants-from-smallJ}. The larger $\frac{n-1}{\mathsf{a}(n,\mathbf{d})}$ is, the more complicated the computation is. At this stage we only have a closed formula when $\frac{n-1}{\mathsf{a}(n,\mathbf{d})}=2$.
\begin{theorem}\label{thm-F20-a(n,d)=(n-1)/2}
Suppose $\frac{n-1}{\mathsf{a}(n,\mathbf{d})}=2$. Then
\begin{equation*}
	F^{(2)}(0)=\frac{\prod_{i=1}^r d_i!(d_i-1)!}{2}.
\end{equation*}
\end{theorem}
The proof is given in Appendix \ref{sec:proof-F2(0)-closedFormula}.

\subsection{Reconstruction theorem II for \texorpdfstring{$X_n(3)$}{Xn(3)}  and \texorpdfstring{$X_n(2,2)$}{Xn(2)}}\label{sec:reconstruction-FanoIndex=(n-1)}

Once we have computed $F^{(2)}(0)$, we can try to use these  equations to solve $F^{(k)}(0)$ for $k\geq 3$ 
from lower order terms (in the order defined in Definition \ref{def-order-gwInvariant}). In the following theorem, we will  use Theorem  \ref{thm-F(2)(0)-cubicHypersurfaces} for cubic hypersurfaces, which will be proved in Section \ref{sec:FanoVarietyOfLines} and Section \ref{sec:genusReduction}.

\begin{theorem}\phantomsection\label{thm-reconstructcubicandquadric}
\begin{itemize}
\item[(i)] For the cubic threefold $X$, $F$ can be reconstructed by (\ref{eq-system1-even}) when $n$ is even (resp. (\ref{eq-system1-odd}) when $n$ is odd) and (\ref{eq-Dim}) from $F^{(0)}$ and  $F^{(2)}(0)$, $F^{(4)}(0)$.
\item[(ii)] For cubic hypersurfaces $X$ with $\dim X\geq 4$, $F$ can be reconstructed by (\ref{eq-system1-even}) when $n$ is even (resp. (\ref{eq-system1-odd}) when $n$ is odd) and (\ref{eq-Dim}) from $F^{(0)}$ and  $F^{(2)}(0)$.
\item[(iii)] For odd-dimensional intersections $X$ of two quadrics, with $\dim X> 2$, $F$ can be reconstructed by (\ref{eq-system1-even}) when $n$ is even (resp. (\ref{eq-system1-odd}) when $n$ is odd) and (\ref{eq-Dim}) from $F^{(0)}$.
\end{itemize}
\end{theorem}
\begin{proof} In this proof we use the symbol $\sim$ to denote that the two sides differ  by  invariants of lower orders in the sense of Definition \ref{def-order-gwInvariant}. 

Let $t^0,\dots,t^n$ be the basis dual to $1,\sfh,\dots,\sfh_n$. Denote $F_e=\frac{\partial F}{\partial t^e}$.  
For $k\geq 2$, by (\ref{eq-higher25}) ,
\begin{eqnarray*}
F_{ e}^{(1)}g^{ef}F_{ f}^{(k+1)}+2kF^{(2)}F^{(k+1)}\sim 0.
\end{eqnarray*}
By (\ref{eq-F1-leadingTerm-FanoIndex=n-1}) and  Theorem \ref{thm-F(2)(0)-cubicHypersurfaces}, we get
\begin{eqnarray}\label{eq-thm-reconstructcubicandquadric-1}
\frac{1}{3}F_{n}^{(k+1)}(0)-2F_{1}^{(k+1)}(0)+2kF^{(k+1)}(0)\sim 0.
\end{eqnarray}
By (\ref{eq-Dim}), $F^{(k+1)}(0)\neq 0$ forces
\begin{eqnarray*}
2(k+1)\cdot \frac{n}{2}=n-3+2(k+1)+\beta\cdot (n-1),
\end{eqnarray*}
i.e.,
\begin{eqnarray*}
\beta=k-\frac{k-1}{n-1}.
\end{eqnarray*}
Thus by the divisor equation (see the observation below (\ref{eq-EV-0}))
\begin{eqnarray}\label{eq-thm-reconstructcubicandquadric-2}
F_{1}^{(k+1)}(0)=\big(k-\frac{k-1}{n-1}\big)F^{(k+1)}(0).
\end{eqnarray}
By (\ref{eq-M-FanoIndex=n-1}), 
\begin{eqnarray}\label{eq-smallQuantumProduct-h}
	&&\sfh\sqp \sfh_{n-1}=\sfh\sqp(\tsfh_{n-1}-\elld)
	=\tsfh_n-\elld\sfh\nn\\
	&&=\sfh_n-\big(\elld-\sfbd\big)\sfh-\elld\sfh
	=\sfh_n+\big(\sfbd-2\elld\big)\sfh.
\end{eqnarray}
By (\ref{eq-smallQuantumProduct-h})  one finds that the leading terms of LHS of (\ref{eq-wdvv23expand}) (\ref{eq-higher24}) is
\begin{eqnarray}\label{eq-thm-reconstructcubicandquadric-leadingTerms}
&&F_{n}^{(k+1)}(0)+\big(\sfbd-2\elld\big)F_{1}^{(k+1)}(0)
+2k F^{(1)}_{1,n-1}(0) F^{(k+1)}(0)\nn\\
&&-F_{1}^{(k+1)}(0) F_{n-1}^{(1)}(0)-F_{1}^{(1)}(0)F_{n-1}^{(k+1)}(0).
\end{eqnarray}
Then by (\ref{eq-F1-leadingTerm-FanoIndex=n-1}) we get
\begin{eqnarray}\label{eq-thm-reconstructcubicandquadric-3}
F_{n}^{(k+1)}(0)\sim 12k F^{(k+1)}(0)-21F_{1}^{(k+1)}(0).
\end{eqnarray}
Inserting  (\ref{eq-thm-reconstructcubicandquadric-2}) and (\ref{eq-thm-reconstructcubicandquadric-3}) into (\ref{eq-thm-reconstructcubicandquadric-1}), we obtain
\begin{eqnarray*}
0&\sim& \frac{1}{3}\Big(12 k F^{(k+1)}(0)-21F_{1}^{(k+1)}(0)\Big)-2F_{1}^{(k+1)}(0)+2kF^{(k+1)}(0)\\
&\sim& 6k F^{(k+1)}(0)-9F_{1}^{(k+1)}(0)=\Big(\frac{9(k-1)}{n-1}-3k\Big)F^{(k+1)}(0).
\end{eqnarray*}
The coefficient $\frac{9(k-1)}{n-1}-3k=0$ if and only if $n=3$ and $k=4$. Combining with Theorem \ref{thm-reconstruction-I}, we complete the proof of (i) and (ii).

(iii) By  Theorem \ref{thm-reconstruction-I} we only need to compute $F^{(k)}(0)$ for $k\geq 2$. By (\ref{eq-Dim}), $F^{(k)}(0)=0$ unless
\begin{equation}\label{eq-dimConstraint-odd(2,2)}
\beta(k):=\frac{(n-2)k-n+3}{\mathsf{a}(n,\mathbf{d})}=k-1-\frac{k-2}{n-1}\in \mathbb{Z}.
\end{equation}
 But the dimension $n$ is odd and $\dim H^{n}(X)=n+1$ (see Corollary \ref{cor-dim-prim}). Then by (\ref{eq-Sym}), it follows that $F^{(k)}(0)=0$ for $k>\frac{n+1}{2}$. The only $k$ in this range satisfying (\ref{eq-dimConstraint-odd(2,2)}) is $k=2$. 
By Theorem \ref{thm-higher10.1} We have $F^{(2)}(0)=1$. So we are done.
\end{proof}

\section{Cohomology ring of the Fano varieties of lines in  cubic hypersurfaces}\label{sec:FanoVarietyOfLines}
The computation of the invariants of degree 1 involving primitive classes is closely related to the geometry of the Fano variety of lines. In this section we 
study the cohomology ring of the Fano variety of lines in the cubic hypersurfaces, and in passing obtain $F^{(2)}(0)$.
\subsection{\texorpdfstring{$F^{(2)}(0)$}{F(2)(0)} and Fano varieties of lines}
Let $X=X_n(3)$ be a smooth cubic hypersurface in $\mathbb{P}^{n+1}$. 
For any smooth projective variety $Y$ and $\alpha_{1},\dots, \alpha_k\in H^*(Y)$, we define
\begin{eqnarray*}
(\alpha_{1},\dots, \alpha_k)_Y=\int_Y \alpha_{1}\cup\cdots\cup\alpha_k.
\end{eqnarray*}
Let $\alpha_{1}, \alpha_2,\alpha_3,\alpha_4$ be primitive cohomology classes of  $X$. 
By Proposition \ref{prop-initialValues-meaning}, $F^{(2)}(0)$ is defined via
\begin{eqnarray}\label{eq-Fano1}
&&F^{(2)}(0)\cdot\Big((\alpha_1,\alpha_2)_X\cdot (\alpha_3,\alpha_4)_X+(\alpha_1,\alpha_3)_X\cdot (\alpha_4,\alpha_2)_X+(\alpha_1,\alpha_4)_X\cdot (\alpha_2,\alpha_3)_X\Big)\nn\\
&=&\langle\alpha_1,\alpha_2,\alpha_3,\alpha_4 \rangle_{0,4,1}.
\end{eqnarray}
The \emph{Fano variety of lines} in $X$ is another name of $\overline{\mathcal{M}}_{0,0}(X,1)$. Following  \cite[Chapter 13]{Lew99}, we denote $\overline{\mathcal{M}}_{0,0}(X,1)$ by $\Omega_X$. By \cite[Theorem 7.2]{CG72}, $\Omega_X$ is smooth of  pure dimension $2n-4$. The universal family of lines  is denoted by $P(X)$. 
Let $G_2(\mathbb{C}^{n+2})$ be the Grassmannian parametrizing the 2-dimensional subspaces of $\mathbb{C}^{n+2}$, and let $E$ be the universal rank 2 vector bundle over $G_2(\mathbb{C}^{n+2})$. There is a natural embedding 
\begin{equation*}
	\iota_X:\Omega_X\hookrightarrow G_2(\mathbb{C}^{n+2}).
\end{equation*}
Then $P(X)$ is naturally identified with the projectivization of $\iota_X^* E$. Denote the projection $P(X)\rightarrow \Omega_X$ by $\rho_X$, and the universal morphism $P(X)\rightarrow X$ by $\pi_X$. Let 
\[
\Psi=\rho_{X*} \pi_X^*: H^{n}(X)\rightarrow H^{n-2}(\Omega_X).
\]
The following lemma can be regarded as a variant of the last equation in the proof of \cite[Proposition 5]{Bea95}.
\begin{lemma}\label{lem-4pointInvariant-FanoVarLine}
\begin{equation}\label{eq-4pointInvariant-FanoVarLine}
	\langle\alpha_1,\alpha_2,\alpha_3,\alpha_4 \rangle_{0,4,1}=(\Psi\alpha_1, \Psi\alpha_2,\Psi\alpha_3, \Psi\alpha_4)_{\Omega_X}
\end{equation}
\end{lemma}
\begin{proof}
The forgetful map $\varphi:\overline{\mathcal{M}}_{0,4}(X,1)\rightarrow \overline{\mathcal{M}}_{0,0}(X,1)=\Omega_X$ is flat, with all geometric fibers isomorphic to $\overline{\mathcal{M}}_{0,4}(\mathbb{P}^1,1)$, which is smooth (\cite[Theorem 2]{FP97}). Thus $\varphi$ is smooth  (\cite[Corollaire 17.5.2]{Gro67}). Since $\Omega_X$ is smooth as we have seen above, $\overline{\mathcal{M}}_{0,4}(X,1)$ is also smooth and has the expected dimension. Then $[\overline{\mathcal{M}}_{0,4}(X,1)]^{\mathrm{vir}}=[\overline{\mathcal{M}}_{0,4}(X,1)]$, and
\begin{equation}\label{eq-lem-4pointInvariant-FanoVarLine-1}
	\langle\alpha_1,\alpha_2,\alpha_3,\alpha_4 \rangle_{0,4,1}
	=\int_{[\overline{\mathcal{M}}_{0,4}(X,1)]}\mathrm{ev}_1^* \alpha_1\cup \mathrm{ev}_2^* \alpha_2\cup \mathrm{ev}_3^* \alpha_3\cup \mathrm{ev}_4^* \alpha_4.
\end{equation}
There is a birational morphism
\begin{eqnarray*}
f: \overline{\mathcal{M}}_{0,4}(X,1)\rightarrow P(X)\times_{\Omega_X}P(X)\times_{\Omega_X}P(X)\times_{\Omega_X}P(X).
\end{eqnarray*}
We denote the composition
\[
P(X)\times_{\Omega_X}P(X)\times_{\Omega_X}P(X)\times_{\Omega_X}P(X)\xrightarrow{\mathrm{pr}_i}P(X)\xrightarrow{\pi_X}X
\]
by $\pi_{X,i}$, for $1\leq i\leq 4$. Then $\mathrm{ev}_i=\pi_{X,i}\circ f$. So by (\ref{eq-lem-4pointInvariant-FanoVarLine-1}) and the projection formula, we have
\begin{eqnarray}\label{eq-lem-4pointInvariant-FanoVarLine-2}
&&\langle\alpha_1,\alpha_2,\alpha_3,\alpha_4 \rangle_{0,4,1}\nn\\
&=&\int_{P(X)\times_{\Omega_X}P(X)\times_{\Omega_X}P(X)\times_{\Omega_X}P(X)}
\pi_{X,1}^{*}\alpha_1\cup \pi_{X,2}^{*}\alpha_2\cup \pi_{X,3}^{*}\alpha_3\cup \pi_{X,4}^{*}\alpha_4.
\end{eqnarray}
As we saw above, $P(X)$ is the projectivization of a rank 2 vector bundle over $\Omega_X$. Let $h=c_1(\mathcal{O}_{P(X)}(1))$. The cohomology group $H^*(P(X))$ is generated by $1$ and $h$ as a module over $H^*(\Omega_X)$. So for $1\leq i\leq 4$ there exists $\xi_i$ such that
\[
\pi_{X}^* \alpha_i=\xi_i+\big(\rho_{X*}\pi_{X}^* \alpha_i\big)\cup h=\xi_i+\big(\Psi(\alpha_i)\big)\cup h.
\]
Then integrations along fibers yield
\begin{eqnarray}\label{eq-lem-4pointInvariant-FanoVarLine-3}
&&\int_{P(X)\times_{\Omega_X}P(X)\times_{\Omega_X}P(X)\times_{\Omega_X}P(X)}
\pi_{X,1}^{*}\alpha_1\cup \pi_{X,2}^{*}\alpha_2\cup \pi_{X,3}^{*}\alpha_3\cup \pi_{X,4}^{*}\alpha_4\nn\\
&=&(\Psi\alpha_1, \Psi\alpha_2,\Psi\alpha_3, \Psi\alpha_4)_{\Omega_X}.
\end{eqnarray}
From (\ref{eq-lem-4pointInvariant-FanoVarLine-1}) and (\ref{eq-lem-4pointInvariant-FanoVarLine-3}) we obtain (\ref{eq-4pointInvariant-FanoVarLine}).
\end{proof}

We denote by $\sigma_1$ the hyperplane class on Grassmannians induced by the Plücker embeddings. 
We denote the pullback of $\sigma_1$ via $\iota$ to $\Omega_X$  also by $\sigma_1$. 
\begin{theorem}\label{thm-cubic-FanoVar-F2-identity}
For every smooth cubic hypersurface $X$ of dimension $n$, and for 
every $\alpha_i \in H_{\mathrm{prim}}^{n}(X)$, and $\beta_i=\Psi(\alpha_i)$, $1\leq i\leq 4$,  we have
\begin{eqnarray}\label{eq-cubic-FanoVar-F2-identity}
&&F^{(2)}(0)\cdot \Big((\beta_1,\beta_2,\sigma_1^{n-2})_{\Omega_X}\cdot (\beta_3,\beta_4,\sigma_1^{n-2})_{\Omega_X}
+(\beta_1,\beta_3,\sigma_1^{n-2})_{\Omega_X}\cdot (\beta_4,\beta_2,\sigma_1^{n-2})_{\Omega_X}\nn\\
&&+(\beta_1,\beta_4,\sigma_1^{n-2})_{\Omega_X}\cdot (\beta_2,\beta_3,\sigma_1^{n-2})_{\Omega_X}\Big)=
36(\beta_1,\beta_2,\beta_3,\beta_4 )_{\Omega_X}.
\end{eqnarray}
\end{theorem}

For the proof, we need to recall the following construction in \cite[Page 203-204]{Lew99}.
\begin{construction}\label{cons-FanoVar-Cubic}
Let $Z$ be a general smooth cubic hypersurface in $\mathbb{P}^{n+2}$, and $X$ the intersection of $Z$ with a general hyperplane. 
We have a natural embedding
\begin{eqnarray*}
\iota_Z:\Omega_Z\hookrightarrow G_2(\mathbb{C}^{n+3}).
\end{eqnarray*}
Taking the intersection of $n-2$ general hyperplane sections of $G_2(\mathbb{C}^{n+3})$ (induced by its Plücker embedding) with $\Omega_Z$, we obtain a smooth subvariety $\Omega'_Z$ of dimension $n$, and the intersection of $\Omega'_Z$ with $\Omega_{X}$ is a smooth subvariety $\Omega'_X$ of dimension $n-2$.  The restriction of $P(Z)$ (resp. $P(X)$)
to $\Omega'_Z$ (resp. $\Omega'_X$) is denoted by $P(Z)'$ (resp. $P(X)'$). Denote by $j:X\hookrightarrow Z$  and $j_2:\Omega'_X\hookrightarrow \Omega'_Z$ the embeddings. Denote by $\rho'_Z$ the restriction of $\rho_Z$ to $P(Z)'$, and similarly $\pi'_Z$, $\rho'_X$, $\pi'_X$. We borrow the following commutative diagram from \cite[(13.2)]{Lew99}.
\begin{eqnarray*}
\xymatrix{
 P(X)^{\prime} \ar[ddd]^{\rho_{X}^{\prime}} \ar[rr]^{\pi_{X}^{\prime}} \ar[dr]^{j_0}& & X \ar[dr]^{j} &\\
& \widetilde{X} \ar[ur]^{\pi} \ar[dr]^{j_1} \ar[ddr]_{\rho} & & Z\\
&& P(Z)^{\prime} \ar[ur]^{\pi_Z^{\prime}} \ar[d]^{\rho_Z^{\prime}}&\\
 \Omega_{X}^{\prime} \ar@{^{(}->}[rr]^{j_2}    &    & \Omega_{Z}^{\prime}  &  .}
\end{eqnarray*}
Here $ \widetilde{X}=X\times_Z  P(Z)'$. The morphism $j_0$ is induced by $\pi'_X$ and the embedding $P(X)'\hookrightarrow P(Z)'$ induced by $j_2$.
 By \cite[Page 204]{Lew99}, when $X$ is sufficiently general, $\widetilde{X}$ is the blow-up of $\Omega_{Z}^{\prime}$ along $\Omega_{X}^{\prime}$, with the  exceptional divisor equal to $P(X)^{\prime} $. Moreover, $ \widetilde{X} $ and $X$ have the same dimension $n$, and
the morphism $\pi$ is surjective.
\end{construction}

\begin{proof}[Proof of Theorem \ref{thm-cubic-FanoVar-F2-identity}] Let $\mathcal{W}$ be the family of all smooth cubic hypersurfaces in $\mathbb{P}^{n+1}$. By \cite[Theorem 4.2]{AK77}, the Fano varieties of lines in the fibers of $\mathcal{W}$ form a smooth family over $\mathcal{W}$. So to show (\ref{eq-cubic-FanoVar-F2-identity}) for all smooth cubic hypersurfaces in $\mathbb{P}^{n+1}$, it suffices to show it for one. In particular, we can assume that $X$ is a general cubic hypersurface such that Construction \ref{cons-FanoVar-Cubic} is available. 

 By \cite[Remark 13.11]{Lew99},  $\mathrm{deg}(\pi)=3!=6$. The \emph{cylinder homomorphism} $\Phi_{*}= \pi_{X*}^{\prime}\rho_X^{\prime*}: H^{n-2}(\Omega_X^{\prime})\rightarrow H^{n}(X)$,  by \cite[Proposition 13.16]{Lew99},  is surjective. 
Since $\pi$ is surjective of finite degree, it follows that for any primitive class $\alpha\in H_{\mathrm{prim}}^{n}(X)$, there exists 
$\beta\in H^{n-2}(\Omega_X^{\prime})$ such that 
\begin{eqnarray}\label{eq-thm-cubic-FanoVar-F2-identity-1}
j_{0*}\rho_X^{\prime*}\beta=\pi^* \alpha.
\end{eqnarray}
Since $j_0$ is the embedding of the exceptional divisor $P(X)'$ in $\widetilde{X}$, we have
\begin{eqnarray}\label{eq-thm-cubic-FanoVar-F2-identity-2}
\rho_{X*}^{\prime}j_{0}^* j_{0*}\rho_X^{\prime*}\beta=\rho_{X*}\big(c_1(\mathcal{O}_{P(X)'}(-1))\cup \rho_{X}^{\prime*}\beta)
\big)=-\beta.
\end{eqnarray}
For $i=1,2$, suppose $\alpha_i\in H_{\mathrm{prim}}^{n}(X)$, $\beta_i\in H^{n-2}(\Omega_X^{\prime})$ such that $j_{0*}\rho_X^{\prime*}\beta_i=\pi^* \alpha_i$. Then
\begin{eqnarray*}
	&&\int_{\widetilde{X}}\pi^* \alpha_1\cup \pi^* \alpha_2=\int_{\widetilde{X}}\big(j_{0*}\rho_X^{\prime*}\beta_1\big)\cup \big(j_{0*}\rho_X^{\prime*}\beta_2\big)\nn\\
	&=& \int_{\widetilde{X}}j_{0*}\big(\rho_X^{\prime*}\beta_1\cup j_{0}^*j_{0*}\rho_X^{\prime*}\beta_2\big)\nn\\
	&=& \int_{P(X)'}\rho_X^{\prime*}\beta_1\cup j_{0}^*j_{0*}\rho_X^{\prime*}\beta_2\nn\\
	&=& \int_{P(X)'}\rho_X^{\prime*}\beta_1\cup c_1\big(\mathcal{O}_{P(X)'}(-1)\big)\cup\rho_X^{\prime*}\beta_2\nn\\
	&=&-\int_{\Omega_X^{\prime}}\beta_1\cup\beta_2.
\end{eqnarray*}
On the other hand, since $\mathrm{deg}(\pi)=6$, we have
\begin{equation*}
	\int_{\widetilde{X}}\pi^* \alpha_1\cup \pi^* \alpha_2=6\int_{X}\alpha_1\cup\alpha_2.
\end{equation*}
So
\begin{eqnarray}\label{eq-thm-cubic-FanoVar-F2-identity-3}
6\int_{X}\alpha_1\cup\alpha_2=-\int_{\Omega_X^{\prime}}\beta_1\cup\beta_2.
\end{eqnarray}

By (\ref{eq-thm-cubic-FanoVar-F2-identity-1}) and (\ref{eq-thm-cubic-FanoVar-F2-identity-2}),
\begin{eqnarray}\label{eq-thm-cubic-FanoVar-F2-identity-4}
\rho_{X*}^{\prime}\pi_{X}^{\prime*}\alpha=\rho_{X*}^{\prime} j_{0}^{*}\pi^{*}\alpha=\rho_{X*}^{\prime}j_{0}^*(j_{0*}\rho_X^{\prime*}\beta)=-\beta.
\end{eqnarray}
Consider the following commutative diagram with a cartesian square,
\begin{eqnarray*}
\xymatrix{\ar @{} [dr] |{\square}
 P(X)^{\prime} \ar[d]^{\rho_{X}^{\prime}} \ar[r] ^{i^{\prime}} \ar @/^1.5pc/[rr]^{\pi^{\prime}_X} & P(X) \ar[d]^{\rho_{X}} \ar[r]^{\pi_X}  & X \\
 \Omega_{X}^{\prime} \ar[r]^{i}  & \Omega_X & . }
\end{eqnarray*}
We have
\begin{eqnarray}\label{eq-thm-cubic-FanoVar-F2-identity-5}
-\beta=\rho_{X*}^{\prime}\pi_{X}^{\prime*}\alpha=\rho_{X*}^{\prime}i^{\prime*}\pi_{X}^{ *}\alpha=i^{*}(\rho_{X*}\pi_{X}^{ *}\alpha)=i^{*}\Psi(\alpha).
\end{eqnarray}
Now let $\alpha_i$ and $\beta_i$ be as the assumption in Theorem \ref{thm-cubic-FanoVar-F2-identity}, for $1\leq i\leq 4$. Then by (\ref{eq-thm-cubic-FanoVar-F2-identity-3}) and (\ref{eq-thm-cubic-FanoVar-F2-identity-5}),
\begin{multline}\label{eq-thm-cubic-FanoVar-F2-identity-6}
	6(\alpha_i,\alpha_j)_X=-\big(i^*\Psi(\alpha_i),i^*\Psi(\alpha_j)\big)_{\Omega'_X}\\
	=-\big(\Psi(\alpha_i),\Psi(\alpha_j),\sigma_1^{n-2}\big)_{\Omega_X}
	=-\big(\beta_i,\beta_j,\sigma_1^{n-2}\big)_{\Omega_X}.
\end{multline}
From (\ref{eq-Fano1}), (\ref{eq-4pointInvariant-FanoVarLine}) and 
(\ref{eq-thm-cubic-FanoVar-F2-identity-6}), we obtain (\ref{eq-cubic-FanoVar-F2-identity}).
\end{proof}

\begin{remark}
We need this statement for every smooth cubic hypersurface because we will use the monodromy argument for $\Omega_X$ in the next section.
\end{remark} 
Before performing the computation for the smooth cubic hypersurfaces in all dimensions, let us first see two examples.

\begin{example}
When $n=3$, $\Omega_X$ is the Fano surface. A detailed study of the intersection ring of $\Omega_X$ is given in \cite{CG72}. They showed that there are divisors $D_s$ on $\Omega_X$ subject to  linear equivalences
\begin{eqnarray*}
\sigma_1=[\Omega_X^{\prime}]\sim D_{s_1}+D_{s_2}+D_{s_3}.
\end{eqnarray*}
Moreover, they showed that there is a basis of $H^{1}(\Omega_X)$, denoted by $\chi$, $\delta$, $\eta_{k}$, $1\leq k\leq 8$ such that all the nonzero integrations of the products 
of these bases over $\Omega_X$ and $D_s$ are
\begin{eqnarray*}
\int_{D_{s}}\chi\cup \delta=\int_{D_s}\eta_k\cup\eta_{k+4}=2, & k=1,\dots,4,\\
\end{eqnarray*} 
and
\begin{eqnarray*}
\int_{\Omega_X}\chi\cup \delta\cup\eta_k\cup\eta_{k+4}=1, & k=1,\dots,4,\\
\int_{\Omega_X}\eta_k\cup\eta_{k+4}\cup\eta_l\cup\eta_{l+4}=1, & 1\leq k<l\leq 4.
\end{eqnarray*}
From this one can check that (\ref{eq-cubic-FanoVar-F2-identity}) is true with $F^{(2)}(0)=1$.
\pqed
\end{example}

\begin{example}
When $n=4$, by \cite{BD85}, $\Omega_X$ is deformation equivalent to $S^{[2]}$, the Hilbert scheme of 2 points on some special K3 surface $S$. We recall the construction of $S$ in
 \cite{BD85} and compute $F^{(2)}(0)$ as follows.  
Let $V$ be a complex vector space of dimension 6, with a basis $e_1,\dots,e_6$, and let $e_1^*,\dots,e_6^*$ be the dual basis of $V^*$. Let $\{a_{ij}\}_{1\leq i<j\leq 6}$ be the dual basis of 
$\{e_i\wedge e_j\}_{1\leq i<j\leq 6}$. Consider the skew-symmetric matrix $N=(a_{ij})_{1\leq i,j\leq 6}$ where $a_{ji}=-a_{ij}$ for $j>i$ and $a_{ii}=0$. The cubic Pfaffian 
is a subvariety of $\mathbb{P}(\bigwedge^2 V^*)$ defined by
\begin{eqnarray*}
\mathrm{Pf}(N)=\sum_{\sigma\in S_6}\mathrm{sign}(\sigma)a_{\sigma(1)\sigma(2)}a_{\sigma(3)\sigma(4)}a_{\sigma(5)\sigma(6)}.
\end{eqnarray*}
It is straightforward to check that the singular locus of $\mathrm{Pf}(N)$ is of dimension 8. So take a generic 8-plane $L$ of $\mathbb{P}(\bigwedge^2 V)$, the intersection of $\mathrm{Pf}(N)$ with the 5-plane $L^{\perp}$ is a smooth cubic 4-fold, denoted by $X$. The intersection of $G$ and $L$ is a K3 surface, denoted by $S$. The hyperplane class of $S$ is denoted by $l$. The Hilbert scheme $S^{[2]}$ is obtained by blowing up the symmetric product $S^{(2)}$ along the diagonal, and a half of the exceptional divisor is denoted by $\delta$. There is a canonical isomorphism 
\begin{eqnarray}\label{eq-4fold0}
H^{2}(\Omega_X)\cong H^{2}(S)\oplus \mathbb{C}\delta,
\end{eqnarray}
 and via this isomorphism, we have $\sigma_1=2l-5\delta.$
 The intersections of these classes on $S^{[2]}$ is given in \cite[proof of Proposition 6]{BD85}. For $\gamma_1,\dots,\gamma_4\in H^{2}(S)$, 
\begin{eqnarray*}
 &&(\gamma_1, \gamma_2, \gamma_3, \gamma_4)_{S^{[2]}}\\
 &=&\frac{1}{2}\int_{S\times S}\prod_{i=1}^{4}(\mathrm{pr}_1^* \gamma_i+\mathrm{pr}_2^* \gamma_i) \\
 &=&(\gamma_1\cdot \gamma_2)(\gamma_3\cdot \gamma_4)
 +(\gamma_1\cdot \gamma_3)(\gamma_2\cdot \gamma_4)+(\gamma_1\cdot \gamma_4)(\gamma_2\cdot \gamma_3).
 \end{eqnarray*}
 So 
\begin{eqnarray}\label{eq-4fold1}
&&(\gamma_1+a_1\delta, \gamma_2+a_2\delta, \gamma_3+a_3\delta, \gamma_4+a_4\delta)_{S^{[2]}}\nn\\
&=&(\gamma_1\cdot \gamma_2)(\gamma_3\cdot \gamma_4)
 +(\gamma_1\cdot \gamma_3)(\gamma_2\cdot \gamma_4)+(\gamma_1\cdot \gamma_4)(\gamma_2\cdot \gamma_3)\nn\\
 &&-(2a_1 a_2 \gamma_3 \cdot \gamma_4+2a_3 a_4 \gamma_1 \cdot\gamma_2+2a_1 a_3 \gamma_2\cdot \gamma_4\nn\\
&&+2a_2 a_4 \gamma_1\cdot\gamma_3+2a_2 a_3 \gamma_1 \cdot\gamma_4+2a_1 a_4 \gamma_2 \cdot\gamma_3)+12a_1a_2a_3a_4.
\end{eqnarray}
On the other hand,
\begin{eqnarray}\label{eq-4fold2}
(\sigma_1,\sigma_1, \gamma_1+a_1\delta,\gamma_2+a_2\delta)_{S^{[2]}}= 6 (\gamma_1\cdot\gamma_2-2a_1a_2).
\end{eqnarray}
Substituting (\ref{eq-4fold1}) and (\ref{eq-4fold2}) into (\ref{eq-cubic-FanoVar-F2-identity}), we obtain $F^{(2)}(0)=1$.\pqed
\end{example}

\subsection{Cohomology ring of Fano variety of lines in cubic hypersurfaces}

In this section for the cubic hypersurface $X=X_n(3)$ of dimension $n\geq 3$, we study the cohomology ring of $\Omega_X$, via the Schubert calculus on the Grassmannian 
$G_2(\mathbb{C}^{n+2})$ and monodromy arguments, together with the  result \cite{GS14} on the Betti number of 
$\Omega_X$. We obtain a complete description of the ring structure of $H^{*}(\Omega_X)$, and by the way obtain $F^{(2)}(0)=1$. 

For the Schubert calculus we adhere to the notations in \cite[chap.14]{Ful98}. To each partition 
$\lambda=(\lambda_0,\lambda_1)$ with
 $n\geq \lambda_0\geq \lambda_1\geq 0$ is associated a Schubert class on $G_2(\mathbb{C}^{n+2})$, denoted by $\{\lambda_0,\lambda_1\}$. In particular, $\sigma_{i}=\{i,0\}$ for $0\leq i\leq n$ is the $i$-th Chern class of the universal quotient bundle.  The cohomology ring $H^{*}(G_2(\mathbb{C}^{n+2}))$ is generated by $\sigma_1,\sigma_2$, and the Schubert classes form an additive basis of 
 $H^{*}(G_2(\mathbb{C}^{n+2}))$. Thus
 \begin{eqnarray*}
\mathrm{rk}\ H^{2i}(G_2(\mathbb{C}^{n+2}))=\left\{
\begin{array}{cc}
\lfloor \frac{i}{2}\rfloor+1, & i\leq n,\\
n+1-\lceil \frac{i}{2}\rceil, & i\geq n.
\end{array}
\right.
\end{eqnarray*}
The pull-backs of the Schubert classes to $\Omega_X$ are still denoted by the same symbols when no confusion arises. The Hodge structure of 
$\Omega_X$ is given in \cite{GS14}. Let 
\begin{eqnarray*}
\delta_{i,j}^{\mathrm{mod}\ 2}=\left\{
\begin{array}{cc}
1, & \mathrm{if}\ i-j\in 2\mathbb{Z},\\
0, & \mathrm{if}\ i-j\in 2\mathbb{Z}+1.
\end{array}\right.
\end{eqnarray*} 
Then as a corollary of \cite[Theorem 6.1]{GS14} we obtain
\begin{eqnarray}\label{eq-cubic0}
&&\mathrm{rk}\ H^{i}(\Omega_X)-\mathrm{rk}\ H^{i}(G_2(\mathbb{C}^{n+2}))\nn\\
&=&\left\{
\begin{array}{cc}
0, & i<n-2,\\
\delta_{0,i-n}^{\mathrm{mod}\ 2}\cdot\mathrm{rk}\ H_{\mathrm{prim}}(X), & n-2\leq i<2n-4,\\
\delta_{0,i-n}^{\mathrm{mod}\ 2}\cdot\mathrm{rk}\ H_{\mathrm{prim}}(X)+\mathrm{rk}\ \mathrm{Sym}^2 H_{\mathrm{prim}}(X)-1, & i=2n-4,\\
\delta_{0,i-n}^{\mathrm{mod}\ 2}\cdot\mathrm{rk}\ H_{\mathrm{prim}}(X)-\delta_{0,i}^{\mathrm{mod}\ 2}, & 2n-4< i\leq 2n-2,\\
\delta_{0,i-n}^{\mathrm{mod}\ 2}\cdot\mathrm{rk}\ H_{\mathrm{prim}}(X)-2\delta_{0,i}^{\mathrm{mod}\ 2}, & 2n-2< i\leq 3n-6,\\
-2\delta_{0,i}^{\mathrm{mod}\ 2}, & \max\{2n-2,3n-6\}<i\leq 4n-8.
\end{array}
\right.
\end{eqnarray}

In the the following proposition we use Schubert calculus to obtain estimates of the rank of the image of the pull-back homomorphism
$H^{2i}(G_2(\mathbb{C}^{n+2}))\rightarrow H^{2i}(\Omega_X)$.
\begin{proposition}\label{prop-cubic2}\quad\\
(i) $H^{2i}(G_2(\mathbb{C}^{n+2}))\rightarrow H^{2i}(\Omega_X)$ is injective for $0\leq i\leq n-2$.\\
(ii) $\mathrm{rk}\ \mathrm{Im}\big(H^{2n-2}(G_2(\mathbb{C}^{n+2}))\rightarrow H^{2n-2}(\Omega_X)\big)\geq
 \lceil\frac{n}{2}\rceil-1$.\\
(iii) $\mathrm{rk}\ \mathrm{Im}\big(H^{2i}(G_2(\mathbb{C}^{n+2}))\rightarrow H^{2i}(\Omega_X)\big)\geq 
n-1-\lceil \frac{i}{2}\rceil$ for  $n\leq i\leq 2n-4$.
\end{proposition}
\begin{proof} The fundamental class of $\Omega_X$ in $H^*(G_{2}(\mathbb{C}^{n+2}))$ is
\begin{eqnarray}\label{eq-cubic3}
9(3\sigma_1^4-4\sigma_1^2\sigma_2+\sigma_2^2).
\end{eqnarray}
(see e.g., \cite[proof of prop.1.6]{AK77}, \cite[Example 14.7.13]{Ful98}.)
 We use the convention $\{k_1,k_2\}=0$ for $k_1\geq n+1$. Then by Pieri's rule \cite[Page 271]{Ful98}, in $H^*(G_2(\mathbb{C}^{n+2}))$, 
\begin{eqnarray}\label{eq-cubic4}
&&\{k_1,k_2\}\cdot (3\sigma_1^4-4\sigma_1^2\sigma_2+\sigma_2^2)\nn\\
&=&\left\{
\begin{array}{lc}
2\{k_1+3,k_2+1\}+5\{k_1+2,k_2+2\}+2\{k_1+1,k_2+3\}, & k_1-k_2\geq 2,\\
2\{k_1+3,k_2+1\}+5\{k_1+2,k_2+2\}, & k_1-k_2=1,\\
2\{k_1+3,k_2+1\}+3\{k_1+2,k_2+2\}, & k_1-k_2=0.
\end{array}\right.
\end{eqnarray}
(i) Since $H^*(G_2(\mathbb{C}^{n+2}))$ is generated by $\sigma_1$ and $\sigma_2$, to show 
$H^{2i}(G_2(\mathbb{C}^{n+2}))\rightarrow H^{2i}(\Omega_X)$ is injective for $0\leq i\leq n-2$, it suffices to show this for $i=n-2$.
 Let $n_0=\lfloor \frac{n}{2}\rfloor$. Then $\{n-2-k,k\}$ for $0\leq k\leq n_0-1$ form a basis of $H^{2n-4}(G_2(\mathbb{C}^{n+2}))$. Suppose 
\begin{eqnarray*}
\sum_{k=0}^{n_0-1}y_{k}\{n-2-k,k\}\cdot (3\sigma_1^4-4\sigma_1^2\sigma_2+\sigma_2^2)=0
\end{eqnarray*} in $H^{2n-4}(\Omega_X)$, then in 
$H^{2n+4}(G_2(\mathbb{C}^{n+2}))$ we have
\begin{eqnarray*}
\sum_{k=0}^{n_0-1}y_{k}\{n-2-k,k\}\cdot (3\sigma_1^4-4\sigma_1^2\sigma_2+\sigma_2^2)=0.
\end{eqnarray*}
Then by (\ref{eq-cubic4}) we have
\begin{eqnarray}\label{eq-prop-cubic2-linearSystem-1}
\left\{\begin{array}{c}
5y_0+2y_1=0,\\
2y_0+5y_1+2y_2=0,\\
\dots,\\
2y_{n_0-3}+5y_{n_0-2}+2y_{n_0-1}=0,\\
\left\{\begin{array}{cc}
2y_{n_0-2}+3y_{n_0-1}=0, & \mbox{if}\ n=2n_0,\\
2y_{n_0-2}+5y_{n_0-1}=0, & \mbox{if}\ n=2n_0+1.
\end{array}\right.
\end{array}\right.
\end{eqnarray}
Let $d_{n_0}$ be the determinant of this linear system. Then $d_{k}=5d_{k-1}-4d_{k-2}$, and 
\begin{equation*}
	\begin{cases}
	d_1=3,\ d_2=11,& \mbox{if}\ n=2n_0,\\
	d_1=5,\ d_2=21, & \mbox{if}\ n=2n_0+1. \\
	\end{cases}
\end{equation*}
By an elementary linear recursion, one computes
\begin{equation*}
	\begin{cases}
	d_{n_0}=\frac{5\cdot 4^n-17}{3},& \mbox{if}\ n=2n_0\geq 4,\\
	d_{n_0}=\frac{2\cdot 4^n+1}{3}, & \mbox{if}\ n=2n_0+1\geq 3.
	\end{cases}
\end{equation*}
So in either case there is no nontrivial  solutions.

(ii) It is equivalent to showing that the kernel of $H^{2n-2}(G_2(\mathbb{C}^{n+2}))\rightarrow H^{2n-2}(\Omega_X)$ has rank $\leq 1$. 
Similarly as in (i), suppose $n_0=\lceil \frac{n}{2}\rceil$, then by (\ref{eq-cubic4})
$$\sum_{k=0}^{n_0-1}y_{k}\{n-1-k,k\}\cdot (3\sigma_1^4-4\sigma_1^2\sigma_2+\sigma_2^2)=0$$
in $H^{2n+6}(G_2(\mathbb{C}^{n+2}))$ if and only if
\begin{eqnarray}\label{eq-prop-cubic2-linearSystem-2}
\left\{\begin{array}{c}
2y_0+5y_1+2y_2=0,\\
2y_1+5y_2+2y_3=0,\\
\dots,\\
2y_{n_0-3}+5y_{n_0-2}+2y_{n_0-1}=0,\\
\left\{\begin{array}{cc}
2y_{n_0-2}+5y_{n_0-1}=0, & \mathrm{if}\ n=2n_0,\\
2y_{n_0-2}+3y_{n_0-1}=0, & \mathrm{if}\ n=2n_0-1.
\end{array}\right.
\end{array}\right.
\end{eqnarray}
The matrix of (\ref{eq-prop-cubic2-linearSystem-2}) is obtained by dropping the first row of that of (\ref{eq-prop-cubic2-linearSystem-1}). 
So in either case $n=2n_0$ or $n=2n_0-1$, there is a unique solution up to a common factor.

(iii) It is equivalent to show that the kernel of $H^{2i}(G_2(\mathbb{C}^{n+2}))\rightarrow H^{2i}(\Omega_X)$ has rank $\leq 2$ for $n\leq i\leq 2n-4$.
For $0\leq l\leq n-4$, suppose $n_0=\lfloor \frac{n+l}{2}\rfloor$, then
$$\sum_{k=l}^{n_0}y_{k}\{n+l-k,k\}\cdot (3\sigma_1^4-4\sigma_1^2\sigma_2+\sigma_2^2)=0$$
in $G_2(\mathbb{C}^{n+2})$ if and only if
\begin{eqnarray}\label{eq-prop-cubic2-linearSystem-3}
\left\{\begin{array}{c}
2y_{l+1}+5y_{l+2}+2y_{l+3}=0,\\
2y_{l+2}+5y_{l+3}+2y_{l+4}=0,\\
\dots,\\
2y_{n_0-2}+5y_{n_0-1}+2y_{n_0}=0,\\
\left\{\begin{array}{cc}
2y_{n_0-1}+3y_{n_0}=0, & \mathrm{if}\ n+l=2n_0,\\
2y_{n_0-1}+5y_{n_0}=0, & \mathrm{if}\ n+l=2n_0+1.
\end{array}\right.
\end{array}\right.
\end{eqnarray}
The matrix of (\ref{eq-prop-cubic2-linearSystem-3}) is the last $n_0-l-1$ rows of that of (\ref{eq-prop-cubic2-linearSystem-1}), thus has rank $=n_0-l-1$. 
So in either case  $n=2n_0$ or $n=2n_0-1$, there are only two linearly independent solutions.
\end{proof}

\begin{proposition}\label{prop-FanoVar-SchubertCalculus}
Let $\alpha\in H_{\mathrm{prim}}^{n}(X)$, and $\beta=\Psi(\alpha)$. Then in $H^*(\Omega_X)$, 
\begin{eqnarray}\label{eq-FanoVar-SchubertCalculus}
\beta\cup \sigma_{1}^2=\beta\cup\sigma_2.
\end{eqnarray}
\end{proposition}
\begin{proof} Choose a general hyperplane of $\mathbb{P}^{n+1}$, such that its intersection with $X$ is a smooth cubic hypersurface $Y$ of dimension $n-1$. Then we have the following commutative diagram, where the lower square is cartesian.
\begin{eqnarray*}
\xymatrix{
Y\ar[r]^{i_1}  & X \\
\ar @{} [dr] |{\square}
 P(Y) \ar[d]_{\rho_Y} \ar[r]^{i_2}  \ar[u]^{\pi_Y}  & P(X)\ar[d]^{\rho_X} \ar[u]_{\pi_X} \\
 \Omega_{Y} \ar[r]^{i_3}   & \Omega_X & }
\end{eqnarray*}
By Pieri's rule,
\[
\sigma_{1}^2-\sigma_2=\{1,1\}.
\]
The Schubert class $\{1,1\}$ is represented by a particular  Schubert variety, the sub-Grassmannian 
$G_2(\mathbb{C}^{n+1})\hookrightarrow G_2(\mathbb{C}^{n+2}) $ induced by any  hyperplane in $\mathbb{P}^{n+1}$. Thus 
$\{1,1\}\cap [\Omega_X]=i_{3*}[\Omega_Y]$. Since $i_1^{*}\alpha=0$ for a primitive class $\alpha$ on $X$, and that the lower square a cartesian one, we have $i_3^{*}\rho_{X*}\pi_{X}^{*}\alpha=\rho_{Y*}i_2^{*}\pi_{X}^{*}\alpha=\rho_{Y*}\pi_Y^{*}i_1^{*}\alpha=0$, i.e., $\beta\cup (\sigma_{1}^2-\sigma_2)=0$.
\end{proof}

Let $\alpha_{1},\dots,\alpha_{m}$ be a basis of $H_{\mathrm{prim}}^{n}(X)$, $g_{ij}=\int_X \alpha_i\cup\alpha_j$, and define $\beta_{i}=\Psi(\alpha_i)$ for $1\leq i\leq m$.

\begin{theorem}\label{thm-cubic5}\quad\\
(i) For $1\leq i,j\leq m$,
\begin{eqnarray}\label{cubic6}
\int_{\Omega_X}\beta_i\cup\beta_j\cup \sigma_1^{n-2}=-6g_{ij}.
\end{eqnarray}
(ii) For $\gamma\in \mathrm{Im}\big(H^{*}(G_2(\mathbb{C}^{n+2}))\rightarrow H^{*}(\Omega_X)\big)$, 
\begin{eqnarray}\label{eq-cubic6.5}
\int_{\Omega_X}\beta_i\cup\gamma=0, & \int_{\Omega_X}\beta_i\cup\beta_j\cup\beta_k\cup\gamma=0,
\end{eqnarray}
for $1\leq i,j,k\leq m$.\\
(iii) \begin{eqnarray}\label{eq-cubic7}
\sum_{i,j=1}^{m}\beta_{i}g^{ij}\beta_{j}
=\sum_{k=0}^{\lfloor\frac{n}{2}\rfloor-1}\frac{(-2)^{n+1-k}-(-2)^{k+2}}{9}\{n-2-k,k\}.
\end{eqnarray}
(iv) For $1\leq i_1,i_2,i_3,i_4\leq m$,
\begin{eqnarray}\label{eq-cubic14}
\int_{\Omega_X}\beta_{i_1}\cup\beta_{i_2}\cup\beta_{i_3}\cup\beta_{i_4}
=g_{i_1i_2}g_{i_3 i_4}+g_{i_1 i_4}g_{i_2 i_3}+g_{i_1 i_3}g_{i_4 i_2}.
\end{eqnarray}
Equivalently, $F^{(2)}(0)=1$.
\end{theorem}
\begin{proof} (i) This is a restatement of (\ref{eq-thm-cubic-FanoVar-F2-identity-3}).

(ii) When $X$ deforms in the whole family of smooth cubic hypersurfaces in $\mathbb{P}^{n+1}$, $\Omega_X$ deforms in the same Grassmannian, and thus $\gamma$ is an invariant class on $\Omega_X$, and the integrals in (\ref{eq-cubic6.5}) transforms in the same way as $\alpha_i$ or $\alpha_i,\alpha_j,\alpha_k$. So by the monodromy reason as in Section \ref{sec:monodromyGroup}, these integrals vanish.

(iii) By (\ref{eq-cubic0}) and Proposition \ref{prop-cubic2} (i), when $n$ is odd, there is  a relation among the classes 
$\{n-2,0\},\dots,\{n-1-\lfloor\frac{n}{2}\rfloor,\lfloor\frac{n}{2}\rfloor-1\}$ and  $\{\beta_i\cup\beta_j\}_{i,j=1}^{m}$, and when $n$ is even, there is  a relation among the classes 
$\{n-2,0\},\dots,\{n-1-\lfloor\frac{n}{2}\rfloor,\lfloor\frac{n}{2}\rfloor-1\}$,  $\{\beta_i\cup\sigma_1^{\frac{n-2}{2}}\}_{i=1}^{m}$, and $\{\beta_i\cup\beta_j\}_{i,j=1}^{m}$. In the latter case, suppose 
\begin{eqnarray*}
\sum_{k=0}^{\lfloor\frac{n}{2}\rfloor-1}a_{k}\{n-2-k,k\}
+\sum_{k=1}^{m}b_k\beta_k\cup\sigma_1^{\frac{n-2}{2}}
+\sum_{i,j=1}^{m}c_{ij}\beta_i\cup\beta_j=0,
\end{eqnarray*}
where $a_k,b_k,c_{ij}\in\mathbb{C}$. Then taking the cup product with $\beta_{j}\sigma_{1}^{\frac{n-2}{2}}$ for $j=1,\dots,m$, by (i) and (ii) we have $b_k=0$ for $1\leq k\leq m$. So in either case, there is a relation of the form
\begin{eqnarray*}
\sum_{k=0}^{\lfloor\frac{n}{2}\rfloor-1}a_{k}\{n-2-k,k\}
+\sum_{i,j=1}^{m}c_{ij}\beta_i\cup\beta_j=0.
\end{eqnarray*}
By the monodromy argument as in (ii) and in Section \ref{sec:monodromyGroup}, the relation must take the following form
\begin{eqnarray}\label{eq-cubic7.9}
\sum_{k=0}^{\lfloor\frac{n}{2}\rfloor-1}a_{k}\{n-2-k,k\}
+c\sum_{i,j=1}^{m}\beta_i g^{ij}\beta_j=0.
\end{eqnarray}
Again by Proposition \ref{prop-cubic2} (i), $c\neq 0$. Let $n_0=\lfloor\frac{n}{2}\rfloor$. Suppose
\begin{eqnarray*}
\sum_{i,j=1}^{m}\beta_{i}g^{ij}\beta_{j}
=\sum_{k=0}^{n_0-1}x_k\{n-2-k,k\}.
\end{eqnarray*}
Then by Proposition \ref{prop-FanoVar-SchubertCalculus},
\begin{eqnarray}\label{eq-cubic8}
(\sigma_1^2-\sigma_2)\sum_{k=0}^{n_0-1}x_{k}\{n-2-k,k\}\cdot (3\sigma_1^4-4\sigma_1^2\sigma_2+\sigma_2^2)=0.
\end{eqnarray}
On the other hand, by Pieri's rule,
\begin{eqnarray}\label{eq-cubic9}
&&\{k_1,k_2\}\cdot (3\sigma_1^4-4\sigma_1^2\sigma_2+\sigma_2^2)\cdot(\sigma_1^2-\sigma_2)\nn\\
&=&\left\{
\begin{array}{lc}
2\{k_1+4,k_2+2\}+5\{k_1+3,k_2+3\}+2\{k_1+2,k_2+4\}, & k_1-k_2\geq 2,\\
2\{k_1+4,k_2+2\}+5\{k_1+3,k_2+3\}, & k_1-k_2=1,\\
2\{k_1+4,k_2+2\}+3\{k_1+3,k_2+3\}, & k_1-k_2=0.
\end{array}\right.
\end{eqnarray}
Comparing (\ref{eq-cubic8}) and (\ref{eq-cubic9}) we obtain a system of linear equations on $x_0,\dots, x_{n_0-1}$, as in the proof of Proposition \ref{prop-cubic2}:
\begin{eqnarray}\label{eq-cubic9.5}
\left\{\begin{array}{c}
2x_0+5x_1+2x_2=0,\\
2x_1+5x_2+2x_3=0,\\
\dots,\\
2x_{n_0-3}+5x_{n_0-2}+2x_{m-1}=0,\\
\left\{\begin{array}{cc}
2x_{n_0-2}+5x_{n_0-1}=0, & \mathrm{if}\ n=2n_0+1,\\
2x_{n_0-2}+3x_{n_0-1}=0, & \mathrm{if}\ n=2n_0.
\end{array}\right.
\end{array}\right.
\end{eqnarray}
Solving this system, we obtain, up to a common factor,
\begin{eqnarray}\label{eq-cubic10}
x_k=
\left\{
\begin{array}{cc}
\frac{2}{3}\Big((-\frac{1}{2})^{n_0-k}-(-2)^{n_0-k}\Big), & n=2n_0+1,\\
-\frac{2}{3}\Big((-\frac{1}{2})^{n_0-k}-(-2)^{n_0-k-1}\Big), & n=2n_0.
\end{array}\right.
\end{eqnarray}
To fix the common factor, we  compute
\begin{eqnarray*}
\sum_{k=0}^{n_0-1}\Bigg(x_k\int_{\Omega_X}\sigma_1^{n-2}\cup\{n-2-k,k\}\Bigg).
\end{eqnarray*}
in two ways. As in \cite[\S 14.7]{Ful98}, we will use $(n-\lambda_0,n+1-\lambda_1):=\{\lambda_0,\lambda_1\}\cap[G_2(\mathbb{C}^{n+2})]$ to denote the Schubert variety. Then when $n=2n_0+1$,
\begin{eqnarray}\label{eq-cubic11}
&&\sum_{k=0}^{n_0-1}\Bigg(x_k\int_{\Omega_X}\sigma_1^{n-2}\cup\{n-2-k,k\}\Bigg)\nn\\
&=& 9x_0\Big(5\deg(0,n-1)+2\deg(1,n-2)\Big)\nn\\
&&+9\sum_{k=1}^{n_0-2}x_k\Big(2\deg(k-1,n-k)+5\deg(k,n-k-1)+2\deg(k+1,n-k-2)\Big)\nn\\
&&+9x_{n_0-1}\Big(2\deg(n_0-2,n-n_0+1)+5\deg(n_0-1,n-n_0)\Big)\nn\\
&=& 9(2x_{1}+5x_0)\deg(0,n-1)+9\sum_{k=1}^{n_0-2}(2x_{k+1}+5x_k+2x_{k-1})\deg(k,n-k-1)\nn\\
&&+9(5x_{n_0-1}+2x_{n_0-2})\deg(n_0-1,n-n_0)\nn\\
&=& 9(2x_{1}+5x_0)\deg(0,n-1)\nn\\
&=&3(-1)^{n_0+1}(2^{n_0+3}-2^{1-n_0}).
\end{eqnarray}
For the first equality we use (\ref{eq-cubic4}), for the third equality we use (\ref{eq-cubic9.5}), and for the last one we use the degree formula of Schubert varieties, see e.g., \cite[Example 14.7.11]{Ful98}. Similarly, when $n=2n_0$,
\begin{eqnarray}\label{eq-cubic12}
&&\sum_{k=0}^{n_0-1}\Bigg(x_k\int_{\Omega_X}\sigma_1^{n-2}\cup\{n-2-k,k\}\Bigg)\nn\\
&=&3(-1)^{n_0+1}(2^{n_0+2}+2^{1-n_0}).
\end{eqnarray}
On the other hand, by (\ref{cubic6}), 
\begin{eqnarray}\label{eq-cubic13}
&& \int_{\Omega_X}\sigma_1^{n-2}\cup \sum_{i,j}\beta_{i}g^{ij}\beta_{j}\nn\\
&=& -6 (\chi(X)-n-1)=(-2)^{n+3}-4.
\end{eqnarray}
Comparing (\ref{eq-cubic11}) and (\ref{eq-cubic12}) to (\ref{eq-cubic13}), we obtain, from (\ref{eq-cubic10}),
\begin{eqnarray*}
 \sum_{i,j}\beta_{i}g^{ij}\beta_{j}=\sum_{k=0}^{n_0-1}\frac{(-2)^{n+1-k}-(-2)^{k+2}}{9}\{n-2-k,k\}.
\end{eqnarray*}
(iv) By Theorem \ref{thm-cubic-FanoVar-F2-identity}, there exists $c(n)$ such that
\begin{eqnarray}\label{eq-cubic15}
\int_{\Omega_X}\beta_{i_1}\cup\beta_{i_2}\cup\beta_{i_3}\cup\beta_{i_4}
=c(n)\big(g_{i_1i_2}g_{i_3 i_4}+g_{i_1 i_4}g_{i_2 i_3}+g_{i_1 i_3}g_{i_4 i_2}\big).
\end{eqnarray}
By taking summations on both sides, one finds that
\begin{eqnarray}\label{eq-cubic16}
&&\int_{\Omega_X}\sum_{i,j=1}^{m}\beta_{i}g^{ij}\beta_{j}\cup\sum_{i,j=1}^{m}\beta_{i}g^{ij}\beta_{j}\nn\\
&=&c(n)(m^2+2m)=c(n)\big((\chi(X)-n)^2-1\big).
\end{eqnarray}
Write the result of (iii) as $\sum_{i,j=1}^{m}\beta_{i}g^{ij}\beta_j=\sum_{k=0}^{n_0-1}z_k\{n-2-k,k\}$. Then $\{z_k\}$ is a solution of (\ref{eq-cubic9.5}). We compute
\begin{eqnarray*}
&&\int_{\Omega_X}\sum_{k=0}^{n_0-1}z_{k}\{n-2-k,k\}\cup\sum_{k=0}^{n_0-1}z_{k}\{n-2-k,k\}\\
&=&\int_{G_{2}(\mathbb{C}^{n+2})}\sum_{k=0}^{n_0-1}z_{k}\{n-2-k,k\}\cup\sum_{k=0}^{n_0-1}z_{k}\{n-2-k,k\}\cup 9(3\sigma_1^4-4\sigma_1^2\sigma_2+\sigma_2^2)\\
&=& \int_{G_{2}(\mathbb{C}^{n+2})}\sum_{k=0}^{n_0-1}z_{k}\{n-2-k,k\}\cup 9(5z_0+2z_1)\{n,2\}\\
&=&\int_{G_{2}(\mathbb{C}^{n+2})}9z_0(5z_0+2z_1)\{n,n\}\\
&=&9z_0(5z_0+2z_1).
\end{eqnarray*}
For the second equality we use (\ref{eq-cubic4}), and then (\ref{eq-cubic9.5}) for $z_k$; the third equality follows from the duality theorem of Schubert classes, see e.g. \cite[Page 271]{Ful98}. Now applying (\ref{eq-cubic7}) we get
\begin{eqnarray}\label{eq-cubic17}
&&\int_{\Omega_X}\sum_{i,j=1}^{m}\beta_{i}g^{ij}\beta_{j}\cup\sum_{i,j=1}^{m}\beta_{i}g^{ij}\beta_{j}=
9 \frac{(-2)^{n+1}-4}{9}\cdot \frac{5\big((-2)^{n+1}-4\big)+2\big((-2)^n+8\big)}{9}\nn\\
&=& \frac{16(-2)^{2n}+40(-2)^n+16}{9}.
\end{eqnarray}
Using (\ref{eq-EulerChar}),  we obtain $c(n)=1$ from (\ref{eq-cubic15}) and (\ref{eq-cubic17}). 
\end{proof}

\begin{remark}
One easily checks that for  cubic threefolds (\ref{eq-cubic7}) coincides with \cite[Lemma 11.27]{CG72}.
\end{remark}

\begin{theorem}\label{thm-cubic}\quad\\
(i) For any basis $\{\alpha_i\}_{i=1}^{m}$ of $H_{\mathrm{prim}}^{n}(X)$, let $\beta_{i}=\Psi(\alpha_i)$, the cohomology ring $H^{*}(\Omega_X)$ is generated by $\sigma_1$, $\sigma_2$, and $\beta_1,\dots,\beta_m$.\\
(ii) The equalities in Proposition \ref{prop-cubic2} (ii), (iii) hold.\\
(iii) The ring structure of $H^{*}(\Omega_X)$ is given by Proposition \ref{prop-FanoVar-SchubertCalculus} and Theorem \ref{thm-cubic5}, together with Pieri's rule on 
$G_{2}(\mathbb{C}^{n+2})$.
\end{theorem}
\begin{proof} Denote the homomorphism $H^{k}(G_{2}(\mathbb{C}^{n+2}))\rightarrow H^{k}(\Omega_X)$ by $j_k^*$.
By (\ref{eq-cubic0}) and Proposition \ref{prop-cubic2} (i), for $0\leq k<n-2$ or $n-2\leq k<2n-4$ and $k-n\neq 0\ \mathrm{mod}\ 2$, $j_k^*$ is an isomorphism. 
For $n-2\leq k<2n-4$ and $i-n= 0\ \mathrm{mod}\ 2$, by the monodromy argument and Theorem \ref{thm-cubic5} (i), (ii), there is no nontrivial relations among
$\mathrm{Im}(j_k^*)$ and $\{\beta_i\cup\sigma_{1}^{\frac{k-n+2}{2}}\}_{i=1}^{m}$, thus by (\ref{eq-cubic0}) they form a basis of $H^{k}(\Omega_X)$. 
For $k=2n-4$, by the proof of  Theorem \ref{thm-cubic5} (iii), it suffices to show that there is only one relation of the form (\ref{eq-cubic7.9}). But two distinct such relations would lead to a nontrivial relation among $\{n-2-k,k\}$,  $0\leq k\leq \lfloor\frac{n}{2}\rfloor$, which contradicts Proposition \ref{prop-cubic2} (i). For $2n-4<k\leq 4n-8$ the conclusion follows easily by similar monodromy arguments and Proposition \ref{prop-cubic2} (ii), (iii), and Theorem \ref{thm-cubic5}. Part (iii) is an immediate consequence of (i), (ii), and Theorem \ref{thm-cubic5}.
\end{proof}

\begin{remark}
Equip $\Omega_X$ with the polarization induced from $G_2(\mathbb{C}^{n+2})$. Then Theorem \ref{thm-cubic5} and Theorem \ref{thm-cubic} implies that  the homomorphism $\Psi$ induces an isomorphism  $H^{n}_{\mathrm{prim}}(X)\rightarrow H_{\mathrm{prim}}^{n-2}(\Omega_X)$, which recovers results in \cite{Iza99} and \cite{Shi90}.
\end{remark}

\begin{remark}
One can use the same method to study the Fano variety of lines on $X_n(2,2)$, the smooth complete intersections of two quadrics. When $n$ is odd, it is hopeful to recover $F^{(2)}(0)=1$ (Theorem \ref{thm-reconstructcubicandquadric}). When $n$ is even, there are two constants to be determined, rather than only one constant like in (\ref{eq-cubic15}),  as we will see at the beginning of  Section \ref{sec:4points-Invariants}. The method in this section can only produce one equation for the two constants.
To compute the two constants we need some additional data. By the $D_{n+3}$-symmetry of WDVV equations and some integrality reason, we compute these constants in \cite[Theorem 1.1]{Hu21}.
I expect that this will be enough to determine the cohomology ring of the Fano variety of lines on $X_n(2,2)$.
\end{remark}

\section{
Some genus 0 GW invariants
via reduced genus one invariants}\label{sec:genusReduction}
This section is devoted to computing certain genus 0 invariants of cubic hypersurfaces and intersections of two quadrics, which are not computed by the monodromy group method of Section \ref{sec:ReconstructionII-F(2)}. More precisely, we do the following:
\begin{enumerate}
        \item We compute $F^{(2)}(0)$ of cubic hypersurfaces, thus give another proof of Theorem \ref{thm-F(2)(0)-cubicHypersurfaces}. In the process we also compute $F^{(2)}(0)$ of odd dimension intersection of two quadrics, giving another proof of Theorem \ref{thm-higher10.1} in this case.
        \item We obtain a partial result on 4-point genus 0 GW invariants with primitive insertions of even dimension intersection of two quadrics, which are exceptional complete intersections according to Definition \ref{def-exceptional}. This result (Proposition \ref{prop-4points-sum}) is used in \cite{Hu21} to give a nearly complete reconstruction theorem in this case.
        \item We compute $F^{(4)}(0)$ of cubic 3-folds. In view of Theorem \ref{thm-reconstructcubicandquadric} (i), this yields  reconstruction for genus 0 GW invariants of cubic 3-folds.
  \end{enumerate}  

Recall that the genus reduction axiom (see e.g., \cite{CK99}) says that, the map $\phi: \overline{\mathcal{M}}_{g-1,n+2}\rightarrow \overline{\mathcal{M}}_{g,n}$ gluing the last two marked points induces the following identity of Gromov-Witten classes
\begin{eqnarray*}
\phi^{*}I_{g,n,\beta}(\alpha_1,\dots,\alpha_n)=\sum_{i,j}I_{g-1,n+2,\beta}(\alpha_1,\dots,\alpha_n, \alpha_i, g^{ij}\alpha_j).
\end{eqnarray*}
In fact one of our original motivations to study the Gromov-Witten invariants involving primitive classes is that such invariants inevitably occur even if one wants to compute the invariants in higher genera with only ambient insertions when using various tautological relations. \\

In this section we go in the 
inverse direction. By Proposition \ref{prop-initialValues-meaning}, and certain \emph{contraction of (permanent) Pfaffians}, the invariant $F^{(k)}(0)$ is determined by 
\begin{eqnarray}\label{eq-gr1}
\sum_{i_1,j_1}\cdots \sum_{i_k,j_k}\langle \gamma_{i_1}, g^{i_1 j_1}\gamma_{j_1},\dots, \gamma_{i_k}, g^{i_k j_k}\gamma_{j_k}\rangle_{0},
\end{eqnarray}
where the summations are taken over the primitive classes. So  one can compute $F^{(k)}(0)$ by topological recursion relations in genus $g>0$ and invariants in genus not greater than $g$ with less than $2k$ primitive insertions. If we have sufficiently many relations we can try to reduce (\ref{eq-gr1}) to higher genus invariants with only ambient insertions. But in general the higher genus invariants are very hard to compute. However, on some occasions the involved higher genus invariants are of lower degrees so that the computation is possible. As we will see, this is the case for the cubic hypersurfaces, where we can prove some vanishing theorems on the so called \emph{reduced genus 1 GW invariants} defined by \cite{Zin09}, and then use Zinger's \emph{standard versus reduced} formula in genus 1. For an application of the trick of this section, see \cite{Ke18}.

\subsection{Reduced genus 1 invariants and the Standard versus Reduced formula}
Let $X$ be a smooth projective scheme over $\mathbb{C}$ of dimension $n$. 
Let $k\in \mathbb{Z}_{\geq 0}$, and $\beta\in H_2(X;\mathbb{Z})$.  Let $\Mbar^0_{1,k}(X,\beta)$ be the stack of genus one stable maps $f:C\rightarrow X$ of degree $\beta$ into $X$ satisfying one of  the following two  conditions:
\begin{enumerate}
	\item[(i)] There is no subcurve of arithmetic genus one contracted to a point by $f$;
	\item[(ii)]  There exists a subcurve of arithmetic genus one contracted to a point by $f$. We denote by $C_0$ the unique maximal connected contracted subcurve of arithmetic genus one. There is a unique way to write $C$ as
	\begin{equation}\label{eq-def-reducedGenusOne-decompositionOfCurve}
	 	C=C_0\cup \bigcup_{i=1}^m C_m
	 \end{equation} 
	 such that $C_1,\dots,C_m$ are connected subcurves satisfying that $C_i\cap C_j=\emptyset$ for $1\leq i\neq j\leq m$, and $C_0\cap C_m=\{p_m\}$. Then $f(p_1)=\dots=f(p_m)=x$ for some $x\in X$. Select a nonzero tangent vector $\partial_i$ of $C_i$ at $p_i$, for $1\leq i\leq m$. The condition is that $f_*(\partial_1),\dots,f_*(\partial_m)$ are \emph{linearly dependent} at $x$. In other words, 
	 \begin{equation}\label{eq-reducedGenusOne-tangentCondition}
	 \mathrm{dim}\ \mathrm{span}\{f_*(T_{p_1}C_1),\dots,f_*(T_{p_m}C_m)\}\leq m-1.
	 \end{equation}
\end{enumerate}
Then $\Mbar^0_{1,k}(X,\beta)$ is a closed substack  of $\Mbar_{1,k}(X,\beta)$.
By \cite[Corollary 1.6]{Zin09}, $\Mbar^0_{1,k}(X,\beta)$ carries a virtual fundamental class $[\Mbar^0_{1,k}(X,\beta)]^{\mathrm{vir}}$ of  the same dimension as $[\Mbar_{1,k}(X,\beta)]^{\mathrm{vir}}$. The \emph{reduced genus 1} Gromov-Witten invariants of $X$ are defined using this class:
\begin{equation*}
	\langle \psi_1^{a_1}\gamma_1,\dots,\psi_k^{a_k}\gamma_k\rangle^0_{1,k,\beta}:=\int_{[\Mbar_{1,k}^0(X,\beta)]^{\mathrm{vir}}}\prod_{i=1}^{k}\psi_i^{a_i}\mathrm{ev}_i^*\gamma_i.
\end{equation*}

In \cite{Zin08} Zinger proved a comparison formula between the standard and the reduced genus 1 GW invariants, the so called \emph{SvR formula} for short. To state his formula, we need to introduce more ingredients.
\subsubsection{A variant of genus 0 GW invariants}
Let $m\in \mathbb{Z}_{>0}$ and $J$ a finite set, and $\beta\in H_2(X;\mathbb{Z})$. Define $\Mbar_{(m,J)}(X,\beta)$ by the cartesian diagram
\begin{equation*}
	\xymatrix{
	\Mbar_{(m,J)}(X,\beta) \ar[r] \ar[d] & \bigsqcup_{\begin{subarray}{c}
	\beta_1+\dots+\beta_m=\beta\\
	J_1\sqcup\dots\sqcup J_m=J\end{subarray}}\Big(\prod_{i=1}^m \Mbar_{0,0\sqcup J_i}(X,\beta_i)\Big) \ar[d]^{\prod_i \mathrm{ev}_0}\\
	X \ar[r]^{\Delta} & X^m
		}
\end{equation*}
where $\Mbar_{0,\{0\}\sqcup J_i}(X,\beta_i)$ means the stack of  stable maps of genus 0 of degree $\beta_i$ with marked smooth points marked by the set $J_i$, and a distinguished smooth point marked by $0$. One can think of a point of $\Mbar_{(m,J)}(X,\beta)$ as representing a variant of genus 0 stable maps, the domain curve allowing more than two components to meet at exactly one point. The virtual fundamental cycle on $\Mbar_{(m,J)}(X,\beta)$ is defined via Gysin pullback:
\begin{equation}\label{eq-virtualFunClass-(m,J)}
	[\Mbar_{(m,J)}(X,\beta)]^{\mathrm{vir}}=\Delta^{!}
	\sum_{\begin{subarray}{c}
	\beta_1+\dots+\beta_m=\beta\\
	J_1\sqcup\dots\sqcup J_m=J\end{subarray}}\Big(\prod_{i=1}^m [\Mbar_{0,0\sqcup J_i}(X,\beta_i)]^{\mathrm{vir}}\Big).
\end{equation}
The dimension of $[\Mbar_{(m,J)}(X,\beta)]^{\mathrm{vir}}$ is
\begin{equation*}
	2\big(n-2m+|J|+c_1(T_X)\cap \beta\big).
\end{equation*}
Let $\mathcal{C}_{(m,J)}(X,\beta)$ be the universal curve over $\Mbar_{(m,J)}(X,\beta)$, and $f: \mathcal{C}_{(m,J)}(X,\beta)\rightarrow X$ the universal stable morphism. With each $j\in J$ there is an associated section $\sigma_j:\Mbar_{(m,J)}(X,\beta)\rightarrow \mathcal{C}_{(m,J)}(X,\beta)$. Let $\psi_j$ be the pullback of the relative cotangent line bundle via $\sigma_j$. For $1\leq i\leq m$, the $i$-th projection
\[
\prod_{i=1}^m \Mbar_{0,0\sqcup J_i}(X,\beta_i)\rightarrow \Mbar_{0,0\sqcup J_i}(X,\beta_i)
\]
induces a projection
\begin{equation*}
	\Mbar_{(m,J)}(X,\beta)\xrightarrow{\pi_i} \bigsqcup_{\beta_i\leq \beta} \bigsqcup_{J_i\subset J}\Mbar_{0,0\sqcup J_i}(X,\beta_i),
\end{equation*}
where $\beta_i\leq \beta$ means that $\beta- \beta_i$ is an effective curve class. Let $\eta_p\in H^{2p}\big(\Mbar_{(m,J)}(X,\beta)\big)$ be the degree $2p$ term of 
\[
\prod_{i=1}^{i=m}\pi_i^*\frac{1}{1-\psi_0}.
\]
For $j\in J$, define the evaluation map $\mathrm{ev}_j=f\circ \sigma_j$. The $0$-th evaluation map is defined to be $\mathrm{ev}_0=\mathrm{ev}_0\circ \pi_i$, which is independent of the choice of $i$, where $1\leq i\leq m$. 
For 
\[
\mu=(a_1,\dots,a_k;\mu_1,\dots,\mu_k)\in \mathbb{Z}_{\geq 0}^k\times H^*(X;\mathbb{Q})^k,
\]
and $J\subset [k]=\{1,\dots,k\}$, define 
\[
\mu_J=\prod_{j\in J}\mu_j.
\]
Now for given $\mu$ as above, and $\mu_0\in H^*(X;\mathbb{Q})$, and $\eta_p$, Zinger's variant of genus 0 GW invariant is defined as
\begin{equation}\label{eq-def-invariants-(m,J)}
	\mathrm{GW}_{m,J}^{\beta}(\eta_p,\mu_0;\mu):=\frac{1}{m!}\int_{[\Mbar_{(m,[k]-J)}(X,\beta)]^{\mathrm{vir}}}
	\eta_p \mathrm{ev}_0^*(\mu_0 \mu_J)\prod_{j\not\in J}\psi_j^{a_j}\mathrm{ev}_j^* \mu_j.
\end{equation}
Note that
\begin{itemize}
	\item the relevant moduli space are stable maps with marked points indexed by  $[k]-J$;
	\item when some $\mu_j$ have odd degrees, one needs to take care of the signs. 
\end{itemize}
 
By (\ref{eq-virtualFunClass-(m,J)}), such invariants can be expressed as the standard genus zero Gromov-Witten invariants. In particular,
\begin{equation}\label{eq-def-invariants-(1,J)}
	\mathrm{GW}_{1,J}^{\beta}(\eta_p,\mu_0;\mu)
	=\langle \psi^p \mu_0 \prod_{j\in J}\mu_j, \prod_{j\not\in J}\psi_j^{a_j}\mathrm{ev}_j^* \mu_j\rangle_{0,1+[k]-J,\beta}.
\end{equation}
\subsubsection{Intersection numbers on blowups of \texorpdfstring{$\Mbar_{1,k}$}{Mbar1,k}}
The second is the intersection numbers on certain blowing up $\widetilde{\mathcal{M}}_{1,([m],J)}$ of $\Mbar_{1,[m]\sqcup J}$, defined in \cite{Zin07}. 

For a non-negative integer $m$ and a finite set $J$ satisfying $[m]\cap J=\emptyset$ and $m+|J|\geq 1$, $\widetilde{\mathcal{M}}_{1,([m],J)}$ is a smooth proper Deligne-Mumford stack
obtained by successively blowing-up  of $\Mbar_{1,[m]\sqcup J}$ along certain natural substacks. Denote by $\pi$ the morphism $\widetilde{\mathcal{M}}_{1,([m],J)}\rightarrow \Mbar_{1,[m]\sqcup J}$. There is a distinguished line bundle $\tilde{\mathbb{E}}$ on $\widetilde{\mathcal{M}}_{1,([m],J)}$, which is obtained by certain twisting of the pullback of the Hodge bundle $\mathbb{E}$ on $\Mbar_{1,[m]\sqcup J}$ by certain exceptional divisors. For non-negative integers $\tilde{a}$ and $a_1,\dots,a_{|J|}$, define
\begin{equation*}
	(\tilde{a};\{a_j\}_{j\in J})_{m,J}=\int_{\widetilde{\mathcal{M}}_{1,([m],J)}}
	\big(c_1(\tilde{\mathbb{E}})\big)^{\tilde{a}}\cup \prod_{j\in J}\pi^* \psi_j^{a_j}.
\end{equation*}
 We do not recall the precise construction of $\widetilde{\mathcal{M}}_{1,([m],J)}$, because we need only know that these intersection numbers  are determined by the following recursions (R0)-(R3). In (R1)-(R3), we assume $m+|J|\geq 2$.
\begin{enumerate}
	\item[(R0)] If $\tilde{a}+\sum_{j\in J} a_j\neq m+|J|$, then
	\[
	(\tilde{a};\{a_j\}_{j\in J})_{m,J}=0.
	\]
	Moreover
	\begin{gather}\label{eq-R0}\tag{R0}
	(1;)_{1,\emptyset}=(0;1)_{0,\{*\}}=\frac{1}{24}.
	\end{gather}
	\item[(R1)] If  $m>0$ and $a_j>0$ for all $j\in J$, then
	\begin{gather}\label{eq-R1}\tag{R1}
	(\tilde{a};\{a_j\}_{j\in J})_{m,J}=(\tilde{a};\{a_j\}_{j\in J},0)_{m-1,|J|\sqcup \{*\}}.
	\end{gather}
	\item[(R2)] If $a_{j^*}=1$ for some  $j\in J$, then
	\begin{gather}\label{eq-R2}\tag{R2}
	(\tilde{a};\{a_j\}_{j\in J})_{m,J}=(m+|J|-1)\cdot(\tilde{a};\{a_j\}_{j\in J-\{j^*\}})_{m,J-\{j^*\}}.
	\end{gather}
	\item[(R3)] If $a_{j^*}=0$ for some $j^*\in J$, then
	\begin{gather}\label{eq-R3}\tag{R3}
	(\tilde{a};\{a_j\}_{j\in J})_{m,J}\\
	=m\cdot (\tilde{a}-1;\{a_j\}_{j\in J-\{j^*\}})_{m,J-\{j^*\}}
	+\sum_{j\in J-\{j^*\}}(\tilde{a};a_j-1,\{a_j\}_{j\in J-\{j^*,j\}})_{m,J-\{j^*\}},\nn
	\end{gather}
	where by convention $(\tilde{a};\{a_j\}_{j\in J})_{m,J}=0$ if $\tilde{a}<0$ or $a_j<0$ for some $j\in J$.
\end{enumerate}
In this paper we will use only (\cite[Corollary 1.2]{Zin07}),
\begin{equation}\label{eq-intersectionNum-M(1,m,J)-1}
(m+|J|;0,\dots,0)_{m,J}=\frac{m^{|J|}\cdot (m-1)!}{24},
\end{equation}
 and
\begin{eqnarray}\label{eq-intersectionNum-M(1,m,J)-2}
&&(m+|J|-1;1,0,\dots,0)_{m,J}\nn\\
&\stackrel{\mbox{by (R2)}}{=}& (m+|J|-1)\cdot (m+|J|-1;0,\dots,0)_{m,|J|-1}\nn\\
&=&\frac{(m+|J|-1)\cdot m^{|J|-1}\cdot (m-1)!}{24}.
\end{eqnarray}

\subsubsection{The SvR formula}

For 
\[
\mu=(a_1,\dots,a_k;\mu_1,\dots,\mu_k)\in \mathbb{Z}_{\geq 0}^k\times H^*(X;\mathbb{Q})^k,
\]
let 
\[
p_{J}(\mu)=\sum_{j\in J}a_j,\
d_{m,J}(\mu)=n-2m-|J|+p_J(\mu).
\]
The following is Zinger's SvR formula (\cite[Theorem 1A]{Zin08}).
\begin{theorem}\label{thm-SvR}
\begin{gather}
 \langle \psi^{a_1}\mu_1,\dots,\psi^{a_k} \mu_k\rangle_{1,\beta}
-\langle \psi^{a_1}\mu_1,\dots,\psi^{a_k} \mu_k\rangle_{1,\beta}^0\nn\\
= \sum_{m=1}^{\infty}\sum_{J\subset [k]}\Big(
(-1)^{m+|J|-p_J(\mu)}\big(m+|J|-p_J(\mu);(a_j)_{j\in J}\big)_{[m],J}\nn\\
 \times \sum_{p=0}^{d_{m,J}(\mu)}
\mathrm{GW}_{(m,J)}^{\beta} \big(\eta_p,c_{d_{m,J}(\mu)-p}(TX);\psi^{a_1}\mu_1,\dots,\psi^{a_k} \mu_k\big)
\Big).\label{eq-SvR}\tag{SvR}
\end{gather}
\end{theorem}

\subsection{Vanishing of certain lower degree reduced genus 1 invariants}\label{sec:vanishing-reducedGenus1}
The following trivial vanishing result enables us to compute $F^{(2)}(0)$ of cubic hypersurfaces.
\begin{proposition}\label{prop-vanishing-degreeOne}
Let $X$ be a smooth closed subscheme of $\mathbb{P}^N$. Let $\beta\in H_2(X;\mathbb{Z})$ such that $\sfh\cdot \beta=1$, where $\sfh$ is the hyperplane class restricted to $X$. Then any reduced genus one invariant of degree $\beta$ is 0.
\end{proposition}
\begin{proof}
Suppose that $f:C\rightarrow X$ is a genus one stable map that  lies in $\Mbar^0_{1,k}(X,\beta)$ for some $k$. Since any non-constant map from a curve of arithmetic genus 1 into $\mathbb{P}^N$ has degree $>1$, $f$ must be the type (ii) in the definition of $\Mbar^0_{1,k}(X,\beta)$. Moreover in the decomposition (\ref{eq-def-reducedGenusOne-decompositionOfCurve}), $m$ must be 1, and $C_1$ maps to $X$ with degree 1. Let $C_1^{\circ}$ be the irreducible component of $C_1$ intersecting $C_0$ at $p_1$. Then $f|_{C_1^{\circ}}$ is an embedding of $C_1^{\circ}$ as a line in $\mathbb{P}^N$. But such a map does not satisfy the condition (\ref{eq-reducedGenusOne-tangentCondition}). Hence $\Mbar^0_{1,k}(X,\beta)$ is empty, and the conclusion follows.
\end{proof}
According to Theorem \ref{thm-reconstructcubicandquadric} we need also to compute $F^{(4)}(0)$ of cubic 3-folds. For this, we need some vanishing results of degree 2 reduced genus 1 invariants.
\begin{lemma}\label{lem-genusone-degreetwo-stablemaps}
Let $f:C\rightarrow \mathbb{P}^N$ be a genus $1$ stable map  of degree $2$ with marked points, satisfying one of the following two conditions:
\begin{enumerate}
   \item[(i)] no connected subcurve of arithmetic genus one is contracted;
   \item[(ii)] there is a contracted subcurve $C_0$ of arithmetic genus one, and  the images of the tangent vectors of the non-contracted rational components  at the attaching nodes are linearly dependent.
 \end{enumerate} 
  Then $h^0(C, f^* \mathcal{O}(1))=2$, $h^1(C, f^* \mathcal{O}(1))=0$, and the image of $f$ is a line.
\end{lemma}
\begin{proof}
Case (i):
let $C_0\subset C$ be the minimal connected subcurve of arithmetic genus one. Then the degree of   $f^* \mathcal{O}(1)|_{C_0}$ is $2$. If $C_0$ is a smooth curve of genus one, by Riemann-Roch and Serre duality we have $h^0(C_0, f^* \mathcal{O}(1))=2$. If the irreducible components of $C_0$ are rational curves, we contract the components of $C_0$ that is contracted by $f$, thus obtain a semistable curve $C'_0$. Then $C'_0$ is either an irreducible rational nodal curve of arithmetic genus 1, or a rational nodal curve with two components. In either case it is straightforward to get $h^0(C_0, f^* \mathcal{O}(1))=h^0(C'_0, f^* \mathcal{O}(1))=2$. Finally by Riemann-Roch for semistable curves, $h^1(C, f^* \mathcal{O}(1))=0$.

Case (ii): Let $C_0$ be the maximal connected contracted subcurve with $p_a(C_0)=1$. There are two subcases. In the first subcase, there exist two rational subcurves $C_1$ and $C_2$, such that $C=C_0\cup C_1\cup C_2$, and  $x_1=C_0\cap C_1$ and $x_2=C_0\cap C_2$ are two distinct points of $C_0$, and $\mathrm{deg}\big(f^*\mathcal{O}(1)|_{C_1}\big)=\mathrm{deg}\big(f^*\mathcal{O}(1)|_{C_2}\big)=1$. Thus $C_i$ maps onto to a line in $\mathbb{P}^n$, for $i=1$ or $2$. Let $\vartheta_i$ be a nonzero tangent vector of $C_i$ at $x_i$, for $i=1,2$. Then $f(x_1)=f(x_2)$, and the two tangent vectors $(df)_{x_i}(\vartheta_i)$ and $(df)_{x_2}(\vartheta_2)$ are linearly dependent, i.e. parallel. So $f(C_1)$ and $f(C_2)$ are the same line in $\mathbb{P}^n$.

In the second subcase, there is a rational subcurve $C_1$, such that $C=C_0\cup C_1$, $x_1=C_0\cap C_1$, and $\mathrm{deg}\big(f^*\mathcal{O}(1)|_{C_1}\big)=2$. Let $C_2\subset C_2$ be the irreducible component of $C_1$ that contains $x_1$. Let $\vartheta$ be a nonzero tangent vector of $C_1$ at $x_1$. Then $(df)_{x_1}(\vartheta)=0$. So $\mathrm{deg}\big(f^*\mathcal{O}(1)|_{C_2}\big)=2$, otherwise $C_2$ will be a line imbedded into $\mathbb{P}^N$, contradicting with $(df)_{x_1}(\vartheta)=0$.  Now if $f(C_2)$ is a degree $2$ curve in $\mathbb{P}^n$,  then $p_a(f(C_2))=2$ and $C_2$ is a normalization of $f(C_2)$. Thus $p_a(C_2)=p_a(f(C_2))$ implies that $C_2\cong f(C_2)$. This contradicts that $(df)_{x_1}(\vartheta)=0$. So $f(C_2)$ must be a line in $\mathbb{P}^n$.

Hence we have shown that in case (ii) $f(C)$ is a line. It follows that $h^0(C, f^* \mathcal{O}(1))=2$, and thus $h^1(C, f^* \mathcal{O}(1))=0$ by Riemann-Roch.
\end{proof}

\begin{proposition}\label{prop-degree2-reducedGenusOne-factorThru}
Let $X\subset \mathbb{P}^N$ be a smooth closed subscheme. Let $k\geq 1$. Set
\begin{equation*}
  \overline{\mathcal{M}}_{0,[k]}(X,1):=
  \underbrace{\overline{\mathcal{M}}_{0,1}(X,1)\times_{\overline{\mathcal{M}}_{0,0}(X,1)}\times\cdots\times_{\overline{\mathcal{M}}_{0,0}(X,1)} \overline{\mathcal{M}}_{0,1}(X,1)}_{k\ \mbox{factors}}.
\end{equation*}
 Then we have a factorization of the evaluation maps $\mathrm{ev}_{[k]}=\mathrm{ev}_1\times\cdots\times \mathrm{ev}_k$
\begin{equation*}
  \xymatrix{
  \overline{\mathcal{M}}^0_{1,k}(X,2) \ar[r]^>>>>>>{\mathrm{ev}_{[k]}} \ar[d]_{\Phi_k} & X^k \\
  \overline{\mathcal{M}}_{0,[k]}(X,1) 
  \ar[ur]_{\mathrm{ev}_{[k]}} 
  &
  }
\end{equation*}
\end{proposition}
\begin{proof}
Let $S\rightarrow \overline{\mathcal{M}}^0_{1,k}(X,2)$ be a morphism, which is induced by a family of stable maps, 
\begin{equation}\label{graph-reducedgenusone-stablemap}
  \xymatrix{
  \mathcal{C} \ar[r]^{f} \ar[d]_{\pi} & X\\
  S
  }
\end{equation}
where $\pi$ is flat, satisfying that every fiber map is of the form in the case (i) or (ii) of Lemma \ref{lem-genusone-degreetwo-stablemaps}.
Then by the conclusion of this lemma and the semicontinuity theorem, we have $R^1\pi_* f^* \mathcal{O}(1)=0$ and that $R^0 \pi_* f^* \mathcal{O}(1)$ is locally free of rank 2. Then the projective bundle $\mathcal{D}=\mathrm{Proj}\Big(S^\bullet\big(R^0 \pi_* f^* \mathcal{O}(1)\big)\Big)$ is a flat family of lines in $X$ over $S$. This construction is natural in (\ref{graph-reducedgenusone-stablemap}). So we obtain morphisms of moduli stacks of stable maps
\[
\overline{\mathcal{M}}^0_{1,k}(X,2)\xrightarrow{\mathrm{forget}} \overline{\mathcal{M}}^0_{1,0}(X,2)\xrightarrow{\varphi_0} \overline{\mathcal{M}}_{0,0}(X,1), 
\]
and $\mathcal{D}$ is the image of $\mathcal{C}$ by $f$. For $k=1$, a section $\sigma:S\rightarrow \mathcal{C}$ induces a section $S\rightarrow \mathcal{D}$, so we obtain a morphism of moduli stacks of stable maps
\[
\overline{\mathcal{M}}^0_{1,1}(X,2)\xrightarrow{\varphi_1} \overline{\mathcal{M}}_{0,1}(X,1).
\]
One can also regard $\varphi_1$ as the morphism of universal curves over $\varphi_0$. From  the commutative diagram of evaluation maps
\[
\xymatrix{
  \overline{\mathcal{M}}^0_{1,k}(X,2) \ar[r]^<<<<<{\mathrm{ev}_i} \ar[d]_{\mathrm{forget}_{\hat{i}}} & X \\
  \overline{\mathcal{M}}^0_{1,k}(X,2) \ar[ur]_{\mathrm{ev}_1}
}
\]
where $\mathrm{forget}_{\hat{i}}$ is the morphism forgetting the marked points except the $i$-th one, the conclusion follows.
\end{proof}

\begin{remark}
There exists a natural lifting of $\Phi_k$ to a morphism $\widetilde{\Phi}_k:\overline{\mathcal{M}}^0_{1,k}(X,2)\rightarrow \overline{\mathcal{M}}_{0,k}(X,1)$. We sketch a construction of $\widetilde{\Phi}_k$. 
Let  $q:S\rightarrow\overline{\mathcal{M}}^0_{1,k}(X,2)$ be an \'{e}tale chart. The sections $\sigma_1,\dots,\sigma_{k}$ over $S$ induces  sections $\varsigma_1,\dots,\varsigma_{k}$ on $\mathcal{D}\rightarrow S$. For $1\leq i\neq j\leq k$, the locus $\{\varsigma_i=\varsigma_j\}$ has codimension one in $S$,  and thus the locus $\varsigma_i\cap \varsigma_j$ has codimension 2 in $\mathcal{D}$, for if two marked points $x_i$ and $x_j$ in $C$ maps to the same point $x \in f(C)$, we can deform this stable map by just moving $x_j$ while fixing $x_i$. Then blow up this locus $\varsigma_i\cap \varsigma_j$ in $\mathcal{D}$, for all unordered pairs $\{i,j\}$, $1\leq i\neq j\leq k$, we obtain a family of degree one and genus one stable maps $\widetilde{\mathcal{D}}\rightarrow X$ over $S$. We are left to show $\widetilde{\mathcal{D}}\rightarrow S$  is flat. This can be done by a local computation, but it has essentially been done in the construction of $\overline{M}_{0,k}$ and its modular interpretation.
\end{remark}

\begin{theorem}\label{thm-vanishing-deg2-cubics}
Let $X$ be a cubic hypersurface in $\mathbb{P}^N$. Let $\alpha_1,\dots,\alpha_k\in H^*(X)$. Then 
\begin{equation}\label{eq-vanishing-deg2-cubics}
\langle \alpha_1,\dots,\alpha_k\rangle_{1,2}^0=0=\langle \psi \alpha_1,\alpha_2,\dots,\alpha_k\rangle_{1,2}^0.
\end{equation}
\end{theorem}
\begin{proof}
By \cite[Theorem 7.2]{CG72}, the Fano variety of lines $\overline{\mathcal{M}}_{0,0}(X,1)$ is smooth and its  dimension is $2N-6$.  So $\dim\ \overline{\mathcal{M}}_{0,k}(X,1)=2N-6+k $. But $\mathrm{vir.dim}\ \overline{\mathcal{M}}^0_{1,k}(X,2)=2N-4+k$. We use the morphism $\Phi_k$ defined in Proposition \ref{prop-degree2-reducedGenusOne-factorThru}. For the dimension reason one gets
\[
\Phi_{k*}[\overline{\mathcal{M}}^0_{1,k}(X,2)]^{\mathrm{vir}}=0=
\Phi_{k*}\big(\psi\cap [\overline{\mathcal{M}}^0_{1,k}(X,2)]^{\mathrm{vir}}).
\]
Thus by the projection formula,
\begin{eqnarray*}
&&\langle \alpha_1,\dots,\alpha_k\rangle_{1,2}^0\\
&=& \mathrm{deg}\big(\mathrm{ev}_{[k]}^*(p_1^*\alpha_1\cup\cdots\cup p_k^* \alpha_k)\cap [\overline{\mathcal{M}}^0_{1,k}(X,2)]^{\mathrm{vir}}\big)\\
&=& \mathrm{deg}\big(\Phi_k^*\mathrm{ev}_{[k]}^*(p_1^*\alpha_1\cup\cdots\cup p_k^* \alpha_k)\cap [\overline{\mathcal{M}}^0_{1,k}(X,2)]^{\mathrm{vir}}\big)\\
&=& \mathrm{deg}\big(\mathrm{ev}_{[k]}^*(p_1^*\alpha_1\cup\cdots\cup p_k^* \alpha_k)\cap \Phi_{k*}[\overline{\mathcal{M}}^0_{1,k}(X,2)]^{\mathrm{vir}}\big)\\
&=&0,
\end{eqnarray*}
and similarly
\[
\langle \psi \alpha_1,\alpha_2,\dots,\alpha_k\rangle_{1,2}^0
=\mathrm{deg}\big(\mathrm{ev}_{[k]}^*(p_1^*\alpha_1\cup\cdots\cup p_k^* \alpha_k)\cap \Phi_{k*}(\psi\cap [\overline{\mathcal{M}}^0_{1,k}(X,2)]^{\mathrm{vir}})\big)=0.
\]
\end{proof}

\subsection{Correlators of length 4 with only primitive insertions when Fano index \texorpdfstring{$=n-1$}{=n-1}}\label{sec:4points-Invariants}
In this subsection we consider  smooth complete intersections $X$ of dimension $n\geq 3$ and Fano index $\mathsf{a}(n,\mathbf{d})=n-1$. There are three kinds:
\begin{enumerate}
	\item[(i)] Cubic hypersurfaces;
	\item[(ii)] Odd dimensional intersection of two quadrics;
	\item[(iii)] Even dimensional intersection of two quadrics.
\end{enumerate}
According to Definition \ref{def-exceptional}, cases (i) and (ii) are non-exceptional, while the third is exceptional. When $X$ is non-exceptional the genus zero  4-point invariants with only primitive insertions are determined by $F^{(2)}(0)$ by Proposition \ref{prop-initialValues-meaning}. 

The exceptional case, i.e. $X=X_n(2,2)$ with $n$ even, is more complicated. For in this case the monodromy group is a finite group $D_{n+3}$,  and has a larger invariant subring in the ring of the variables $t^{n+1},\dots,t^{n+m}$ dual to a suitable orthonormal chose basis of $H^*_{\mathrm{prim}}(X)$ . This subring  is generated by
\begin{equation}\label{eq-invariantsOf-typeD-1}
      s_{i}=\frac{1}{(2i)!}\sum_{j=n+1}^{2n+3}(t^j)^{2i},\ \mbox{for}\ 1\leq i\leq n+2,
\end{equation}
and
\begin{equation}\label{eq-invariantsOf-typeD-2}
      s_{n+3}=\prod_{j=n+1}^{2n+3}t^j.
\end{equation}
In particular, in this case the genus zero 4-point invariants with only primitive insertions are determined by 
\begin{equation}\label{eq-initialValues-even(2,2)}
\frac{\partial^2 F}{(\partial s_1)^2}(0)\ \mbox{and}\ \frac{\partial F}{\partial s_2}(0).
\end{equation}
(For more details we refer the reader to \cite[Section 2.2]{Hu21}; in the following we will make no explicit use of this paragraph.)

To compute $F^{(2)}(0)$ in the cases (i) and (ii), and the values (\ref{eq-initialValues-even(2,2)}) in the case (iii), are equivalent to compute the 4-point invariants with primitive insertions. 
The case (ii) is done in Theorem \ref{thm-higher10.1} and the case (i) is done in Theorem \ref{thm-cubic5}. 
In the following of this subsection, using (\ref{eq-SvR}), we show an identity in Proposition \ref{prop-4points-sum}  on the correlators of length 4 in the cases (i), (ii) and (iii) in a uniform way. This will determine $F^{(2)}(0)$ in the cases (i) and (ii) and thus all the 4-point invariants, while the exceptional case (iii) still needs some additional ad hoc treatment which will be completed in a separate paper \cite{Hu21}.

\subsubsection{From genus 1 to genus 0}
Let $\gamma_i=\sfh_i$ the $i$-th power of the hyperplane class for $0\leq i\leq n$, 
and  $\gamma_{n+1},\dots,\gamma_{n+m}$ a basis of $H^*_{\mathrm{prim}}(X)$. The Poincar\'{e} pairing is still denoted by $g_{ij}=(\gamma_i,\gamma_j)=\int_X \gamma_i\cup \gamma_j$. 
By the genus 1 topological recursion relation (\ref{eq-TRR1}), for $n+1\leq b,c\leq n+m$,
\begin{eqnarray}\label{eq-apply-TRR1-FanoIndex=n-1-1}
&&\langle \psi\gamma_b,\gamma_c\rangle_{1,1}\nn\\
&=&\frac{1}{\prod_{i=1}^r d_i}\sum_{e=0}^{n+m}\sum_{f=0}^{n+m}\langle \gamma_b,\gamma_c,\gamma_e\rangle_{0,1}g^{ef}\langle \gamma_f \rangle_{1,0}
+\frac{1}{\prod_{i=1}^r d_i}\sum_{e=0}^{n+m}\sum_{f=0}^{n+m}\langle \gamma_b,\gamma_c,\gamma_e\rangle_{0,0} g^{ef} \langle  \gamma_f \rangle_{1,1}\nn\\
&&+\frac{1}{24}\sum_{e=0}^{n+m}\sum_{f=0}^{n+m}\langle \gamma_b, \gamma_e, g^{ef}\gamma_f,\gamma_c\rangle_{0,1}.
\end{eqnarray}
By (\ref{eq-Dim}), (\ref{eq-Deg0}), and Corollary \ref{cor-mono3}, most of the invariants on the RHS are zero. Collecting the possibly nonzero ones, we get
\begin{eqnarray}\label{eq-apply-TRR1-FanoIndex=n-1-2}
&&\langle \psi\gamma_b,\gamma_c\rangle_{1,1}\nn\\
&=&\frac{1}{\prod_{i=1}^r d_i}\langle \gamma_b,\gamma_c,\sfh_{n-1}\rangle_{0,1}\langle \sfh \rangle_{1,0}
+\frac{1}{\prod_{i=1}^r d_i}\langle \gamma_b,\gamma_c,1\rangle_{0,3,0}\langle \sfh_n \rangle_{1,1,1}\nn\\
&&+\frac{1}{24}\sum_{e=0}^{n+m}\sum_{f=0}^{n+m}\langle \gamma_b, \gamma_e, g^{ef}\gamma_f,\gamma_c\rangle_{0,1}.
\end{eqnarray}
In the following applications of (\ref{eq-TRR1}), such arguments will be abbreviated and we will display directly the resulted equations.
\subsubsection{Applications of SvR}\label{sec:application-of-SvR-4points}
We are going to use (\ref{eq-SvR}) to compute 
$\langle \psi\gamma_b,\gamma_c\rangle_{1,1}$ and $\langle \sfh_n \rangle_{1,1}$. 
By Proposition \ref{prop-vanishing-degreeOne}, the reduced genus 1 invariants are 0, so  (\ref{eq-SvR}) yields, using (\ref{eq-def-invariants-(1,J)}), 
\begin{eqnarray*}
&&\langle \psi\gamma_b, \gamma_c\rangle_{1,1}=\sum_{J\subset [2]}\Big((-1)^{1+|J|-p_J}(1+|J|-p_J;(c_j)_{j\in J})_{[1],J}\\
&&\times\sum_{p=0}^{n-2-|J|+p_J}
\langle \psi^p \mathrm{ev}_0^*(c_{n-2-|J|+p_J-p}(TX)\mu_{J}),\prod_{j\not\in J}\psi_j^{c_j}\mathrm{ev}_{j}^* (\mu_j)\rangle_{0,1+[2]-|J|,1}\Big).
\end{eqnarray*}
Using (\ref{eq-intersectionNum-M(1,m,J)-1}) and (\ref{eq-intersectionNum-M(1,m,J)-2}), and that $\sfh\cup \gamma_b=0$, we obtain
\begin{eqnarray}\label{eq-psi-gammab-gammac-SvR}
&&\langle \psi\gamma_b, \gamma_c\rangle_{1,1}\nn\\
&=&-\frac{1}{24}\sum_{p=0}^{n-2}
\langle \psi^p c_{n-2-p}(T_X),\psi\gamma_b, \gamma_c\rangle_{0,1}
-\frac{1}{24}\sum_{p=0}^{n-2}
\langle \psi^p c_{n-2-p}(T_X)\gamma_b, \gamma_c\rangle_{0,1}\nn\\
&&+(-1)^{|\gamma_b|\cdot|\gamma_c|}\frac{1}{24}\sum_{p=0}^{n-3}
\langle \psi^p c_{n-3-p}(T_X)\gamma_c,\psi\gamma_b\rangle_{0,1}
+\frac{1}{12}\sum_{p=0}^{n-3}
\langle \psi^p c_{n-3-p}(T_X)\gamma_b\gamma_c\rangle_{0,1}\nn\\
&=&-\frac{1}{24}\sum_{p=0}^{n-2}
\langle \psi^p c_{n-2-p}(T_X),\psi\gamma_b, \gamma_c\rangle_{0,1}
-\frac{1}{24}\langle \psi^{n-2}\gamma_b, \gamma_c\rangle_{0,1}\nn\\
&&+(-1)^{|\gamma_b|\cdot|\gamma_c|}\langle \psi^{n-3}\gamma_c,\psi\gamma_b\rangle_{0,1}
+\frac{g_{bc}}{12\prod_{i=1}^r d_i}\langle \psi^{n-3}\sfh_{n}\rangle_{0,1}.
\end{eqnarray}

\begin{lemma}\phantomsection\label{lem-someGenus0Vanishing}
\begin{enumerate}
	\item[(i)] For $\gamma_b,\gamma_c\in \mathrm{H}^*_{\mathrm{prim}}(X)$,
	\begin{equation}\label{eq-someGenus0Vanishing-1}
	\langle \psi^{n-2}\gamma_b, \gamma_c\rangle_{0}=0.
	\end{equation}
	\item[(ii)] For $\alpha_1,\alpha_2\in \mathrm{H}^*(X)$, and $\gamma\in \mathrm{H}^*_{\mathrm{prim}}(X)$, and $a_1,a_2\in \mathbb{Z}_{\geq 0}$,
	\begin{equation}\label{eq-someGenus0Vanishing-2}
	\langle \psi \gamma,\psi^{a_1}\alpha_1,\psi^{a_2}\alpha_2\rangle_{0,1}=0.
	\end{equation}
	\item[(iii)] For $\gamma_b,\gamma_c\in \mathrm{H}^*_{\mathrm{prim}}(X)$, and $a\in \mathbb{Z}_{\geq 0}$, 
	\begin{equation}\label{eq-someGenus0Vanishing-3}
	\langle \psi\gamma_b, \psi^{a}\gamma_c\rangle_{0}=0.
	\end{equation}
\end{enumerate}
\end{lemma}
\begin{proof}
(i) Since $\sfh\cup \gamma_{b}=0$, by (\ref{eq-Div}) we have
\begin{eqnarray*}
\langle \gamma_{b}\psi^{n-2}, \gamma_c, \sfh\rangle_{0,1}=\langle \gamma_{b} \psi^{n-2}, \gamma_c\rangle_{0,1}.
\end{eqnarray*}
By (\ref{eq-TRR0}), 
\begin{equation}\label{eq-someGenus0Vanishing-1-TRR0}
\langle \gamma_{b}\psi^{n-2}, \gamma_c, \sfh\rangle_{0,1}
=\sum_{\mu=0}^{n+m}\sum_{\nu=0}^{n+m}\langle \gamma_{b}\psi^{n-3}, \gamma_{\mu}\rangle_{0} g^{\mu\nu}
\langle \gamma_{\nu}, \gamma_c, \sfh\rangle_{0}.
\end{equation}
By Corollary \ref{cor-mono3}, 
$\langle \gamma_\nu,\gamma_c\rangle_{0}=0$ when $\mu\leq n$, i.e. when $\gamma_{\nu}$ is an ambient class. If $\gamma_{\nu}$ is a primitive class, $\langle \gamma_\nu,\gamma_c\rangle_{0}=0$ for the degree reason because the Fano index is $n-1>1$. 
Then by (\ref{eq-Div}), $\langle\gamma_{\nu}, \gamma_c, \sfh\rangle_{0}=0$, and thus (\ref{eq-someGenus0Vanishing-1}) follows from (\ref{eq-someGenus0Vanishing-1-TRR0}).

(ii) By (\ref{eq-TRR0}),
\begin{equation}\label{eq-someGenus0Vanishing-2-TRR0}
	\langle \psi \gamma,\psi^{a_1}\alpha_1,\psi^{a_2}\alpha_2\rangle_{0}
	=\sum_{\mu=0}^{n+m}\sum_{\nu=0}^{n+m}\langle \gamma,\gamma_\mu\rangle_{0}g^{\mu \nu}\langle \gamma_\nu,\psi^{a_1}\alpha_1,\psi^{a_2}\alpha_2\rangle_{0}.
\end{equation}
As we see above, $\langle \gamma,\gamma_\mu\rangle_{0}=0$, so (\ref{eq-someGenus0Vanishing-2}) follows.

(iii) Since $\sfh\cup \gamma_b=\sfh\cup \gamma_c=0$, (\ref{eq-someGenus0Vanishing-3}) follows by (\ref{eq-Div}) from
\begin{equation}\label{eq-someGenus0Vanishing-3-1}
	\langle \psi\gamma_b, \psi^{a}\gamma_c,\sfh\rangle_{0}=0.
	\end{equation}
which is a special case of (ii).	
\end{proof}
\begin{lemma}
\begin{eqnarray}\label{eq-h(n)-psi(n-3)}
\langle \psi^{n-3}\sfh_{n}\rangle_{0,1}=\prod_{i=1}^r (d_i!)d_i.
\end{eqnarray}
\end{lemma}
\begin{proof}
This is extracted from the mirror theorem (\ref{eq-mirrorFormula-a(n,d)>1}).
\end{proof}

\begin{corollary}\label{cor-psi-gammab-gammac}
\begin{equation}\label{eq-psi-gammab-gammac}
\langle \psi\gamma_b, \gamma_c\rangle_{1,1}=\frac{\elld g_{bc}}{12}.
\end{equation}
\end{corollary}
\begin{proof}
This follows from (\ref{eq-psi-gammab-gammac-SvR}), Lemma \ref{lem-someGenus0Vanishing} and (\ref{eq-h(n)-psi(n-3)}).
\end{proof}

Similarly by Proposition \ref{prop-vanishing-degreeOne}, (\ref{eq-SvR}) and (\ref{eq-def-invariants-(1,J)}) we get
\begin{eqnarray*}
&&\langle \sfh_n\rangle_{1,1}\\
&=&-(1;)_{1,\emptyset}\sum_{p=0}^{n-2} \mathrm{GW}^1_{(1,\emptyset)}(\eta_p,c_{n-2-p}(T_X);\sfh_n)
+(2;0)_{1,1}\sum_{p=0}^{n-3} \mathrm{GW}^1_{(1,1)}(\eta_p,c_{n-3-p}(T_X);\sfh_n)\\
&=&-(1;)_{1,\emptyset}\sum_{p=0}^{n-2}\langle \psi^p c_{n-2-p}(T_X), \sfh_n\rangle_{0,1}
+(2;0)_{1,1}\sum_{p=0}^{n-3}\langle \psi^p c_{n-3-p}(T_X)\sfh_n\rangle_{0,1}.
\end{eqnarray*}
Using (\ref{eq-intersectionNum-M(1,m,J)-1})  and that $\sfh_n\cup \sfh_i=0$ for $i>0$, we obtain
\begin{equation}\label{eq-hn-(1,1,1)-SvR}
\langle \sfh_n\rangle_{1,1}=-\frac{1}{24}\sum_{p=0}^{n-2}
\langle \psi^{p} c_{n-2-p}(T_X),\sfh_n\rangle_{0,1}
+\frac{1}{24}\langle \psi^{n-3}\sfh_n\rangle_{0,1}.
\end{equation}
\begin{lemma}\label{lem-length2-invariant-1}
For $0\leq i,j\leq n$ and $i+j\leq 2n-2$, 
\begin{multline}\label{eq-length2-invariant-1}
\langle \psi^{2n-2-i-j}\sfh_i, \sfh_j\rangle_{0,1}=(-1)^{n-i}\prod_{i=1}^r d_i\Big[\binom{2n-2-i-j}{n-i}\ell_0\\
-\binom{2n-2-i-j}{n-1-i}\ell_1+\binom{2n-2-i-j}{n-j}\ell_0\Big],
\end{multline}
where by convention $\binom{x}{k}=0$ for $k\in\mathbb{Z}_{<0}$.
\end{lemma}
\begin{proof} By (\ref{eq-Div}) and (\ref{eq-TRR0}),
\begin{eqnarray*}
&&\langle \psi^{2n-2-i-j}\sfh_i, \sfh_j\rangle_{0,1}\\
&=& \langle \psi^{2n-2-i-j}\sfh_i, \sfh_j, \sfh\rangle_{0,1}-\langle \psi^{2n-3-i-j}\sfh_{i+1}, \sfh_j\rangle_{0,1}\\
&=&\sum_{\mu}\sum_{\nu} \langle \psi^{2n-3-i-j}\sfh_i, \gamma_{\mu}\rangle_{0,1}g^{\mu\nu}\langle \gamma_{\nu}, \sfh_j, \sfh\rangle_{0,0}
-\langle \psi^{2n-3-i-j}\sfh_{i+1}, \sfh_j\rangle_{0,1}\\
&=&\langle \psi^{2n-3-i-j}\sfh_{i}, \sfh_{j+1}\rangle_{0,1}-\langle \psi^{2n-3-i-j}\sfh_{i+1}, \sfh_j\rangle_{0,1}.
\end{eqnarray*}
Since $\langle \sfh_{n}, \sfh_{n-2}\rangle_{0,1}=\ell_0\prod_{i=1}^r d_i$, and $\langle \sfh_{n-1}, \sfh_{n-1}\rangle_{0,1}=\ell_1\prod_{i=1}^r d_i$, by an easy induction on $i+j$ from $2n-2$ to $0$ we obtain (\ref{eq-length2-invariant-1}).
\end{proof}
\subsubsection{Some calculations of residues}
\begin{lemma}\label{lem-hn-(1,1,1)}
\begin{eqnarray}\label{eq-hn-(1,1,1)}
\langle \sfh_n\rangle_{1,1}=
\begin{cases}
\frac{-(-2)^{n+2}-9n^2-3n+58}{72},&
 \mbox{if}\ \mathbf{d}=3;\\
\frac{(-1)^{n+1}(2n+1)-2n^2-2n+17}{24},&
 \mbox{if}\ \mathbf{d}=(2,2).
\end{cases} 
\end{eqnarray}
\end{lemma}
\begin{proof}
By the $j=n$ case of Lemma \ref{lem-length2-invariant-1},
\begin{eqnarray*}
&&\sum_{p=0}^{n-2}
\langle \psi^{p} c_{n-2-p}(T_X),\sfh_n\rangle_{0,1}\\
&=& \ell_0\prod_{i=1}^r d_i\sum_{p=0}^{n-2}(-1)^p \mathrm{Coeff}_{x^{n-2-p}}\Big(\frac{(1+x)^{n+r+1}}{\prod_{i=1}^r (1+d_i x)}\Big)\\
&=&\ell_0\prod_{i=1}^r d_i\cdot \mathrm{Coeff}_{x^{n-2}}\Big(\frac{1}{1+x}\cdot\frac{(1+x)^{n+r+1}}{\prod_{i=1}^r (1+d_i x)}\Big)\\
&=&\ell_0\prod_{i=1}^r d_i\cdot \mathrm{Res}_{x=0}\Big(\frac{(1+x)^{n+r}}{x^{n-1}\prod_{i=1}^r (1+d_i x)}\Big),
\end{eqnarray*}
where the notation $\mathrm{Coeff}_{x^i}\big(f(x)\big)$ denotes the coefficient of $x^i$ in the Laurent expansion of $f(x)$ at $0$. By the residue theorem,
\begin{eqnarray*}
&&\mathrm{Res}_{x=0}\Big(\frac{(1+x)^{n+r}}{x^{n-1}\prod_{i=1}^r (1+d_i x)}\Big)\\
&=&-\sum_{i=1}^{r}\mathrm{Res}_{x=-\frac{1}{d_i}}\Big(\frac{(1+x)^{n+r}}{x^{n-1}\prod_{i=1}^r (1+d_i x)}\Big)-\mathrm{Res}_{x=\infty}\Big(\frac{(1+x)^{n+r}}{x^{n-1}\prod_{i=1}^r (1+d_i x)}\Big).
\end{eqnarray*}
We compute the residues respectively:
\begin{eqnarray}\label{eq-residue-1}
&&\mathrm{Res}_{x=\infty}\Big(\frac{(1+x)^{n+r}}{x^{n-1}\prod_{i=1}^r (1+d_i x)}\Big)\nn\\
&=&-\frac{1}{\prod_{i=1}^r d_i}\Big(\binom{n+r}{2}-(n+r)\sum_{i=1}^r\frac{1}{d_i}
+\sum_{i=1}^r \frac{1}{d_i^2}+\sum_{1\leq i\neq j\leq r}\frac{1}{d_i d_j}
\Big)\nn\\
&=& \begin{cases}
-\frac{3n^2+n-2}{18}-\frac{1}{27}, & \mbox{if}\ \mathbf{d}=3;\\
-\frac{2n^2+2n-1}{16}, & \mbox{if}\ \mathbf{d}=(2,2),
\end{cases}
\end{eqnarray}
and
\begin{eqnarray}\label{eq-residue-3}
\sum_{i=1}^{r}\mathrm{Res}_{x=-\frac{1}{d_i}}\Big(\frac{(1+x)^{n+r}}{x^{n-1}\prod_{i=1}^r (1+d_i x)}\Big)
=\begin{cases}
\frac{(-2)^{n+1}}{27},&
 \mbox{if}\ \mathbf{d}=3;\\
\frac{(-1)^{n-1}(2n+1)}{16}, &
 \mbox{if}\ \mathbf{d}=(2,2).
\end{cases}  
\end{eqnarray}

So
\begin{equation}\label{eq-hn-(1,1,1)-SvR-firstGroupOfTerms}
\sum_{p=0}^{n-2}
\langle \psi^{p} c_{n-2-p}(T_X),\sfh_n\rangle_{0,1}
=\begin{cases}
\frac{(-2)^{n+2}}{3}+3n^2+n-\frac{4}{3},&
 \mbox{if}\ \mathbf{d}=3;\\
(-1)^n(2n+1)+2n^2+2n-1,&
 \mbox{if}\ \mathbf{d}=(2,2).
\end{cases}
\end{equation}
Putting (\ref{eq-h(n)-psi(n-3)}), (\ref{eq-hn-(1,1,1)-SvR-firstGroupOfTerms}) into (\ref{eq-hn-(1,1,1)-SvR}) we obtain (\ref{eq-hn-(1,1,1)}).
\end{proof}

\begin{lemma}\label{lem-h-(1,1,0)}
\begin{equation}\label{eq-h-(1,1,0)}
\langle \sfh \rangle_{1,0}
=\begin{cases}
\frac{1}{24}\big(\frac{(-2)^{n+2}-1}{9}-\frac{3n^2+7n+2}{6}\big),&
 \mbox{if}\ \mathbf{d}=3;\\
\frac{1}{24}\big(\frac{(-1)^n(2n+3)}{4}-\frac{2n^2+6n+3}{4}\big),
 &
 \mbox{if}\ \mathbf{d}=(2,2).
\end{cases} 
\end{equation}
\end{lemma}
\begin{proof}
By (\ref{eq-Deg0}) and the residue theorem,
\begin{eqnarray*}
&&\langle \sfh \rangle_{1,0}=-\frac{1}{24}\int_{X}\sfh\cup c_{n-1}(T_X)
=-\frac{\prod_{i=1}^r d_i}{24}\cdot \mathrm{Coeff}_{x^{n-1}}\Big(\frac{(1+x)^{n+r+1}}{\prod_{i=1}^r (1+d_i x)}\Big)\\
&=&-\frac{\prod_{i=1}^r d_i}{24}\cdot \mathrm{Res}_{x=0}\big(\frac{(1+x)^{n+r+1}}{x^n\prod_{i=1}^r (1+d_i x)}\big)\\
&=& \frac{\prod_{i=1}^r d_i}{24}\Big(
\sum_{i=1}^r \mathrm{Res}_{x=-\frac{1}{d_i}}\big(\frac{(1+x)^{n+r+1}}{x^n\prod_{i=1}^r (1+d_i x)}\big)+\mathrm{Res}_{x=\infty}\big(\frac{(1+x)^{n+r+1}}{x^n\prod_{i=1}^r (1+d_i x)}\big)\Big).
\end{eqnarray*}
Replacing $n$ by $n+1$ in the formulae of residues (\ref{eq-residue-1}) and (\ref{eq-residue-1}) we get
\begin{eqnarray*}
&&\sum_{i=1}^r \mathrm{Res}_{x=-\frac{1}{d_i}}\big(\frac{(1+x)^{n+r+1}}{x^n\prod_{i=1}^r (1+d_i x)}\big)+\mathrm{Res}_{x=\infty}\big(\frac{(1+x)^{n+r+1}}{x^n\prod_{i=1}^r (1+d_i x)}\big)\\
&=&\begin{cases}
\frac{(-2)^{n+2}}{27}-\frac{3(n+1)^2+(n+1)-2}{18}-\frac{1}{27},&
 \mbox{if}\ \mathbf{d}=3;\\
\frac{(-1)^{n}(2n+3)}{16}-\frac{2(n+1)^2+2(n+1)-1}{16}, &
 \mbox{if}\ \mathbf{d}=(2,2).
\end{cases} 
\end{eqnarray*}
So (\ref{eq-h-(1,1,0)}) follows.
\end{proof}

\begin{proposition}\label{prop-4points-sum}
Suppose $n+1\leq b,c\leq n+m$. Then
\begin{equation}\label{eq-4points-sum}
\sum_{e=n+1}^{n+m}\sum_{f=n+1}^{n+m}\langle \gamma_b, \gamma_e, g^{ef}\gamma_f,\gamma_c\rangle_{0,1}
=g_{bc}\cdot\begin{cases}
\frac{(-2)^{n+2}+8}{3},&
 \mbox{if}\ \mathbf{d}=3;\\
 (-1)^n(n+1)+2,
 &
 \mbox{if}\ \mathbf{d}=(2,2)
\end{cases}
\end{equation}
\end{proposition}
\begin{proof}
By (\ref{eq-psi-gammab-gammac}) and (\ref{eq-F1-leadingTerm-FanoIndex=n-1}), 
(\ref{eq-apply-TRR1-FanoIndex=n-1-2}) reads
\begin{equation*}\label{eq-apply-TRR1-FanoIndex=n-1-3}
\frac{\elld g_{bc}}{12}
=\frac{-\elld g_{bc}}{\prod_{i=1}^r d_i}\langle \sfh \rangle_{1,0}
+\frac{g_{bc}}{\prod_{i=1}^r d_i}\langle \sfh_n \rangle_{1,1}
+\frac{1}{24}\sum_{e=0}^{n+m}\sum_{f=0}^{n+m}\langle \gamma_b, \gamma_e, g^{ef}\gamma_f,\gamma_c\rangle_{0,1}.
\end{equation*}
Then by (\ref{eq-hn-(1,1,1)})  and (\ref{eq-h-(1,1,0)}),
\begin{eqnarray*}
&& \sum_{e=0}^{n+m}\sum_{f=0}^{n+m}\langle \gamma_b, \gamma_e, g^{ef}\gamma_f,\gamma_c\rangle_{0,1}\\
&=& 24g_{bc}\Big(\frac{\elld}{12}
+\frac{\elld}{\prod_{i=1}^r d_i}\langle \sfh \rangle_{1,0}
-\frac{1}{\prod_{i=1}^r d_i}\langle \sfh_n \rangle_{1,1}\Big)\\
&=&g_{bc}\Big(2\elld
+\frac{\elld}{\prod_{i=1}^r d_i}\cdot
\begin{cases}
\frac{(-2)^{n+2}-1}{9}-\frac{3n^2+7n+2}{6},&
 \mbox{if}\ \mathbf{d}=3;\\
\frac{(-1)^n(2n+3)}{4}-\frac{2n^2+6n+3}{4},
 &
 \mbox{if}\ \mathbf{d}=(2,2)
\end{cases}\\
&&-\frac{1}{\prod_{i=1}^r d_i}\cdot\begin{cases}
\frac{-(-2)^{n+2}-9n^2-3n+58}{3},&
 \mbox{if}\ \mathbf{d}=3;\\
(-1)^{n+1}(2n+1)-2n^2-2n+17,&
 \mbox{if}\ \mathbf{d}=(2,2)
\end{cases}\Big)\\
&=&g_{bc}\cdot\begin{cases}
\frac{(-2)^{n+2}}{3}-2n+\frac{14}{3},&
 \mbox{if}\ \mathbf{d}=3;\\
 (-1)^n(n+1)-n+3,
 &
 \mbox{if}\ \mathbf{d}=(2,2).
\end{cases}
\end{eqnarray*}
On the other hand, by (\ref{eq-F1-leadingTerm-FanoIndex=n-1}),
\begin{eqnarray*}
&& \sum_{e=0}^{n+m}\sum_{f=0}^{n+m}\langle \gamma_b, \gamma_e, g^{ef}\gamma_f,\gamma_c\rangle_{0,1}\\
&=& \sum_{e=0}^{n}\sum_{f=0}^{n}\langle \gamma_b, \gamma_e, g^{ef}\gamma_f,\gamma_c\rangle_{0,1}+\sum_{e=n+1}^{n+m}\sum_{f=n+1}^{n+m}\langle \gamma_b, \gamma_e, g^{ef}\gamma_f,\gamma_c\rangle_{0,1}\\
&=& \frac{1}{\prod_{i=1}^r d_i}\sum_{i=1}^{n-1}\langle \gamma_b, \sfh_{i}, \sfh_{n-i},\gamma_c\rangle_{0,4,1}+\sum_{e=n+1}^{n+m}\sum_{f=n+1}^{n+m}\langle \gamma_b, \gamma_e, g^{ef}\gamma_f,\gamma_c\rangle_{0,1}\\
&=& -\frac{(n-1)\elld g_{bc}}{\prod_{i=1}^r d_i}+\sum_{e=n+1}^{n+m}\sum_{f=n+1}^{n+m}\langle \gamma_b, \gamma_e, g^{ef}\gamma_f,\gamma_c\rangle_{0,1}.
\end{eqnarray*}
So (\ref{eq-4points-sum}) follows.
\end{proof}

\begin{lemma}\label{lem-EulerChar}
\begin{equation}\label{eq-EulerChar}
\chi(X_n(\mathbf{d}))=
\begin{cases}
\frac{(-2)^{n+2}+5}{3}+n,&
 \mbox{if}\ \mathbf{d}=3;\\
 \big((-1)^n+1\big)(n+2),
 &
 \mbox{if}\ \mathbf{d}=(2,2).
 \end{cases}
\end{equation}
\end{lemma}
\begin{proof}
\begin{eqnarray*}
&& \chi(X) =\int_X c_n(T_X)=\prod_{i=1}^r d_i\cdot \mathrm{Coeff}_{x^n}\big(\frac{(1+x)^{n+r+1}}{\prod_{i=1}^r (1+d_i x)}\big)\\
&=&\prod_{i=1}^r d_i\cdot\mathrm{Res}_{x=0}\big(\frac{(1+x)^{n+r+1}}{x^{n+1}\prod_{i=1}^r (1+d_i x)}\big)\\
&=& \prod_{i=1}^r d_i\cdot\Big(-\sum_{i=1}^r\mathrm{Res}_{x=-\frac{1}{d_i}}\big(\frac{(1+x)^{n+r+1}}{x^{n+1}\prod_{i=1}^r (1+d_i x)}\big)-\mathrm{Res}_{x=\infty}\big(\frac{(1+x)^{n+r+1}}{x^{n+1}\prod_{i=1}^r (1+d_i x)}\big)\Big).
\end{eqnarray*}
We compute residues respectively: 
\begin{eqnarray}\label{eq-residue-4}
&&\mathrm{Res}_{x=\infty}\Big(\frac{(1+x)^{n+r+1}}{x^{n+1}\prod_{i=1}^r (1+d_i x)}\Big)\nn\\
&=&-\frac{1}{\prod_{i=1}^r d_i}\big(n+r+1-\sum_{i=1}^r \frac{1}{d_i}\big)\nn\\
&=& \begin{cases}
-\frac{n}{3}-\frac{5}{9}, & \mbox{if}\ \mathbf{d}=3;\\
-\frac{n+2}{4}, & \mbox{if}\ \mathbf{d}=(2,2),
\end{cases}
\end{eqnarray}
and
\begin{eqnarray}\label{eq-residue-6}
\sum_{i=1}^{r}\mathrm{Res}_{x=-\frac{1}{d_i}}\Big(\frac{(1+x)^{n+r+1}}{x^{n+1}\prod_{i=1}^r (1+d_i x)}\Big)
=\begin{cases}
\frac{(-1)^{n+1}2^{n+2}}{9},&
 \mbox{if}\ \mathbf{d}=3;\\
\frac{(-1)^{n+1}(n+2)}{4}, &
 \mbox{if}\ \mathbf{d}=(2,2).
\end{cases}  
\end{eqnarray}
Then (\ref{eq-EulerChar}) follows
\end{proof}

\begin{corollary}\label{cor-dim-prim}
\begin{equation}\label{eq-dim-prim}
	\mathrm{rank}\ H_{\mathrm{prim}}^n(X_n(\mathbf{d}))=\begin{cases}
	\frac{2^{n+2}+2}{3},& \mbox{if}\ \mathbf{d}=3,\ n\ \mbox{is even};\\
	\frac{2^{n+2}-2}{3},& \mbox{if}\ \mathbf{d}=3,\ n\ \mbox{is odd};\\
	n+3,& \mbox{if}\ \mathbf{d}=(2,2),\ n\ \mbox{is even};\\
	n+1,& \mbox{if}\ \mathbf{d}=(2,2),\ n\ \mbox{is odd}.
	\end{cases}
\end{equation}
\end{corollary}

\subsubsection{The initial values}
\begin{theorem}\label{thm-4points-fanoIndex-(n-1)}
For $X_n(3)$  of dimension $\geq 3$, and $X_n(2,2)$  of odd dimension $n\geq 3$, $F^{(2)}(0)=1$.
\end{theorem}
\begin{proof}
By Proposition \ref{prop-initialValues-meaning} and Lemma \ref{lem-contraction-Pfaffian},
\begin{equation}\label{eq-4points-sum-1}
\sum_{e=n+1}^{n+m}\sum_{f=n+1}^{n+m}\langle \gamma_b, \gamma_e, g^{ef}\gamma_f,\gamma_c\rangle_{0,1}
=\big(2+(-1)^n m\big)F^{(2)}(0)g_{bc}.
\end{equation}
On the other hand by  (\ref{eq-4points-sum}) and (\ref{eq-dim-prim}) one finds
\begin{equation}\label{eq-4points-sum-2}
\sum_{e=n+1}^{n+m}\sum_{f=n+1}^{n+m}\langle \gamma_b, \gamma_e, g^{ef}\gamma_f,\gamma_c\rangle_{0,1}
=\big(2+(-1)^n m\big)g_{bc}.
\end{equation}
Comparing (\ref{eq-4points-sum-1}) and (\ref{eq-4points-sum-2}) we get $F^{(2)}(0)=1$.
\end{proof}

\subsection{\texorpdfstring{$F^{(4)}(0)$}{F(4)(0)} of cubic 3-folds}\label{sec:8points-Invariants-Cubic3Fold}
In this section we consider a smooth cubic 3-fold $X=X_n(3)$. Recall
\begin{equation}\label{eq-chernClass-cubic3fold}
	c_1(T_X)=2\sfh_1,\ c_2(T_X)=4\sfh_2,\ c_3(T_X)=-2\sfh_3,
\end{equation}
and
\begin{equation*}
	\mathrm{rank}\ H^*_{\mathrm{prim}}(X)=10.
\end{equation*}
Let $\Gamma_i=\sfh_i$, $0\leq i\leq 3$, and $\Gamma_4,\dots,\Gamma_{13}$ be a basis of $H^{3}_{\mathrm{prim}}(X)=H^3(X)$.
Let $t^0,\dots,t^{13}$ be the dual basis. 
 Let $g_{i,j}=(\Gamma_i,\Gamma_j)$, and $\Gamma^a=\sum_{e}g^{a,e}\gamma_e$.

As we have seen in Theorem \ref{thm-reconstructcubicandquadric} (i), to complete the reconstruction theorem of cubic 3-fold,  using Theorem \ref{thm-4points-fanoIndex-(n-1)} (or Theorem \ref{thm-cubic5}) we are left to determine $F^{(4)}(0)$. The aim of this subsection is to show:
\begin{theorem}\label{thm-8points-cubic3fold}
\begin{equation}\label{eq-8points-cubic3fold}
  F^{(4)}(0)=-\frac{2}{3}\frac{\partial F^{(3)}}{\partial t^2}(0).
\end{equation}
\end{theorem}

Now running
\[
\mbox{\texttt{correlatorInTauCoord}}\ \{3,\{3\},3,\{0,0,1,0\}\}
\]
in the Macaulay2 package \texttt{QuantumCohomologyFanoCompleteIntersection} returns the value of $\frac{\partial F^{(3)}}{\partial \tau^2}(0)$:
\[
-2\,{z}_{2}^{2}+18\,{z}_{2}-16,
\]
where $z_2:=F^{(2)}(0)$. By Theorem \ref{thm-cubic5} or Theorem \ref{thm-4points-fanoIndex-(n-1)}, $z_2=1$. 
So by (\ref{eq-tauTot-FanoIndex=n-1}),
\begin{equation*}
      \frac{\partial F^{(3)}}{\partial t^2}(0)=\frac{\partial F^{(3)}}{\partial \tau^2}(0)=0.
\end{equation*}
Then using (\ref{eq-8points-cubic3fold}) we get
\begin{theorem}\label{thm-8points-cubic3fold-final}
$F^{(4)}(0)=0$.
\end{theorem}

The strategy of the proof of Theorem \ref{thm-8points-cubic3fold} is similar to that of Section \ref{sec:4points-Invariants}; the main difference is that now we need to use Theorem \ref{thm-vanishing-deg2-cubics}. 
Take arbitrary $\gamma_0,\gamma_1,\dots,\gamma_5\in H^3(X)$.  We apply (\ref{eq-TRR1}) to $\langle \psi\gamma_0,\gamma_1,\gamma_2,\gamma_3,\gamma_4,\gamma_5\rangle_{1,2}$. By (\ref{eq-Dim}) and Theorem \ref{thm-monodromythm}, only the following terms are left:
\begin{eqnarray}\label{eq-TRR-genus1-deg2}
&&\langle \psi\gamma_0,\gamma_1,\gamma_2,\gamma_3,\gamma_4,\gamma_5\rangle_{1,2}\nonumber\\
&=&\frac{1}{3}\sum_{i=1}^{5}(\pm)\langle \gamma_0,\gamma_i,1\rangle_{0,3,0}
 \langle  \sfh_3,\dots,\hat{\gamma_i},\dots\rangle_{1,5,2}
 +\frac{1}{3}\sum_{i=1}^{5}(\pm)\langle \gamma_0,\gamma_i,\sfh_2\rangle_{0,1}
 \langle  \sfh,\dots,\hat{\gamma_i},\dots\rangle_{1,5,1}\nonumber\\
 &&+\frac{1}{3}\sum_{\{i,j,k\}\subset [5]}(\pm)\langle \gamma_0,\gamma_i,\gamma_j,\gamma_k,\sfh\rangle_{0,5,1}
 \langle  \sfh_2,\dots,\hat{\gamma_i},\hat{\gamma_j},\hat{\gamma_k},\dots\rangle_{1,3,1}\nn\\
 &&+\frac{1}{3}\langle \gamma_0,\gamma_1,\gamma_2,\gamma_3,\gamma_4,\gamma_5,\sfh_2\rangle_{0,7,2}
 \langle \sfh\rangle_{1,1,0}+\frac{1}{24}\sum_{a=0}^{13}\langle \gamma_0,\Gamma_a,\Gamma^a,\gamma_1,\gamma_2,\gamma_3,\gamma_4,\gamma_5\rangle_{0,8,2},
\end{eqnarray}
where $(\pm)$ means the sign arising from permutations of $\gamma_i$'s. By (\ref{eq-Deg0}) and (\ref{eq-chernClass-cubic3fold}),
\begin{equation}\label{eq-h110}
	\langle \sfh\rangle_{1,1,0}=-\frac{1}{2}.
\end{equation}
In the next section we use (\ref{eq-SvR}) and Theorem \ref{thm-vanishing-deg2-cubics} to compute the other genus 1 invariants in both sides of (\ref{eq-TRR-genus1-deg2}).

\subsubsection{Applications of SvR}
\begin{lemma}\label{lem-8points-applicationOfSvR-1}
\begin{gather}\label{eq-8points-applicationOfSvR-1}
\langle \psi\gamma_0,\gamma_1,\gamma_2,\gamma_3,\gamma_4,\gamma_5\rangle_{1,2}
\nn\\= -\frac{1}{6} \langle\psi \gamma_0, \gamma_1,\gamma_2,\gamma_3,\gamma_4,\gamma_5\rangle_{0,6,2}
+\frac{1}{12}\sum_{j=1}^{5}(-1)^{j-1}\langle  \gamma_0 \gamma_j,\prod_{\begin{subarray}{c}
1\leq i\leq 5\\
i\neq j\end{subarray}} \mathrm{ev}_i^* \gamma_i \rangle_{0,5,2}.
\end{gather}
\end{lemma}
\begin{proof}
Let 
\[
\mu=(a_0,a_1,a_2,a_3,a_4,a_5;\gamma_0,\gamma_1,\gamma_2,\gamma_3,\gamma_4,\gamma_5)
=(1,0,0,0,0,0;\gamma_0,\gamma_1,\gamma_2,\gamma_3,\gamma_4,\gamma_5).
\] 
By (\ref{eq-SvR}),
\begin{eqnarray}\label{eq-lem-8points-applicationOfSvR-1}
&&\langle \psi\gamma_0,\gamma_1,\gamma_2,\gamma_3,\gamma_4,\gamma_5\rangle_{1,2}
-\langle \psi\gamma_0,\gamma_1,\gamma_2,\gamma_3,\gamma_4,\gamma_5\rangle_{1,2}^0\nn\\
&=&\sum_{J\subset \{0,\dots,5\}}\Big(
(-1)^{1+|J|-p_J(\mu)}\big(1+|J|-p_J(\mu);(a_j)_{j\in J}\big)_{[1],J}\nn\\
&& \times \sum_{p=0}^{d_{1,J}(\mu)}
\mathrm{GW}_{(1,J)}^{2}
\big(\eta_p,c_{d_{1,J}(\mu)-p}(T_X);\psi\gamma_0,\gamma_1,\gamma_2,\gamma_3,\gamma_4,\gamma_5\big)\Big)\nn\\
&&+\sum_{J\subset \{0,\dots,5\}}\Big(
(-1)^{2+|J|-p_J(\mu)}\big(2+|J|-p_J(\mu);(a_j)_{j\in J}\big)_{[2],J}\nn\\
&& \times \sum_{p=0}^{d_{2,J}(\mu)}
\mathrm{GW}_{(2,J)}^2
\big(\eta_p,c_{d_{2,J}(\mu)-p}(T_X);\psi\gamma_0,\gamma_1,\gamma_2,\gamma_3,\gamma_4,\gamma_5\big)
\Big).
\end{eqnarray}
First by (\ref{eq-vanishing-deg2-cubics}), 
\begin{equation}\label{eq-8points-applicationOfSvR-1-1}
\langle \psi\gamma_0,\gamma_1,\gamma_2,\gamma_3,\gamma_4,\gamma_5\rangle_{1,2}^0=0.
\end{equation}
We compute the ingredients of (\ref{eq-lem-8points-applicationOfSvR-1}) separately.
If $0\not\in J$, then by (\ref{eq-intersectionNum-M(1,m,J)-1}),
\begin{equation}\label{eq-8points-applicationOfSvR-1-2}
\big(1+|J|-p_J(\mu);(a_j)_{j\in J}\big)_{[1],J}
=\big(1+|J|;0\big)_{[1],J}
=\frac{1}{24},
\end{equation}
and
\begin{equation}\label{eq-8points-applicationOfSvR-1-3}
\big(2+|J|-p_J(\mu);(a_j)_{j\in J}\big)_{[2],J}
=\big(2+|J|;0\big)_{[2],J}
= \frac{2^{|J|}}{24}.
\end{equation}
If $0\in J$, then by (\ref{eq-intersectionNum-M(1,m,J)-2}),
\begin{equation}\label{eq-8points-applicationOfSvR-1-4}
\big(1+|J|-p_J(\mu);(a_j)_{j\in J}\big)_{[1],J}
=\big(|J|;1,0^{|J|-1}\big)_{[1],J}=\frac{|J|}{24},
\end{equation}
and
\begin{equation}\label{eq-8points-applicationOfSvR-1-5}
\big(2+|J|-p_J(\mu);(a_j)_{j\in J}\big)_{[2],J}
= \big(1+|J|;1,0^{|J|-1}\big)_{[2],J}
=\frac{2^{|J|-1}(|J|+1)}{24}.
\end{equation}

If $0\not\in J$, then
\[
p_{J}(\mu)=0,\
d_{1,J}(\mu)=1-|J|,
\]
and by (\ref{eq-def-invariants-(1,J)}) and (\ref{eq-chernClass-cubic3fold}),
\begin{eqnarray*}
&&\sum_{p=0}^{d_{1,J}(\mu)}\mathrm{GW}_{(1,J)}^{2}
\big(\eta_p, c_{d_{1,J}(\mu)-p}(T_X);\psi\gamma_0,\gamma_1,\gamma_2,\gamma_3,\gamma_4,\gamma_5\big)\\
&=&\pm \langle c_{1-|J|}(T_X)\prod_{j\in J}\gamma_j,\psi \gamma_0, \prod_{j\not\in J\cup\{0\}} \mathrm{ev}_j^* \gamma_j\rangle_{0,7-|J|,2}\\
&&\pm\langle \psi c_{-|J|}(T_X)\prod_{j\in J}\gamma_j,\psi \gamma_0,\prod_{j\not\in J\cup\{0\}} \mathrm{ev}_j^* \gamma_j\rangle_{0,7-|J|,2}\\
&=& \begin{cases}
0, & \mbox{if } |J|\geq 2,\\
(-1)^j\langle\gamma_j,\psi \gamma_0, \prod_{k\neq j} \mathrm{ev}_k^* \gamma_k\rangle_{0,6,2}, & \mbox{if } |J|=\{j\},\ 1\leq j\leq 5,\\
\langle 2\sfh_1,\psi \gamma_0, \gamma_1,\dots,\gamma_5\rangle_{0,7,2}
+\langle \psi ,\psi \gamma_0,\gamma_1,\dots,\gamma_5\rangle_{0,7,2},
& \mbox{if } J=\emptyset.
\end{cases}
\end{eqnarray*}
Since $\sfh_1\cup \gamma_i=0$, by (\ref{eq-Div}) and (\ref{eq-Dilaton}) we get, in the case $0\not\in J$,
\begin{eqnarray}\label{eq-8points-applicationOfSvR-1-6}
&&\sum_{p=0}^{d_{1,J}(\mu)}\mathrm{GW}_{(1,J)}^{2}
\big(\eta_p, c_{d_{1,J}(\mu)-p}(T_X);\psi\gamma_0,\gamma_1,\gamma_2,\gamma_3,\gamma_4,\gamma_5\big)\nn\\
&=& \begin{cases}
0, & \mbox{if } |J|\geq 2,\\
(-1)^j
\langle \gamma_j,\psi \gamma_0, \prod_{j\not\in J\cup\{0\}} \mathrm{ev}_j^* \gamma_j\rangle_{0,6,2}, & \mbox{if } |J|=\{j\},\ 1\leq j\leq 5,\\
8\langle \psi \gamma_0, \gamma_1,\dots,\gamma_5\rangle_{0,6,2},& \mbox{if } J=\emptyset.
\end{cases}
\end{eqnarray}

If $0\in J$, then
\[
p_{J}(\mu)=1,
d_{1,J}(\mu)=2-|J|,\ 
\]
and by by (\ref{eq-def-invariants-(1,J)}),
\begin{eqnarray*}
&&\sum_{p=0}^{d_{1,J}(\mu)}\mathrm{GW}_{(1,J)}^{2}
\big(\eta_p, c_{d_{1,J}(\mu)-p}(T_X);\psi\gamma_0,\gamma_1,\gamma_2,\gamma_3,\gamma_4,\gamma_5\big)\\
&=&\pm \sum_{p=0}^{2-|J|}
\langle\psi^p c_{2-|J|-p}(T_X) \prod_{j\in J}\gamma_j,\prod_{j\not\in J} \gamma_j \rangle_{0,7-|J|,2}.
\end{eqnarray*}
Since $\sfh_i\cdot \gamma_0=0$ for $i>0$, $c_{2-|J|-p}(T_X)\cup \gamma_0=0$ when $p\neq 2-|J|$. So
\begin{eqnarray}\label{eq-8points-applicationOfSvR-1-7}
&&\sum_{p=0}^{d_{1,J}(\mu)}\mathrm{GW}_{(1,J)}^{2}
\big(\eta_p, c_{d_{1,J}(\mu)-p}(T_X);\psi\gamma_0,\gamma_1,\gamma_2,\gamma_3,\gamma_4,\gamma_5\big)\nn\\
&=& \pm\langle\psi^{2-|J|}  \prod_{j\in J}\gamma_j,\prod_{j\not\in J} \gamma_j \rangle_{0,7-|J|,2}\nn\\
&=&\begin{cases}
\langle \psi \gamma_0,\gamma_1,\dots,\gamma_5 \rangle_{0,6,2}, &
\mbox{if}\ J=\{0\},\\
(-1)^{j-1}\langle  \gamma_0 \gamma_j,\prod_{k\not\in J} \gamma_k \rangle_{0,5,2},&
\mbox{if}\ J=\{0,j\},\ 1\leq j\leq 5,\\
0,& \mbox{if}\ 0\in J\ \mbox{and}\ |J|\geq 3.
\end{cases}
\end{eqnarray}

Whether $0\in J$ or not we have $d_{2,J}(\mu)<0$, 
so
\begin{equation}\label{eq-8points-applicationOfSvR-1-8}
\mathrm{GW}_{(2,J)}^2
\big(\eta_p,c_{d_{2,J}(\mu)-p}(T_X);\psi\gamma_0,\gamma_1,\gamma_2,\gamma_3,\gamma_4,\gamma_5\big)=0.
\end{equation}
Combining (\ref{eq-8points-applicationOfSvR-1-1})-(\ref{eq-8points-applicationOfSvR-1-8}), we obtain
\begin{eqnarray*}
&&\langle \psi\gamma_0,\gamma_1,\gamma_2,\gamma_3,\gamma_4,\gamma_5\rangle_{1,2}\\
&=& \frac{1}{24} \sum_{j=1}^{5}(-1)^j
\langle \gamma_j,\psi \gamma_0, \prod_{i\not\in \{0,j\}} \mathrm{ev}_i^* \gamma_i\rangle_{0,6,2}
- \frac{1}{24}\cdot 8\langle \psi \gamma_0, \gamma_1,\dots,\gamma_5\rangle_{0,6,2}\\
&&- \frac{1}{24}\langle \psi \gamma_0,\gamma_1,\dots,\gamma_5 \rangle_{0,6,2}
+\frac{1}{12}\sum_{j=1}^{5}(-1)^{j-1}\langle  \gamma_0 \gamma_j,\prod_{i\not\in \{0,j\}} \mathrm{ev}_i^* \gamma_i \rangle_{0,5,2}\\
&=& \frac{1}{24} \cdot 5
\langle\psi \gamma_0, \gamma_1,\dots,\gamma_5\rangle_{0,6,2}
- \frac{1}{24}\cdot 8\langle \psi \gamma_0, \gamma_1,\dots,\gamma_5\rangle_{0,6,2}\\
&&- \frac{1}{24}\langle \psi \gamma_0,\gamma_1,\dots,\gamma_5 \rangle_{0,6,2}
+\frac{1}{12}\sum_{j=1}^{5}(-1)^{j-1}\langle  \gamma_0 \gamma_j,\prod_{i\not\in \{0,j\}} \mathrm{ev}_i^* \gamma_i \rangle_{0,5,2}\\
&=& -\frac{1}{6} \langle\psi \gamma_0, \gamma_1,\gamma_2,\gamma_3,\gamma_4,\gamma_5\rangle_{0,6,2}
+\frac{1}{12}\sum_{j=1}^{5}(-1)^{j-1}\langle  \gamma_0 \gamma_j,\prod_{i\not\in \{0,j\}} \mathrm{ev}_i^* \gamma_i \rangle_{0,5,2}.
\end{eqnarray*}
\end{proof}

\begin{lemma}\label{lem-8points-applicationOfSvR-2}
\begin{equation}\label{eq-8points-applicationOfSvR-2}
\langle \sfh_3,\gamma_1,\gamma_2,\gamma_3,\gamma_4\rangle_{1,5,2}
=-\frac{1}{12}\langle \sfh_3,\gamma_1,\gamma_2,\gamma_3,\gamma_4\rangle_{0,5,2}.
\end{equation}
\end{lemma}
\begin{proof}
Let 
\[
\mu=(0,0,0,0,0;\sfh_3,\gamma_1,\gamma_2,\gamma_3,\gamma_4).
\] 
Then 
\[
p_J(\mu)=0,\ d_{l,J}(\mu)=3-2l-|J|.
\]
By (\ref{eq-SvR}), using (\ref{eq-intersectionNum-M(1,m,J)-1}),
\begin{eqnarray*}
&&\langle \sfh_3,\gamma_1,\gamma_2,\gamma_3,\gamma_4\rangle_{1,5,2}
-\langle \sfh_3,\gamma_1,\gamma_2,\gamma_3,\gamma_4\rangle_{1,5,2}^0\\
&=&\sum_{J\subset \{1,\dots,5\}}\Big(
(-1)^{1+|J|}\big(1+|J|;0^{J}\big)_{[1],J}\\
&& \times \sum_{p=0}^{d_{1,J}(\mu)}
\mathrm{GW}_{(1,J)}^{2}
\big(\eta_p,c_{d_{1,J}(\mu)-p}(T_X);\sfh_3,\gamma_1,\gamma_2,\gamma_3,\gamma_4\big)\Big)\\
&=& -\frac{1}{24}\sum_{p=0}^{1}
\mathrm{GW}_{(1,\emptyset)}^{2}
\big(\eta_p,c_{1-p}(T_X);\sfh_3,\gamma_1,\gamma_2,\gamma_3,\gamma_4\big)\\
&&+\sum_{J\subset \{1,\dots,5\}, |J|=1}\Big(\frac{1}{24}
\mathrm{GW}_{(1,J)}^{2}
\big(\eta_0,1;\sfh_3,\gamma_1,\gamma_2,\gamma_3,\gamma_4\big)\Big)\\
&=& -\frac{1}{24}\langle 2\sfh_1,\sfh_3,\gamma_1,\gamma_2,\gamma_3,\gamma_4\rangle_{0,6,2}
-\frac{1}{24}\langle \psi,\sfh_3,\gamma_1,\gamma_2,\gamma_3,\gamma_4\rangle_{0,6,2}\\
&&+\frac{1}{24}\langle \sfh_3,\gamma_1,\gamma_2,\gamma_3,\gamma_4\rangle_{0,5,2}
+\frac{1}{24}\sum_{j=1}^{4}(-1)^{j-1}\langle \gamma_j,\sfh_3,\prod_{k\neq j} \mathrm{ev}_k^* \gamma_k\rangle_{0,5,2}.
\end{eqnarray*}
By (\ref{eq-vanishing-deg2-cubics}),
\[
\langle \sfh_3,\gamma_1,\gamma_2,\gamma_3,\gamma_4\rangle_{1,5,2}^0=0.
\]
Then by (\ref{eq-Div}) and (\ref{eq-Dilaton}) we get
\begin{eqnarray*}
&&\langle \sfh_3,\gamma_1,\gamma_2,\gamma_3,\gamma_4\rangle_{1,5,2}\\
&=& -\frac{1}{6}\langle \sfh_3,\gamma_1,\gamma_2,\gamma_3,\gamma_4\rangle_{0,5,2}
-\frac{1}{8}\langle \sfh_3,\gamma_1,\gamma_2,\gamma_3,\gamma_4\rangle_{0,5,2}\\
&&+\frac{1}{24}\langle \sfh_3,\gamma_1,\gamma_2,\gamma_3,\gamma_4\rangle_{0,5,2}
+\frac{1}{24}\sum_{j=1}^{4}(-1)^{j-1}\langle \gamma_j,\sfh_3,\prod_{k\neq j} \mathrm{ev}_k^* \gamma_k\rangle_{0,5,2}\\
&=& -\frac{1}{6}\langle \sfh_3,\gamma_1,\gamma_2,\gamma_3,\gamma_4\rangle_{0,5,2}
-\frac{1}{8}\langle \sfh_3,\gamma_1,\gamma_2,\gamma_3,\gamma_4\rangle_{0,5,2}\\
&&+\frac{1}{24}\langle \sfh_3,\gamma_1,\gamma_2,\gamma_3,\gamma_4\rangle_{0,5,2}
+\frac{1}{24}\cdot 4\langle \sfh_3,\gamma_1,\gamma_2,\gamma_3,\gamma_4\rangle_{0,5,2}\\
&=&-\frac{1}{12}\langle \sfh_3,\gamma_1,\gamma_2,\gamma_3,\gamma_4\rangle_{0,5,2}.
\end{eqnarray*}
\end{proof}

\begin{lemma}\label{lem-8points-applicationOfSvR-3}
\begin{equation}\label{eq-8points-applicationOfSvR-3}
	\langle \gamma_1,\gamma_2,\gamma_3,\gamma_4\rangle_{1,4,1}=0.
\end{equation}
\end{lemma}
\begin{proof}
Let 
\[
\mu=(0,0,0,0;\gamma_1,\gamma_2,\gamma_3,\gamma_4).
\] 
Then 
\[
p_J(\mu)=0,\ d_{l,J}(\mu)=3-2l-|J|.
\]
By (\ref{eq-SvR}) and (\ref{eq-intersectionNum-M(1,m,J)-1}),
\begin{eqnarray*}
&&\langle \gamma_1,\gamma_2,\gamma_3,\gamma_4\rangle_{1,4,1}
-\langle \gamma_1,\gamma_2,\gamma_3,\gamma_4\rangle_{1,4,1}^0\\
&=&\sum_{J\subset \{1,\dots,4\}}\Big(
(-1)^{1+|J|}\big(1+|J|;0^{J}\big)_{[1],J}\\
&& \times \sum_{p=0}^{d_{1,J}(\mu)}
\mathrm{GW}_{(1,J)}^{1}
\big(\eta_p,c_{d_{1,J}(\mu)-p}(T_X);\gamma_1,\gamma_2,\gamma_3,\gamma_4\big)\Big)\\
&=& -\frac{1}{24}\sum_{p=0}^{1}
\mathrm{GW}_{(1,\emptyset)}^{1}
\big(\eta_p,c_{1-p}(T_X);\gamma_1,\gamma_2,\gamma_3,\gamma_4\big)\\
&&+\sum_{J\subset \{1,\dots,4\}, |J|=1}\Big(\frac{1}{24}
\mathrm{GW}_{(1,J)}^{1}
\big(\eta_0,1;\gamma_1,\gamma_2,\gamma_3,\gamma_4\big)\Big)\\
&=& -\frac{1}{24}\langle 2\sfh_1,\gamma_1,\gamma_2,\gamma_3,\gamma_4\rangle_{0,5,1}
-\frac{1}{24}\langle \psi,\gamma_1,\gamma_2,\gamma_3,\gamma_4\rangle_{0,5,1}\\
&&+\frac{1}{24}\sum_{j=1}^{4}(\pm)\langle \mu_j,\prod_{k\neq j}\mathrm{ev}_k^{*} \mu_k\rangle_{0,4,1}\\
&=& -\frac{1}{24}\langle 2\sfh_1,\gamma_1,\gamma_2,\gamma_3,\gamma_4\rangle_{0,5,1}
-\frac{1}{24}\langle \psi,\gamma_1,\gamma_2,\gamma_3,\gamma_4\rangle_{0,5,1}\\
&&+\frac{1}{24}\cdot 4\langle\gamma_1,\gamma_2,\gamma_3,\gamma_4\rangle_{0,4,1}.
\end{eqnarray*}
By (\ref{eq-vanishing-deg2-cubics}),
\[
\langle \gamma_1,\gamma_2,\gamma_3,\gamma_4\rangle_{1,4,1}^0=0. 
\]
So by (\ref{eq-Div}) and (\ref{eq-Dilaton}) we get
\begin{eqnarray*}
\langle \gamma_1,\gamma_2,\gamma_3,\gamma_4\rangle_{1,4,1}
= -\frac{2}{24}\langle\gamma_1,\gamma_2,\gamma_3,\gamma_4\rangle_{0,4,1}
-\frac{2}{24}\langle \gamma_1,\gamma_2,\gamma_3,\gamma_4\rangle_{0,4,1}\\
+\frac{4}{24}\langle\gamma_1,\gamma_2,\gamma_3,\gamma_4\rangle_{0,4,1}=0.
\end{eqnarray*}
\end{proof}

\begin{lemma}\label{lem-8points-applicationOfSvR-4}
\begin{equation}\label{eq-8points-applicationOfSvR-4}
\langle \sfh_2,\gamma_1,\gamma_2\rangle_{1,3,1}
= 0.
\end{equation}
\end{lemma}
\begin{proof}
Let 
\[
\mu=(0,0,0;\sfh_2,\gamma_1,\gamma_2).
\] 
Then 
\[
p_J(\mu)=0,\ d_{l,J}(\mu)=3-2l-|J|.
\]
By (\ref{eq-SvR}) and (\ref{eq-intersectionNum-M(1,m,J)-1}),
\begin{eqnarray*}
&&\langle \sfh_2,\gamma_1,\gamma_2\rangle_{1,3,1}
-\langle \sfh_2, \gamma_1,\gamma_2\rangle_{1,3,1}^0\\
&=&\sum_{J\subset \{1,2,3\}}\Big(
(-1)^{1+|J|}\big(1+|J|;0^{J}\big)_{[1],J}\\
&& \times \sum_{p=0}^{d_{1,J}(\mu)}
\mathrm{GW}_{(1,J)}^{1}
\big(\eta_p,c_{d_{1,J}(\mu)-p}(T_X);\sfh_2,\gamma_1,\gamma_2\big)\Big)\\
&=& -\frac{1}{24}\sum_{p=0}^{1}
\mathrm{GW}_{(1,\emptyset)}^{1}
\big(\eta_p,c_{1-p}(T_X);\sfh_2,\gamma_1,\gamma_2\big)\\
&&+\sum_{J\subset \{1,2,3\}, |J|=1}\Big(\frac{1}{24}
\mathrm{GW}_{(1,J)}^{1}
\big(\eta_0,1;\sfh_2,\gamma_1,\gamma_2\big)\Big)\\
&=& -\frac{1}{24}\langle 2\sfh_1,\sfh_2,\gamma_1,\gamma_2\rangle_{0,4,1}
-\frac{1}{24}\langle \psi,\sfh_2,\gamma_1,\gamma_2\rangle_{0,4,1}\\
&&+\frac{1}{24}\sum_{j=1}^{3}(\pm)\langle \mu_j,\prod_{k\neq j}\mathrm{ev}_k^{*} \mu_k\rangle_{0,3,1}\\
&=& -\frac{1}{24}\langle 2\sfh_1,\sfh_2,\gamma_1,\gamma_2\rangle_{0,4,1}
-\frac{1}{24}\langle \psi,\sfh_2,\gamma_1,\gamma_2\rangle_{0,4,1}\\
&&+\frac{1}{24}\cdot 3\langle \sfh_2,\gamma_1,\gamma_2\rangle_{0,3,1}.
\end{eqnarray*}
By (\ref{eq-vanishing-deg2-cubics}),
\[
\langle \sfh_2, \gamma_1,\gamma_2\rangle_{1,3,1}^0=0.
\]
So by (\ref{eq-Div}) and (\ref{eq-Dilaton}),
\begin{equation*}
\langle \sfh_2, \gamma_1,\gamma_2\rangle_{1,3,1}
= -\frac{2}{24}\langle \sfh_2,\gamma_1,\gamma_2\rangle_{0,3,1}
-\frac{1}{24}\langle \sfh_2,\gamma_1,\gamma_2\rangle_{0,3,1}
+\frac{3}{24}\langle \sfh_2,\gamma_1,\gamma_2\rangle_{0,3,1}= 0.
\end{equation*}
\end{proof}

\subsubsection{Further simplifications}
By (\ref{eq-h110}) and Lemmas \ref{lem-8points-applicationOfSvR-1} to
\ref{lem-8points-applicationOfSvR-4}, (\ref{eq-TRR-genus1-deg2}) becomes
\begin{eqnarray}\label{eq-TRR-genus1-deg2-1}
&& -\frac{1}{6} \langle\psi \gamma_0, \gamma_1,\gamma_2,\gamma_3,\gamma_4,\gamma_5\rangle_{0,6,2}
+\frac{1}{12}\sum_{j=1}^{5}(-1)^{j-1}\langle  \gamma_0 \gamma_j,\prod_{\begin{subarray}{c}
1\leq i\leq 5\\
i\neq j\end{subarray}} \mathrm{ev}_i^* \gamma_i \rangle_{0,5,2}\nonumber\\
&=&-\frac{1}{36}\sum_{i=1}^{5}(-1)^{i-1}\langle\gamma_0,\gamma_i,1\rangle_{0,3,0}
 \langle  \sfh_3,\dots,\hat{\gamma_i},\dots\rangle_{0,5,2}
 -\frac{1}{6}\langle \sfh_2,\gamma_0,\gamma_1,\gamma_2,\gamma_3,\gamma_4,\gamma_5\rangle_{0,7,2}\nn\\
&& +\frac{1}{24}\sum_{a=0}^{13}\langle \Gamma_a,\Gamma^a,\gamma_0,\gamma_1,\gamma_2,\gamma_3,\gamma_4,\gamma_5\rangle_{0,8,2}.
\end{eqnarray}
By (\ref{eq-FCA}) and (\ref{eq-Div}),
\begin{eqnarray}\label{eq-TRR-genus1-deg2-simplification-1}
&&\sum_{a=0}^{13}\langle \Gamma_a,\Gamma^a,\gamma_0,\gamma_1,\gamma_2,\gamma_3,\gamma_4,\gamma_5\rangle_{0,8,2}\nn\\
&=& \sum_{\Gamma_a\in H^3(X)}\langle \Gamma_a,\Gamma^a,\gamma_0,\gamma_1,\gamma_2,\gamma_3,\gamma_4,\gamma_5\rangle_{0,8,2}\nn\\
&&+\frac{1}{3}\langle 1,\sfh_3,\gamma_0,\gamma_1,\gamma_2,\gamma_3,\gamma_4,\gamma_5\rangle_{0,8,2}+\frac{1}{3}\langle \sfh_1,\sfh_2, \gamma_0,\gamma_1,\gamma_2,\gamma_3,\gamma_4,\gamma_5\rangle_{0,8,2}\nn\\
&&+\frac{1}{3}\langle \sfh_2,\sfh_1,\gamma_0,\gamma_1,\gamma_2,\gamma_3,\gamma_4,\gamma_5\rangle_{0,8,2}
+\frac{1}{3}\langle \sfh_3,1, \gamma_0,\gamma_1,\gamma_2,\gamma_3,\gamma_4,\gamma_5\rangle_{0,8,2}\nn\\
&=& \sum_{\Gamma_a\in H^3(X)}\langle \Gamma_a,\Gamma^a,\gamma_0,\gamma_1,\gamma_2,\gamma_3,\gamma_4,\gamma_5\rangle_{0,8,2}
+\frac{4}{3}\langle \sfh_2,\gamma_0,\gamma_1,\gamma_2,\gamma_3,\gamma_4,\gamma_5\rangle_{0,7,2}.
\end{eqnarray}
Note that
\[
\gamma_i\cup\gamma_j=\frac{1}{3}\langle \gamma_i,\gamma_j,1\rangle_{0,3,0}\sfh_3.
\]
So from (\ref{eq-TRR-genus1-deg2-1}) and (\ref{eq-TRR-genus1-deg2-simplification-1}) we get
\begin{eqnarray}\label{eq-TRR-genus1-deg2-3}
&&\frac{1}{24}\sum_{\Gamma_a\in H^3}\langle \gamma_0,\gamma_a,\gamma^a,\gamma_1,\gamma_2,\gamma_3,\gamma_4,\gamma_5\rangle_{0,8,2}\nonumber\\
&=&  -\frac{1}{6} \langle\psi \gamma_0, \gamma_1,\gamma_2,\gamma_3,\gamma_4,\gamma_5\rangle_{0,6,2}
+\frac{1}{18}\sum_{i=1}^{5}(-1)^{i-1}\langle \gamma_0,\gamma_i,1\rangle_{0,3,0}
 \langle  \sfh_3,\dots,\hat{\gamma_i},\dots\rangle_{0,5,2}\nonumber\\
&&+\frac{1}{9}\langle H_2,\gamma_0,\gamma_1,\gamma_2,\gamma_3,\gamma_4,\gamma_5\rangle_{0,7,2}.
\end{eqnarray}

\begin{lemma}\label{lem-TRR-genus1-deg2-simplification-2}
\begin{equation}\label{eq-TRR-genus1-deg2-simplification-2}
	   -\frac{1}{6} \langle\psi \gamma_0, \gamma_1,\gamma_2,\gamma_3,\gamma_4,\gamma_5\rangle_{0,6,2}
+\frac{1}{18}\sum_{i=1}^{5}(-1)^{i-1}\langle \gamma_0,\gamma_i,1\rangle_{0,3,0}
 \langle  \sfh_3,\dots,\hat{\gamma_i},\dots\rangle_{0,5,2}=0.
\end{equation}
\end{lemma}
\begin{proof}
Denote
\[
c_{i,j}=(\gamma_i,\gamma_j)=\langle \gamma_i,\gamma_j,1\rangle_{0,3,0}.
\]
We apply (\ref{eq-TRR0}) to $\langle \psi \gamma_0, \gamma_1,\dots,\gamma_5\rangle_{0,6,2}$. For brevity we use the Einstein summation convention for $\Gamma_a$'s, with $a$ running through $0\leq a\leq 13$. So
\begin{align}
& \langle \psi \gamma_0, \gamma_1,\dots,\gamma_5\rangle_{0,6,2}\nn\\
=& \langle \gamma_0, \Gamma_e\rangle_{0,2,1}g^{ef}\langle \Gamma_f, \gamma_1,\dots,\gamma_5\rangle_{0,6,1}+\langle \gamma_0, \Gamma_e\rangle_{0,2,2}g^{ef}\langle \Gamma_f, \gamma_1,\dots,\gamma_5\rangle_{0,6,0}\label{eq-groupOfTerms-S1}\tag{S1}\\
&+\sum_{\beta_1+\beta_2=2}\big(\langle \gamma_0, \gamma_3, \Gamma_e\rangle_{0,3,\beta_1}g^{ef}\langle \Gamma_f, \gamma_1,\gamma_2, \gamma_4,\gamma_5\rangle_{0,5,\beta_2}\nn\\
&-\langle \gamma_0, \gamma_4, \Gamma_e\rangle_{0,3,\beta_1}g^{ef}\langle \Gamma_f, \gamma_1,\gamma_2, \gamma_3,\gamma_5\rangle_{0,5,\beta_2}\nn\\
&+\langle \gamma_0, \gamma_5, \Gamma_e\rangle_{0,3,\beta_1}g^{ef}\langle \Gamma_f, \gamma_1,\gamma_2, \gamma_3,\gamma_4\rangle_{0,5,\beta_2}\big)\label{eq-groupOfTerms-S2}\tag{S2}\\
&+\sum_{\beta_1+\beta_2=2}\big(\langle \gamma_0, \gamma_3, \gamma_4, \Gamma_e\rangle_{0,4,\beta_1}g^{ef}\langle \Gamma_f, \gamma_1,\gamma_2, \gamma_5\rangle_{0,4,\beta_2}\nn\\
&-\langle \gamma_0, \gamma_3, \gamma_5, \Gamma_e\rangle_{0,4,\beta_1}g^{ef}\langle \Gamma_f, \gamma_1,\gamma_2, \gamma_4\rangle_{0,4,\beta_2}\nn\\
&+\langle \gamma_0, \gamma_4, \gamma_5, \Gamma_e\rangle_{0,4,\beta_1}g^{ef}\langle \Gamma_f, \gamma_1,\gamma_2, \gamma_3,\rangle_{0,4,\beta_2}\big)\label{eq-groupOfTerms-S3}\tag{S3}\\
&+\sum_{\beta_1+\beta_2=2}\langle \gamma_0, \gamma_3, \gamma_4,\gamma_5, \Gamma_e\rangle_{0,5,\beta_1}g^{ef}\langle \Gamma_f, \gamma_1,\gamma_2\rangle_{0,3,\beta_2}.\label{eq-groupOfTerms-S4}\tag{S4}
\end{align}

We evaluate the groups of terms (\ref{eq-groupOfTerms-S1})-(\ref{eq-groupOfTerms-S4})
separately. By (\ref{eq-Dim}), 
\begin{equation*}
\langle \gamma_0, \Gamma_e\rangle_{0,2,1}=\langle \gamma_0, \Gamma_e\rangle_{0,2,2}=0,
\end{equation*}
\begin{equation*}
	\langle \gamma_i,\gamma_j,\Gamma_e\rangle_{0,3,\beta}=0\ \mbox{unless}\ \Gamma_e=1\ \mbox{or}\ \sfh_2,
\end{equation*}
and
\begin{equation*}
	\langle \gamma_i,\gamma_j, \gamma_k,\Gamma_e\rangle_{0,4,\beta}=0\ \mbox{unless}\ \Gamma_e\in H^3(X).
\end{equation*}
so (\ref{eq-groupOfTerms-S1})=0, and
\begin{eqnarray*}
(\ref{eq-groupOfTerms-S2})&=&
\frac{1}{3}\big(c_{0,3}\langle \sfh_3, \gamma_1,\gamma_2, \gamma_4,\gamma_5\rangle_{0,5,2}
-c_{0,4}\langle \sfh_3, \gamma_1,\gamma_2, \gamma_3,\gamma_5\rangle_{0,5,2}
+c_{0,5}\langle \sfh_3, \gamma_1,\gamma_2, \gamma_3,\gamma_4\rangle_{0,5,2}\big)\\
&&+\frac{1}{3}\big(\langle \gamma_0, \gamma_3, \sfh_2\rangle_{0,3,1}\langle \sfh_1, \gamma_1,\gamma_2, \gamma_4,\gamma_5\rangle_{0,5,1}
-\langle \gamma_0, \gamma_4, \sfh_2\rangle_{0,3,1}\langle \sfh_1, \gamma_1,\gamma_2, \gamma_3,\gamma_5\rangle_{0,5,1}\\
&&+\langle \gamma_0, \gamma_5, \sfh_2\rangle_{0,3,1}\langle \sfh_1, \gamma_1,\gamma_2, \gamma_3,\gamma_4\rangle_{0,5,1}\big),
\end{eqnarray*}
\begin{eqnarray*}
(\ref{eq-groupOfTerms-S3})&=&
\sum_{\Gamma_e, \Gamma_f\in H^3(X)}
\big(\langle \gamma_0, \gamma_3, \gamma_4, \gamma_e\rangle_{0,4,1}g^{ef}\langle \gamma_f, \gamma_1,\gamma_2, \gamma_5\rangle_{0,4,1}\\
&&-\langle \gamma_0, \gamma_3, \gamma_5, \gamma_e\rangle_{0,4,1}g^{ef}\langle \gamma_f, \gamma_1,\gamma_2, \gamma_4\rangle_{0,4,1}\\
&&+\langle \gamma_0, \gamma_4, \gamma_5, \gamma_e\rangle_{0,4,1}g^{ef}\langle \gamma_f, \gamma_1,\gamma_2, \gamma_3,\rangle_{0,4,1}\big),
\end{eqnarray*}
and
\begin{eqnarray*}
(\ref{eq-groupOfTerms-S4})&=&
\frac{1}{3}\langle \gamma_0, \gamma_3, \gamma_4,\gamma_5, \sfh_1\rangle_{0,5,1}
\langle \sfh_2, \gamma_1,\gamma_2\rangle_{0,3,1}
+\frac{1}{3}\langle \gamma_0, \gamma_3, \gamma_4,\gamma_5, \sfh_3\rangle_{0,5,2}\cdot c_{1,2}.
\end{eqnarray*}
From (\ref{eq-F1-leadingTerm-FanoIndex=n-1}) one finds
\begin{equation}\label{eq-lem-TRR-genus1-deg2-simplification-2-3point}
\langle \gamma_i,\gamma_j,\sfh_2\rangle_{0,3,1}=-6c_{i,j}.
\end{equation}
By Theorem \ref{thm-4points-fanoIndex-(n-1)},  $F^{(2)}(0)=1$, so by (\ref{eq-initialValues-meaning-odd})  we have
\begin{equation}\label{eq-lem-TRR-genus1-deg2-simplification-2-4point}
\langle \gamma_{i},\gamma_{j},\gamma_{k},\gamma_{l}\rangle_{0,4,1}
= c_{i,j}c_{k,l}-c_{i,k}c_{j,l}+c_{i,l}c_{j,k}.
\end{equation}
Moreover by (\ref{eq-lowerTerms-Cubic-F(2)}) one finds
\begin{equation}\label{eq-lem-TRR-genus1-deg2-simplification-2-5point}
  \langle \sfh_3,\gamma_i,\gamma_j,\gamma_k,\gamma_l\rangle_{0,5,2}=
  3(c_{i,j}c_{k,l}-c_{i,k}c_{j,l}+c_{i,l}c_{j,k}).
\end{equation}
Using (\ref{eq-lem-TRR-genus1-deg2-simplification-2-3point}), (\ref{eq-lem-TRR-genus1-deg2-simplification-2-4point}), (\ref{eq-lem-TRR-genus1-deg2-simplification-2-5point}), and (\ref{eq-Div}), we get
\begin{eqnarray}\label{eq-groupOfTerms-S2-1}
(\ref{eq-groupOfTerms-S2})&=&
\frac{1}{3}\big(c_{0,3}\langle \sfh_3, \gamma_1,\gamma_2, \gamma_4,\gamma_5\rangle_{0,5,2}
-c_{0,4}\langle \sfh_3, \gamma_1,\gamma_2, \gamma_3,\gamma_5\rangle_{0,5,2}
+c_{0,5}\langle \sfh_3, \gamma_1,\gamma_2, \gamma_3,\gamma_4\rangle_{0,5,2}\big)\nn\\
&&+\frac{1}{3}\big(\langle \gamma_0, \gamma_3, \sfh_2\rangle_{0,3,1}\langle \gamma_1,\gamma_2, \gamma_4,\gamma_5\rangle_{0,5,1}
-\langle \gamma_0, \gamma_4, \sfh_2\rangle_{0,3,1}\langle \gamma_1,\gamma_2, \gamma_3,\gamma_5\rangle_{0,5,1}\nn\\
&&+\langle \gamma_0, \gamma_5, \sfh_2\rangle_{0,3,1}\langle  \gamma_1,\gamma_2, \gamma_3,\gamma_4\rangle_{0,5,1}\big)\nn\\
&=&-c_{0,3}(c_{1,2}c_{4,5}-c_{1,4}c_{2,5}+c_{1,5}c_{2,4})
+c_{0,4}(c_{1,2}c_{3,5}-c_{1,3}c_{2,5}+c_{1,5}c_{2,3})\nn\\
&&-c_{0,5}(c_{1,2}c_{3,4}-c_{1,3}c_{2,4}+c_{1,4}c_{2,3}),
\end{eqnarray}
\begin{eqnarray}\label{eq-groupOfTerms-S3-1}
(\ref{eq-groupOfTerms-S3})
&=&\sum_{\Gamma_e, \Gamma_f\in H^3(X)}\Big(\big(c_{0,3}\cdot(\gamma_4,\Gamma_e)-c_{0,4}\cdot(\gamma_3,\Gamma_e)+(\gamma_0,\Gamma_e)\cdot c_{3,4}\big)\nn\\
&&g^{ef}\big((\Gamma_f,\gamma_1)\cdot c_{2,5}-(\Gamma_f,\gamma_2)\cdot c_{1,5}+(\Gamma_f,\gamma_5)\cdot c_{1,2}\big)\nn\\
&&-\big(c_{0,3}\cdot(\gamma_5,\Gamma_e)-c_{0,5}\cdot(\gamma_3,\Gamma_e)+(\gamma_0,\Gamma_e)\cdot c_{3,5}\big)\nn\\
&&g^{ef}\big((\Gamma_f,\gamma_1)\cdot c_{2,4}-(\Gamma_f,\gamma_2)\cdot c_{1,4}+(\Gamma_f,\gamma_4)\cdot c_{1,2}\big)\nn\\
&&+\big(c_{0,4}\cdot(\gamma_5,\Gamma_e)-c_{0,5}\cdot(\gamma_4,\Gamma_e)+(\gamma_0,\Gamma_e)\cdot c_{4,5}\big)\nn\\
&&g^{ef}\big(\Gamma_f,\gamma_1)\cdot c_{2,3}-(\Gamma_f,\gamma_2)\cdot c_{1,3}+(\Gamma_f,\gamma_3)\cdot c_{1,2}\big)\Big)\nn\\
&=&c_{0,3}(c_{4,1}c_{2,5}-c_{4,2}c_{1,5}+c_{4,5}c_{1,2})
-c_{0,4}(c_{3,1}c_{2,5}-c_{3,2}c_{1,5}+c_{3,5}c_{1,2})\nn\\
&&+c_{3,4}(c_{0,1}c_{2,5}-c_{0,2}c_{1,5}+c_{0,5}c_{1,2})\nn\\
&&-c_{0,3}(c_{5,1}c_{2,4}-c_{5,2}c_{1,4}+c_{5,4}c_{1,2})
+c_{0,5}(c_{3,1}c_{2,4}-c_{3,2}c_{1,4}+c_{3,4}c_{1,2})\nn\\
&&-c_{3,5}(c_{0,1}c_{2,4}-c_{0,2}c_{1,4}+c_{0,4}c_{1,2})\nn\\
&&+c_{0,4}(c_{5,1}c_{2,3}-c_{5,2}c_{1,3}+c_{5,3}c_{1,2})
-c_{0,5}(c_{4,1}c_{2,3}-c_{4,2}c_{1,3}+c_{4,3}c_{1,2})\nn\\
&&+c_{4,5}(c_{0,1}c_{2,3}-c_{0,2}c_{1,3}+c_{0,3}c_{1,2}),
\end{eqnarray}
and
\begin{eqnarray}\label{eq-groupOfTerms-S4-1}
(\ref{eq-groupOfTerms-S4})&=&
\frac{1}{3}\langle \gamma_0, \gamma_3, \gamma_4,\gamma_5\rangle_{0,5,1}
\langle \sfh_2, \gamma_1,\gamma_2\rangle_{0,3,1}
+\frac{1}{3}\langle \gamma_0, \gamma_3, \gamma_4,\gamma_5, \sfh_3\rangle_{0,5,2} c_{1,2}\nn\\
&=&-c_{1,2}(c_{0,3}c_{4,5}-c_{0,4}c_{3,5}+c_{0,5}c_{3,4}).
\end{eqnarray}
Denote by $c$ the skew-symmetric matrix $(c_{i,j})_{0\leq i,j\leq 5}$. Summing (\ref{eq-groupOfTerms-S2-1})-(\ref{eq-groupOfTerms-S4-1}) one gets
\begin{equation}\label{eq-sum-S2-S3-S4}
	\langle \psi \gamma_0, \gamma_1,\dots,\gamma_5\rangle_{0,6,2}=(\ref{eq-groupOfTerms-S2})+(\ref{eq-groupOfTerms-S3})+(\ref{eq-groupOfTerms-S4})
	=\mathrm{Pf}(c).
\end{equation}
On the other hand by (\ref{eq-lem-TRR-genus1-deg2-simplification-2-5point}) one also finds
\begin{eqnarray}\label{eq-TRR-genus1-deg2-simplification-2-2}
&&\sum_{i=1}^{5}(-1)^{i-1}\langle \gamma_0,\gamma_i,1\rangle_{0,3,0}
 \langle  \sfh_3,\dots,\hat{\gamma_i},\dots\rangle_{0,5,2}\nn\\
&=&3\big(c_{0,1}(c_{2,3}c_{4,5}-c_{2,4}c_{3,5}+c_{2,5}c_{3,4})
-c_{0,2}(c_{1,3}c_{4,5}-c_{1,4}c_{3,5}+c_{1,5}c_{3,4})\nn\\
&&c_{0,3}(c_{1,2}c_{4,5}-c_{1,4}c_{2,5}+c_{1,5}c_{2,4})
-c_{0,4}(c_{1,2}c_{3,5}-c_{1,3}c_{2,5}+c_{1,5}c_{2,3})\nn\\
&&+c_{0,5}(c_{1,2}c_{3,4}-c_{1,3}c_{2,4}+c_{1,4}c_{2,3})\big)\nn\\
&=& 3\mathrm{Pf}(c).
\end{eqnarray}
Hence (\ref{eq-TRR-genus1-deg2-simplification-2}).
\end{proof}

\begin{proof}[Proof of Theorem \ref{thm-8points-cubic3fold}]
As a consequence of (\ref{eq-TRR-genus1-deg2-3})  and (\ref{eq-TRR-genus1-deg2-simplification-2}),
\begin{equation}\label{eq-TRR-genus1-deg2-4}
\frac{1}{24}\sum_{\Gamma_a\in H^3(X)}\langle \Gamma_a,\Gamma^a,\gamma_0,\gamma_1,\gamma_2,\gamma_3,\gamma_4,\gamma_5\rangle_{0,8,2}
=\frac{1}{9}\langle \sfh_2,\gamma_0,\gamma_1,\gamma_2,\gamma_3,\gamma_4,\gamma_5\rangle_{0,7,2}.
\end{equation}
Recall that $m=\mathrm{rank}\ H^*_{\mathrm{prim}}(X)=10$. So by Proposition \ref{prop-initialValues-meaning} and Lemma \ref{lem-contraction-Pfaffian},
\begin{equation}\label{eq-thm-8points-cubic3fold-1}
	\sum_{\Gamma_a\in H^3(X)}\langle \Gamma_a,\Gamma^a,\gamma_0,\gamma_1,\gamma_2,\gamma_3,\gamma_4,\gamma_5\rangle_{0,8,2}=-4 F^{(4)}(0) \mathrm{Pf}(c),
\end{equation}
where $c=\big((\gamma_i,\gamma_j)_{0\leq i,j\leq 5}\big)$ is the matrix used in the above proof. 
As a variant of Proposition \ref{prop-initialValues-meaning}, one has
\begin{equation}\label{eq-thm-8points-cubic3fold-2}
	\langle \sfh_2,\gamma_0,\gamma_1,\gamma_2,\gamma_3,\gamma_4,\gamma_5\rangle_{0,7,2}
	=\frac{\partial F^{(3)}}{\partial t^2}(0) \mathrm{Pf}(c)
\end{equation}
As we can take $\gamma_0,\dots,\gamma_5$ to be a part of a symplectic basis of $H^3(X)$ so that $\mathrm{Pf}(c)\neq 0$, 
(\ref{eq-8points-cubic3fold}) follows from (\ref{eq-TRR-genus1-deg2-4}), (\ref{eq-thm-8points-cubic3fold-1}) and (\ref{eq-thm-8points-cubic3fold-2}).
\end{proof}

\section{Higher order constant terms: examples and conjectures}\label{sec:HigherOrderLeadingTerms}
\subsection{Double root recursion}\label{sec:sqrtRecursion}

Recall (\ref{eq-wdvv24expand}) written in $\tau$-coordinates,  for $l\geq 2$, ($2\leq l\leq \frac{m}{2}$ when $n$ is odd)
\begin{eqnarray}\label{eq-wdvv24expand-sqrt}
\sum_{k=1}^{l}\sum_{a=0}^n\sum_{b=0}^n
\frac{\partial_{\tau^a}F^{(k)}\eta^{ab}\partial_{\tau^b}F^{(l+1-k)}}{(k-1)!(l-k)!}
+2\sum_{k=2}^{l}\frac{F^{(k)}F^{(l+2-k)}}{(k-2)!(l-k)!}=0.
\end{eqnarray}
For $I=(i_1,\dots,i_r)$ and $J=(j_1,\dots,j_r)\in \mathbb{Z}^{r}$, we say $J\leq I$ if $j_k\leq i_k$ for $1\leq k\leq r$. Define
\begin{equation*}
	I-J=(i_1-j_1,\dots,i_r-j_r),
\end{equation*}
and
\begin{equation}\label{eq-binomOfLists}
\binom{I}{J}=\prod_{k=1}^r \binom{i_k}{j_k}.
\end{equation}
We denote $(0,\dots,0)\in \mathbb{Z}^r$ by 0, when no confusion arises in the context.
For $I=(i_0,i_1,\dots,i_{n})\in \mathbb{Z}_{\geq 0}^{n+1}$, we define
\begin{equation}\label{eq-def-tauI}
\partial_{\tau^I}:=(\partial_{\tau^0})^{i_0}\circ\dots\circ (\partial_{\tau^{n}})^{i_{n}}.
\end{equation}

Let $I=(i_0,i_1,\dots,i_{n})\in \mathbb{Z}_{\geq 0}^{n+1}$ be given. 
We multiply (\ref{eq-wdvv24expand-sqrt}) by $(l-1)!$, then apply the operator $\partial_{\tau^I}$, and then take the constant term of the resulted equation, we get 
\begin{eqnarray}\label{eq-constantTerm-I-0}
&&\sum_{k=1}^{l}\sum_{0\leq J\leq I}\sum_{a=0}^n\sum_{b=0}^n
	\binom{l-1}{k-1}\binom{I}{J}
\partial_{\tau^I}\partial_{\tau^a}F^{(k)}(0)\eta^{ab}\partial_{\tau^{I-J}}\partial_{\tau^b}F^{(l+1-k)}(0)\nn\\
&&+2(l-1)\sum_{k=2}^{l}\sum_{0\leq J\leq I}
\binom{l-2}{k-2}\binom{I}{J} F^{(k)}(0)F^{(l+2-k)}(0)=0.
\end{eqnarray}
By the proof of Theorem \ref{thm-reconstruction-I}, the LHS of (\ref{eq-constantTerm-I-0}) can be written as 
a polynomial of 
\begin{equation*}
	F^{(2)}(0),F^{(3)}(0),\dots,F^{(l)}(0).
\end{equation*}
One can see this from Lemma \ref{lem-recursion-F(1)-WDVV-tau} and \ref{lem-recursion-F(l)-WDVV-tau} manifestly. 
In Appendix \ref{sec:algorithm}, we provide an algorithm to compute the equation \eqref{eq-constantTerm-I-0} for $F^{(i)}(0)$ for $i\geq 2$.
We denote (the LHS of) the  resulted equation by $\mathrm{Eqc}(n,\mathbf{d},l,I)$ (see \ref{sub:equations_of_higher_order_leading_terms}).

\begin{definition}
Let $g(z)=a z^2+b z+c$ be a quadratic polynomial with the variable $z$, where $a,b,c\in \mathbb{C}$. We say $g(z)$ is a \emph{complete square in $z$} if $a\neq 0$ and the two roots of $g(z)$ are equal.
\end{definition}

\begin{example}\label{example-sqrtRecursion-(3,(2,2,2))}
Let $n=3$, $\mathbf{d}=(2,2,2)$. Then $m=\mathrm{rank}\ H_{\mathrm{prim}}^3(X)=28$. For brevity we denote $z_l=F^{(l)}(0)$, which is also the notation in our package \texttt{QuantumCohomologyFanoCompleteIntersection}
in Macaulay2. Running the command
\begin{equation*}
      \mbox{\texttt{equationOfConstTerm}}\ \{3,\{2,2,2\},2,\{0,0,0,0\}\}
\end{equation*}
returns
\begin{equation*}\label{eq-eqc-(3,(2,2,2))-2}
      2\,{z}_{2}^{2}-16\,{z}_{2}+32.
\end{equation*}
In the above notation this means 
\begin{equation*}
      \mathrm{Eqc}\big(3,(2,2,2),2,(0,0,0,0)\big)=2\,{z}_{2}^{2}-16\,{z}_{2}+32=2(z_2-4)^2.
\end{equation*}
i.e. a complete square in $z_2$. This is an example for our general computation (\ref{eq-F2-case2-final}). We thus get $z_2=4$. Then running 
\begin{equation*}
      \mbox{\texttt{equationOfConstTerm}}\ \{3,\{2,2,2\},3,\{0,0,0,0\}\}
\end{equation*}
returns
\begin{equation}\label{eq-eqc-(3,(2,2,2))-3}
      112\,{z}_{2}^{2}+8\,{z}_{2}{z}_{3}-832\,{z}_{2}-32\,{z}_{3}+1536.
\end{equation}
Substituting $z_2=4$ into (\ref{eq-eqc-(3,(2,2,2))-3}), we get a trivial equation $0=0$. Similarly we have
\begin{eqnarray}\label{eq-eqc-(3,(2,2,2))-4}
      \mathrm{Eqc}\big(3,(2,2,2),4,(0,0,0,0)\big)&=&24\,{z}_{2}^{3}+21016\,{z}_{2}^{2}+936\,{z}_{2}{z}_{3}+12\,{z}_{3}^{2}+12\,{z}_{2}{z}_{4}\nn\\
      &&-162176\,{z}_{2}-3552\,{z}_{3}-48\,{z}_{4}+311680.
\end{eqnarray}
Substituting $z_2=4$ into (\ref{eq-eqc-(3,(2,2,2))-4}), we get
\begin{equation}\label{eq-eqc-(3,(2,2,2))-5}
     12 \left({{z}_{3}+8}\right)^{2},
\end{equation}
again a complete square in $z_3$! We thus obtain $z_3=-8$. Continuing this process, we have
\begin{eqnarray*}
 &&\mathrm{Eqc}\big(3,(2,2,2),5,(0,0,0,0)\big)=31540\,{z}_{2}^{3}+516\,{z}_{2}^{2}{z}_{3}+5680656\,{z}_{2}^{2}+231408\,{z}_{2}{z}_{3}+2256\,{z}_{3}^{2}\\
 &&+2240\,{z}_{2}{z}_{4}+48\,{z}_{3}{z}_{4}+16\,{z}_{2}{z}_{5}-45143360\,{z}_{2}-899328\,{z}_{3}
 -8576\,{z}_{4}-64\,{z}_{5}+87796480,\\
 &&\mathrm{Eqc}\big(3,(2,2,2),6,(0,0,0,0)\big)=13328\,{z}_{2}^{4}+30226526\,{z}_{2}^{3}+777199\,{z}_{2}^{2}{z}_{3}+3930\,{z}_{2}{z}_{3}^{2}\\
 &&+1920\,{z}_{2}^{2}{z}_{4}+2783397920\,{z}_{2}^{2}+79826776\,{z}_{2}{z}_{3}+727400\,{z}_{3}^{2}+724872\,{z}_{2}{z
      }_{4}+12880\,{z}_{3}{z}_{4}\\
      &&+60\,{z}_{4}^{2}+4280\,{z}_{2}{z}_{5}+80\,{z}_{3}{z}_{5}+20\,{z}_{2}{z}_{6}-23056231136\,{z}_{2}-320248528\,{z}_{3}-2831008\,{z}_{4}\\
      &&-16480\,{z}_{5}-80\,{z}_{6}+45797100032,\\
&& \dots.
\end{eqnarray*}
Since $n=3$ is odd, the allowed $l$ is allowed to be at most $\frac{m}{2}=14$. So  we stop at 
\[
\mathrm{Eqc}\big(3,(2,2,2),14,(0,0,0,0)\big)=
141803916446616765088\,{z}_{2}^{8}+\frac{243995302522352123860857835}{128}\,{z}_{2}^{7}+\dots
\]
which is too long to be spelt out. 
Substituting the solutions of  $z_i$ inductively in the above equations $\mathrm{Eqc}\big(3,(2,2,2),l,(0,0,0,0)\big)$ from $l=5$ to $l=15$, we get 
\begin{eqnarray}\label{eq-eqc-(3,(2,2,2))-6}
    &0,\  60 \left({{z}_{4}-32}\right)^{2},\ 0,\ 280\left({{z}_{5}+200}\right)^{2},\ 0,\ 1260\left({{z}_{6}-1728}\right)^{2},\nn\\
     &0,\ 5544 \left({{z}_{7}+19208}\right)^{2},\ 0,\ 24024\left({{z}_{8}-262144}\right)^{2}.
\end{eqnarray}
Up to now we take only $I=(0,0,0,0)$, i.e. we use only the equation arisen from the constant term of (\ref{eq-wdvv24expand-sqrt}). We can also take $I>(0,0,0,0)$. For instance:
\begin{eqnarray*}
      \mathrm{Eqc}\big(3,(2,2,2),3,(0,0,2,2)\big)&=&51275366400\,{z}_{2}^{3}-241350486261760\,{z}_{2}^{2}-10787741499392\,{z}_{2}{z}_{3}\\
      &&+1842040740511744\,{z}_{2}+43150965997568\,{z}_{3}-3509836805308416,
\end{eqnarray*}
substituting $z_2=4$ we get 0. 
\begin{eqnarray*}
      &&\mathrm{Eqc}\big(3,(2,2,2),4,(0,0,2,2)\big)=-93158375424\,{z}_{2}^{4}+1194036432470016\,{z}_{2}^{3}+14849800667136\,{z}_{2}^{2}{z}_{3}\\
      &&+12959514301038592\,{z}_{2}^{2}-2942643071877120\,{z}_{2}{z}_{3}-25577071312896\,{z}_{3}^{2}-25577071312896\,{z}_{2}{z}_{4}\\
      &&-182738305673068544\,{z}_{2}+11123742335827968\,{z}_{3}+102308285251584\,{z}_{4}+445569578177658880,
\end{eqnarray*}
substituting $z_2=4$ and $z_3=-8$ we get again 0. Examples shows that the equation $\mathrm{Eqc}\big(3,(2,2,2),l,I\big)$ becomes trivial after substituting the values of $z_i$ for $2\leq i\leq \lfloor \frac{l+2}{2}\rfloor$.
\end{example}

For more examples see Section \ref{sec:sqrtRecursion-examples}. We encapsulate such phenomena into a precise conjecture.

\begin{conjecture}[Double root recursion]\label{conj-sqrtRecursion}
Let $X=X_n(\mathbf{d})$ be a non-exceptional smooth complete intersection, with $n\geq 3$, and $\mathbf{d}\neq (3)$. Let $m=\mathrm{rank}\ H^n_{\mathrm{prim}}(X)$. Suppose $n$ is even (resp. $n$ is odd). Denote the unknown $F^{(i)}(0)$ in $\mathrm{Eqc}(n,\mathbf{d},l,0)$ by $z_i$. Then
\begin{enumerate}
      \item[(i)] For even $l\geq 2$ (resp. $2\leq l\leq \frac{m}{2}$), the polynomial $\mathrm{Eqc}(n,\mathbf{d},l,0)$ of $z_2,\dots,z_l$, after substituting the values of $z_i$ for  $2\leq i\leq \frac{l}{2}$, becomes a complete square in $z_{\frac{l+2}{2}}$;
      \item[(ii)] For odd $l\geq 2$ (resp. $2\leq l\leq \frac{m}{2}$), the polynomial $\mathrm{Eqc}(n,\mathbf{d},l,0)$ of $z_2,\dots,z_l$, after substituting the values of $z_i$ for  $2\leq i\leq \frac{l+1}{2}$, becomes 0;
      \item[(iii)] For $l\geq 2$  (resp. $2\leq l\leq \frac{m}{2}$) and $I> 0$, the polynomial $\mathrm{Eqc}(n,\mathbf{d},l,I)$ of $z_2,\dots,z_l$, after substituting the values of $z_i$ for  $2\leq i\leq \lfloor\frac{l+2}{2}\rfloor$, becomes 0.
\end{enumerate}
\end{conjecture}

\begin{remark}\label{rem:sqrtRecursion}
The significance of Conjecture \ref{conj-sqrtRecursion} is twofold.
\begin{enumerate}
      \item For a given dimension $n$ and a given multidegree $\mathbf{d}$,  one can check the statement (i) step by step, starting from $l=2$. If (i) is true for $l$, one gets the value of $z_{\frac{l}{2}+1}=F^{(\frac{l}{2}+1)}(0)$. Suppose  (i) is true at for all  $l$ in the allowed range, then  when $n$ is even  the full generating function $F$ of genus zero Gromov-Witten invariants can be reconstructed, and 
      when $n$ is odd, one gets $F^{(k)}$ for $k\leq \lfloor \frac{m}{4}\rfloor+1$.
      \item The parts (ii) and (iii) imply non-trivial relations among genus 0 invariants with both ambient and primitive insertions. As a consequence one can find a \emph{closed formula} of $F^{(k)}$ in terms of $F^{(0)}$ and $F^{(1)}$. We will illustrate this by the example $F^{(2)}$ in Section \ref{sec:closedFormula}. 
\end{enumerate}
\end{remark}

\subsection{Speculations in odd dimensions}\label{sec:sqrtRecursion-oddDim}
As we see in the statement of Conjecture \ref{conj-sqrtRecursion}, when the dimension $n$ is odd, the monodromy reduced WDVV equations (Theorem \ref{thm-wdvvoddthm}) do not supply enough information to compute $F^{(k)}(0)$ when $k>\frac{m}{4}+1$. In this section we make some attempts on this issue.

Our first attempt is to assume that the equations (\ref{eq-wdvv23odd}) and (\ref{eq-wdvv24odd}) hold at all orders of $s$, i.e. omit the restriction \textquotedblleft mod $s^{\frac{m}{2}}$\textquotedblright. Note that $F^{(k)}(0)=0$ for $k>\frac{m}{2}$. Then examples show that one cannot get consistent solutions of the resulted system of equations. So this naive attempt cannot be valid.

Now we propose a much bolder attempt. We assume again that the equations (\ref{eq-wdvv23odd}) and (\ref{eq-wdvv24odd}) hold at all orders of $s$. But we do not use $F^{(k)}(0)=0$ for $k>\frac{m}{2}$. 

\begin{example}\label{example-sqrtRecursion-(3,(2,2,2))-hypothetical}
Using this \emph{hypothetical} approach, we continue the computation as Example \ref{example-sqrtRecursion-(3,(2,2,2))}. We get 
\begin{multline}\label{eq-sqrtRecursion-(3,(2,2,2))-hypothetical-1}
      \mathrm{Eqc}\big(3,(2,2,2),15,(0,0,0,0)\big)=\frac{397233094426766697927933005}{64}\,{z}_{2}^{8}\\
      +\frac{5030242236275270025855639}{128}\,{z}_{2}^{7}{z}_{3}+\frac{1585053987640659484569374404693}{64}\,{z}_{2}^{7}+\dots
\end{multline}
\begin{multline}\label{eq-sqrtRecursion-(3,(2,2,2))-hypothetical-2}
      \mathrm{Eqc}\big(3,(2,2,2),16,(0,0,0,0)\big)=8901138728521104194857128\,{z}_{2}^{9}\\
      +\frac{79657522496696989359751556310147}{512}\,{z}_{2}^{8}+\frac{904662082772496883754648382693}{512}\,{z}_{2}^{7}{z}_{3}+\dots
\end{multline}
Substituting, from (\ref{eq-eqc-(3,(2,2,2))-5}) and (\ref{eq-eqc-(3,(2,2,2))-6}),
\begin{equation*}
      z_2=4,z_3=-8,z_4=32,z_5=-200,z_6=1728,z_7=-19208,z_8=262144,
\end{equation*}
(\ref{eq-sqrtRecursion-(3,(2,2,2))-hypothetical-1}) becomes 0, and (\ref{eq-sqrtRecursion-(3,(2,2,2))-hypothetical-2}) becomes
\begin{equation*}
      102960 \left({{z}_{9}+4251528}\right)^{2}.
\end{equation*}
Continuing this process, we find
\begin{eqnarray*}
    &  z_{10}=80000000,\ z_{11}=-1714871048,\ z_{12}= 41278242816,\nn\\
    & z_{13}=-1102867934792, z_{14}= 32396521357312,\ z_{15}=-1037970703125000,\ \dots
\end{eqnarray*}
Then we conjecture that $F^{(k)}(0)=z_l$ for $l\leq 14$. This can be written as a uniform formula
\begin{equation}\label{eq-sqrtRecursion-(3,(2,2,2))-hypothetical-5}
      F^{(k)}(0)=8(-1)^k k^{k-3},\ 1\leq k\leq 14.
\end{equation}
Our computation in Example \ref{example-sqrtRecursion-(3,(2,2,2))} shows the validity of (\ref{eq-sqrtRecursion-(3,(2,2,2))-hypothetical-5}) for $k\leq 8$. \pqed
\end{example}

We encapsulate this hypothetical approach as the following conjecture.
\begin{conjecture}\label{conj-sqrtRecursion-oddDim}
Let $X=X_n(\mathbf{d})$ be a non-exceptional smooth Fano complete intersection, with odd $n\geq 3$, and $\mathbf{d}\neq (3)$ or $(2,2)$. 
Let $m=\mathrm{rank}\ H^n(X)$. Let $\gamma_0,\dots,\gamma_n$ be a basis of $H^*_{\mathrm{amb}}(X)$, and $t^0,\dots,t^n$ be the dual basis. Let $g_{e,f}=(\gamma_e,\gamma_f)$ be the Poincaré pairing, and $(g^{e,f})_{0\leq e,f\leq n}$ the dual matrix of $(g_{e,f})$. Let $E$ be the Euler vector field. Let $G$ be a series in $t^0,\dots,t^n$ and $s$. Let 
\begin{equation*}
      G^{(i)}(t^0,\dots,t^n):= \Big(\big(\frac{\partial}{\partial s}\big)^i G\Big)|_{s=0}.
\end{equation*}
Let $E$ be the symmetric reducted Euler vector field (\ref{eq-EV-symmetricReducted}) and $c$ be the triple form (\ref{eq-tripleIntersectionForm}). 
Then the system 
\begin{subequations}\label{eq-system1-odd-G}
\begin{align}[left ={\empheqlbrace}]
& G^{(0)}=F^{(0)},\\
&\sum_{e=0}^n \sum_{f=0}^n \frac{\partial^3 G}{\partial t^a\partial t^b \partial t^e} g^{ef}\frac{\partial^2 G}{\partial s \partial t^f}+2s \frac{\partial^3 G}{\partial s\partial t^a \partial t^b} \frac{\partial^2 G}{\partial s\partial s} =\frac{\partial^2 G}{\partial s\partial t^a} \frac{\partial^2 G}{\partial s\partial t^b},\quad \mbox{for}\ 0\leq a,b\leq n,\\
&\sum_{e=0}^n \sum_{f=0}^n \frac{\partial^2 G}{\partial s\partial t^e}g^{ef}\frac{\partial^2 G}{\partial s\partial t^f}+2s \big(\frac{\partial^2 G}{\partial s\partial s})^2=0, \\
&EG=(3-n)G+\mathsf{a}(n,\mathbf{d})\frac{\partial}{\partial t^1}c, \\
&G^{(k)}(0)=0\ \mbox{unless}\ \frac{(n-2)k-n+3}{\mathsf{a}(n,\mathbf{d})}\in \mathbb{Z},
\end{align}
\end{subequations}
has a unique solution, and the solution can be given by the double root recursion as Conjecture \ref{conj-sqrtRecursion}. Moreover, the properties (i)-(iii) in Conjecture \ref{conj-sqrtRecursion} holds without any upper bound on $k$. Finally, denote by $\widetilde{F}$ the unique solution to the above system. Then 
\begin{equation}\label{eq-sqrtRecursion-OddDim}
      F=\sum_{k=0}^{\frac{m}{2}}\frac{s^k}{k!}\widetilde{F}^{(k)}.
\end{equation}
\end{conjecture}

I have no direct verifications for the statement (\ref{eq-sqrtRecursion-OddDim}) in Conjecture \ref{conj-sqrtRecursion-oddDim}.  I  can only take the uniformity of the formula (\ref{eq-conj-dim3-(2,2,2)}) as  indirect evidence. 
\begin{question}\label{ques-series-G}
Assume that Conjecture \ref{conj-sqrtRecursion-oddDim} is true. Then what is the geometry underlying the series $\widetilde{F}$? Does it encode new invariants?
\end{question}

\begin{remark}\label{rem:sqrtRecursion-hypergeometric-chow}
The equation  (\ref{eq-sqrtRecursion-OddDim}) is a geometric statement. 
Recall that the coefficients of $F^{(l)}$ can be interpreted as ratios (Proposition \ref{prop-F(l)-ratio}). This suggests us to directly define the coefficients of $\widetilde{F}^{(l)}$ as ratios; if this can be done  there would be no conflict with the anti-commutativity of Gromov-Witten invariants with insertions of odd real degrees.


The statements in Conjecture \ref{conj-sqrtRecursion}, and Conjecture \ref{conj-sqrtRecursion-oddDim} except  (\ref{eq-sqrtRecursion-OddDim}), are formal properties of the relevant system of equations. Note that our algorithm, which we describe in detail in Appendix \ref{sec:algorithm}, starts from the small $J$-function, which is expressed as a hypergeometric series by the mirror formula. So I would like to  regard the double root recursion for the system (\ref{eq-system1-odd-G}), both in even and odd dimensions,  as a deeply-hidden formal property of hypergeometric series.
\end{remark}

\subsection{Examples of \texorpdfstring{$F^{(k)}(0)$}{F(k)(0)}}\label{sec:sqrtRecursion-examples}
In this section, we provide examples of $F^{(k)}(0)$ of non-exceptional complete intersections and make some observations.
First recall that
\begin{eqnarray*}
      m&=&\mathrm{rank}\ H_{\mathrm{prim}}^n(X_n(\mathbf{d}))=(-1)^n \big(\int_X c_n(T_X)-(n+1)\big)\nn\\
      &=&\prod_{i=1}^r d_i\cdot \mathrm{Coeff}_{h^n}\frac{(1+h)^{n+r+1}}{\prod_{i=1}^r(1+d_i h)}-(n+1).
\end{eqnarray*}
This is an important number when $n$ is odd. For by (\ref{eq-Sym}), $F^{(k)}(0)=0$ if $k>\frac{m}{2}$. 
On the other hand by (\ref{eq-Dim}), $F^{(k)}(0)=0$ unless 
\begin{equation}\label{eq-degreeBeta}
\beta(k):=\frac{(n-2)k-n+3}{\mathsf{a}(n,\mathbf{d})}
\end{equation}
is an integer. This imposes conditions on $(n,\mathbf{d})$. For example, if $\mathbf{d}=(4)$, then $\beta$  takes integer values only when $n=3$.

\subsubsection{Cubic hypersurfaces}
The case $\mathbf{d}=(3)$, i.e. the cubic hypersurfaces, is excluded in the assumption of  Conjecture \ref{conj-sqrtRecursion}. In fact by Theorem \ref{thm-reconstructcubicandquadric} $F^{(k)}(0)$ can be computed by an essentially linear recursion on the leading terms. By Theorem \ref{thm-cubic5} or Theorem \ref{thm-4points-fanoIndex-(n-1)} we have  $F^{(2)}(0)=1$ in all dimensions $n\geq 3$.  
\begin{example}
The case $n=3$ is an atypical case in Theorem \ref{thm-reconstructcubicandquadric}. We have  $\mathrm{rank}\ H^3(X)=10$. So the only possible $3\leq k\leq 5$ such that $\beta(k)\in \mathbb{Z}$ is $4$. By Theorem \ref{thm-8points-cubic3fold}, we have
\[
F^{(4)}(0)=0.
\]
\end{example}
For $n>3$, using an algorithm based on the proof of Theorem \ref{thm-reconstructcubicandquadric}, we get the following results.
\begin{example}
$n=4$.
\[
F^{(5)}(0)=0,\ F^{(8)}(0)=1,\ F^{(11)}=-20,\ F^{(14)}(0)=6363.
\]
\end{example}

\begin{example}
$n=5$.
\[
F^{(6)}(0)=0,\ F^{(10)}(0)=33,\ F^{(14)}(0)=72912.
\]
\end{example}

\begin{example}
$n=6$.
\[
F^{(7)}(0)=0,\ F^{(12)}(0)=3231,\ F^{(17)}(0)= -539677008.
\]
\end{example}

\begin{example}
$n=7$.
\[
F^{(8)}(0)=0,\ F^{(14)}(0)=547335.
\]
\end{example}

In these examples one finds that $F^{(n+1)}(0)=0$ for $n$-dimensional cubic hypersurfaces. But we recall that the methods for $n=3$ and $n>3$ cases on this computation are different.

\subsubsection{Non-exceptional complete intersections other than the cubic hypersurfaces}
We define a function  
\begin{equation*}
      \mbox{\texttt{sqrtRecursion}}
\end{equation*}
in our Macaulay2 package \texttt{QuantumCohomologyFanoCompleteIntersection}. With the input $\{n,\mathbf{d},k\}$ it automates the verification of Conjecture \ref{conj-sqrtRecursion} (i) for $l=2k-2$. For example,
\begin{equation*}
        \mbox{\texttt{sqrtRecursion}}\ \{3,\{2,2,2\},2\}
 \end{equation*} 
 returns
 \begin{equation*}
       \left\{2\,{z}_{2}^{2}-16\,{z}_{2}+32,\,4\right\}
 \end{equation*}
 where $4$ is the unique solution to the complete square of the first entry. Similarly, 
 \begin{equation*}
        \mbox{\texttt{sqrtRecursion}}\ \{3,\{2,2,2\},3\}
 \end{equation*} 
 returns
 \begin{equation*}
    \left\{12\,{z}_{3}^{2}+192\,{z}_{3}+768,\,-8\right\}.    
 \end{equation*}
In our design, if  Conjecture \ref{conj-sqrtRecursion} (i) fails for $2k-2$, the command 
\begin{equation*}
      \mbox{\texttt{sqrtRecursion}}\ \{n,\mathbf{d},2k-2\}
\end{equation*}
will return the failing information and the first $i$ such that Conjecture \ref{conj-sqrtRecursion} (i) fails  for $2i-2$. In all our computations up to now, no failing information shows up for the non-exceptional complete intersections other than the cubic hypersurfaces. In the following we display some examples. We begin by repeating the computations for $(n,\mathbf{d})=(3,(2,2,2))$. It turns out that in dimension 3 we have always a simple conjectural formula for $F^{(k)}(0)$.
\begin{example}
$n=3$, $\mathbf{d}=(2,2,2)$.  $\mathrm{rank}\ H_{\mathrm{prim}}^3(X)=28$. Using Theorem \ref{thm-estimate-rank-primCoh} it is easily seen that this is the smallest dimension of primitive cohomology among the non-exceptional smooth complete intersections
besides the cases $X_n(3)$ ($n\geq 3$) and $X_n(2,2)$ (odd $n\geq 3$).
 The results indicated in blue are hypothetical.
\[
F^{(2)}(0)=4=2^2,\ F^{(3)}(0)=-8=-2^3,\ F^{(4)}=32=2^5,\ F^{(5)}(0)=-200=-2^{3}5^{2}, 
\]
\[
F^{(6)}(0)=1728=2^{6}3^{3},\
F^{(7)}(0)=-19208=-2^{3}7^{4},\ F^{(8)}(0)=262144=2^{18},
\]
\textcolor{blue}
{
\[
F^{(9)}(0)=-4251528=-2^{3}3^{12},\ F^{(10)}=80000000=2^{10}5^{7},
\]
\[
F^{(11)}(0)=-1714871048=-2^{3}11^{8},\ F^{(12)}(0)=41278242816=2^{21}3^{9},
\]
\[
F^{(13)}(0)=-1102867934792=-2^{3}13^{10},\ F^{(14)}(0)=32396521357312=2^{14}7^{11}.
\]
}
\end{example}
\begin{conjecture}\label{conj-dim3-(2,2,2)}
When $n=3$, $\mathbf{d}=(2,2,2)$,
\begin{equation}\label{eq-conj-dim3-(2,2,2)}
	F^{(k)}(0)=8(-1)^k k^{k-3},\ \mbox{for}\ 1\leq k\leq 14.
\end{equation}
\end{conjecture}

\begin{example}
$n=3$, $d=4$. $\mathrm{rank}\ H_{\mathrm{prim}}^3(X)=60$.
\begin{equation*}
F^{(2)}(0)=72=2^{3}3^{2},\
      F^{(3)}(0)=-864=-2^{5}3^{3},\  F^{(4)}(0)=20736=2^{8}3^{4},
\end{equation*}
\[
F^{(5)}(0)=-777600=-2^{7}3^{5}5^{2},\ F^{(6)}(0)=40310784=2^{11}3^{9},
\]
\[
F^{(7)}(0)=-2688505344=-2^{9}3^{7}7^{4},\ F^{(8)}(0)=220150628352=2^{25}3^{8}.
\]

\end{example}

\begin{conjecture}
When $n=3$, $d=4$,
\[
F^{(k)}(0)=(-1)^k 2^{k+2}3^k k^{k-3},\ \mbox{for}\ 1\leq k\leq 30.
\]
\end{conjecture}

\begin{example}
$n=3$, $\mathbf{d}=(2,3)$. $\mathrm{rank}\ H_{\mathrm{prim}}^3(X)=40$.
\[
F^{(2)}(0)=12=2^{2}3,\ F^{(3)}(0)=-48=-2^{4}3,\ F^{(4)}(0)=384=2^{7}3,
\]
\[
 F^{(5)}(0)=-4800=-2^{6}3\cdot 5^{2},\ F^{(6)}(0)=82944=2^{10}3^{4},
\]
\[
F^{(7)}(0)=-1843968=-2^{8}3\cdot 7^{4},\ F^{(8)}(0)=50331648=2^{24}3.
\]

\end{example}
\begin{conjecture}
When $n=3$, $\mathbf{d}=(2,3)$,
\[
F^{(k)}(0)=(-1)^k 2^{k+1}3\cdot k^{k-3},\ \mbox{for}\ 1\leq k\leq 20.
\]
\end{conjecture}

\begin{example}
$n=4$, $d=5$. 
\[
F^{(2)}(0)=2088000=2^{6}3^{2}5^{3}29,\
F^{(3)}(0)=-413985600000=-2^{9}3^{3}5^{5}7\cdot 37^{2},
\]
\[
F^{(4)}(0)=199423892160000000=2^{12}3^{4}5^{7}1877\cdot 4099,
\]
\[
F^{(5)}(0)=-161117359277760000000000=-2^{15}3^{5}5^{10}359\cdot 1117\cdot 5167.
\]
\end{example}

\begin{example}
$n=5$, $d=5$. $F^{(l)}(0)=0$ when $2\nmid l$.
\[
F^{(2)}(0)=1440=2^{5}3^{2}5,\ F^{(4)}(0)=26077593600=2^{15}3^{5}5^{2}131,
\]
\[
F^{(6)}(0)=5549953864826880000=2^{20}3^{6}5^{4}11616677.
\]
\end{example}

\begin{example}
$n=6$, $d=5$. $F^{(l)}(0)=0$ when $3\nmid l$.
\[
F^{(3)}(0)=-69120=-2^{9}3^{3}5,\ F^{(6)}(0)=1315624550400000=2^{18}3^{6}5^{5}2203.
\]
\end{example}

\begin{example}
$n=6$, $\mathbf{d}=(3,3)$. $F^{(l)}(0)=0$ when $3\nmid l$.
\[
F^{(3)}(0)=-576,\ 
\]
\end{example}

\begin{example}
$n=5$, $d=6$. 
\[
F^{(2)}(0)=20558229235200,
\]
\[
F^{(3)}(0)=-15278161374523225276416000,
\]
\[
F^{(4)}(0)=30403236725336498691688248933089280000.
\]
\end{example}

\begin{example}
$n=4$, $\mathbf{d}=(2,2,2,2)$. 
\[
F^{(2)}(0)=896=2^{7}7,\ F^{(3)}(0)=-600832=-2^{8}2347,\ F^{(4)}(0)=982757376=2^{12}3^{2}53\cdot 503,
\]
\[
F^{(5)}(0)=-2699323002880=-2^{12}5\cdot 7\cdot 18828983,
\]
\[
F^{(6)}(0)=10606611483492352=2^{15}211\cdot 337\cdot 881\cdot 5167,
\]
\[
F^{(7)}(0)=-54646751402746904576=-2^{16}7\cdot 5227\cdot 22789451819.
\]
\end{example}

\begin{example}
$n=7$, $d=6$. $F^{(l)}(0)=0$ when $l\not\equiv 2 \mod 3$.

\[
F^{(2)}(0)=43200,\
F^{(5)}(0)=-21381767820096307200000,
\]
\[
F^{(8)}(0)=735166560884689201501389038026752000000000.
\]
\end{example}

\begin{example}
$n=9$, $d=7$. $F^{(l)}(0)=0$ when $l\not\equiv 2 \mod 4$.
\[
F^{(2)}(0)=1814400.
\]
\end{example}

\subsection{Integrality and positivity}\label{sec:integrality-positivity}
From the above examples we find some patterns.
\begin{conjecture}\label{conj-integrality-positivity}
When $\mathbf{d}= (3)$, $F^{(n+1)}(0)=0$. When $\mathbf{d}\neq (3)$, or $\mathbf{d}=(3)$ and $l\neq n+1$, $F^{(l)}(0)$ is a positive integer if $l$ is even, and is a negative integer if $l$ is odd.
\end{conjecture}

\begin{remark}\label{rem:integrality}
Let $L$ be the lattice $H^n_{\mathrm{prim}}(X)\cap H^n(X;\mathbb{Z})$. For even dimensions $n$ (resp. odd dimensions $n$), let $I_k$ be the ideal of $\mathbb{Z}$ generated by all the permanent Pfaffians (resp. Pfaffians)
\begin{equation*}
	\mathrm{P}\big(G(\alpha_1,\dots,\alpha_{2k})\big),\ (\mbox{resp.}\ \mathrm{Pf}\big(G(\alpha_1,\dots,\alpha_{2k})\big))
\end{equation*}
for $\alpha_i$ running over $L$. By Proposition \ref{prop-initialValues-meaning}  and the integrality of genus 0 Gromov-Witten invariants of semipositive symplectic manifolds (\cite[Theorem A]{Ruan96}, and also \cite[Theorem 7.1.1]{MS12}), we have
\begin{equation*}
	F^{(k)}(0)\cdot I_k\subset \mathbb{Z}.
\end{equation*}
In the odd dimensions,  $L=H^n(X;\mathbb{Z})$ and thus is a unimodular skew-symmetric lattice. By \cite[P. 79, Théorème 1]{Bou59}, such a lattice is a direct sum of the lattice $\mathbb{Z}\{e_1,e_2\}$ with $(e_1,e_2)=1$. Then $I_k$ is the unit ideal for $1\leq k\leq m$. Hence the integrality of $F^{(k)}$ follows.

The even dimensional cases are more complicated. By \cite[Theorem 2.1]{LW81} $L$ is an even lattice, i.e. $(x,x)\in 2 \mathbb{Z}$ for all $x\in L$. Then one can show that for  $k>m$, $I_k\subset 2 \mathbb{Z}$. So one cannot show the integrality of $F^{(k)}(0)$ for all $k$ in this way. 
But note that we have not exhausted the full strength of  the integrality of $F^{(k)}(0)$ from the integrality of genus 0 Gromov-Witten invariants of semipositive symplectic manifolds; for there exist integral vectors in $H^n(X;\mathbb{Z})$ that cannot be written as an integral combination of $\sfh_{\frac{n}{2}}$ and a vector in $L$. We will not pursue this further in this paper.


The positivity and negativity seem  mysterious and  I can say nothing at present.
\end{remark}

\subsection{Conjectural closed formula of \texorpdfstring{$F^{(2)}$}{F(2)}}\label{sec:closedFormula}
In this section we elaborate Remark \ref{rem:sqrtRecursion} (ii), by giving a conjectural formula of $F^{(2)}$ in terms of $F^{(0)}$ and $F^{(1)}$. The  method applies to arbitrary $F^{(k)}$. 
Let $t^0,\dots,t^n$ be the basis dual to $1,\sfh,\dots,\sfh_n$. To save space we often use the notation $F_a$ to denote the partial derivative $\partial_{t^a}F$, and use Einstein's summation convention, where the indices run over $0,\dots,n$.

Recall  (\ref{eq-wdvv23expand}), for $k\geq 1$,
\begin{eqnarray*}
&&F_{ abe}^{(0)}g^{ef}F_{ f}^{(k+1)}+2k F_{ ab}^{(1)}F^{(k+1)}-F_{ a}^{(k+1)}F_{ b}^{(1)}-F_{ a}^{(1)}F_{ b}^{(k+1)}\nn\\
&=&\sum_{j=2}^{k}\binom{k}{j-1}F_{ a}^{(j)}F_{b}^{(k-j+2)}-\sum_{j=1}^{k}\binom{k}{j}F_{ abe}^{(j)}g^{ef}F_{ f}^{(k-j+1)}\nn\\
&&-2k\sum_{j=2}^{k}\binom{k-1}{j-1}F_{ ab}^{(j)}F^{(k-j+2)}.
\end{eqnarray*}
Taking $a=1$ and $1\leq b\leq n$ we get, for $k\geq 2$,
\begin{eqnarray}\label{eq-linearSystemOf-F(k)b}
      &&\sum_{i=1}^n (\sum_{e=0}^{n}F_{1be}^{(0)}g^{ei}-\delta_{i,1}F_b^{(1)}-\delta_{i,b}F_1^{(1)})F_i^{(k)}\nn\\
      &=&-2(k-1) F_{1b}^{(1)}F^{(k)}+\sum_{j=2}^{k-1}\binom{k-1}{j-1}F_{1}^{(j)}F_{b}^{(k-j+1)}-\sum_{j=1}^{k-1}\binom{k-1}{j}F_{1be}^{(j)}g^{ef}F_{ f}^{(k-j)}\nn\\
&&-2(k-1)\sum_{j=2}^{k-1}\binom{k-2}{j-1}F_{1b}^{(j)}F^{(k-j+1)}.
\end{eqnarray}

The Euler vector field
\begin{equation*}
      E= \sum_{i=0}^{n}(1-i)t^{i}\frac{\partial}{\partial t^i}+(2-n)s\frac{\partial }{\partial s}
+\mathsf{a}(n,\mathbf{d})\frac{\partial}{\partial t^1}
\end{equation*}
yields, for $k\geq 1$,
\begin{equation}\label{eq-linearSystemOf-F(k)1}
\mathsf{a}(n,\mathbf{d}) F_1^{(k)}
+\sum_{i=2}^{n}(1-i)t^{i}F_i^{(k)}=(nk-n-2k+3)F^{(k)}.
\end{equation}
For a given $k\geq 2$, we regard (\ref{eq-linearSystemOf-F(k)1}) and (\ref{eq-linearSystemOf-F(k)b}) as a linear system of $F_i^{(k)}$ for $1\leq i\leq n$. 
Denote by $\Phi$ be the matrix of the coefficients of this linear system. Namely, $\Phi$ is the $n\times n$ matrix with entries
\begin{eqnarray*}
      \Phi_{j}^i&=&\begin{cases}
  \mathsf{a}(n,\mathbf{d}),& \mbox{if}\ j=1,i=1,\\
  (1-i)t^i,                         & \mbox{if}\ j=1,i=\geq 2,\\
\sum_{e=0}^{n}F_{1,j-1,e}^{(0)}g^{ei}-\delta_{i,1}F_{j-1}^{(1)}-\delta_{i,j-1}F_1^{(1)}, & \mbox{if}\ 2\leq j\leq n
      \end{cases}\\
      &=&\begin{cases}
  \mathsf{a}(n,\mathbf{d}),& \mbox{if}\ j=1,i=1,\\
  (1-i)t^i,                         & \mbox{if}\ j=1,i=\geq 2,\\
\frac{1}{\prod_{i=1}^r d_i}F_{1,j-1,n-i}^{(0)}-\delta_{i,1}F_{j-1}^{(1)}-\delta_{i,j-1}F_1^{(1)}, & \mbox{if}\ 2\leq j\leq n.
      \end{cases}
\end{eqnarray*}
Theorem \ref{thm-reconstruction-I} essentially states that the constant matrix $\Phi(0)$ is  nonsingular. 
Then
\begin{eqnarray}\label{eq-linearSystemOf-F(k)i-solution}
      &&\begin{pmatrix}
      \partial_{t^1}F^{(k)}\\ \dots \\ \partial_{t^n}F^{(k)}
      \end{pmatrix}
      = \Phi^{-1}\bigg(\begin{pmatrix}
      nk-n-2k+3\\ (2-2k) F_{1,1}^{(1)} \\ \dots \\ (2-2k) F_{1,n-1}^{(1)}
      \end{pmatrix}F^{(k)}\nn\\
      &&+\begin{pmatrix}
      0\\
      \sum_{j=2}^{k-1}\binom{k-1}{j-1}F_{1}^{(j)}F_{1}^{(k-j+1)}-\sum_{j=1}^{k-1}\binom{k-1}{j}F_{1,1,e}^{(j)}g^{ef}F_{f}^{(k-j)}-2(k-1)\sum_{j=2}^{k-1}\binom{k-2}{j-1}F_{1,1}^{(j)}F^{(k-j+1)}\\
      \dots\\
      \sum_{j=2}^{k-1}\binom{k-1}{j-1}F_{1}^{(j)}F_{n-1}^{(k-j+1)}-\sum_{j=1}^{k-1}\binom{k-1}{j}F_{1,n-1,e}^{(j)}g^{ef}F_{ f}^{(k-j)}-2(k-1)\sum_{j=2}^{k-1}\binom{k-2}{j-1}F_{1,n-1}^{(j)}F^{(k-j+1)}
      \end{pmatrix}\bigg).\nn\\
\end{eqnarray}
In particular,
\begin{equation}\label{eq-linearSystemOf-F(2)i-solution}
      \begin{pmatrix}
      \partial_{t^1}F^{(2)}\\ \dots \\ \partial_{t^n}F^{(2)}
      \end{pmatrix}
      = \Phi^{-1} \begin{pmatrix}
      n-1\\ -2F_{1,1}^{(1)} \\ \dots \\ -2F_{1,n-1}^{(1)}
      \end{pmatrix}F^{(2)}
      -\Phi^{-1}\begin{pmatrix} 
      0\\ F_{1,1,e}^{(1)}g^{ef}F_f^{(1)} \\ \dots\\ F_{1,n-1,e}^{(1)}g^{ef}F_f^{(1)}
      \end{pmatrix}.
\end{equation}
Now we recall (\ref{eq-wdvv24expand}), for $k\geq 2$,
\begin{eqnarray}\label{eq-quadraticSystemOf-F(i)}
&&F_{ e}^{(1)}g^{ef}F_{ f}^{(k+1)}+2kF^{(2)}F^{(k+1)}\nn\\
&=&-\frac{1}{2}\sum_{j=2}^{k}\binom{k}{j-1}F_{ e}^{(j)}g^{ef}F_{ f}^{(k+2-j)}
-k\sum_{j=3}^{k}\binom{k-1}{j-2}F^{(j)}F^{(k+3-j)}.
\end{eqnarray}
Substituting (\ref{eq-linearSystemOf-F(k)i-solution}) and (\ref{eq-linearSystemOf-F(2)i-solution}) into (\ref{eq-quadraticSystemOf-F(i)}) we get
\begin{eqnarray*}
      &&\left(\frac{1}{\prod_{i=1}^r d_i} (\partial_{t^{n-1}}F^{(1)},\dots,\partial_{t^{0}}F^{(1)})
      \Phi^{-1} \begin{pmatrix}
      (n-2)(k-1)+1\\ (2-2k) F_{1,1}^{(1)} \\ \dots \\ (2-2k) F_{1,n-1}^{(1)}
      \end{pmatrix}+2(k-1)F^{(2)}\right)F^{(k)}\nn\\
&=&-\frac{1}{\prod_{i=1}^r d_i} (\partial_{t^{n-1}}F^{(1)},\dots,\partial_{t^{0}}F^{(1)})
      \Phi^{-1}\nn\\
&&\cdot      \begin{pmatrix}
      0\\
      \sum_{j=2}^{k-1}\binom{k-1}{j-1}F_{1}^{(j)}F_{1}^{(k-j+1)}-\sum_{j=1}^{k-1}\binom{k-1}{j}F_{1,1,e}^{(j)}g^{ef}F_{f}^{(k-j)}-2(k-1)\sum_{j=2}^{k-1}\binom{k-2}{j-1}F_{1,1}^{(j)}F^{(k-j+1)}\\
      \dots\\
      \sum_{j=2}^{k-1}\binom{k-1}{j-1}F_{1}^{(j)}F_{n-1}^{(k-j+1)}-\sum_{j=1}^{k-1}\binom{k-1}{j}F_{1,n-1,e}^{(j)}g^{ef}F_{ f}^{(k-j)}-2(k-1)\sum_{j=2}^{k-1}\binom{k-2}{j-1}F_{1,n-1}^{(j)}F^{(k-j+1)}
      \end{pmatrix}\nn\\
&&-\frac{1}{2}\sum_{j=2}^{k-1}\binom{k-1}{j-1}F_{ e}^{(j)}g^{ef}F_{ f}^{(k+1-j)}
-(k-1)\sum_{j=3}^{k-1}\binom{k-2}{j-2}F^{(j)}F^{(k+2-j)}.   
\end{eqnarray*}
The vanishing statements in Conjecture \ref{conj-sqrtRecursion} imply that, for $k\geq 3$, the coefficient of $F^{(k)}$ vanishes. So we get the following identities.

\begin{conjecture}\label{conj-formula-F(2)}
Let $X=X_n(\mathbf{d})$ be an $n$-dimensional smooth complete intersection of multidegree $\mathbf{d}$, with $n\geq 3$ and $\mathbf{d}\neq (3)$. Then
\begin{equation*}
      F^{(2)}=\frac{1}{\prod_{i=1}^r d_i} (\partial_{t^{n-1}}F^{(1)},\dots,\partial_{t^{0}}F^{(1)})
      \Phi^{-1} \begin{pmatrix}
      0\\ \partial_{t^1}\partial_{t^1} F^{(1)}  \\ \dots \\ \partial_{t^1}\partial_{t^{n-1}} F^{(1)}
      \end{pmatrix},
\end{equation*}
and
\begin{equation*}
       (\partial_{t^{n-1}}F^{(1)},\dots,\partial_{t^{0}}F^{(1)})
      \Phi^{-1} \begin{pmatrix}
      1\\ 0 \\ \dots \\ 0
      \end{pmatrix}=0.
\end{equation*}
\end{conjecture}

By analogy and numerical experiments, we have a conjecture also in the case $\mathbf{d}=(3)$.
\begin{conjecture}\label{conj-formula-F(2)-d=3}
For cubic hypersurfaces of dimension $n\geq 3$,
\begin{equation*}
      F^{(2)}=\frac{1}{3} (\partial_{t^{n-1}}F^{(1)},\dots,\partial_{t^{0}}F^{(1)})
      \Phi^{-1} \begin{pmatrix}
      \vspace{0.1cm}
      -\frac{n-1}{3}\\ \partial_{t^1}\partial_{t^{1}}  F^{(1)} \\ \dots \\ \partial_{t^1}\partial_{t^{n-1}} F^{(1)}
      \end{pmatrix}.
\end{equation*}
\end{conjecture}

We have checked Conjecture \ref{conj-formula-F(2)} and \ref{conj-formula-F(2)-d=3} in low $\mathbf{t}$-degrees in some of the examples in Section \ref{sec:sqrtRecursion-examples}.

\newpage

\begin{appendix}

\section{An estimate for  the rank of the primitive cohomology}\label{sec:estimate-dimPrim}

Suppose $n\geq 1$, $r\geq 1$, $\mathbf{d}=(d_{1},\dots,d_{r})\in \mathbb{Z}^{r}$, where $2\leq d_{1}\leq d_{2}\leq \cdots \leq d_{r}$.\
In this appendix we show an estimate on the rank of $H^{n}_{\mathrm{prim}}(X)$ for a smooth complete intersection of dimension $n$ and multidegree $\mathbf{d}$ in $\mathbb{P}^{n+r}$. This estimate is crucial in performing the symmetric reduction in Section \ref{sec:red}, and is also used in Example \ref{example-sqrtRecursion-(3,(2,2,2))}.
By Riemann-Roch we have
\begin{equation*}
 \mathrm{rank}\ H^{n}_{\mathrm{prim}}(X)=(-1)^{n}f(n,\mathbf{d}),
\end{equation*}
  where
\begin{equation*}
f(n,\mathbf{d}):=\prod_{i=1}^{r}d_{i}\cdot \mathrm{Coeff}_{x^{n}}\Bigg(\frac{(1+x)^{n+r+1}}{\prod_{i=1}^{r}(1+d_{i}x)}\Bigg)-(n+1).
\end{equation*}
\begin{lemma}\label{lem-induciton-rank-primCoh}
Let  $1\leq i\leq r$. Then
\begin{equation}\label{eq-induciton-rank-primCoh}
	f(n,\mathbf{d})+(d_{i}-1)f(n-1,\mathbf{d})=
	\begin{cases}
	d_1-1,& \mbox{if}\ r=1;\\
	d_i f\big(n,(d_1,\dots,\hat{d_i},\dots,d_r)\big)+d_i-1,& \mbox{if}\ r\geq 2.
	\end{cases}
\end{equation}
\end{lemma}
\begin{proof}
\begin{eqnarray*}
&& \mathrm{Coeff}_{x^n}\Big(\frac{(1+x)^{n+r+1}}{\prod_{i=1}^r (1+d_i x)}\Big)
+(d_i-1)\mathrm{Coeff}_{x^{n-1}}\Big(\frac{(1+x)^{n+1}}{\prod_{i=1}^r (1+d_i x)}\Big)\\
&=&  \mathrm{Coeff}_{x^n}\Big(\frac{(1+x)^{n+r+1}+(d_i-1)x(1+x)^{n+r}}{\prod_{i=1}^r (1+d_i x)}\Big)\\
&=&\mathrm{Coeff}_{x^n}\big(\frac{(1+x)^{n+r}}{\prod_{j\neq i}(1+d_j x)}\big).
\end{eqnarray*}
So 
\begin{eqnarray*}
	&&f(n,\mathbf{d})+(d_i-1)f(n-1,\mathbf{d})=
	\prod_{i=1}^r d_i\cdot \mathrm{Coeff}_{x^n}\big(\frac{(1+x)^{n+r}}{\prod_{j\neq i}(1+d_j x)}\big)
	-(n+1)-(d_i-1)n\\
	&=& d_i\Big(f\big(n,(d_1,\dots,\hat{d_i},\dots,d_r)\big)+n+1\Big)-d_i n-1
	=\mbox{RHS of (\ref{eq-induciton-rank-primCoh})}.
\end{eqnarray*}
\end{proof}

\begin{theorem}\label{thm-estimate-rank-primCoh}
If $\mathbf{d}\neq (2)$, 
\begin{eqnarray}\label{eq-estimate-rank-primCoh}
\mathrm{rk}\ H^{n}_{\mathrm{prim}}(X_{n}(\mathbf{d}))=(-1)^{n}f(n,\mathbf{d})\geq (d_r-1)^{n-1}\Big(\prod_{i=1}^{r}d_{i}\cdot\big(\sum_{i=1}^{r}d_{i}-r-2\big)+2\Big)+n-1.
\end{eqnarray}
\end{theorem}
\begin{proof}
When $n=1$,
\begin{eqnarray*}
-f(1,\mathbf{d})&=&-\prod_{i=1}^{r}d_{i}\cdot \mathrm{Coeff}_x\Bigg(\frac{(1+x)^{r+2}}{\prod_{i=1}^{r}(1+d_{i}x)}\Bigg)+2\\
&=&\prod_{i=1}^{r}d_{i}\cdot\Big(\sum_{i=1}^{r}d_{i}-r-2\Big)+2.
\end{eqnarray*}
Suppose  $r=1$. We show (\ref{eq-estimate-rank-primCoh}) by induction on $n$. Suppose it holds for $n-1$. Then by (\ref{eq-induciton-rank-primCoh}),
\begin{eqnarray*}
&&(-1)^n f(n,d)=(d-1)\big((-1)^{n-1}f(n-1,d)+(-1)^n\big)\\
&\geq& (d-1)\Big((d-1)^{n-1}\big(d(d-r-2)+2\big)+n-2+(-1)^n\Big)\\
&\geq & \begin{cases}
(d-1)^n\big(d(d-r-2)+2\big)+(d-1)(n-1),& \mbox{if $n$ is even};\\
(d-1)^n\big(d(d-r-2)+2\big)+(d-1)(n-3),& \mbox{if $n$ is odd}.
\end{cases}
\end{eqnarray*}
Since $d\neq 2$, we get (\ref{eq-estimate-rank-primCoh}) when $n$ is even or $n\geq 5$. We show the $n=3$ case directly:
\begin{eqnarray*}
&&-f(3,d)-\big((d-1)^3(d-2)+2\big)\\
&=& (d^4-5d^3+10d^2-10d+4)-\big((d-1)^3(d-2)+2\big)=d(d-3)\geq 0.
\end{eqnarray*}
Hence the $r=1$ is done. Then by induction on $r$, and (\ref{eq-induciton-rank-primCoh}), we have
\begin{equation}\label{eq-thm-estimate-rank-primCoh-1}
	(-1)^{n}f(n,\mathbf{d})\geq 2
\end{equation}
and
\begin{eqnarray}\label{eq-thm-estimate-rank-primCoh-2}
(-1)^{n}\big(f(n,\mathbf{d})+(d_{r}-1)f(n-1,\mathbf{d})\big)>0.
\end{eqnarray}
Then for fixed $r\geq 2$, (\ref{eq-estimate-rank-primCoh}) from (\ref{eq-thm-estimate-rank-primCoh-2}) and an induction on $n$.
\end{proof}

\begin{corollary}\label{cor-estimate-rank-primCoh}
Suppose $\mathbf{d}\neq (2)$. 
\begin{itemize}
\item[(i)] If $n$ is even,  $\mathrm{rk}\ H^{n}_{\mathrm{prim}}(X_{n}\big(\mathbf{d})\big)\geq 3$.
\item[(ii)] 
If $n\geq 3$ is odd, $\mathrm{rk}\ H^{n}_{\mathrm{prim}}(X_{n}\big(\mathbf{d})\big)\geq 4$.
\end{itemize}
\end{corollary} 

\section{Symmetric reduction of  the big \texorpdfstring{$J$}{J}-function}
\label{sec:symmetricReduction-JFunction}
Let $X$ be a non-exceptional smooth complete intersection of dimension $n$. Let $\gamma_0,\dots,\gamma_n$ be a basis of $H^*_{\mathrm{amb}}(X)$ and $\gamma_{n+1},\dots,\gamma_{n+m}$ an orthonormal basis of $H^n_{\mathrm{prim}}(X)$. Let $t^0,\dots,t^{n+m}$ be the dual basis. 
Recall that the big $J$-function of $X$ is defined to be
\begin{eqnarray}\label{eq-Jfunction0}
&&\mathcal{J}(t^0,\dots, t^{n+m},z)
:=z+\sum_{a=0}^{n+m}t^{a}\gamma_{a}\nn\\
&&+\sum_{n\geq 0}\sum_{d}\sum_{a=0}^{n+m}\frac{1}{n!}\big\langle 
\sum_{b=0}^{n+m}t^{b}\gamma_{b},\dots,\sum_{b=0}^{n+m}t^{b}\gamma_{b}, \frac{\gamma_{a}}{z-\psi}\big\rangle_{0,n+1,d}\gamma^{a}.
\end{eqnarray}
For later convenience we define, for $ 0\leq a\leq n+m$,
\begin{eqnarray}\label{eq-Jfunction1}
\mathcal{J}_a(t^0,\dots,t^{n+m},z):=\sum_{c=0}^{n+m}g_{ac}t^{c}+\sum_{n\geq 0}\sum_{d}\sum_{a=0}^{n+m}\frac{1}{n!}\big\langle 
\sum_{b=0}^{n+m}t^{b}\gamma_{b},\dots,\sum_{b=0}^{n+m}t^{b}\gamma_{b}, \frac{\gamma_{a}}{z-\psi}\big\rangle_{0,n+1,d}.
\end{eqnarray}
Then the \emph{quantum differential equations} read, for $0\leq i,j\leq n+m$,
\begin{eqnarray}\label{eq-qde0}
z \frac{\partial^2 \mathcal{J}_a}{\partial t^i \partial t^j}
=\sum_{b=0}^{n+m}\sum_{c=0}^{n+m}\frac{\partial^3 F}{\partial t^i \partial t^j \partial t^b}g^{bc}\frac{\partial \mathcal{J}_a}{ \partial t^c}.
\end{eqnarray}
We define \emph{ambient big $J$-function}  to be
\begin{eqnarray}\label{eq-Jfunctionamb}
\mathcal{J}_{\mathrm{amb}}(t^0,\dots, t^{n},z)=z+\sum_{a=0}^{n}\gamma^{a}\mathcal{J}_a|_{t^{n+1}=\dots =t^{n+m}=0}.
\end{eqnarray}
By Theorem \ref{thm-monodromythm}, for $0\leq a\leq n$, $\mathcal{J}_a(t^0,\dots,t^{n+m},z)$ is a series in $t^0,\dots,t^{n},s$. 
And moreover, for $n+1\leq a\leq n+m$, there exists $\tilde{\mathcal{J}}(t^0,\dots,t^{n},s,z)$ such that
 \begin{eqnarray}\label{eq-Jfunction2}
 \mathcal{J}_a(t^0,\dots,t^{n+m},z)=t^a \tilde{\mathcal{J}}(t^0,\dots,t^{n},s,z).
 \end{eqnarray}
When the dimension $n$ is odd, the degree of $s$ in $\tilde{\mathcal{J}}(t^0,\dots,t^{n},s,z)$ is not greater than $\frac{m}{2}-1$. 
In the following we adopt Einstein's summation convention, where the indices run  from 
$0$ to $n$.
\begin{proposition}\label{eq-qdered}
Let $X$ be a non-exceptional complete intersection with dimension $\geq 3$. 
The system of quantum differential equations (\ref{eq-qde0}) is equivalent to the collection of the following systems. 
(i) If $n$ is even,
\begin{eqnarray}\label{eq-qde4}
z\frac{\partial \mathcal{J}_{a}}{\partial s}=\frac{\partial^2 F}{\partial s\partial t^b}g^{bc}\frac{\partial \mathcal{J}_{a}}{\partial t^c}+2s \frac{\partial^2F}{\partial s \partial s}\frac{\partial \mathcal{J}_{a}}{\partial s}, & 0\leq a\leq n,
\end{eqnarray}
\begin{eqnarray}\label{eq-qde7}
z\frac{\partial \tilde{\mathcal{J}}}{\partial t^{i}}=\frac{\partial^2 F}{\partial s\partial t^i}\tilde{\mathcal{J}}, & 0\leq i\leq n,
\end{eqnarray}
\begin{eqnarray}\label{eq-qde9}
z\frac{\partial \tilde{\mathcal{J}}}{\partial s}=\frac{\partial^2F}{\partial s \partial s}\tilde{\mathcal{J}}.
\end{eqnarray}
(ii) If $n$ is odd,
\begin{eqnarray}\label{eq-qde4odd}
z\frac{\partial \mathcal{J}_{a}}{\partial s}=\frac{\partial^2 F}{\partial s\partial t^b}g^{bc}\frac{\partial \mathcal{J}_{a}}{\partial t^c}+2s \frac{\partial^2F}{\partial s \partial s}\frac{\partial \mathcal{J}_{a}}{\partial s} \mod s^{\frac{m}{2}}, & 0\leq a\leq n,
\end{eqnarray}
\begin{eqnarray}\label{eq-qde7odd}
z\frac{\partial \tilde{\mathcal{J}}}{\partial t^{i}}=\frac{\partial^2 F}{\partial s\partial t^i}\tilde{\mathcal{J}} \mod s^{\frac{m}{2}}, & 0\leq i\leq n,
\end{eqnarray}
\begin{eqnarray}\label{eq-qde9odd}
z\frac{\partial \tilde{\mathcal{J}}}{\partial s}=\frac{\partial^2F}{\partial s \partial s}\tilde{\mathcal{J}} \mod s^{\frac{m}{2}-1}.
\end{eqnarray}
\end{proposition}
\begin{proof}: The proof is similar to the proof of Theorem \ref{thm-wdvveventhm} and Theorem \ref{thm-wdvvoddthm}, and we omit it.
\end{proof}

\begin{theorem}\label{thm-Jreconstruction1}
Let $X$ be a non-exceptional complete intersection with dimension $\geq 3$. 
\begin{itemize}
\item[(i)]  For $0\leq a\leq n$, $\mathcal{J}_{a}$ can be reconstructed by (\ref{eq-qde4}) (resp., (\ref{eq-qde4odd}) when $n$ is odd) from $\mathcal{J}_{\mathrm{amb}}$.
\item[(ii)] For $n+1\leq a\leq n+m$,
\begin{equation*}
\mathcal{J}_{a}(t^0,\dots,t^{n+m},z)=t^{a}\exp\left(\frac{1}{z}\frac{\partial F}{\partial s}\right).
\end{equation*}
\end{itemize}
\end{theorem}
\begin{proof}: By the dimension constraint (\ref{eq-Dim}) the genus zero Gromov-Witten invariants for non-Fano complete intersections are trivial. So we assume $X$ is Fano.\\
(i) For $0\leq a\leq n$, we expand $\mathcal{J}_a$ in $s$ as
\begin{eqnarray*}
\mathcal{J}_a=\mathcal{J}_a^{(0)}+s\mathcal{J}_a^{(1)}+\frac{s^2}{2}\mathcal{J}_a^{(2)}+\dots
\end{eqnarray*}
Then $\mathcal{J}_a^{(0)}$ is the coefficient of $\gamma^a$ in $\mathcal{J}_{\mathrm{amb}}$, and by (\ref{eq-qde4}),  for $k\geq 0$,
\begin{eqnarray*}
\mathcal{J}_a^{(k+1)}=\frac{1}{z}\sum_{i=0}^{k} \binom{k}{i} \frac{\partial F}{\partial t^b}^{(i+1)}g^{bc}\frac{\partial \mathcal{J}_{a}^{(k-i)}}{\partial t^c}+\frac{2k}{z}\sum_{i=0}^{k-1}\binom{k-1}{i}F^{(i+2)}\mathcal{J}_a^{(k-i)}.
\end{eqnarray*}
(ii) The general solution for  (\ref{eq-qde7}) (resp., (\ref{eq-qde7odd}) when $n$ is odd) and (\ref{eq-qde9}) (resp., (\ref{eq-qde9odd}) when $n$ is odd) is
\begin{eqnarray*}
\mathcal{J}(t^0,\dots,t^n, s,z)=C(1/z)\exp\left(\frac{1}{z}\frac{\partial F}{\partial s}\right),
\end{eqnarray*}
where 
\begin{eqnarray*}
C(1/z)=1+\frac{c_{1}}{z}+\frac{c_2}{z^2}+\frac{c_3}{z^3}+\dots
\end{eqnarray*}
for some $c_i\in \mathbb{C}$, $i\geq 1$. These coefficients can be determined by $F^{(1)}(0)$ and 
\begin{eqnarray*}
\sum_{a=n+1}^{n+m}\sum_{b=n+1}^{n+m}g^{ab}\langle \gamma_{a}\psi^{k}, \gamma_b \rangle_{0,2}, & k\geq 1.
\end{eqnarray*}
Since $\sfh\cup \gamma_{a}=0$ for $\gamma_{a}\in H_{\mathrm{prim}}(X)$, by (\ref{eq-Div}) we have
\begin{eqnarray*}
\langle \gamma_{a}\psi^{k}, \gamma_b, \sfh\rangle_{0,3}=\frac{k+1}{\mathsf{a}(n,\mathbf{d})}\langle \gamma_{a} \psi^{k}, \gamma_b\rangle_{0,2}.
\end{eqnarray*}
On the other hand (\ref{eq-TRR0}) implies
\begin{eqnarray}\label{eq-trr}
\langle \gamma_{a}\psi^{k}, \gamma_b, \sfh\rangle_{0,3}
=\sum_{\mu=0}^{n+m}\sum_{\nu=0}^{n+m}\langle \gamma_{a}\psi^{k-1}, \gamma_{\mu}\rangle_{0,2} g^{\mu\nu}
\langle \gamma_{\nu}, \gamma_b, \sfh\rangle_{0,3}.
\end{eqnarray}
If the Fano index $\mathsf{a}(n,\mathbf{d})>1$, the proof of Lemma \ref{lem-someGenus0Vanishing} shows that $\langle \gamma_{\nu}, \gamma_b, \sfh\rangle_{0,3}=0$, 
 so  we have
\begin{eqnarray*}
\langle \gamma_{a} \psi^{k}, \gamma_b\rangle_{0,2}=0,\ \mbox{for}\ k\geq 0.
\end{eqnarray*}
If the Fano index $\mathsf{a}(n,\mathbf{d})=1$, (\ref{eq-trr}) implies
\begin{eqnarray*}
&&(k+1)\langle \gamma_{a}\psi^{k}, \gamma_b\rangle_{0,2,k+1}\\
&=&\sum_{\mu=n+1}^{n+m}\sum_{\nu=n+1}^{n+m}\langle \gamma_{a}\psi^{k-1}, \gamma_{\mu}\rangle_{0,2,k} g^{\mu\nu}
\langle \gamma_{\nu}, \gamma_b\rangle_{0,2,1}\\
&=&\sum_{\mu=n+1}^{n+m}\sum_{\nu=n+1}^{n+m} \langle \gamma_{a}\psi^{k-1}, \gamma_{\mu}\rangle_{0,2,k} g^{\mu\nu}F^{(1)}(0)g_{\nu b}\\
&=&F^{(1)}(0)\langle \gamma_{a}\psi^{k-1}, \gamma_{b}\rangle_{0,2,k}.
\end{eqnarray*}
By Proposition \ref{prop-Theta-inSmallQuantumCohomology}, $F^{(1)}(0)=-\elld$.
So by induction on $k$ we get
\begin{eqnarray*}
\langle \gamma_{a}\psi^{k}, \gamma_b\rangle_{0,2,k+1}=\frac{(-\elld)^{k+1}}{(k+1)!}g_{ab}.
\end{eqnarray*}
Therefore in either case we obtain
\begin{eqnarray*}
C(1/z)=1.
\end{eqnarray*}
\end{proof}

\section{An identity of contractions of (permanent) Pfaffians}
 Let $\gamma_1,\dots,\gamma_{m}$ be  a basis of  $H^*_{\mathrm{prim}}(X)$. Let $G(\gamma_1,\dots,\gamma_m)$ to be the matrix $G=(g_{i,j})_{1\leq i,j\leq m}$ with $g_{i,j}=(\alpha_i,\alpha_j)$. Then when $\dim X=n$ is even $G$ is symmetric, and when $\dim X=n$ is odd $G$ is skew-symmetric. Let $\gamma^a=\sum_{e=1}^{m}g^{a,e}\gamma_e$. Recall the notations introduced in Section \ref{sec:meaningOfF(l)(0)}.

\begin{lemma}\label{lem-contraction-Pfaffian}
When $n$ is even,
\begin{equation}\label{eq-contraction-Pfaffian-even}
  \sum_{a=1}^{m}\mathrm{P}\big(G(\gamma_a,\gamma^a,\gamma_1,\dots,\gamma_{2l})\big)
  =(2l+m)\mathrm{P}\big(G(\gamma_1,\dots,\gamma_{2l})\big).
\end{equation}
When $n$ is odd,
\begin{equation}\label{eq-contraction-Pfaffian-odd}
  \sum_{a=1}^{m}\mathrm{Pf}\big(G(\gamma_a,\gamma^a,\gamma_1,\dots,\gamma_{2l})\big)
  =(2l-m)\mathrm{Pf}\big(G(\gamma_1,\dots,\gamma_{2l})\big)
\end{equation}
\end{lemma}
\begin{proof}
We show the odd dimensional case (\ref{eq-contraction-Pfaffian-odd}).
Recall the set $A_{2l}$ defined as (\ref{eq-def-A2l}).
We define maps $\phi_{k}:A_{2l}\rightarrow A_{2l+2},\ 0\leq k\leq 2l$ as follows:
\begin{equation*}
\phi_0\Big(\big((i_1,j_1),(i_2,j_2),\dots,(i_l,j_l)\big)\Big)= 
\big((1,2),(i_1+2,j_1+2),(i_2+2,j_2+2),\dots,(i_l+2,j_l+2)\big);
\end{equation*}
and for $1\leq k\leq l$,
\begin{multline*}
\phi_{k}\Big(\big((i_1,j_1),(i_2,j_2),\dots,(i_l,j_l)\big)\Big)\\=
\big((1,i_k+2),(2,j_k+2),(i_1+2,j_1+2),\dots,(i_k+2,j_k+2)^{\land},\dots,
(i_l+2,j_l+2)\big);
\end{multline*}
and for $l+1\leq k\leq 2l$,
\begin{multline*}
\phi_{k}\Big(\big((i_1,j_1),(i_2,j_2),\dots,(i_l,j_l)\big)\Big)\\=
\big((1,j_k+2),(2,i_k+2),(i_1+2,j_1+2),\dots,(i_k+2,j_k+2)^{\land},\dots,
(i_l+2,j_l+2)\big),
\end{multline*}
where the superscript $\land$ means to delete this pair. Then $A_{2l+2}$ is the disjoint union of the images of $\phi_k$ for $0\leq k\leq 2l$. By (\ref{eq-pfaffian-def-odd}) we have
\begin{eqnarray}\label{eq-lem-contraction-Pfaffian-1}
&&\sum_{a=1}^{m}\sum_{\sigma\in \mathrm{Im}(\phi_0)}\mathrm{Pf}_{\sigma}\big(G(\gamma_a,\gamma^a,\gamma_1,\dots,\gamma_{2l})\big)\nonumber\\
&=&\sum_{a=1}^m \langle \gamma_a, \gamma^a\rangle 
\sum_{\sigma\in A_{2l}}\mathrm{Pf}_{\sigma}\big(G(\gamma_1,\dots,\gamma_{2l})\big)\nonumber\\
&=& -m \mathrm{Pf}\big(G(\gamma_1,\dots,\gamma_{2l})\big).
\end{eqnarray}
For $1\leq k\leq l$, from 
\begin{equation}\label{eq-lem-contraction-Pfaffian-1.1}
\sum_{a=1}^m \langle \gamma_a, \gamma_{i}\rangle \langle \gamma^a, \gamma_{j}\rangle=
-\langle \gamma_i,\gamma_j\rangle
\end{equation}
and 
\begin{equation}\label{eq-lem-contraction-Pfaffian-1.2}
\mathrm{sgn}(\phi_k(\sigma))=-\mathrm{sgn}(\sigma)
\end{equation}
it follows that
\begin{equation}\label{eq-lem-contraction-Pfaffian-2}
\sum_{a=1}^{m}\sum_{\sigma\in \mathrm{Im}(\phi_k)}\mathrm{Pf}_{\sigma}\big(G(\gamma_a,\gamma^a,\gamma_1,\dots,\gamma_{2l})\big)=
 \mathrm{Pf}\big(G(\gamma_1,\dots,\gamma_{2l})\big).
\end{equation}
Similarly for $l+1\leq k\leq 2l$ we have
\begin{equation}\label{eq-lem-contraction-Pfaffian-3}
\sum_{a=1}^{m}\sum_{\sigma\in \mathrm{Im}(\phi_k)}\mathrm{Pf}_{\sigma}\big(G(\gamma_a,\gamma^a,\gamma_1,\dots,\gamma_{2l})\big)=
 \mathrm{Pf}\big(G(\gamma_1,\dots,\gamma_{2l})\big).
\end{equation}
Combining (\ref{eq-lem-contraction-Pfaffian-1}), (\ref{eq-lem-contraction-Pfaffian-2}) and (\ref{eq-lem-contraction-Pfaffian-3}) we get (\ref{eq-contraction-Pfaffian-odd}). 

A similar proof as above, neglecting the signs in (\ref{eq-lem-contraction-Pfaffian-1})-(\ref{eq-lem-contraction-Pfaffian-1.2}), shows the  even dimensional case (\ref{eq-contraction-Pfaffian-even}).
\end{proof}

\newpage

\section{Algorithm}\label{sec:algorithm}
Let $X$ be a non-exceptional smooth complete intersection of dimension $n$ and multidegree $d=(d_1,\dots,d_r)$. Let  $t^0,\dots,t^{n+m}$ be the basis dual to  the basis $1,\sfh,\dots,\sfh_n$ of $H^*_{\mathrm{amb}}(X)$. 
Let  $\tau^0,\dots,\tau^{n+m}$ be the basis dual to  the basis $1,\tsfh,\dots,\tsfh_n$ of $H^*_{\mathrm{amb}}(X)$.

In this appendix we describe an algorithm to the compute the $F^{(0)}$ based on the mirror formula, and to compute $F^{(k)}$ with $F^{(k)}(0)$ as  an unknown based on Theorem \ref{thm-reconstruction-I}. Intermediate formulae are computed explicitly so that they match the corresponding functions in our package  \texttt{QuantumCohomologyFanoCompleteIntersection} in

\url{https://github.com/huxw06/Quantum-cohomology-of-Fano-complete-intersections}

Examples of computations of ambients invariants  via this package
have been verified by comparison with the Julia package

\url{https://github.com/mgemath/AtiyahBott.jl}

\noindent pertinent to \cite{MS21}, which is based on a direct Atiyah-Bott localization computation and summations over graphs.

\subsection{Summary}\label{sec:algorithm-summary}
The main novelty of our algorithm is the systematic use of the $\tau$-coordinates. For this, we need first compute the transition matrices $W$ and $M$ between the $t$-coordinates and $\tau$-coordinates. So our algorithm consists of the following steps.
\begin{enumerate}
      \item We extract correlators of length 1 and length 2 from the mirror formula of the small $J$-function. Some combinatorial and numerical tricks are involved.
      \item Compute the  matrices $W$ and $M$. This is necessary for the expression of the Euler field.
      \item Use the WDVV recursion for $F^{(0)}$. The advantage of the use of $\tau$-coordinates will show up in this step.
\end{enumerate}

\subsection{Small and big \texorpdfstring{$J$}{J}-functions}
\label{sec:computation-twoPointInvariants-from-smallJ}
Define the small $J$-function by
\begin{equation*}
      J(t^1,z):=\mathcal{J}(t^0,\dots,t^{n+m},z)|_{t^i=0\ \mathrm{for}\ i\neq 1}.
\end{equation*}
When the Fano index $\mathsf{a}(n,\mathbf{d})\geq 1$, the mirror formula  \cite{Giv96} is
\begin{equation}\label{eq-mirrorFormulaOfSmallJfunction}
      J(t^1,z)=
      \begin{cases}
      ze^{\frac{t^1 \sfh}{z}}\sum_{\beta\geq 0}\frac{\prod_{i=1}^r \prod_{k=1}^{d_i \beta}(d_i \sfh+kz)}{\prod_{j=1}^{\beta}(\sfh+jz)^{n+r+1}} \mathsf{q}^{\beta}e^{\beta t^1},& \mbox{when}\ \mathsf{a}(n,\mathbf{d})\geq 2,\\
       ze^{\frac{t^1 \sfh}{z}-\elld\frac{e^{t^1}}{z}}\sum_{\beta\geq 0}\frac{\prod_{i=1}^r \prod_{k=1}^{d_i \beta}(d_i \sfh+kz)}{\prod_{j=1}^{\beta}(\sfh+jz)^{n+r+1}} \mathsf{q}^{\beta}e^{\beta t^1},& \mbox{when}\ \mathsf{a}(n,\mathbf{d})=1.      
      \end{cases}
\end{equation}
Recall the big $J$-function (\ref{eq-Jfunction0})\footnote{One can also use the ambient big $J$-function $\mathcal{J}_{\mathrm{amb}}$ in this appendix.}.
By (\ref{eq-TRR0}),
\begin{equation*}
	z^2\frac{\partial^2}{\partial t^i \partial t^j} \mathcal{J}(t,z)=\mathcal{C}_{i,j}^k
	z\frac{\partial}{\partial t^k} \mathcal{J}(t,z),
\end{equation*}
where $\mathcal{C}_{i,j}^k$ is the structure coefficient of the big quantum product. By induction on $l$, there exist integers $\nu_{i_1,\dots,i_p}$ such that
\begin{eqnarray}\label{eq-derivative-J}
&&z^l\big(\frac{\partial}{\partial t^1}\big)^l \mathcal{J}\nn\\
&=& \sum_{\begin{subarray}{c}i_1+\dots+i_p+p=l-1\\
p\geq 1, i_1,\dots,i_p\geq 0\end{subarray}} \Big(\nu_{i_1,\dots,i_p}
\sum_{\begin{subarray}{c}0\leq a_j\leq n\\ \mbox{for } 1\leq j\leq p\end{subarray}}
(z^{i_1}\frac{\partial^{i_1}}{(\partial t^1)^{i_1}}\mathcal{C}_{1,1}^{a_1})
(z^{i_2}\frac{\partial^{i_2}}{(\partial t^1)^{i_2}}\mathcal{C}_{1,a_1}^{a_2})\nn\\
&&\cdots
(z^{i_p}\frac{\partial^{i_p}}{(\partial t^1)^{i_p}}\mathcal{C}_{1,a_{p-1}}^{a_p})z\frac{\partial \mathcal{J}}{\partial t^{a_p}}\Big).
\end{eqnarray}
These coefficients are determined by  the recursion
\begin{equation}\label{eq-recursion-coefficient-derivativeOfJ}
\nu_{i_1,\dots,i_p}=\sum_{j=1}^{p}\nu_{i_1,\dots,i_{j-1},i_j-1,i_{j+1},\dots,i_{p}}+\delta_{i_p,0}\nu_{i_1,\dots,i_{p-1}},
\end{equation}
and the initial value
\begin{equation*}
	\nu_0=1.
\end{equation*}

\begin{lemma}\label{lem-coefficient-derivativeOfJ}
\begin{eqnarray}\label{eq-coefficient-derivativeOfJ}
&& \nu_{i_1,\dots,i_p}=\binom{i_1+\dots+i_p+p-1}{i_1}\binom{i_2+\dots+i_p+p-2}{i_2}\dots\binom{i_{p-1}+i_p+1}{i_{p-1}}.
\end{eqnarray}
\end{lemma}
\begin{proof} 
It suffices to verify that $\nu_{i_1,\dots,i_p}$'s defined as (\ref{eq-coefficient-derivativeOfJ}) satisfy (\ref{eq-recursion-coefficient-derivativeOfJ}). We use repeadedly the identity
\[
\binom{a+b}{b}-\binom{a+b-1}{b-1}=\binom{a+b-1}{b}.
\]
Then
\begin{eqnarray*}
\nu_{i_1,i_2,\dots,i_p}-\nu_{i_1-1,i_2,\dots,i_p}
&=&\binom{i_1+(i_2-1)+i_3+\dots+i_p+p-1}{i_1}\\
&&\cdot\binom{i_2+\dots+i_p+p-2}{i_2}\dots\binom{i_{p-1}+i_p+1}{i_{p-1}},
\end{eqnarray*}
\begin{multline*}
\nu_{i_1,i_2,\dots,i_p}-\nu_{i_1-1,i_2,\dots,i_p}-\nu_{i_1,i_2-1,i_3,\dots,i_p}
=\binom{i_1+i_2+(i_3-1)+\dots+i_p+p-1}{i_1}\\
\binom{i_2+(i_3-1)+i_p+p-2}{i_2}
\binom{i_1+\dots+i_p}{i_3}
\dots\binom{i_{p-1}+i_p+1}{i_{p-1}},
\end{multline*}
and so on. 
\end{proof}

\subsection{Correlators of length 1}
In this section we use the mirror formula (\ref{eq-mirrorFormulaOfSmallJfunction}) to compute the length 1  genus 0 invariants.
We define a function
\begin{equation}\label{eq-def-functionMu}
      \mu(d,k)=
      \begin{cases}
      1, & \mbox{if}\ k=0;\\
      \sum_{1\leq i_1<\dots<i_k\leq d}\frac{1}{i_1\cdots i_k},& \mbox{if } k\geq 1.
      \end{cases}
\end{equation}
Then
\begin{equation}\label{eq-identity-muFunction-1}
      \prod_{k=1}^{d}(dx+k)=d!\sum_{j=0}^d d^j \mu(d,j)x^j.
\end{equation}
This function will be used repeatedly in this and the next appendix.
By (\ref{eq-Dim}),
\begin{eqnarray*}
&&\langle \sfh_i\rangle_{0}
=\begin{cases}
\langle \sfh_i\rangle_{0,1},& \mbox{if}\ i=n-2+\mathsf{a}(n,\mathbf{d}),\\
\langle \sfh_i\rangle_{0,2},& \mbox{if}\ i=n-2+2\mathsf{a}(n,\mathbf{d})
\end{cases}\\
&=& \begin{cases}
 \delta_{i,n-1}\langle \sfh_{n-1}\rangle_{0,1},& \mbox{if}\ \mathsf{a}(n,\mathbf{d})=1,\\
\delta_{i,n}\langle \sfh_{n}\rangle_{0,1} ,& \mbox{if}\ \mathsf{a}(n,\mathbf{d})=2,\\
\delta_{i,n} \langle \sfh_n\rangle_{0,2} ,& \mbox{if}\ \mathsf{a}(n,\mathbf{d})=1.
\end{cases}
\end{eqnarray*}
When $\mathsf{a}(n,\mathbf{d})=2$,
\begin{eqnarray*}
&& \langle \sfh_{n}\rangle_{0,1}=\prod_{i=1}^r d_i\cdot
\mathrm{Coeff}_{\frac{1}{z}}\Big(ze^{t^1 \sfh/z}\sum_{\beta\geq 0}\frac{\prod_{i=1}^r \prod_{k=1}^{d_i \beta}(d_i \sfh+kz)}{\prod_{j=1}^{\beta}(\sfh+jz)^{n+r+1}} \mathsf{q}^{\beta}e^{\beta t^1}\Big)\\
&=& \prod_{i=1}^r d_i\cdot [\frac{1}{z}]\Big(ze^{t^1 \sfh/z}\frac{\prod_{i=1}^r \prod_{k=1}^{d_i }(d_i \sfh+kz)}{(\sfh+z)^{n+r+1}}\Big)
= \prod_{i=1}^r d_i\cdot \elld.
\end{eqnarray*}

When $\mathsf{a}(n,\mathbf{d})=1$,
\begin{eqnarray*}
&& \langle \sfh_{n-1}\rangle_{0,1}=\prod_{i=1}^r d_i\cdot\sum_{\xi\geq 0} [\frac{\sfh}{z}] \big(ze^{t^1(\sfh/z+\xi)-\elld e^{t^1}/z}\frac{\prod_{i=1}^r \prod_{k=1}^{d_i \xi}(d_i \sfh+k z)}{\prod_{j=1}^{\xi}(\sfh+j z)^{n+r+1}}\big)\\
&=& \prod_{i=1}^r d_i\cdot \mathrm{Coeff}_{\frac{\sfh}{z}} \big(e^{t^1(\sfh/z+1)-\elld e^{t^1}/z}\frac{\prod_{i=1}^r \prod_{k=1}^{d_i }(d_i \sfh/z+k)}{(\sfh/z+ 1)^{n+r+1}}\big)\\
&=& \prod_{i=1}^r d_i\cdot \mathrm{Coeff}_{x} \big(\frac{\prod_{i=1}^r \prod_{k=1}^{d_i }(d_i x+k)}{(x+ 1)^{n+r+1}}\big)\\
&=& \prod_{i=1}^r d_i\cdot\elld\big(\sum_{i=1}^r d_i \mu(d_i,1)-n-r-1\big),
\end{eqnarray*}
and
\begin{eqnarray*}
&& \langle \sfh_{n}\rangle_{0,2}=\prod_{i=1}^r d_i\cdot\sum_{\xi\geq 0} \mathrm{Coeff}_{\frac{1}{z}} \big(ze^{t^1(\sfh/z+\xi)-\elld e^{t^1}/z}\frac{\prod_{i=1}^r \prod_{k=1}^{d_i \xi}(d_i \sfh+k z)}{\prod_{j=1}^{\xi}(\sfh+j z)^{n+r+1}}\big)\\
&=& \prod_{i=1}^r d_i\cdot\sum_{\xi\geq 0} \mathrm{Coeff}_{\frac{1}{z^2}} \big(e^{-\elld/z}\frac{\prod_{i=1}^r \prod_{k=1}^{d_i \xi}k z}{\prod_{j=1}^{\xi}(j z)^{n+r+1}}\big)\\
&=& \prod_{i=1}^r d_i\cdot\elld^2\big(-\frac{1}{2}+2^{-n-r-1}\prod_{i=1}^r\binom{2d_i}{d_i}\big).
\end{eqnarray*}

\subsection{Recursion for genus 0 GW invariants of length two  with a \texorpdfstring{$\psi$-class}{psi}}
 We are going to compute
\begin{equation*}
      \langle \sfh_a,\sfh_b\psi^c\rangle_{0,\beta}
\end{equation*}
recursively (on the degree $\beta$) by (\ref{eq-derivative-J}). Define
\begin{equation}\label{eq-def-functionA}
\mathsf{A}_N(x_1,\dots,x_p)=\sum_{\begin{subarray}{c}i_1,\dots,i_p\geq 0\\
i_1+\dots+i_p=N\end{subarray}}
\nu_{i_1,\dots,i_p}x_1^{i_1}\cdots x_p^{i_p}.
\end{equation}
In particular,
\begin{equation*}
      \mathsf{A}_0(x_1,\dots,x_p)\equiv 1.
\end{equation*}

\begin{lemma}
Suppose $a,b,c\in \mathbb{Z}_{\geq 0}$, $\beta\in \mathbb{Z}_{\geq 1}$, and
\[
a+b+c=n-1+\beta\cdot \mathsf{a}(n,\mathbf{d}).
\]
Then
\begin{eqnarray}\label{eq-recusion-twoPointInv-FanoIndexOne}
&&\langle \sfh_a,\sfh_b\psi^c\rangle_{0,\beta}\nn\\
&=& \prod_{i=1}^r d_i\cdot [\frac{\sfh_{n-b}}{z^{c}}]\big((z\partial/\partial t^1)^a J(t^1,z)\big)\nn\\
&&-\sum_{\begin{subarray}{c}
      1\leq p\leq a-1\\
      \alpha_0=1,\ 
      0\leq \alpha_j\leq 1+\alpha_{j-1},\ \frac{1+\alpha_{p-1}- \alpha_p}{\mathsf{a}(n,\mathbf{d})}\in \mathbb{Z},\ \mathrm{for}\  1\leq j\leq p\\
      (p,\alpha_1,\dots,\alpha_p)\neq (a-1,2,3,\dots,a)
      \end{subarray}}
      \Big( \mathsf{A}_{a-1-p}\big(\frac{2- \alpha_1}{\mathsf{a}(n,\mathbf{d})},\dots,\frac{1+\alpha_{p-1}- \alpha_p}{\mathsf{a}(n,\mathbf{d})}\big)\nn\\
&&\cdot     \prod_{j=1}^{p}C_{1,\alpha_{j-1}}^{\alpha_j}
      \cdot \langle \sfh_{\alpha_p},\sfh_{b}\psi^{a+c-p-1}\rangle_{0,\beta-\frac{p+1- \alpha_p}{\mathsf{a}(n,\mathbf{d})}}\Big),   
\end{eqnarray}
where $C_{i,j}^k$ is the structure coefficient of the small quantum product $\sfh_i\sqp \sfh_j=C_{i,j}^k \sfh_k$. More precisely,
\[
C_{1,\alpha_{j-1}}^{\alpha_j}=
\begin{cases}
1, &  \mbox{if}\ \alpha_j=1+\alpha_{j-1},\\
\frac{1}{\prod_{i=1}^r d_i}\cdot \frac{1+\alpha_{j-1}-\alpha_{j}}{\mathsf{a}(n,\mathbf{d})}
\langle \sfh_{\alpha_{j-1}},\sfh_{n- \alpha_{j}}\rangle_{0,\frac{1+\alpha_{j-1}-\alpha_{j}}{\mathsf{a}(n,\mathbf{d})}}
, & \mbox{if}\ \alpha_j-1- \alpha_{j-1}>0,
\end{cases}
\]
\end{lemma}
\begin{proof}
By (\ref{eq-Jfunction0}) and (\ref{eq-derivative-J}),
\begin{eqnarray*}
&&\langle \sfh_a,\sfh_b\psi^c\rangle_{0,\beta}\\
&=& \prod_{i=1}^r{d_i}\cdot \mathrm{Coeff}_{\frac{\sfh_{n-b}}{z^{c}}}\bigg((z\frac{\partial}{\partial t^1})^a \mathcal{J}\\
&& -  \sum_{\begin{subarray}{c}i_1+\dots+i_p+p=a-1,\ p\geq 1, i_1,\dots,i_p\geq 0\\
0\leq \alpha_j\leq n\ \mbox{for } 1\leq j\leq p\\
(p,\alpha_1,\dots,\alpha_p)\neq (a-1,2,3,\dots,a)
\end{subarray}}
 \nu_{i_1,\dots,i_p}
(z^{i_1}\frac{\partial^{i_1}}{(\partial t^1)^{i_1}}\mathcal{C}_{1,1}^{\alpha_1})
(z^{i_2}\frac{\partial^{i_2}}{(\partial t^1)^{i_2}}\mathcal{C}_{1,\alpha_1}^{\alpha_2})\nn\\
&&\cdots
(z^{i_p}\frac{\partial^{i_p}}{(\partial t^1)^{i_p}}\mathcal{C}_{1,\alpha_{p-1}}^{\alpha_p})z\frac{\partial \mathcal{J}}{\partial t^{\alpha_p}}\bigg).
\end{eqnarray*}
We have
\begin{equation*}
      \mathrm{Coeff}_{\frac{\sfh_{n-b}}{z^{c}}}\big((z\frac{\partial}{\partial t^1})^a \mathcal{J} \big)
      =\mathrm{Coeff}_{\frac{\sfh_{n-b}}{z^{c}}}\big((z\frac{\partial}{\partial t^1})^a J(t^1,z) \big).
\end{equation*}
By (\ref{eq-Div}), 
\begin{equation*}
      \big((\frac{\partial}{\partial t^1})^i\mathcal{C}_{1,j}^k\big)|_{\mathrm{all}\ t^p=0}
      =(\frac{1+j-k}{\mathsf{a}(n,\mathbf{d})})^i C_{1,j}^k,
\end{equation*}
where we adopt the convention 
\[
0^0=1.
\]
It follows that 
\begin{eqnarray*}
&& \mathrm{Coeff}_{\frac{\sfh_{n-b}}{z^{c}}}\bigg(\sum_{\begin{subarray}{c}i_1+\dots+i_p+p=a-1,\ p\geq 1, i_1,\dots,i_p\geq 0\\
0\leq \alpha_j\leq n\ \mbox{for } 1\leq j\leq p\\
(p,\alpha_1,\dots,\alpha_p)\neq (a-1,2,3,\dots,a)
\end{subarray}}
 \nu_{i_1,\dots,i_p}
(z^{i_1}\frac{\partial^{i_1}}{(\partial t^1)^{i_1}}\mathcal{C}_{1,1}^{\alpha_1})
(z^{i_2}\frac{\partial^{i_2}}{(\partial t^1)^{i_2}}\mathcal{C}_{1,\alpha_1}^{\alpha_2})\nn\\
&&\cdots
(z^{i_p}\frac{\partial^{i_p}}{(\partial t^1)^{i_p}}\mathcal{C}_{1,\alpha_{p-1}}^{\alpha_p})z\frac{\partial \mathcal{J}}{\partial t^{\alpha_p}}\bigg)\\
&=&\sum_{\begin{subarray}{c}i_1+\dots+i_p+p=a-1,\ p\geq 1, i_1,\dots,i_p\geq 0\\
\alpha_0=1,\ 
0\leq \alpha_j\leq 1+\alpha_{j-1},\ \mbox{for } 1\leq j\leq p\\
(p,\alpha_1,\dots,\alpha_p)\neq (a-1,2,3,\dots,a)
\end{subarray}}
\Big(\nu_{i_1,\dots,i_p} \\
&&\cdot \big(\prod_{j=1}^{p}(\frac{1+\alpha_{j-1}-\alpha_{j}}{\mathsf{a}(n,\mathbf{d})})^{i_j}C_{1,\alpha_{j-1}}^{\alpha_j}\big)
\cdot \langle \sfh_{\alpha_p},\sfh_{b}\psi^{a+c-p-1}\rangle_{0,\beta-\frac{p+1- \alpha_p}{\mathsf{a}(n,\mathbf{d})}}\Big)\\
&=&\sum_{\begin{subarray}{c}
      1\leq p\leq a-1\\
      \alpha_0=1,\ 
      0\leq \alpha_j\leq 1+\alpha_{j-1},\ \frac{1+\alpha_{p-1}- \alpha_p}{\mathsf{a}(n,\mathbf{d})}\in \mathbb{Z},\ \mathrm{for}\  1\leq j\leq p\\
      (p,\alpha_1,\dots,\alpha_p)\neq (a-1,2,3,\dots,a)
      \end{subarray}}
      \Big( \mathsf{A}_{a-1-p}\big(\frac{2- \alpha_1}{\mathsf{a}(n,\mathbf{d})},\dots,\frac{1+\alpha_{p-1}- \alpha_p}{\mathsf{a}(n,\mathbf{d})}\big)\nn\\
&&\cdot     \prod_{j=1}^{p}C_{1,\alpha_{j-1}}^{\alpha_j}
      \cdot \langle \sfh_{\alpha_p},\sfh_{b}\psi^{a+c-p-1}\rangle_{0,\beta-\frac{p+1- \alpha_p}{\mathsf{a}(n,\mathbf{d})}}\Big).       
\end{eqnarray*}

\end{proof}

\begin{proposition}\label{prop-recusion-twoPointInv}
Suppose the Fano index $\mathsf{a}(n,\mathbf{d})\geq 2$. Then
\begin{eqnarray}\label{eq-recusion-twoPointInv}
&&\langle \sfh_a,\sfh_b\psi^c\rangle_{0,\beta}\nn\\
&=& \prod_{i=1}^r\big(d_i\cdot (d_i \beta)!\big)\cdot \frac{\beta^{a}}{(\beta!)^{n+r+1}}
      \sum_{\begin{subarray}{c}j_1+\dots+j_{\beta}+k_1+\dots+k_r=n- b\\
      j_1,\dots,j_{\beta},k_1,\dots,k_r\geq 0\end{subarray}}  
      \bigg(\binom{a-n-r-1}{j_\beta}\beta^{-j_\beta}\nn\\
&&    \cdot       \prod_{l=1}^{\beta-1}\binom{-n-r-1}{j_l}l^{-j_l}\prod_{i=1}^r  d_i^{k_i}\mu(d_i \beta,k_i)\bigg)\nn\\
&&-\sum_{\begin{subarray}{c}
      1\leq p\leq a-1\\
      \alpha_0=1,\ 
      0\leq \alpha_j\leq 1+\alpha_{j-1},\ \frac{1+\alpha_{p-1}- \alpha_p}{\mathsf{a}(n,\mathbf{d})}\in \mathbb{Z},\ \mathrm{for}\  1\leq j\leq p\\
      (p,\alpha_1,\dots,\alpha_p)\neq (a-1,2,3,\dots,a)
      \end{subarray}}
      \Big( \mathsf{A}_{a-1-p}\big(\frac{2- \alpha_1}{\mathsf{a}(n,\mathbf{d})},\dots,\frac{1+\alpha_{p-1}- \alpha_p}{\mathsf{a}(n,\mathbf{d})}\big)\nn\\
&&\cdot     \prod_{j=1}^{p}C_{1,\alpha_{j-1}}^{\alpha_j}
      \cdot \langle \sfh_{\alpha_p},\sfh_{b}\psi^{a+c-p-1}\rangle_{0,\beta-\frac{p+1- \alpha_p}{\mathsf{a}(n,\mathbf{d})}}\Big).    
\end{eqnarray}
In particular
\begin{equation}\label{eq-degOne-invariant}
\langle \sfh_a, \sfh_b\rangle_{0,1}
= \prod_{i=1}^r\big(d_i\cdot (d_i)!\big)\cdot
      \sum_{\begin{subarray}{c}j+k_1+\dots+k_r=n-b\\
      j,k_1,\dots,k_r\geq 0\end{subarray}}  
      \bigg(\binom{a-n-r-1}{j}\prod_{i=1}^r  d_i^{k_i}\mu(d_            i,k_i)\bigg).
\end{equation}
\end{proposition}
\begin{proof}
It suffices to evaluate the first term on RHS of (\ref{eq-recusion-twoPointInv-FanoIndexOne}). By (\ref{eq-mirrorFormulaOfSmallJfunction}),
\begin{eqnarray*}
&&\mathrm{Coeff}_{\frac{\sfh_{n-b}}{z^{c}}}\big((z\frac{\partial}{\partial t^1})^a J(t^1,z) \big)\\
&=&  \mathrm{Coeff}_{x^{n-b}}\Big((\partial_{t^1})^a \big(e^{t^1(x+\beta)}\frac{\prod_{i=1}^r \prod_{k=1}^{d_i \beta}(d_i x+k)}{\prod_{j=1}^{\beta}(x+j)^{n+r+1}}\big)
\Big)\\
&=&  \mathrm{Coeff}_{x^{n-b}}\Big(\frac{\prod_{i=1}^r \prod_{k=1}^{d_i \beta}(d_i x+k)}{(x+\beta)^{n+r+1-a}\prod_{j=1}^{\beta-1}(x+j)^{n+r+1}}\big)
\Big)\\
&=& 
      \sum_{\begin{subarray}{c}j_1+\dots+j_{\beta}+k_1+\dots+k_r=n- b\\
      j_1,\dots,j_{\beta},k_1,\dots,k_r\geq 0\end{subarray}}  
      \Big(\beta^{a-n-r-1}((\beta-1)!)^{-n-r-1}\binom{a-n-r-1}{j_\beta}\beta^{-j_\beta}\\
&&\cdot     \prod_{l=1}^{\beta-1}\binom{-n-r-1}{j_l}l^{-j_l}\prod_{i=1}^r ((d_i \beta)! d_i^{k_i}\mu(d_i \beta,k_i))\Big).
\end{eqnarray*}
\end{proof}

\subsubsection{The case Fano index =1}
In the Fano index =1 case, the computation of the first term on RHS of (\ref{eq-recusion-twoPointInv-FanoIndexOne}) is slightly more complicated. 
\begin{lemma}
Let $f_{i,j}(x)$ be polynomials of $x$ defined by
\begin{equation}\label{eq-def-function-fij}
\sum_{i=0}^{\infty}f_{i,j}t^i=\frac{t^j}{\prod_{k=0}^{j}\big(1-(x+k)t\big)}.
\end{equation}
Let 
\[
\sigma_{i,j}(x,y)=f_{i,j}(x)y^j.
\]
Then
\begin{equation}\label{eq-diff-exp}
      (\partial_{t})^i(e^{xt+ye^t})
      =\sum_{j=0}^i \sigma_{i,j}(x,y)e^{(x+j)t+ye^t}.
\end{equation}
\end{lemma}
\begin{proof}
Let $A$ and $B$ be quantities independent of $t$. Then for $i\geq 0$ there exists polynomials $\sigma_{i,j}(A,B)$ of $A$, $B$ for $0\leq j\leq i$ such that
\begin{equation*}
      (\partial_{t})^i(e^{At+Be^t})
      =\sum_{j=0}^i \sigma_{i,j}(A,B)e^{(A+j)t+Be^t}.
\end{equation*}
Since
\begin{eqnarray*}
&& \partial_t\big(\sum_{j=0}^i \sigma_{i,j}(A,B)e^{(A+j)t+Be^t}\big)\\
&=&\sum_{j=0}^i \sigma_{i,j}(A,B)(A+j+Be^t)e^{(A+j)t+Be^t},
\end{eqnarray*}
The polynomials $\sigma_{i,j}(x,y)$ are determined by $\sigma_{0,0}=1$
and 
\begin{equation*}
      \sigma_{i,j}=(x+j)\sigma_{i-1,j}+y \sigma_{i-1,j-1}.
\end{equation*}
By induction on $j$ one sees that there exists polynomials $f_{i,j}(x)$ such that 
\[
\sigma_{i,j}=f_{i,j}(x)y^j.
\]
The polynomials $f_{i,j}$ are determined by $f_{0,0}=1$
and
\[
f_{i,j}=(x+j)f_{i-1,j}+f_{i-1,j-1}.
\]
Solving the generating functions $\sum_{i=0}^{\infty}f_{i,j}t^i$  inductively yields (\ref{eq-def-function-fij}).
\end{proof}

Suppose  $a,b,c\in \mathbb{Z}_{\geq 0}$, $\beta\in \mathbb{Z}_{\geq 1}$, and
\begin{equation*}
      a+b+c=n-1+\beta.
\end{equation*}

Then by (\ref{eq-mirrorFormulaOfSmallJfunction}) and (\ref{eq-diff-exp}),
\begin{eqnarray*}
&&\mathrm{Coeff}_{\frac{\sfh_{n-b}}{z^{c}}}\big((z\frac{\partial}{\partial t^1})^a J(t^1,z) \big)\\
&=& \sum_{\xi\geq 0} \mathrm{Coeff}_{\frac{\sfh_{n-b}}{z^{c}}}\Big(z^a(\partial_{t^1})^a \big(ze^{t^1(\frac{\sfh}{z}+\xi)-\elld\frac{e^{t^1}}{z}}\frac{\prod_{i=1}^r \prod_{k=1}^{d_i \xi}(d_i \sfh+k z)}{\prod_{j=1}^{\xi}(\sfh+j z)^{n+r+1}}\big)\Big)\\
&=& \sum_{\xi\geq 0} \mathrm{Coeff}_{\frac{\sfh_{n-b}}{z^{c}}}\Big(
      z^{a+1}\sum_{l=0}^a f_{a,l}(\frac{\sfh}{z}+\xi)(-\frac{\elld}{z})^{l} e^{(\frac{\sfh}{z}+\xi+l)t^1-\frac{\elld}{z}e^{t^1}}
\cdot\frac{\prod_{i=1}^r \prod_{k=1}^{d_i \xi}(d_i \sfh+k z)}{\prod_{j=1}^{\xi}(\sfh+j z)^{n+r+1}}\Big)\\
&=& \sum_{\xi\geq 0} \mathrm{Coeff}_{\frac{\sfh_{n-b}}{z^{c}}}\Big(
      z^{a+1}\sum_{l=0}^a f_{a,l}(\frac{\sfh}{z}+\xi)(-\frac{\elld}{z})^{l} e^{-\frac{\elld}{z}}
\cdot\frac{\prod_{i=1}^r \prod_{k=1}^{d_i \xi}(d_i \sfh+k z)}{\prod_{j=1}^{\xi}(\sfh+j z)^{n+r+1}}\Big)\\
&=& \sum_{\xi\geq 0} \mathrm{Coeff}_{\frac{1}{z^{b+c-n}}\frac{\sfh_{n-b}}{z^{n-b}}}\Big(
      z^{a+1}\sum_{l=0}^a f_{a,l}(\frac{\sfh}{z}+\xi)(-\frac{\elld}{z})^{l} e^{-\frac{\elld}{z}}
\cdot\frac{\prod_{i=1}^r \prod_{k=1}^{d_i \xi}(d_i \sfh/z+k )}{z^{\xi}\prod_{j=1}^{\xi}(\sfh/z+j )^{n+r+1}}\Big)\\
&=& \sum_{\xi\geq 0} \mathrm{Coeff}_{z^{-\beta}x^{n-b}}\Big(
      \sum_{l=0}^a f_{a,l}(x+\xi)(-\frac{\elld}{z})^{l} e^{-\frac{\elld}{z}}
\cdot\frac{\prod_{i=1}^r \prod_{k=1}^{d_i \xi}(d_i x+k )}{z^{\xi}\prod_{j=1}^{\xi}(x+j )^{n+r+1}}\Big)\\
&=& \sum_{\xi\geq 0} \mathrm{Coeff}_{z^{-\beta}x^{n-b}}\Big(
      \sum_{l=0}^a f_{a,l}(x+\xi)(-\frac{\elld}{z})^{l}  
            \frac{1}{(\beta-l-\xi)!}      (-\frac{\elld}{z})^{\beta-l-\xi}
\cdot\frac{\prod_{i=1}^r \prod_{k=1}^{d_i \xi}(d_i x+k )}{z^{\xi}\prod_{j=1}^{\xi}(x+j )^{n+r+1}}\Big)\\
&=& \sum_{\xi=0}^{\beta} \mathrm{Coeff}_{x^{n-b}}\Big(
      \sum_{l=0}^{\beta- \xi} f_{a,l}(x+\xi)
            \frac{\big(-\elld\big)^{\beta-\xi}}{(\beta-l-\xi)!}    
\cdot\frac{\prod_{i=1}^r \prod_{k=1}^{d_i \xi}(d_i x+k )}{\prod_{j=1}^{\xi}(x+j )^{n+r+1}}\Big)\\
&=& \sum_{\xi=0}^{\beta}\sum_{l=0}^{\beta- \xi}\bigg(\frac{\big(-\elld\big)^{\beta-\xi}}{(\beta-l-\xi)!}\sum_{p=0}^{n-b}\Big(
      \mathrm{Coeff}_{x^p}\big( f_{a,l}(x+\xi)\big)\cdot
       \mathrm{Coeff}_{x^{n-b-p}}\big(\frac{\prod_{i=1}^r \prod_{k=1}^{d_i \xi}(d_i x+k )}{\prod_{j=1}^{\xi}(x+j )^{n+r+1}}\big) 
      \Big)\bigg),
\end{eqnarray*}
where
\begin{eqnarray*}
&&  f_{a,l}(x+\xi)\\
&=& \sum_{i_0+\dots+i_l=a-l}(x+\xi)^{i_0}\cdots(x+\xi+l)^{i_l}\\
&=& \sum_{p\geq 0}\bigg(\sum_{\begin{subarray}{c}i_0+\dots+i_l=a-l\\ j_0+\dots+j_l=p\\
i_0,\dots,i_l,j_0,\dots,j_l\geq 0\end{subarray}} \binom{i_0}{j_0}\xi^{i_0-j_0}\cdots\binom{i_l}{j_l}(\xi+l)^{i_l-j_l}\bigg)x^p,
\end{eqnarray*}
and
\begin{eqnarray*}
&&\mathrm{Coeff}_{x^{n-b-p}}\big(\frac{\prod_{i=1}^r \prod_{k=1}^{d_i \xi}(d_i x+k )}{\prod_{j=1}^{\xi}(x+j )^{n+r+1}}\big)\\
&=&\sum_{\begin{subarray}{c}j_1+\dots+j_{\xi}+k_1+\dots+k_r=n-b-p\\
      j_1,\dots,j_{\xi},k_1,\dots,k_r\geq 0\end{subarray}}  
      \Big((\xi!)^{-n-r-1}\prod_{e=1}^{\xi}\binom{-n-r-1}{j_e}e^{-j_e}\prod_{i=1}^r ((d_i \xi)! d_i^{k_i}\mu(d_i \xi,k_i))\Big)\\
&=&\frac{\prod_{i=1}^r (d_i \xi)! }{(\xi!)^{n+r+1}}\sum_{\begin{subarray}{c}j_1+\dots+j_{\xi}+k_1+\dots+k_r=n-b-p\\
      j_1,\dots,j_{\xi},k_1,\dots,k_r\geq 0\end{subarray}}  
      \Big(\prod_{e=1}^{\xi}\binom{-n-r-1}{j_e}e^{-j_e}\prod_{i=1}^r  d_i^{k_i}\mu(d_i \xi,k_i)\Big).
\end{eqnarray*}

\subsection{Computation of the function \texorpdfstring{$\mathsf{A}_N$}{A N}}
The practical computation of (many occurrences of) the function $\mathsf{A}_N$ on  the RHS of (\ref{eq-recusion-twoPointInv}) turns out to be  rather time-consuming. In this section we provide another way to compute the terms like
\begin{equation}\label{eq-functionA-toCompute}
\mathsf{A}_{N}\big(x_1,\dots,x_p\big),\ \mbox{with nonnegative integer arguments}\  x_i,
\end{equation}
which will make the computations of the correlators of length 2 and therefore the matrices $W$ and $M$ much faster. The readers that do not care about the performance of our Macaulay2 package, but only the structure of the algorithm, can skip this section.\\

First we regard $x_1,\dots,x_p$ as indeterminates and prove identities of rational functions of  $x_1,\dots,x_p$.
\begin{lemma}\label{lem-functionA-identity-2}
\begin{eqnarray}\label{eq-functionA-identity-2}
      &&\sum_{i=0}^{p-1} (-1)^i \frac{(\sum_{j=1}^{p-i}x_j)^{N}}{\prod_{j=1}^{p-i}(\sum_{k=j}^{p-i}x_k)\cdot \prod_{j=p-i+1}^{p}(\sum_{k=p-i+1}^j x_k)}\nn\\
      &=&\begin{cases}
       \frac{(-1)^{p-1}}{ \prod_{j=1}^{p}(\sum_{k=1}^j x_k)},& \mbox{if}\ N=0;\\
      0, & \mbox{if}\ 1\leq N\leq p-1.
      \end{cases}
\end{eqnarray}
\end{lemma}
\begin{proof}
Let $y_i=\sum_{j=1}^{p}x_j$. Then (\ref{eq-functionA-identity-2}) is equivalent to 
\begin{equation*}
\sum_{i=0}^{p-1}  \frac{y_{p-i}^{N}}{y_{p-i}\prod_{j\neq p-i}(y_{p-i}-y_j)}   
=\begin{cases}
       \frac{(-1)^{p-1}}{ \prod_{j=1}^{p}y_j},& \mbox{if}\ N=0;\\
      0, & \mbox{if}\ 1\leq N\leq p-1.
      \end{cases}
\end{equation*}
\end{proof}

\begin{proposition}\label{prop-generatingFuction-nu}
\begin{eqnarray}\label{eq-generatingFuction-nu}
\sum_{\begin{subarray}{c}i_1,\dots,i_p\geq 0\\
i_1+\dots+i_p=N\end{subarray}}
\nu_{i_1,\dots,i_p}x_1^{i_1}\cdots x_p^{i_p}
=\sum_{i=0}^{p-1}  \frac{(-1)^i(\sum_{j=1}^{p-i}x_j)^{N+p-1}}{\prod_{j=2}^{p-i}(\sum_{k=j}^{p-i}x_k)\cdot \prod_{j=p-i+1}^{p}(\sum_{k=p-i+1}^j x_k)}.
\end{eqnarray}
\end{proposition}
\begin{proof}
We do induction on $p$. Suppose (\ref{eq-generatingFuction-nu}) holds for $p-1$. Then by (\ref{eq-coefficient-derivativeOfJ}),
\begin{eqnarray*}
&& \sum_{\begin{subarray}{c}i_1,\dots,i_p\geq 0\\
i_1+\dots+i_p=N\end{subarray}}
\nu_{i_1,\dots,i_p}x_1^{i_1}\cdots x_p^{i_p}\\
&=& \sum_{\begin{subarray}{c}i_1,\dots,i_p\geq 0\\
i_1+\dots+i_p=N\end{subarray}}\binom{i_1+\dots+i_p+p-1}{i_1}\binom{i_2+\dots+i_p+p-2}{i_2}\dots\binom{i_{p-1}+i_p+1}{i_{p-1}}x_1^{i_1}\cdots x_p^{i_p}\\
&=& \sum_{i_1=0}^N \bigg(\binom{i_1+\dots+i_p+p-1}{i_1} x_1^{i_1}\\
 &&\cdot \sum_{\begin{subarray}{c}i_2,\dots,i_p\geq 0\\
i_2+\dots+i_p=N-i_1\end{subarray}}
\binom{i_2+\dots+i_p+p-2}{i_2}\dots\binom{i_{p-1}+i_p+1}{i_{p-1}}x_2^{i_2}\cdots x_p^{i_p}\bigg)\\
&=& \sum_{i_1=0}^N \bigg(\binom{N+p-1}{i_1} x_1^{i_1}
\sum_{i=0}^{p-2}(-1)^i \frac{(\sum_{j=2}^{p-i}x_j)^{N-i_1+p-2}}{\prod_{j=3}^{p-i}(\sum_{k=j}^{p-i}x_k)\cdot \prod_{j=p-i+1}^{p}(\sum_{k=p-i+1}^j x_k)}\bigg)\\
&=& \sum_{i=0}^{p-2}  \bigg(\frac{(-1)^i}{\prod_{j=3}^{p-i}(\sum_{k=j}^{p-i}x_k)\cdot \prod_{j=p-i+1}^{p}(\sum_{k=p-i+1}^j x_k)}\sum_{i_1=0}^N \binom{N+p-1}{i_1} x_1^{i_1}\big(\sum_{j=2}^{p-i}x_j\big)^{N-i_1+p-2}\bigg)\\
&=& \sum_{i=0}^{p-2}  \frac{(-1)^i}{\prod_{j=2}^{p-i}(\sum_{k=j}^{p-i}x_k)\cdot \prod_{j=p-i+1}^{p}(\sum_{k=p-i+1}^j x_k)}\sum_{i_1=0}^N \binom{N+p-1}{i_1} x_1^{i_1}\big(\sum_{j=2}^{p-i}x_j\big)^{N-i_1+p-1}\\
&=& \sum_{i=0}^{p-2}\bigg( \frac{(-1)^i }{\prod_{j=2}^{p-i}(\sum_{k=j}^{p-i}x_k)\cdot \prod_{j=p-i+1}^{p}(\sum_{k=p-i+1}^j x_k)}\\
&&\cdot \Big(\big(\sum_{j=1}^{p-i}x_j\big)^{N+p-1}-
\sum_{i_1=N+1}^{N+p-1} \binom{N+p-1}{i_1} x_1^{i_1}\big(\sum_{j=2}^{p-i}x_j\big)^{N-i_1+p-1}\Big)\bigg)\\
&=& \sum_{i=0}^{p-2} (-1)^i \frac{(\sum_{j=1}^{p-i}x_j)^{N+p-1}}{\prod_{j=2}^{p-i}(\sum_{k=j}^{p-i}x_k)\cdot \prod_{j=p-i+1}^{p}(\sum_{k=p-i+1}^j x_k)}\\
&&- \sum_{i_1=N+1}^{N+p-1} \bigg(\binom{N+p-1}{i_1} x_1^{i_1}\sum_{i=0}^{p-2}\frac{ (-1)^i(\sum_{j=2}^{p-i}x_j)^{N-i_1+p-1} }{\prod_{j=2}^{p-i}(\sum_{k=j}^{p-i}x_k)\cdot \prod_{j=p-i+1}^{p}(\sum_{k=p-i+1}^j x_k)}\bigg).
\end{eqnarray*}
Applying Lemma \ref{lem-functionA-identity-2} to $x_2,\dots,x_p$, the sum
\[
\sum_{i=0}^{p-2}\frac{ (-1)^i(\sum_{j=2}^{p-i}x_j)^{N-i_1+p-1} }{\prod_{j=2}^{p-i}(\sum_{k=j}^{p-i}x_k)\cdot \prod_{j=p-i+1}^{p}(\sum_{k=p-i+1}^j x_k)}
\]
vanishes for $N+1\leq i_1<N+p-1$, and equals $ \frac{(-1)^{p-2}}{ \prod_{j=2}^{p}(\sum_{k=1}^j x_k)}$ when $i_1=N+p-1$. So
\[
\sum_{i_1=N+1}^{N+p-1} \bigg(\binom{N+p-1}{i_1} x_1^{i_1}\sum_{i=0}^{p-2}\frac{ (-1)^i(\sum_{j=2}^{p-i}x_j)^{N-i_1+p-1} }{\prod_{j=2}^{p-i}(\sum_{k=j}^{p-i}x_k)\cdot \prod_{j=p-i+1}^{p}(\sum_{k=p-i+1}^j x_k)}\bigg)
=\frac{(-1)^{p-2}x_1^{N+p-1}}{ \prod_{j=2}^{p}(\sum_{k=1}^j x_k)}.
\]
Hence (\ref{eq-generatingFuction-nu}) holds for $p$.
\end{proof}

To compute (\ref{eq-functionA-toCompute}) by using (\ref{eq-generatingFuction-nu}), there is an obstacle is that when some $x_i$'s are 0, a certain denominator on RHS of  (\ref{eq-generatingFuction-nu}) may  vanish. We bypass this obstacle by \emph{perturbing} the arguments $x_i$. To do this we need an a priori estimate.
\begin{proposition}\label{prop-perturbationComputationOfFunctionA}
Let $K$ be a proper subset of $\{1,\dots,p\}$. Let $C$ be a positive real number. Suppose $0\leq x_i\leq C$ for $1\leq i\leq p$, such that $x_i=0$ for $i\in K$ and $x_i\neq 0$ for $i\not\in K$. Suppose 
\begin{equation}\label{eq-perturbation-epsilon}
     0< \epsilon\leq \frac{1}{|K|N(p-|K|)^{N-1}C^{N-1}\binom{N+p-1}{N}}\Big(1-\frac{N-1}{2N (p-|K|)^NC^N\binom{N+p-1}{N}}\Big).
\end{equation}
Let 
\begin{equation*}
      \tilde{x}_i=\begin{cases}
      \epsilon,& \mbox{if}\ i\in K;\\
      x_i,& \mbox{if}\ i\not \in K.
      \end{cases}
\end{equation*}
Then
\begin{equation}\label{eq-perturbationComputationOfFunctionA}
       0\leq \mathsf{A}_N(\tilde{x}_1,\dots,\tilde{x}_p)-\mathsf{A}_N(x_1,\dots,x_p)<1.
\end{equation}
In particular, if we assume moreover $x_i \in \mathbb{Z}$, then 
\[
\mathsf{A}_N(x_1,\dots,x_p)=\lfloor \mathsf{A}_N(\tilde{x}_1,\dots,\tilde{x}_p)\rfloor.
\]
\end{proposition}
\begin{proof}
The left inequality is obvious because the coefficients $\nu_{i_1,\dots,i_p}$ in (\ref{eq-def-functionA}) are nonnegative. For the right inequality, 
\begin{eqnarray*}
&& \mathsf{A}_N(\tilde{x}_1,\dots,\tilde{x}_p)-\mathsf{A}_N(x_1,\dots,x_p)\\
&=& \big(\sum_{q=1}^{N} \sum_{\begin{subarray}{c}i_1,\dots,i_p\geq 0\\
i_1+\dots+i_p=N\\ \sum_{j\in K}i_j=q
\end{subarray}}
\nu_{i_1,\dots,i_p}x_1^{i_1}\cdots x_p^{i_p}\big)|_{x_i=\epsilon,i\in K}\\
&\leq& \sum_{q=1}^{N} \sum_{\begin{subarray}{c}i_1,\dots,i_p\geq 0\\
i_1+\dots+i_p=N\\ \sum_{j\in K}i_j=q
\end{subarray}}
\nu_{i_1,\dots,i_p} C^{N-q} \epsilon^{q}\\
&=&\sum_{q=1}^{N} \sum_{\begin{subarray}{c}i_1,\dots,i_p\geq 0\\
i_1+\dots+i_p=N\\ \sum_{j\in K}i_j=q
\end{subarray}}
\binom{i_1+\dots+i_p+p-1}{i_1}\binom{i_2+\dots+i_p+p-2}{i_2}\dots\binom{i_{p-1}+i_p+1}{i_{p-1}}C^{N-q} \epsilon^{q}\\
&\leq& \sum_{q=1}^{N} \sum_{\begin{subarray}{c}i_1,\dots,i_p\geq 0\\
i_1+\dots+i_p=N\\ \sum_{j\in K}i_j=q
\end{subarray}} \binom{N+p-1}{i_1,i_2,\dots,i_p} C^{N-q} \epsilon^{q}\\
&=& \sum_{q=1}^{N} \sum_{\begin{subarray}{c}i_1,\dots,i_p\geq 0\\
i_1+\dots+i_p=N\\ \sum_{j\in K}i_j=q
\end{subarray}}
\binom{i_1+\dots+i_p+p-1}{i_1}\binom{i_2+\dots+i_p+p-2}{i_2}\dots\binom{i_{p-1}+i_p+1}{i_{p-1}}C^{N-q} \epsilon^{q}\\
&=& \sum_{q=1}^{N} \sum_{\begin{subarray}{c}i_1,\dots,i_p\geq 0\\
i_1+\dots+i_p=N\\ \sum_{j\in K}i_j=q
\end{subarray}} \binom{N+p-1}{N}\binom{N}{q}\binom{N-q}{\underbrace{i_j,\dots}_{j\not\in K}}\binom{q}{\underbrace{i_k,\dots}_{k\in K}}C^{N-q} \epsilon^{q}\\
&=& \sum_{q=1}^{N}  \binom{N+p-1}{N}\binom{N}{q} (p-|K|)^{N-q} |K|^{q}C^{N-q} \epsilon^{q}\\
&=& \binom{N+p-1}{N}\Big(\big((p-|K|)C+|K|\epsilon\big)^N-(p-|K|)^{N}C^{N} \Big).
\end{eqnarray*}
So $ \mathsf{A}_N(\tilde{x}_1,\dots,\tilde{x}_p)-\mathsf{A}_N(x_1,\dots,x_p)<1$ if 
\begin{equation}\label{eq-estimate-1}
      \big((p-|K|)C+|K|\epsilon\big)^N-(p-|K|)^{N}C^{N} < \binom{N+p-1}{N}^{-1}.
\end{equation}
By the assumption, $|K|<p$. Then (\ref{eq-estimate-1}) is equivalent to
\begin{eqnarray*}
&&|K| \epsilon<\bigg((p-|K|)^NC^N+ \binom{N+p-1}{N}^{-1}\bigg)^{1/N}-(p-|K|)C\\
&=& (p-|K|)C \bigg(1+ \frac{1}{(p-|K|)^NC^N\binom{N+p-1}{N}}\bigg)^{1/N}.
\end{eqnarray*}
Note that for $x>0$ and $N\geq 1$, we have
\[
(1+x)^{1/N}>1+\frac{x}{N}+\binom{1/N}{2}x^2.
\]
Thus
\[
\bigg(1+ \frac{1}{(p-|K|)^NC^N\binom{N+p-1}{N}}\bigg)^{1/N}
>\frac{1}{N(p-|K|)^NC^N\binom{N+p-1}{N}}\big(1-\frac{N-1}{2N (p-|K|)^NC^N\binom{N+p-1}{N}}\big).
\]
This leads to the condition (\ref{eq-perturbation-epsilon}).
\end{proof}

As a consequence we  compute 
\[
\mathsf{A}_{a-1-p}\big(\frac{2- \alpha_1}{\mathsf{a}(n,\mathbf{d})},\dots,\frac{1+\alpha_{p-1}- \alpha_p}{\mathsf{a}(n,\mathbf{d})}\big)
\]
in (\ref{eq-recusion-twoPointInv}) using Proposition \ref{prop-perturbationComputationOfFunctionA}, where we  take
\begin{equation*}
      \epsilon=\frac{1}{N(a-2)^N a^{N-1}2^{a-1}}\leq \frac{1}{N(a-2)^N \lfloor \frac{a}{\mathsf{a}(n,\mathbf{d})}\rfloor^{N-1}2^{a-2}}\cdot\frac{1}{2}.
\end{equation*}

\subsection{Recursions for the matrices \texorpdfstring{$W$}{W} and \texorpdfstring{$M$}{M}}
\label{sec:appendix-computation-MandW}
In this section we give recursive formulae of the entries of the transformation matrix between the basis $\sfh_i$'s and $\tsfh_i$'s, in terms of the length 2 correlators.

\begin{lemma}\label{lem-recursionOfW-FanoIndexAtLeastTwo}
Suppose the Fano index $\mathsf{a}(n,\mathbf{d})\geq 2$. Then
\begin{equation}\label{eq-recursion-W-initialValues-FanoIndexAtLeastTwo}
      W_{0}^j=\delta_{0,j},\ W_1^j=\delta_{1,j},
\end{equation}
and for $i\geq 1$ and $0\leq j\leq n$,
\begin{equation}\label{eq-recursion-W-FanoIndexAtLeastTwo}
      W_{i+1}^j=W_i^{j-1}+\frac{1}{\prod_{i=1}^r d_i}\big(\sum_{\beta=1}^{\lfloor \frac{i-j+1}{\mathsf{a}(n,\mathbf{d})}\rfloor} \beta W_i^{j-1+\beta\cdot \mathsf{a}(n,\mathbf{d})}\langle \sfh_{j-1+\beta\cdot \mathsf{a}(n,\mathbf{d})},\sfh_{n-j}\rangle_{0,\beta}\big).
\end{equation}
\end{lemma}
\begin{proof}
Since $\mathsf{a}(n,\mathbf{d})\geq 2$, so $\tsfh_1=\sfh_1=\sfh$, and thus (\ref{eq-recursion-W-initialValues-FanoIndexAtLeastTwo}) holds. 
By definition of $W$. 
\[
\tsfh_i=\sum_j W_i^j \sfh_j.
\]
So
\begin{eqnarray*}
&& \tsfh_{i+1}=\sum_j W_i^j \sfh\sqp \sfh_j\frac{1}{\prod_{i=1}^r d_i} \sum_j W_i^j \langle \sfh,\sfh_j,\sfh_k\rangle_0 \sfh_{n-k}\\
&=& \frac{1}{\prod_{i=1}^r d_i} \sum_{j=0}^{n} \big( \sum_{k}W_i^k \langle \sfh,\sfh_k,\sfh_{n-j}\rangle_0
\big)\sfh_j.
\end{eqnarray*}
Comparing the coefficients of $\sfh_j$ of both sides and taking into account (\ref{eq-Dim}), we get (\ref{eq-recursion-W-FanoIndexAtLeastTwo}).
\end{proof}

\begin{lemma}\label{lem-recursionOfW-FanoIndexOne}
Suppose the Fano index $\mathsf{a}(n,\mathbf{d})=1$. Then
\begin{equation}\label{eq-recursion-W-initialValues-FanoIndexOne}
      W_{0}^j=\delta_{0,j},\ W_1^0=\elld,\ \mbox{and}\ W_1^j=\delta_{1,j}\ \mbox{for}\ j\geq 1,
\end{equation}
and for $i\geq 1$ and $0\leq j\leq n$,
\begin{equation}\label{eq-recursion-W-FanoIndexOne}
      W_{i+1}^j=W_i^{j-1}+\frac{1}{\prod_{i=1}^r d_i}\big(\sum_{\beta=1}^{i-j+1} \beta W_i^{j-1+\beta}\langle \sfh_{j-1+\beta},\sfh_{n-j}\rangle_{0,\beta}\big)+\elld W_i^j.
\end{equation}
\end{lemma}
\begin{proof}
Since the Fano index $\mathsf{a}(n,\mathbf{d})=1$, $\tsfh_1=\sfh_1+\elld=\sfh+\elld$. 
So
\begin{eqnarray*}
&& \tsfh_{i+1}=\sum_j W_i^j \sfh\sqp \sfh_j+\elld\sum_j W_i^j \sfh_j= \frac{1}{\prod_{i=1}^r d_i} \sum_j W_i^j \langle \sfh,\sfh_j,\sfh_k\rangle_0 \sfh_{n-k}+\elld\sum_j W_i^j \sfh_j\\
&=& \sum_{j=0}^{n} \big(\frac{1}{\prod_{i=1}^r d_i} \sum_{k}W_i^k \langle \sfh,\sfh_k,\sfh_{n-j}\rangle_0+\elld W_i^j\big)\sfh_j.
\end{eqnarray*}
Then (\ref{eq-recursion-W-FanoIndexOne}) follows.
\end{proof}
\begin{lemma}\label{lem-recursionOfM}
Suppose the Fano index $\mathsf{a}(n,\mathbf{d})\geq 1$. Then
\begin{equation}\label{eq-recursionOfM-initial}
      M_i^i=1,\ M_i^{j}=0\ \mbox{for}\ j>i,
\end{equation}
and for $k<i$
\begin{equation}\label{eq-recursionOfM}
      M_i^k=-\sum_{0\leq j<i}\sum_{k=0}^j W_i^j M_j^k.
\end{equation}
\end{lemma}
\begin{proof}
Recall that $W$ is lower triangular with diagonal entries all equal to 1. So we have (\ref{eq-recursionOfM-initial}). Moreover by definition, 
\begin{equation*}
\sfh_i=\tsfh_i-\sum_{j<i} W_i^j \sfh_j=\tsfh_i-\sum_{j<i}\sum_{k=0}^j W_i^j M_j^k \tsfh_k.
\end{equation*}
Hence (\ref{eq-recursionOfM}).
\end{proof}

\subsection{Recursion for \texorpdfstring{$F^{(0)}$}{F(0)}}\label{sec:algorithm-recursionOfAmbientGeneratingFunction}
Recall (\ref{eq-pairing1}), (\ref{eq-pairing2}), and 
\begin{eqnarray}\label{eq-qp1.5-appendix}
\sfF_{\tau^a \tau^b \tau^c}^{(0)}(0)=\left\{
\begin{array}{cc}
(\sfbd)^{\frac{a+b+c-n}{\mathsf{a}(n,\mathbf{d})}} \mathsf{q}^{\frac{a+b+c-n}{\mathsf{a}(n,\mathbf{d})}}\prod_{i=1}^{r}d_i, & \mathrm{if}\ \frac{a+b+c-n}{\mathsf{a}(n,\mathbf{d})}\in \mathbb{Z}_{\geq 0}; \\
0, & \mathrm{otherwise}.
\end{array}
\right.
\end{eqnarray}

\begin{lemma}\label{lem-recursion-F(0)-EulerField}
\begin{eqnarray}\label{eq-recursion-F(0)-EulerField}
&&\partial_{\tau^0}^{p_0}\partial_{\tau^1}^{p_1+1}\cdots\partial_{\tau^n}^{p_n}F^{(0)}(0)\nn\\
&=&\frac{1}{\mathsf{a}(n,\mathbf{d})}\big(-\sum_{i=0}^{n}\sum_{j=0}^n\sum_{k=0}^n
(1-i)W_j^i  M_{i}^k p_j\partial_{\tau^0}^{p_0}\cdots\partial_{\tau^k}^{p_k+1}\cdots
\partial_{\tau^j}^{p_j-1}\cdots 
\partial_{\tau^n}^{p_n}F^{(0)}(0)\nn\\
&&+(3-n)\partial_{\tau^0}^{p_0}\cdots\partial_{\tau^n}^{p_n}F^{(0)}(0)\big).
\end{eqnarray}
In the summation, when $k=j$, the operator $\partial_{\tau^k}^{p_k+1}\partial_{\tau^j}^{p_j-1}$ is understood as $\partial_{\tau^k}^{p_k}$. See (\ref{eq-recursion-F(0)-EulerField-I}) for a more concise expression.
\end{lemma}
\begin{proof}
The ambient part of the Euler field is
\begin{eqnarray*}
E^{(0)}&=& \sum_{i=0}^{n}(1-i)t^{i}\frac{\partial}{\partial t^i}+\mathsf{a}(n,\mathbf{d})\frac{\partial}{\partial t^1}\\
&=&\sum_{i=0}^{n}\sum_{j=0}^n\sum_{k=0}^n
(1-i)W_j^i  M_{i}^k \tau^j\frac{\partial}{\partial \tau^k}+\mathsf{a}(n,\mathbf{d})\sum_{k=0}^n
M_{1}^k
\frac{\partial}{\partial \tau^k}\\
&=&\sum_{i=0}^{n}\sum_{j=0}^n\sum_{k=0}^n
(1-i)W_j^i  M_{i}^k \tau^j\frac{\partial}{\partial \tau^k}
+\mathsf{a}(n,\mathbf{d})\frac{\partial}{\partial \tau^1}
-\delta_{1,\mathsf{a}(n,\mathbf{d})} \elld\frac{\partial}{\partial \tau^0}.
\end{eqnarray*}
Then
\[
E^{(0)} F^{(0)}=(3-n)F^{(0)}+\mathsf{a}(n,\mathbf{d})\partial_{t^1}c|_{s=0}.
\]
So for $p_0+\dots+p_n\geq 3$,
\begin{eqnarray*}
&& (\partial_{\tau^0}^{p_0}\cdots\partial_{\tau^n}^{p_n}\circ E^{(0)}) F^{(0)}=(3-n)\partial_{\tau^0}^{p_0}\cdots\partial_{\tau^n}^{p_n}F^{(0)},
\end{eqnarray*}
i.e.
\begin{eqnarray*}
&& \sum_{i=0}^{n}\sum_{j=0}^n\sum_{k=0}^n
(1-i)W_j^i  M_{i}^k p_j \partial_{\tau^0}^{p_0}\cdots\partial_{\tau^k}^{p_k+1}\cdots
\partial_{\tau^j}^{p_j-1}\cdots 
\partial_{\tau^n}^{p_n}F^{(0)}(0)\\
&&+\mathsf{a}(n,\mathbf{d})\partial_{\tau^0}^{p_0}\partial_{\tau^1}^{p_1+1}\cdots\partial_{\tau^n}^{p_n}F^{(0)}(0)\\
&=&(3-n)\partial_{\tau^0}^{p_0}\cdots\partial_{\tau^n}^{p_n}F^{(0)}(0),
\end{eqnarray*}
thus (\ref{eq-recursion-F(0)-EulerField}) follows.
\end{proof}

For a multiset $S=\{i_1,\dots,i_k\}$ where $0\leq i_j\leq n$, denote by $\partial_{\tau^S}$ the differential operator
\begin{equation*}
\partial_{\tau^{i_1}}\circ\dots\circ \partial_{\tau^{i_k}}.
\end{equation*}

\begin{lemma}\label{lem-recursion-F(0)-WDVV-tau}
For $i\geq 1$, and $0\leq j,k\leq n$,
\begin{eqnarray}\label{eq-recursion-F(0)-WDVV-tau}
&&\partial_{\tau^i}\partial_{\tau^j}\partial_{\tau^k}\partial_{\tau^S}F^{(0)}(0)\nn\\
&=& -\sum_{\begin{subarray}{c}S_1\sqcup S_2=S\\ S_1\neq \emptyset\end{subarray}}\big(\partial_{\tau^{S_1}}\partial_{\tau^1} \partial_{\tau^{i-1}}\partial_{\tau^a}F^{(0)}(0)\big)\eta^{ab}\big(\partial_{\tau^{b}}\partial_{\tau^j}\partial_{\tau^k}
\partial_{\tau^{S_2}}F^{(0)}(0)\big)\nn\\
&&+\sum_{S_1\sqcup S_2=S}\big(\partial_{\tau^{S_1}}\partial_{\tau^1} \partial_{\tau^{j}}\partial_{\tau^a}F^{(0)}(0)\big)\eta^{ab}\big(\partial_{\tau^{b}}\partial_{\tau^{i-1}}\partial_{\tau^k}
\partial_{\tau^{S_2}}F^{(0)}(0)\big).
\end{eqnarray}
\end{lemma}
\begin{proof}
Set 
\begin{equation*}
      \partial_{\tau^{k}}=\begin{cases}
      \partial_{\tau^k}, & \mbox{if}\ 0\leq k\leq n,\\
      (\sfbd)^{1+\lfloor \frac{k-n-1}{\mathsf{a}(n,\mathbf{d})}}\partial_{\tau^{k-\big(1+\lfloor \frac{k-n-1}{\mathsf{a}(n,\mathbf{d})}\big)\mathsf{a}(n,\mathbf{d})}},
      & \mbox{if}\ k>n.
      \end{cases}
\end{equation*}
Then 
\begin{eqnarray}\label{eq-leadingTerm-WDVV}
&& (\partial_{\tau^1}\sqp \partial_{\tau^{i-1}})\circ (\partial_{\tau^j}\circ\partial_{\tau^k})+
(\partial_{\tau^1}\circ \partial_{\tau^{i-1}})\circ (\partial_{\tau^j}\sqp\partial_{\tau^k})\nn\\
&&-(\partial_{\tau^1}\sqp \partial_{\tau^j})\circ (\partial_{\tau^{i-1}}\circ\partial_{\tau^k})
-(\partial_{\tau^1}\circ\partial_{\tau^j})\circ (\partial_{\tau^{i-1}}\sqp\partial_{\tau^{k}})\nn\\
&=& \partial_{\tau^i}\partial_{\tau^j}\partial_{\tau^k}+\partial_{\tau^1}\partial_{\tau^i}\partial_{\tau^{j+k}}
-\partial_{\tau^{i-1}}\partial_{\tau^{j+1}}\partial_{\tau^k}
-\partial_{\tau^1}\partial_{\tau^j}\partial_{\tau^{i+k-1}}.
\end{eqnarray}
We make use of the WDVV
\begin{equation*}
      (\partial_{\tau^1}\partial_{\tau^{i-1}}\partial_{\tau^e}F)\eta^{ef}(\partial_{\tau^f}\partial_{\tau^{j}}\partial_{\tau^k}F)
      = (\partial_{\tau^1}\partial_{\tau^{j}}\partial_{\tau^{e}}F)\eta^{ef} (\partial_{\tau^f}\partial_{\tau^{i-1}}\partial_{\tau^k}F).
\end{equation*}
The leading terms are given by (\ref{eq-leadingTerm-WDVV}). So we get (\ref{eq-recursion-F(0)-WDVV-tau}).
\end{proof}

Lemma \ref{lem-recursion-F(0)-EulerField}, and Lemma \ref{lem-recursion-F(0)-WDVV-tau} applied to $2\leq i\leq n$, compute $F^{(0)}$ recursively from the correlators of length 3.\\

For ease of comparison to our codes in Macaulay2, we adopt the notations
\begin{equation*}
      \binom{I}{J},\ \partial_{\tau^I}
\end{equation*}
for  $I=(i_0,\dots,i_n)\in \mathbb{Z}^{n+1}$  introduced in Section \ref{sec:sqrtRecursion}. Moreover, for $0\leq j\leq n$ let $\mathbf{e}_j$ be the $j$-th unit vector in $\mathbb{Z}^{n+1}$.
Then (\ref{eq-recursion-F(0)-EulerField}) can be rewritten as, for $I=(p_0,\dots,p_n)\in \mathbb{Z}_{\geq 0}^{n+1}$ with $|I|\geq 1$,
\begin{eqnarray}\label{eq-recursion-F(0)-EulerField-I}
\partial_{\tau^{I+\mathbf{e}_1}}F^{(0)}(0)
&=&\frac{1}{\mathsf{a}(n,\mathbf{d})}\big(-\sum_{i=0}^{n}\sum_{j=0}^n\sum_{k=0}^n
(1-i)W_j^i  M_{i}^k p_j\partial_{\tau^{I+\mathbf{e}_k- \mathbf{e}_j}}F^{(0)}(0)\nn\\
&&+(3-n)\partial_{\tau^{I}}F^{(0)}(0)\big),
\end{eqnarray}
and (\ref{eq-recursion-F(0)-WDVV-tau}) as, for $I\in \mathbb{Z}_{\geq 0}^{n+1}$,
\begin{eqnarray}\label{eq-recursion-F(0)-WDVV-tau-I}
&&\partial_{\tau^i}\partial_{\tau^j}\partial_{\tau^k}\partial_{\tau^I}F^{(0)}(0)\nn\\
&=& -\sum_{\begin{subarray}{c}0^{n+1}\leq J\leq I\\J\neq 0^{n+1}\end{subarray}}\binom{I}{J}\sum_{a=0}^n\sum_{b=0}^n
\big(\partial_{\tau^{J}}\partial_{\tau^1} \partial_{\tau^{i-1}}\partial_{\tau^a}F^{(0)}(0)\big)\eta^{ab}\big(\partial_{\tau^{b}}\partial_{\tau^j}\partial_{\tau^k}
\partial_{\tau^{I-J}}F^{(0)}(0)\big)\nn\\
&&+\sum_{0^{n+1}\leq  J\leq I}\binom{I}{J}\sum_{a=0}^n\sum_{b=0}^n\big(\partial_{\tau^{J}}\partial_{\tau^1} \partial_{\tau^{j}}\partial_{\tau^a}F^{(0)}(0)\big)\eta^{ab}\big(\partial_{\tau^{b}}\partial_{\tau^{i-1}}\partial_{\tau^k}
\partial_{\tau^{I-J}}F^{(0)}(0)\big).
\end{eqnarray}
In the following sections we stick to this notation.

\begin{remark}\label{rem:convergence-F(0)}
As a direct consequence of (\ref{eq-recursion-F(0)-EulerField}) and (\ref{eq-recursion-F(0)-WDVV-tau-I}), there exists a constant $C>0$ such that 
\begin{equation*}
      |\partial_{\tau^I}F^{(0)}(0)|\leq (|I|-3)!C^{|I|-3}.
\end{equation*}
for all $I\in \mathbb{Z}^{n+1}_{\geq 0}$. This gives a simple proof of  \cite[Theorem 1]{Zin14} restricted to primary genus 0 invariants, and therefore the ambient Frobenius manifold associated with $X$ is analytic.
\end{remark}

\subsection{Recursion for \texorpdfstring{$F^{(1)}$}{F(1)}}
We rewrite the system of equations (\ref{eq-ss1}) in the $\tau$-coordinates:
\begin{equation}\label{eq-ss1-tau}
      (\partial_{\tau^a}\partial_{\tau^b}\partial_{\tau^e}F^{(0)})\eta^{ef}\partial_{\tau^f}F^{(1)}=(\partial_{\tau^a}F^{(1)})(\partial_{\tau^b}F^{(1)}).
\end{equation}
By Theorem \ref{thm-F(1)} and \ref{thm-F(1)-fanoIndex=1},
\begin{equation}\label{eq-recursion-F(1)-initial}
F^{(1)}(0)=-\elld\delta_{\mathsf{a}(n,\mathbf{d}),1},\
\partial_{\tau^i}F^{(1)}(0)=\delta_{i,0}.
\end{equation}
and we regard these as initial values.

\begin{lemma}\label{lem-recursion-F(1)-EulerField}
For $I=(p_0,\dots,p_n)\in \mathbb{Z}_{\geq 0}^{n+1}$ with $|I|\geq 1$,
\begin{equation}\label{eq-recursion-F(1)-EulerField}
\mathsf{a}(n,\mathbf{d})\partial_{\tau^1}\partial_{\tau^I}F^{(1)}(0)=
- \sum_{i=0}^{n}\sum_{j=0}^n\sum_{k=0}^n
(1-i)W_j^i  M_{i}^k p_j \partial_{\tau^{I+\mathbf{e}_k- \mathbf{e}_j}}F^{(1)}(0)+\partial_{\tau^I}F^{(1)}(0).
\end{equation}
\end{lemma}
\begin{proof}
Let $I=(p_0,\dots,p_n)$.
\begin{eqnarray*}
&&E=E^{(0)}+(2-n)s\partial_{s}\\
&=& \sum_{i=0}^{n}(1-i)t^{i}\frac{\partial}{\partial t^i}+\mathsf{a}(n,\mathbf{d})\frac{\partial}{\partial t^1}+(2-n)s\partial_{s}\\
&=&\sum_{i=0}^{n}\sum_{j=0}^n\sum_{k=0}^n
(1-i)W_j^i  M_{i}^k \tau^j\frac{\partial}{\partial \tau^k}
+\mathsf{a}(n,\mathbf{d})\frac{\partial}{\partial \tau^1}
-\delta_{1,\mathsf{a}(n,\mathbf{d})} \elld\frac{\partial}{\partial \tau^0}
+(2-n)s\partial_{s}.
\end{eqnarray*}
and
\[
E F=(3-n)F+\mathsf{a}(n,\mathbf{d})\partial_{t^1}c.
\]
So
\[
E^{(0)}F^{(1)}+(2-n)F^{(1)}=(3-n)F^{(1)},
\]
i.e.
\[
E^{(0)}F^{(1)}=F^{(1)},
\]
i.e.
\begin{eqnarray*}
&&\mathsf{a}(n,\mathbf{d})\partial_{\tau^1}F^{(1)}
      =-\sum_{i=0}^{n}\sum_{j=0}^n\sum_{k=0}^n
(1-i)W_j^i  M_{i}^k \tau^j\partial_{\tau^k}F^{(1)}
+\delta_{1,\mathsf{a}(n,\mathbf{d})} \elld\partial_{\tau^0}F^{(1)}
+F^{(1)}\\
&=& -\sum_{i=0}^{n}\sum_{j=0}^n\sum_{k=0}^n
(1-i)W_j^i  M_{i}^k \tau^j\partial_{\tau^k}F^{(1)}
+\delta_{1,\mathsf{a}(n,\mathbf{d})} \elld
+F^{(1)}.
\end{eqnarray*}
Applying $\partial_{\tau^0}^{p_0}\partial_{\tau^1}^{p_1}\cdots\partial_{\tau^n}^{p_n}$ to both sides and then evaluating at $\tau=0$, we get (\ref{eq-recursion-F(1)-EulerField}).
\end{proof}

\begin{lemma}\label{lem-recursion-F(1)-WDVV-tau}
For $i\geq 1$ and $I\in \mathbb{Z}_{\geq 0}^{n+1}$, 
\begin{eqnarray}\label{eq-recursion-F(1)-WDVV-tau}
&& \partial_{\tau^I}\partial_{\tau^i}F^{(1)}(0)\nn\\
&=&-\sum_{\begin{subarray}{c}0^{n+1}\leq J\leq I\\J\neq 0^{n+1}\end{subarray}}\binom{I}{J}\sum_{a=0}^n\sum_{b=0}^n
\big(\partial_{\tau^{J}}\partial_{\tau^1} \partial_{\tau^{i-1}}\partial_{\tau^a}F^{(0)}(0)\big)\eta^{ab}\big(\partial_{\tau^{b}}\partial_{\tau^{I-J}}F^{(1)}(0)\big)\nn\\
&&+ \sum_{0^{n+1}\leq  J\leq I}\binom{I}{J}\big(\partial_{\tau^{J}}\partial_{\tau^1} F^{(1)}(0)\big)\big(
\partial_{\tau^{I-J}}\partial_{\tau^{i-1}}F^{(1)}(0)\big).
\end{eqnarray}
\end{lemma}
\begin{proof}
Take the equation (\ref{eq-ss1-tau}) with $a=1$ and $b=i-1$, applying $\partial_{\tau^I}$, and then then evaluating at $\tau=0$, we get (\ref{eq-recursion-F(1)-WDVV-tau}).
\end{proof}

\subsection{Recursion for \texorpdfstring{$F^{(l)}$}{F(l)} for \texorpdfstring{$l\geq 2$}{l>2}}
We regard the unknown
\begin{equation*}
      F^{(l)}(0)
\end{equation*}
as the initial value.

\begin{lemma}\label{lem-recursion-F(l)-EulerField}
For $I=(p_0,\dots,p_n)\in \mathbb{Z}_{\geq 0}^{n+1}$,
\begin{eqnarray}\label{eq-recursion-F(l)-EulerField}
\mathsf{a}(n,\mathbf{d})\partial_{\tau^1}\partial_{\tau^I}F^{(l)}(0)
&=&- \sum_{i=0}^{n}\sum_{j=0}^n\sum_{k=0}^n
(1-i)W_j^i  M_{i}^k p_j \partial_{\tau^{I+\mathbf{e}_k- \mathbf{e}_j}}F^{(l)}(0)\nn\\
&&+(nl-2l-n+3)\partial_{\tau^I}F^{(l)}(0).
\end{eqnarray}
\end{lemma}
\begin{proof}
 For $l\geq 1$,
\[
E^{(0)}F^{(l)}+(2-n)lF^{(l)}=(3-n)F^{(l)},
\]
i.e.
\[
E^{(0)}F^{(l)}=(nl-2l-n+3)F^{(l)}.
\]
When $l\geq 2$,
\begin{eqnarray*}
&&\mathsf{a}(n,\mathbf{d})\partial_{\tau^1}F^{(l)}
      =-\sum_{i=0}^{n}\sum_{j=0}^n\sum_{k=0}^n
(1-i)W_j^i  M_{i}^k \tau^j\partial_{\tau^k}F^{(l)}
+\delta_{1,\mathsf{a}(n,\mathbf{d})} \elld\partial_{\tau^0}F^{(l)}
+(nl-2l-n+3)F^{(1)}\\
&=& -\sum_{i=0}^{n}\sum_{j=0}^n\sum_{k=0}^n
(1-i)W_j^i  M_{i}^k \tau^j\partial_{\tau^k}F^{(l)}
+(nl-2l-n+3)F^{(l)}.
\end{eqnarray*}
Applying $\partial_{\tau^0}^{p_0}\partial_{\tau^1}^{p_1}\cdots\partial_{\tau^n}^{p_n}$ to both sides and then evaluating at $\tau=0$, we get (\ref{eq-recursion-F(l)-EulerField}).
\end{proof}

\begin{lemma}\label{lem-recursion-F(l)-WDVV-tau}
For $i\geq 1$ and $I\in \mathbb{Z}_{\geq 0}^{n+1}$, 
\begin{eqnarray}\label{eq-recursion-F(l)-WDVV-tau}
&& \partial_{\tau^I}\partial_{\tau^i}F^{(l)}(0)\nn\\
&=&-\sum_{\begin{subarray}{c}0^{n+1}\leq J\leq I\\J\neq 0^{n+1}\end{subarray}}\binom{I}{J}\sum_{a=0}^n\sum_{b=0}^n
\big(\partial_{\tau^{J}}\partial_{\tau^1} \partial_{\tau^{i-1}}\partial_{\tau^a}F^{(0)}(0)\big)\eta^{ab}\big(\partial_{\tau^{b}}\partial_{\tau^{I-J}}F^{(l)}(0)\big)\nn\\
&&+\sum_{k=1}^{l}\sum_{0^{n+1}\leq  J\leq I}\binom{l-1}{k-1}\binom{I}{J}
\big(\partial_{\tau^{J}}\partial_{\tau^1} F^{(k)}(0)\big)\big(
\partial_{\tau^{I-J}}\partial_{\tau^{i-1}}F^{(l-k+1)}(0)\big)\nn\\
&&-\sum_{k=1}^{l-1}\sum_{0^{n+1}\leq J\leq I}\binom{l-1}{k}\binom{I}{J}\sum_{a=0}^n\sum_{b=0}^n
\big(\partial_{\tau^{J}}\partial_{\tau^1} \partial_{\tau^{i-1}}\partial_{\tau^a}F^{(k)}(0)\big)\eta^{ab}\big(\partial_{\tau^{b}}\partial_{\tau^{I-J}}F^{(l-k)}(0)\big)\nn\\
&&-2(l-1)\sum_{k=1}^{l-1}\sum_{0^{n+1}\leq  J\leq I}\binom{l-2}{k-1}\binom{I}{J}
\big(\partial_{\tau^{J}}\partial_{\tau^1}\partial_{\tau^{i-1}} F^{(k)}(0)\big)\big(
\partial_{\tau^{I-J}}F^{(l-k+1)}(0)\big).
\end{eqnarray}
\end{lemma}
\begin{proof}
For $l\geq 2$, we rewrite the system of equations (\ref{eq-wdvv23expand}) in the $\tau$-coordinates:
\begin{eqnarray*}
&&\partial_{\tau^a}\partial_{\tau^b}\partial_{\tau^e}F^{(0)}\eta^{ef}\partial_{\tau^f}F^{(l)}\nn\\
&=&\sum_{j=1}^{l}\binom{l-1}{j-1}\partial_{\tau^a}F^{(j)}\partial_{\tau^b}F^{(l-j+1)}-\sum_{j=1}^{l-1}\binom{l-1}{j}\partial_{\tau^a}\partial_{\tau^b}\partial_{\tau^e}F^{(j)}\eta^{ef}\partial_{\tau^f}F^{(l-j)}\nn\\
&&-2(l-1)\sum_{j=1}^{l-1}\binom{l-2}{j-1}\partial_{\tau^a}\partial_{\tau^b}F^{(j)}F^{(l-j+1)}.
\end{eqnarray*}
Take the equation  with $a=1$ and $b=i-1$, applying $\partial_{\tau^I}$, and then evaluating at $\tau=0$, we get (\ref{eq-recursion-F(l)-WDVV-tau}).
\end{proof}

\subsection{Equations of higher order leading terms}
\label{sub:equations_of_higher_order_leading_terms}
As a result of the previous section, we define a function \texttt{correlatorInTauCoord} in our package. Running 
\begin{equation*}
      \mbox{\texttt{correlatorInTauCoord}}\ \{n,\mathbf{d},l,I\}
\end{equation*}
returns 
\begin{equation*}
      \partial_{\tau^I}F^{(l)}
\end{equation*}
in terms of $F^{(k)}(0)$ for $2\leq k\leq l$.
Base on this function we define functions \texttt{equationOfConstTerm} and \texttt{sqrtRecursion} as we mentioned in Section \ref{sec:sqrtRecursion}.

\subsubsection{Cubic hypersurfaces}
The cubic hypersurfaces are  excluded in the double root recursion. 
Their quantum cohomology are reconstructible by Theorem \ref{thm-reconstructcubicandquadric}. We give the explicit recursion as follows.
If $n=3$, $m=10$. We have computed $F^{(2)}(0)$ and $F^{(4)}(0)$ via the Fano variety of lines and reduced genus 1 invariants. The other $F^{(k)}(0)$'s vanish for the dimension reason.
Now suppose $n\geq 4$, and $k\geq 3$. Then
\begin{eqnarray*}
      &&\frac{-3kn+12k+3n-21}{n-1}F^{(k)}(0)\\
&=& -\frac{1}{3}\Big(\sum_{j=2}^{k-1}\binom{k-1}{j-1}\partial_{\tau^1}F^{(j)}(0)\partial_{\tau^{n-1}}F^{(k-j+1)}(0)\\
      &&-\sum_{j=1}^{k-1}\sum_{e=0}^n\sum_{f=0}^n\binom{k-1}{j}\partial_{\tau^1}\partial_{\tau^{n-1}}\partial_{\tau^e}F^{(j)}(0)\eta^{ef}\partial_{\tau^f}F^{(k-j)}(0)\\
&&-2(k-1)\sum_{j=2}^{k-1}\binom{k-2}{j-1}\partial_{\tau^1}\partial_{\tau^{n-1}}F^{(j)}(0)F^{(k-j+1)}(0)\Big)\\
&&-\frac{1}{2}\sum_{j=2}^{k-1}\sum_{e=0}^n\sum_{f=0}^n \binom{k-1}{j-1}\partial_{\tau^e}F^{(j)}(0)\eta^{ef}\partial_{\tau^f}F^{(k+1-j)}(0)\\
&&-(k-1)\sum_{j=3}^{k-1}\binom{k-2}{j-2}F^{(j)}(0)F^{(k+2-j)}(0).
\end{eqnarray*}

\section{Proof of Theorem \ref{thm-F20-a(n,d)=(n-1)/2}}\label{sec:proof-F2(0)-closedFormula}

This section is devoted to showing the following formula. 
\begin{theorem}\label{thm-F20-WM-a(n,d)=(n-1)/2}
When $\mathsf{a}(n,\mathbf{d})=\frac{n-1}{2}$,
\begin{equation}\label{eq-thm-F20}
	-\sum_{j=0}^n j M_{j}^{1}W_{n}^{j}
+\sfbd\sum_{j=0}^n 
j M_{j}^{1}W_{n- \mathsf{a}(n,\mathbf{d})}^j
=\frac{(n-1)(\sum_{i=1}^r d_i!)^2}{4}.
\end{equation}
\end{theorem}

Theorem \ref{thm-F20-a(n,d)=(n-1)/2} follows as a consequence of (\ref{eq-higher22}) and (\ref{eq-thm-F20}).

We deal with the cases the Fano index $\mathsf{a}(n,\mathbf{d})=1$ and $\mathsf{a}(n,\mathbf{d})>1$ separately. There are only finitely many tuples $(n,\mathbf{d})$ such that $\mathsf{a}(n,\mathbf{d})=1$
and $n=3$. We list them below and one can check the validity of (\ref{eq-thm-F20})  in these cases.
\begin{example}
$n=3$, $d=4$, $\mathsf{a}(n,\mathbf{d})=1$.

\[
W=\left(
\begin{array}{cccc}
 1 & 0 & 0 & 0 \\
 24 & 1 & 0 & 0 \\
 4464 & 128 & 1 & 0 \\
 1109376 & 31376 & 232 & 1 \\
\end{array}
\right),\
M=\left(
\begin{array}{cccc}
 1 & 0 & 0 & 0 \\
 -24 & 1 & 0 & 0 \\
 -1392 & -128 & 1 & 0 \\
 -33408 & -1680 & -232 & 1 \\
\end{array}
\right).
\]

\end{example}

\begin{example}
$n=3$, $\mathbf{d}=(2,3)$, $\mathsf{a}(n,\mathbf{d})=1$.

\[
W=\left(
\begin{array}{cccc}
 1 & 0 & 0 & 0 \\
 12 & 1 & 0 & 0 \\
 936 & 54 & 1 & 0 \\
 97632 & 5544 & 96 & 1 \\
\end{array}
\right),\ 
M=\left(
\begin{array}{cccc}
 1 & 0 & 0 & 0 \\
 -12 & 1 & 0 & 0 \\
 -288 & -54 & 1 & 0 \\
 -3456 & -360 & -96 & 1 \\
\end{array}
\right).
\]

\end{example}

\begin{example}
$n=3$, $\mathbf{d}=(2,2,2)$, $\mathsf{a}(n,\mathbf{d})=1$.

\[
W=\left(
\begin{array}{cccc}
 1 & 0 & 0 & 0 \\
 8 & 1 & 0 & 0 \\
 368 & 32 & 1 & 0 \\
 22656 & 1936 & 56 & 1 \\
\end{array}
\right),\
M=\left(
\begin{array}{cccc}
 1 & 0 & 0 & 0 \\
 -8 & 1 & 0 & 0 \\
 -112 & -32 & 1 & 0 \\
 -896 & -144 & -56 & 1 \\
\end{array}
\right).
\]

\end{example}

So we suppose $\frac{n-1}{\mathsf{a}(n,\mathbf{d})}=2$ and $\mathsf{a}(n,\mathbf{d})>1$ from now on.
\subsection{Reduction to the computation of a descendant invariant}
Computing the small quantum powers of $\sfh$ iteratively, one has
\begin{eqnarray*}
\tsfh_i&=&\sfh_i\ \mbox{for}\ 0\leq i\leq \mathsf{a}(n,\mathbf{d})-1,\\
\tsfh_i&=&\sfh_i+\sfh_{i-\mathsf{a}(n,\mathbf{d})}\sum_{j=\mathsf{a}(n,\mathbf{d})-1}^{i-1}\frac{1}{\prod_{i=1}^r d_i}\langle \sfh,\sfh_{j},\sfh_{n-j-1+\mathsf{a}(n,\mathbf{d})}\rangle_{0,1} \
 \mbox{for}\ \mathsf{a}(n,\mathbf{d})\leq i\leq n-2,\\
\tsfh_{n-1}&=&\sfh_{n-1}+ \sfh_{n-1-\mathsf{a}(n,\mathbf{d})}\sum_{j=\mathsf{a}(n,\mathbf{d})-1}^{n-2}\frac{1}{\prod_{i=1}^r d_i}\langle \sfh,\sfh_{j},\sfh_{n-j-1+\mathsf{a}(n,\mathbf{d})}\rangle_{0,1}\\
&&+\frac{1}{\prod_{i=1}^r d_i}\langle \sfh,\sfh_{n-2},\sfh_n\rangle_{0,2}
+\frac{1}{(\prod_{i=1}^r d_i)^2}\langle \sfh,\sfh_{n-2-\mathsf{a}(n,\mathbf{d})},\sfh_{1+2\mathsf{a}(n,\mathbf{d})}\rangle_{0,1}\\
&&\cdot\sum_{j=\mathsf{a}(n,\mathbf{d})-1}^{n-3}\langle \sfh,\sfh_{j},\sfh_{n-j-1+\mathsf{a}(n,\mathbf{d})}\rangle_{0,1},\\
\tsfh_{n}&=& \sfh_{n}+
\sfh_{n-\mathsf{a}(n,\mathbf{d})}\sum_{j=\mathsf{a}(n,\mathbf{d})-1}^{n-1}\frac{1}{\prod_{i=1}^r d_i}\langle \sfh,\sfh_{j},\sfh_{n-j-1+\mathsf{a}(n,\mathbf{d})}\rangle_{0,1}\\
&&+\big(\frac{1}{\prod_{i=1}^r d_i}\langle \sfh,\sfh_{n-1},\sfh_{n-1}\rangle_{0,2}+\frac{1}{\prod_{i=1}^r d_i}\langle \sfh,\sfh_{n-2},\sfh_n\rangle_{0,2}\\
&&+\frac{1}{(\prod_{i=1}^r d_i)^2}\langle \sfh,\sfh_{n-1-\mathsf{a}(n,\mathbf{d})},\sfh_{2\mathsf{a}(n,\mathbf{d})}\rangle_{0,1}
\sum_{j=\mathsf{a}(n,\mathbf{d})-1}^{n-2}\langle \sfh,\sfh_{j},\sfh_{n-j-1+\mathsf{a}(n,\mathbf{d})}\rangle_{0,1}\\
&&+\frac{1}{(\prod_{i=1}^r d_i)^2}\langle \sfh,\sfh_{n-2-\mathsf{a}(n,\mathbf{d})},\sfh_{1+2\mathsf{a}(n,\mathbf{d})}\rangle_{0,1}\sum_{j=\mathsf{a}(n,\mathbf{d})-1}^{n-3}\langle \sfh,\sfh_{j},\sfh_{n-j-1+\mathsf{a}(n,\mathbf{d})}\rangle_{0,1}
\big)\sfh.
\end{eqnarray*}
The inverse transform is then the following:
\begin{eqnarray*}
\sfh_i&=&\tsfh_i\ \mbox{for}\ 0\leq i\leq \mathsf{a}(n,\mathbf{d})-1,\\
\sfh_i&=&\tsfh_i-\tsfh_{i-\mathsf{a}(n,\mathbf{d})}\sum_{j=\mathsf{a}(n,\mathbf{d})-1}^{i-1}\frac{1}{\prod_{i=1}^r d_i}\langle \sfh,\sfh_{j},\sfh_{n-j-1+\mathsf{a}(n,\mathbf{d})}\rangle_{0,1} \
 \mbox{for}\ \mathsf{a}(n,\mathbf{d})\leq i\leq n-2, \\
\sfh_{n-1}&=&\tsfh_{n-1}-\tsfh_{n-1-\mathsf{a}(n,\mathbf{d})} \sum_{j=\mathsf{a}(n,\mathbf{d})-1}^{n-2}\frac{1}{\prod_{i=1}^r d_i}\langle \sfh,\sfh_{j},\sfh_{n-j-1+\mathsf{a}(n,\mathbf{d})}\rangle_{0,1}\\
&&+\big(-\frac{1}{\prod_{i=1}^r d_i}\langle \sfh,\sfh_{n-2},\sfh_n\rangle_{0,2}
+\frac{1}{(\prod_{i=1}^r d_i)^2}\langle \sfh,\sfh_{n-2-\mathsf{a}(n,\mathbf{d})},\sfh_{1+2\mathsf{a}(n,\mathbf{d})}\rangle_{0,1}\langle \sfh,\sfh_{n-2},\sfh_{1+\mathsf{a}(n,\mathbf{d})}\rangle_{0,1}\big),\\
\sfh_n&=&\tsfh_{n}-
\tsfh_{n-\mathsf{a}(n,\mathbf{d})}\sum_{j=\mathsf{a}(n,\mathbf{d})-1}^{n-1}\frac{1}{\prod_{i=1}^r d_i}\langle \sfh,\sfh_{j},\sfh_{n-j-1+\mathsf{a}(n,\mathbf{d})}\rangle_{0,1}\\
&&+\big(-\frac{1}{\prod_{i=1}^r d_i}\langle \sfh,\sfh_{n-1},\sfh_{n-1}\rangle_{0,2}-\frac{1}{\prod_{i=1}^r d_i}\langle \sfh,\sfh_{n-2},\sfh_n\rangle_{0,2}\\
&&+\frac{1}{(\prod_{i=1}^r d_i)^2}\langle \sfh,\sfh_{n-1-\mathsf{a}(n,\mathbf{d})},\sfh_{2\mathsf{a}(n,\mathbf{d})}\rangle_{0,1}
\langle \sfh,\sfh_{n-1},\sfh_{\mathsf{a}(n,\mathbf{d})}\rangle_{0,1}\\
&&+\frac{1}{(\prod_{i=1}^r d_i)^2}\langle \sfh,\sfh_{n-2-\mathsf{a}(n,\mathbf{d})},\sfh_{1+2\mathsf{a}(n,\mathbf{d})}\rangle_{0,1}\sum_{j=n-2}^{n-1}\langle \sfh,\sfh_{j},\sfh_{n-j-1+\mathsf{a}(n,\mathbf{d})}\rangle_{0,1}
\big)\tsfh.
\end{eqnarray*}

\begin{lemma}
Suppose $\mathsf{a}(n,\mathbf{d})\geq 2$. Then
\begin{equation}\label{eq-consequence-mirrorFormula-1}
	\frac{1}{\prod_{i=1}^r d_i}\sum_{j=\mathsf{a}(n,\mathbf{d})-1}^{n}\langle \sfh,\sfh_{j},\sfh_{n-j-1+\mathsf{a}(n,\mathbf{d})}\rangle_{0,1}=\sfbd,
\end{equation}
and
\begin{eqnarray}\label{eq-consequence-mirrorFormula-2}
&&\langle \sfh,\sfh_{n-1},\sfh_{n-1}\rangle_{0,2}+2\langle \sfh,\sfh_{n-2},\sfh_n\rangle_{0,2}
-\frac{2}{\prod_{i=1}^r d_i}\langle \sfh,\sfh_{\mathsf{a}(n,\mathbf{d})+1},\sfh_{n-2}\rangle_{0,1}
\langle \sfh,\sfh_{\mathsf{a}(n,\mathbf{d})-1},\sfh_{n}\rangle_{0,1}\nn\\
&&-\frac{2}{\prod_{i=1}^r d_i}\langle \sfh,\sfh_{\mathsf{a}(n,\mathbf{d})},\sfh_{n-1}\rangle_{0,1}
\langle \sfh,\sfh_{\mathsf{a}(n,\mathbf{d})-1},\sfh_{n}\rangle_{0,1}
-\frac{1}{\prod_{i=1}^r d_i}\langle \sfh,\sfh_{\mathsf{a}(n,\mathbf{d})},\sfh_{n-1}\rangle_{0,1}
\langle \sfh,\sfh_{\mathsf{a}(n,\mathbf{d})},\sfh_{n-1}\rangle_{0,1}\nn\\
&&-\frac{1}{\prod_{i=1}^r d_i}\langle \sfh,\sfh_{\mathsf{a}(n,\mathbf{d})-1},\sfh_{n}\rangle_{0,1}
\langle \sfh,\sfh_{\mathsf{a}(n,\mathbf{d})-1},\sfh_{n}\rangle_{0,1}= 0.
\end{eqnarray}
\end{lemma}
\begin{proof}
Suppose first $\mathsf{a}(n,\mathbf{d})\geq 3$. Then
\begin{eqnarray*}
&& \sfh\sqp \tsfh_n\\
&=& \sfh_{n-\mathsf{a}(n,\mathbf{d})+1}\sum_{j=\mathsf{a}(n,\mathbf{d})-1}^{n}\frac{1}{\prod_{i=1}^r d_i}\langle \sfh,\sfh_{j},\sfh_{n-j-1+\mathsf{a}(n,\mathbf{d})}\rangle_{0,1}\\
&&+\big(\frac{1}{\prod_{i=1}^r d_i}\langle \sfh,\sfh_{n-1},\sfh_{n-1}\rangle_{0,2}
+\frac{2}{\prod_{i=1}^r d_i}\langle \sfh,\sfh_{n-2},\sfh_n\rangle_{0,2}\\
&&+\frac{1}{(\prod_{i=1}^r d_i)^2}\langle \sfh,\sfh_{n-\mathsf{a}(n,\mathbf{d})},\sfh_{2\mathsf{a}(n,\mathbf{d})-1}\rangle_{0,1} \sum_{j=\mathsf{a}(n,\mathbf{d})-1}^{n-1}\langle \sfh,\sfh_{j},\sfh_{n-j-1+\mathsf{a}(n,\mathbf{d})}\rangle_{0,1}\\
&&+\frac{1}{(\prod_{i=1}^r d_i)^2}\langle \sfh,\sfh_{n-1-\mathsf{a}(n,\mathbf{d})},\sfh_{2\mathsf{a}(n,\mathbf{d})}\rangle_{0,1}
\sum_{j=\mathsf{a}(n,\mathbf{d})-1}^{n-2}\langle \sfh,\sfh_{j},\sfh_{n-j-1+\mathsf{a}(n,\mathbf{d})}\rangle_{0,1}\\
&&+\frac{1}{(\prod_{i=1}^r d_i)^2}\langle \sfh,\sfh_{n-2-\mathsf{a}(n,\mathbf{d})},\sfh_{1+2\mathsf{a}(n,\mathbf{d})}\rangle_{0,1}\sum_{j=\mathsf{a}(n,\mathbf{d})-1}^{n-3}\langle \sfh,\sfh_{j},\sfh_{n-j-1+\mathsf{a}(n,\mathbf{d})}\rangle_{0,1}
\big)\sfh_2,
\end{eqnarray*}
and
\begin{eqnarray*}
&&\tsfh_{n-\mathsf{a}(n,\mathbf{d})+1}=\sfh_{n-\mathsf{a}(n,\mathbf{d})+1}
+\sfh_{2}\sum_{j=\mathsf{a}(n,\mathbf{d})-1}^{\mathsf{a}(n,\mathbf{d})+1}\frac{1}{\prod_{i=1}^r d_i}\langle \sfh,\sfh_{j},\sfh_{n-j-1+\mathsf{a}(n,\mathbf{d})}\rangle_{0,1}\\
&=&\sfh_{n-\mathsf{a}(n,\mathbf{d})+1}
+\sfh_{2}\big(\frac{1}{\prod_{i=1}^r d_i}\langle \sfh,\sfh_{\mathsf{a}(n,\mathbf{d})-1},\sfh_{n}\rangle_{0,1}\\
&&+\frac{1}{\prod_{i=1}^r d_i}\langle \sfh,\sfh_{\mathsf{a}(n,\mathbf{d})},\sfh_{n-1}\rangle_{0,1}
+\frac{1}{\prod_{i=1}^r d_i}\langle \sfh,\sfh_{\mathsf{a}(n,\mathbf{d})+1},\sfh_{n-2}\rangle_{0,1}\big).
\end{eqnarray*}
By (\ref{eq-Givental-smallQuantumCohomology}),
\begin{equation}\label{eq-appendix-Givental-smallQuantumCohomology}
\sfh\sqp \tsfh_n=\sfbd\tsfh_{n-\mathsf{a}(n,\mathbf{d})+1}.
\end{equation}
Comparing the coefficients of $\sfh_{n-\mathsf{a}(n,\mathbf{d})+1}$ and $\sfh_2$ in both sides of (\ref{eq-appendix-Givental-smallQuantumCohomology}), 
we get (\ref{eq-consequence-mirrorFormula-1}) and
\begin{eqnarray}\label{eq-consequence-mirrorFormula-2-1}
&&\frac{1}{\prod_{i=1}^r d_i}\langle \sfh,\sfh_{n-1},\sfh_{n-1}\rangle_{0,2}+\frac{2}{\prod_{i=1}^r d_i}\langle \sfh,\sfh_{n-2},\sfh_n\rangle_{0,2}\nn\\
&&+\frac{1}{(\prod_{i=1}^r d_i)^2}\langle \sfh,\sfh_{n-\mathsf{a}(n,\mathbf{d})},\sfh_{2\mathsf{a}(n,\mathbf{d})-1}\rangle_{0,1} \sum_{j=\mathsf{a}(n,\mathbf{d})-1}^{n-1}\langle \sfh,\sfh_{j},\sfh_{n-j-1+\mathsf{a}(n,\mathbf{d})}\rangle_{0,1}\nn\\
&&+\frac{1}{(\prod_{i=1}^r d_i)^2}\langle \sfh,\sfh_{n-1-\mathsf{a}(n,\mathbf{d})},\sfh_{2\mathsf{a}(n,\mathbf{d})}\rangle_{0,1}
\sum_{j=\mathsf{a}(n,\mathbf{d})-1}^{n-2}\langle \sfh,\sfh_{j},\sfh_{n-j-1+\mathsf{a}(n,\mathbf{d})}\rangle_{0,1}\nn\\
&&+\frac{1}{(\prod_{i=1}^r d_i)^2}\langle \sfh,\sfh_{n-2-\mathsf{a}(n,\mathbf{d})},\sfh_{1+2\mathsf{a}(n,\mathbf{d})}\rangle_{0,1}\sum_{j=\mathsf{a}(n,\mathbf{d})-1}^{n-3}\langle \sfh,\sfh_{j},\sfh_{n-j-1+\mathsf{a}(n,\mathbf{d})}\rangle_{0,1}\\
&=& \sfbd\big(\frac{1}{\prod_{i=1}^r d_i}\langle \sfh,\sfh_{\mathsf{a}(n,\mathbf{d})-1},\sfh_{n}\rangle_{0,1}+\frac{1}{\prod_{i=1}^r d_i}\langle \sfh,\sfh_{\mathsf{a}(n,\mathbf{d})},\sfh_{n-1}\rangle_{0,1}
+\frac{1}{\prod_{i=1}^r d_i}\langle \sfh,\sfh_{\mathsf{a}(n,\mathbf{d})+1},\sfh_{n-2}\rangle_{0,1}\big).\nn
\end{eqnarray}
Applying (\ref{eq-consequence-mirrorFormula-1}) to the sums in (\ref{eq-consequence-mirrorFormula-2-1}) we get  (\ref{eq-consequence-mirrorFormula-2}).

Similar computations show that (\ref{eq-consequence-mirrorFormula-1}) and (\ref{eq-consequence-mirrorFormula-2}) hold also in the case $\mathsf{a}(n,\mathbf{d})=2$; in fact the additional nonzero constant  term  does not affect the proof. 
\end{proof}

\begin{lemma}
Suppose $\mathsf{a}(n,\mathbf{d})\geq 2$. Then
\begin{eqnarray}\label{eq-sumOfProductsOfMatrices-simplification-1}
&& -\sum_{j=0}^n j M_{j}^{1}W_{n}^{j}
+\sfbd\sum_{j=0}^n 
j M_{j}^{1}W_{n- \mathsf{a}(n,\mathbf{d})}^j\nn\\
&=&(n-1)\big(-\frac{2}{\prod_{i=1}^r d_i}\langle \sfh_{n-2},\sfh_n\rangle_{0,2}+\frac{1}{(\prod_{i=1}^r d_i)^2}\langle \sfh_{\mathsf{a}(n,\mathbf{d})+1},\sfh_{n-2}\rangle_{0,1}
\langle \sfh_{\mathsf{a}(n,\mathbf{d})-1},\sfh_{n}\rangle_{0,1}\nn\\
&&+\frac{1}{2(\prod_{i=1}^r d_i)^2}\langle \sfh_{\mathsf{a}(n,\mathbf{d})},\sfh_{n-1}\rangle_{0,1}
\langle \sfh_{\mathsf{a}(n,\mathbf{d})-1},\sfh_{n}\rangle_{0,1}
+\frac{1}{2(\prod_{i=1}^r d_i)^2}\langle\sfh_{\mathsf{a}(n,\mathbf{d})-1},\sfh_{n}\rangle_{0,1}
\langle \sfh_{\mathsf{a}(n,\mathbf{d})-1},\sfh_{n}\rangle_{0,1}\big).\nn\\
\end{eqnarray}
\end{lemma}
\begin{proof}
By the transformation matrices between $\sfh_i$'s and $\tsfh_i$'s, 
\begin{eqnarray}\label{eq-sumOfProductsOfMatrices-simplification-0}
&& -\sum_{j=0}^n j M_{j}^{1}W_{n}^{j}
+\sfbd\sum_{j=0}^n 
j M_{j}^{1}W_{n- \mathsf{a}(n,\mathbf{d})}^j\nn\\
&=&-W_n^1-(n- \mathsf{a}(n,\mathbf{d}))M_{n- \mathsf{a}(n,\mathbf{d})}^1 W_n^{n- \mathsf{a}(n,\mathbf{d})}-nM_n^1\nn\\
&&+\sfbd W_{n-\mathsf{a}(n,\mathbf{d})}^1
+\sfbd (n- \mathsf{a}(n,\mathbf{d}))M_{n- \mathsf{a}(n,\mathbf{d})}^1\nn\\
&=& -\big(\frac{1}{\prod_{i=1}^r d_i}\langle \sfh,\sfh_{n-1},\sfh_{n-1}\rangle_{0,2}+\frac{1}{\prod_{i=1}^r d_i}\langle \sfh,\sfh_{n-2},\sfh_n\rangle_{0,2}\nn\\
&&+\frac{1}{(\prod_{i=1}^r d_i)^2}\langle \sfh,\sfh_{n-1-\mathsf{a}(n,\mathbf{d})},\sfh_{2\mathsf{a}(n,\mathbf{d})}\rangle_{0,1}
\sum_{j=\mathsf{a}(n,\mathbf{d})-1}^{n-2}\langle \sfh,\sfh_{j},\sfh_{n-j-1+\mathsf{a}(n,\mathbf{d})}\rangle_{0,1}\nn\\
&&+\frac{1}{(\prod_{i=1}^r d_i)^2}\langle \sfh,\sfh_{n-2-\mathsf{a}(n,\mathbf{d})},\sfh_{1+2\mathsf{a}(n,\mathbf{d})}\rangle_{0,1}\sum_{j=\mathsf{a}(n,\mathbf{d})-1}^{n-3}\langle \sfh,\sfh_{j},\sfh_{n-j-1+\mathsf{a}(n,\mathbf{d})}\rangle_{0,1}
\big)\nn\\
&&-(n- \mathsf{a}(n,\mathbf{d}))\big(-\sum_{j=\mathsf{a}(n,\mathbf{d})-1}^{\mathsf{a}(n,\mathbf{d})}\frac{1}{\prod_{i=1}^r d_i}\langle \sfh,\sfh_{j},\sfh_{n-j-1+\mathsf{a}(n,\mathbf{d})}\rangle_{0,1}\big)\sum_{j=\mathsf{a}(n,\mathbf{d})-1}^{n-1}\frac{1}{\prod_{i=1}^r d_i}\langle \sfh,\sfh_{j},\sfh_{n-j-1+\mathsf{a}(n,\mathbf{d})}\rangle_{0,1}\nn\\
&&-n\big(-\frac{1}{\prod_{i=1}^r d_i}\langle \sfh,\sfh_{n-1},\sfh_{n-1}\rangle_{0,2}-\frac{1}{\prod_{i=1}^r d_i}\langle \sfh,\sfh_{n-2},\sfh_n\rangle_{0,2}\nn\\
&&+\frac{1}{(\prod_{i=1}^r d_i)^2}\langle \sfh,\sfh_{n-1-\mathsf{a}(n,\mathbf{d})},\sfh_{2\mathsf{a}(n,\mathbf{d})}\rangle_{0,1}
\langle \sfh,\sfh_{n-1},\sfh_{\mathsf{a}(n,\mathbf{d})}\rangle_{0,1}\nn\\
&&+\frac{1}{(\prod_{i=1}^r d_i)^2}\langle \sfh,\sfh_{n-2-\mathsf{a}(n,\mathbf{d})},\sfh_{1+2\mathsf{a}(n,\mathbf{d})}\rangle_{0,1}\sum_{j=n-2}^{n-1}\langle \sfh,\sfh_{j},\sfh_{n-j-1+\mathsf{a}(n,\mathbf{d})}\rangle_{0,1}
\big)\nn\\
&&+\sfbd\sum_{j=\mathsf{a}(n,\mathbf{d})-1}^{\mathsf{a}(n,\mathbf{d})}\frac{1}{\prod_{i=1}^r d_i}\langle \sfh,\sfh_{j},\sfh_{n-j-1+\mathsf{a}(n,\mathbf{d})}\rangle_{0,1}\nn\\
&&+\sfbd (n- \mathsf{a}(n,\mathbf{d}))\big(-\sum_{j=\mathsf{a}(n,\mathbf{d})-1}^{\mathsf{a}(n,\mathbf{d})}\frac{1}{\prod_{i=1}^r d_i}\langle \sfh,\sfh_{j},\sfh_{n-j-1+\mathsf{a}(n,\mathbf{d})}\rangle_{0,1}\big).
\end{eqnarray}
Applying (\ref{eq-consequence-mirrorFormula-1}), after some manipulations we get
\begin{eqnarray*}
&& -\sum_{j=0}^n j M_{j}^{1}W_{n}^{j}
+\sfbd\sum_{j=0}^n 
j M_{j}^{1}W_{n- \mathsf{a}(n,\mathbf{d})}^j\\
&=&(n-1)\big(\frac{2}{\prod_{i=1}^r d_i}\langle \sfh_{n-1},\sfh_{n-1}\rangle_{0,2}
+\frac{2}{\prod_{i=1}^r d_i}\langle \sfh_{n-2},\sfh_n\rangle_{0,2}\\
&&-\frac{1}{(\prod_{i=1}^r d_i)^2}\langle \sfh_{\mathsf{a}(n,\mathbf{d})},\sfh_{n-1}\rangle_{0,1}
\langle \sfh_{\mathsf{a}(n,\mathbf{d})},\sfh_{n-1}\rangle_{0,1}
-\frac{3}{2(\prod_{i=1}^r d_i)^2}\langle \sfh_{\mathsf{a}(n,\mathbf{d})},\sfh_{n-1}\rangle_{0,1}
\langle \sfh_{\mathsf{a}(n,\mathbf{d})-1},\sfh_{n}\rangle_{0,1}\\
&&-\frac{1}{2(\prod_{i=1}^r d_i)^2}\langle \sfh_{\mathsf{a}(n,\mathbf{d})-1},\sfh_{n}\rangle_{0,1}
\langle \sfh_{\mathsf{a}(n,\mathbf{d})-1},\sfh_{n}\rangle_{0,1}
-\frac{1}{(\prod_{i=1}^r d_i)^2}\langle \sfh_{\mathsf{a}(n,\mathbf{d})+1},\sfh_{n-2}\rangle_{0,1}
\langle \sfh_{\mathsf{a}(n,\mathbf{d})-1},\sfh_{n}\rangle_{0,1}
\big).
\end{eqnarray*}
Then we use (\ref{eq-consequence-mirrorFormula-2}) to eliminate $\langle \sfh_{n-1},\sfh_{n-1}\rangle_{0,2}$ and obtain (\ref{eq-sumOfProductsOfMatrices-simplification-1}).
\end{proof}

With these preparations we can prove Theorem \ref{thm-F20-WM-a(n,d)=(n-1)/2} in the case $\mathsf{a}(n,\mathbf{d})\geq 2$. 
For brevity of notations we only give  a proof in the case $r=1$, and leave the general case to the reader. So $n$ is an 
odd number $\geq 3$, and 
\[
d=\frac{n+5}{2},\ \mathsf{a}(n,\mathbf{d})=\frac{n-1}{2}.
\]
\begin{lemma}\label{lem-twistedInvariant-descendant-localization-result}
When $r=1$, and $n\geq 5$,
Theorem \ref{thm-F20-WM-a(n,d)=(n-1)/2} is equivalent to
\begin{eqnarray}\label{eq-twistedInvariant-descendant-localization-result}
 \langle \sfh_{n}\psi,\sfh_{n-3}\rangle_{0,2}
&=& \frac{d(d!)^2}{2} \Big(\frac{1}{8}+\frac{1}{4}\sum_{i=1}^{d-1}\frac{d-i}{i}+ \frac{1}{2}\sum_{1\leq i<j\leq d-1}\frac{(d-i)(d-j)}{ij}\nn\\
&&	+ \sum_{1\leq i<j<k\leq d-1}\frac{(d-i)(d-j)(d-k)}{ijk}\Big).
\end{eqnarray}
\end{lemma}
\begin{proof}
By (\ref{eq-degOne-invariant}),
\begin{eqnarray}\label{eq-degOne-n-(n-1)-(n-2)}
&&\langle\sfh_{\mathsf{a}(n,\mathbf{d})-1},\sfh_{n}\rangle_{0,1}= d\cdot d!,\\
&&\langle\sfh_{\mathsf{a}(n,\mathbf{d})},\sfh_{n-1}\rangle_{0,1}
	=d^2\cdot d!\big(\mu(d,1)-1	\big)	
	=d^2\cdot d!\sum_{k=2}^d \frac{1}{k},\nn\\
&&\langle\sfh_{\mathsf{a}(n,\mathbf{d})+1},\sfh_{n-2}\rangle_{0,1}=
	d\cdot d! \Big( \frac{d^2-d}{2}-(d^2-d) \mu(d,1)+d^2 \mu(d,2)\Big).\nn
\end{eqnarray}
So
\begin{eqnarray*}
&&\frac{1}{d^2}\langle \sfh_{\mathsf{a}(n,\mathbf{d})+1},\sfh_{n-2}\rangle_{0,1}
\langle \sfh_{\mathsf{a}(n,\mathbf{d})-1},\sfh_{n}\rangle_{0,1}\\
&&+\frac{1}{2d^2}\langle \sfh_{\mathsf{a}(n,\mathbf{d})},\sfh_{n-1}\rangle_{0,1}
\langle \sfh_{\mathsf{a}(n,\mathbf{d})-1},\sfh_{n}\rangle_{0,1}
+\frac{1}{2d^2}\langle\sfh_{\mathsf{a}(n,\mathbf{d})-1},\sfh_{n}\rangle_{0,1}
\langle \sfh_{\mathsf{a}(n,\mathbf{d})-1},\sfh_{n}\rangle_{0,1}\\
&=&  (d!)^2 \big( \frac{(d-1)^2}{2}+(\frac{3d}{2}-d^2) \mu(d,1)+d^2 \mu(d,2)\big).
\end{eqnarray*}
Thus from (\ref{eq-sumOfProductsOfMatrices-simplification-1})
it follows that (\ref{eq-thm-F20}) is equivalent to
\begin{eqnarray}\label{eq-thm-(n-2,n)-1}
\langle \sfh_{n-2},\sfh_n\rangle_{0,2}
= d(d!)^2\big( \frac{2d^2-4d+1}{8}+\frac{3d-2d^2}{4}\mu(d,1)+\frac{d^2}{2}\mu(d,2)\big).
\end{eqnarray}

By (\ref{eq-TRR0}) and (\ref{eq-Div}),
\begin{eqnarray*}
&&\langle \sfh_n \psi,\sfh_{n-3},\sfh\rangle_{0,2}\\
&=&\langle \sfh_n,\sfh_{n-2}\rangle_{0,2}\frac{1}{d}
\langle \sfh_2,\sfh_{n-3},\sfh\rangle_{0,0}
+\langle \sfh_n,\sfh_{\mathsf{a}(n,d)-1}\rangle_{0,1}\frac{1}{d}
\langle \sfh_{\mathsf{a}(n,d)+2},\sfh_{n-3},\sfh\rangle_{0,1}\\
&=&\langle \sfh_n,\sfh_{n-2}\rangle_{0,2}
+\langle \sfh_n,\sfh_{\mathsf{a}(n,d)-1}\rangle_{0,1}\frac{1}{d}
\langle \sfh_{\mathsf{a}(n,d)+2},\sfh_{n-3}\rangle_{0,1}.
\end{eqnarray*}
On the other hand by (\ref{eq-Div}),
\begin{eqnarray*}
\langle \sfh_n \psi,\sfh_{n-3},\sfh\rangle_{0,2}=
2\langle \sfh_n \psi,\sfh_{n-3}\rangle_{0,2}.
\end{eqnarray*}
So
\begin{equation}\label{eq-degOne-(n,n-2)to(n-3)}
\langle \sfh_n,\sfh_{n-2}\rangle_{0,2}
=	2\langle \sfh_n \psi,\sfh_{n-3}\rangle_{0,2}
 -\frac{1}{d}\langle \sfh_n,\sfh_{\mathsf{a}(n,d)-1}\rangle_{0,1}
\langle \sfh_{\mathsf{a}(n,d)+2},\sfh_{n-3}\rangle_{0,1}.
\end{equation}
By (\ref{eq-degOne-invariant}) again,
\begin{eqnarray}\label{eq-degOne-(n-3)}
\langle\sfh_{\mathsf{a}(n,\mathbf{d})+2},\sfh_{n-3}\rangle_{0,1}
&=& d\cdot d! \Big( \frac{1}{6} \left(-d^3+3 d^2-2 d\right)
	+\frac{d}{2} \left(d^2-3 d+2\right)\mu(d,1)\nn\\
&&	+d^2(2-d)\mu(d,2)
	+d^3 \mu(d,3)
	\Big).	
\end{eqnarray}
Then the equivalence of (\ref{eq-thm-(n-2,n)-1}) and (\ref{eq-twistedInvariant-descendant-localization-result}) follows from (\ref{eq-degOne-n-(n-1)-(n-2)}),
(\ref{eq-degOne-(n,n-2)to(n-3)}), (\ref{eq-degOne-(n-3)}) and the following lemma.
\end{proof}	

\begin{lemma}\label{lem-combinatorialIdentity-1}
\begin{equation}\label{eq-combinatorialIdentity-1}
	\sum_{1\leq i_1<\dots<i_k\leq d}\frac{(d-i_1)\cdots(d-i_k)}{i_1\cdots i_k}
	=\sum_{j=0}^{k}(-1)^j d^{k-j}\binom{d-k+j}{j}\mu(d,k-j).
\end{equation}
\end{lemma}
\begin{proof}
We expand the left handside as polynomials of $d$,
\[
\sum_{j=0}^{k}(-1)^j d^{k-j}\sum_{1\leq i_1<\dots<i_{k-j}\leq d}\frac{C_j}{i_1 \cdots i_{k-j}},
\]
where $C_j$ is the number of ways to enlarge the chain $i_1<\dots<i_{k-j}$ to a chain of the form $i'_1<\dots i'_{k}$. Then $C_j$ is equal to the number of ways to choose $j$ numbers in the complement $\{1,\dots,d\}\backslash\{i_1,\dots,i_k\}$. So (\ref{eq-combinatorialIdentity-1}) follows.
\end{proof}

In principle, for a given $n$, one can use (\ref{eq-derivative-J}) to compute $\langle \sfh_{n-2},\sfh_n\rangle_{0,2}$ or $\langle \sfh_n \psi,\sfh_{n-3}\rangle_{0,2}$, and thus prove Theorem \ref{thm-F20-WM-a(n,d)=(n-1)/2}. But it seems hard to obtain a closed formula in this way. In the following section we compute $\langle \sfh_n \psi,\sfh_{n-3}\rangle_{0,2}$ directly by the virtual torus localization. One can also compute $\langle \sfh_{n-2},\sfh_n\rangle_{0,2}$ by localization and obtains (\ref{eq-thm-(n-2,n)-1}). But it turns out that the appearance of a $\psi$ class makes the contributions of  several types of graphs vanish, and thus greatly simplifies the summation.

\subsection{A localization computation with descendants}
Let $n$ be an odd number $\geq 5$, and $d=\frac{n+5}{2}$. Let $X\subset \mathbb{P}^{n+1}$ be a smooth hypersurface of degree $d$.  As we have seen in Lemma \ref{lem-twistedInvariant-descendant-localization-result} we need to compute $\langle \sfh_{n}\psi,\sfh_{n-3}\rangle_{0,2}^X$. 
By the relation of the virtual fundamental classes $[\Mbar_{0,k}(X,\beta)]^{\mathrm{vir}}$ and $[\Mbar_{0,k}(\mathbb{P}^{n+1},\beta)]^{\mathrm{vir}}$ (\cite[P. 181]{CK99}), we will compute
\begin{equation}\label{eq-twistedInvariant-2}
	\langle \sfh_{n}\psi,\sfh_{n-3}\rangle_{0,2}^X=\int_{[\Mbar_{0,2}(\mathbb{P}^{n+1},2)]^{\mathrm{vir}}} \psi_1 c_{n}\big(\mathrm{ev}_1^*\mathcal{O}(1))c_{n-3}\big(\mathrm{ev}_2^*\mathcal{O}(1)) e(R^0 \pi_* f^* \mathcal{O}(d)).
\end{equation}

\subsubsection{Contributions of graphs}
Let $\mathbb{G}_m^{n+2}$ act on $\mathbb{P}^{n+1}$ with fixed points $P_i$, $0\leq i\leq n+1$, such that the tangent weights at $P_i$ are $\alpha_i- \alpha_k$ for $k\in \{0,\dots,n+1\}\backslash\{i\}$. We linearize $\mathcal{O}(l)$ such that it has weight $l \alpha_i$ at $P_i$. We use torus localization to compute (\ref{eq-twistedInvariant-2}). 

For the virtual localization on  $\Mbar_{0,2}(\mathbb{P}^{n+1},2)$ we  follow the presentation of \cite[Chap. 9]{CK99}. 
The contributions of the $\psi$-class are highlighted in \textcolor{blue}{blue}. Since 
\[
\int_{\Mbar_{0,3}}\psi=0
\]
some contributions vanish.

$$ \Gamma_{i12,j,k}=\xy
(0,0); (10,0), **@{-};(20,0), **@{-};
(0,0);(0,-3),**@{.};(0,0);(-3,0),**@{.};
(-4,0)*+{1};(0,-4)*+{2};
(0,0)*+{\bullet};(10,0)*+{\bullet};(20,0)*+{\bullet};(5,2);
(0,3)*+{i};(10,3)*+{j};(20,3)*+{k};
\endxy,
$$
The contribution is $0$.

$$ \Gamma_{i2,j1,k}=\xy
(0,0); (10,0), **@{-};(20,0), **@{-};
(0,0);(0,-3),**@{.};(0,-4)*+{2};
(10,0);(10,-3),**@{.};(10,-4)*+{1};
(0,0)*+{\bullet};(10,0)*+{\bullet};(20,0)*+{\bullet};(5,2);
(0,3)*+{i};(10,3)*+{j};(20,3)*+{k};
\endxy,
$$
The contribution is $0$.

$$ \Gamma_{i1,j2,k}=\xy
(0,0); (10,0), **@{-};(20,0), **@{-};
(0,0);(0,-3),**@{.};(0,-4)*+{1};
(10,0);(10,-3),**@{.};(10,-4)*+{2};
(0,0)*+{\bullet};(10,0)*+{\bullet};(20,0)*+{\bullet};(5,2);
(0,3)*+{i};(10,3)*+{j};(20,3)*+{k};
\endxy,\quad \mbox{with } i,k\neq j,
$$

\begin{eqnarray}\label{eq-localization-Contribution-1}
	&&\int_{\Mbar_{0,3}}\Big(\textcolor{blue}{(\alpha_j- \alpha_i)}\alpha_i^{n} \alpha_j^{n-3} (d\alpha_j)\prod_{p=0}^{d-1}\big(p\alpha_j+(d-p) \alpha_i\big)
	\prod_{p=0}^{d-1}\big(p\alpha_j+(d-p) \alpha_k\big)\nn\\
	&&\cdot\frac{\prod_{p\neq i}(\alpha_i- \alpha_p)
	\cdot \big(\prod_{p\neq j}(\alpha_j- \alpha_p)\big)^2 
	\cdot \prod_{p\neq k}(\alpha_k- \alpha_p)}
	{(\alpha_j- \alpha_i-\psi)(\alpha_j- \alpha_k-\psi)}\nn\\
	&&\cdot	\frac{1}{\prod_{p\neq i}(\alpha_i- \alpha_p)\cdot
	\prod_{q\neq j}(\alpha_j- \alpha_q)
	\cdot	\prod_{r\neq k,j}(\alpha_k- \alpha_r)}\nn\\
	&&\cdot\frac{1}{(\alpha_i- \alpha_j)^2 \prod_{p\neq i,j}( \alpha_i- \alpha_p)(\alpha_j- \alpha_p)\cdot
	(\alpha_j- \alpha_k)^2 \prod_{p\neq j,k}( \alpha_k- \alpha_p)(\alpha_j- \alpha_p)}\Big)\\
	&=& (\alpha_j- \alpha_i)\alpha_i^{n} \alpha_j^{n-3} (d\alpha_j)\prod_{p=0}^{d-1}\big(p\alpha_j+(d-p) \alpha_i\big)
	\prod_{p=0}^{d-1}\big(p\alpha_j+(d-p) \alpha_k\big)\nn\\
	&&\cdot\frac{\prod_{p\neq i}(\alpha_i- \alpha_p)
	\cdot \big(\prod_{p\neq j}(\alpha_j- \alpha_p)\big)^2 
	\cdot \prod_{p\neq k}(\alpha_k- \alpha_p)}
	{(\alpha_j- \alpha_i)(\alpha_j- \alpha_k)}\nn\\
	&&\cdot	\frac{1}{\prod_{p\neq i}(\alpha_i- \alpha_p)\cdot
	\prod_{q\neq j}(\alpha_j- \alpha_q)
	\cdot	\prod_{r\neq k,j}(\alpha_k- \alpha_r)}\nn\\
	&&\cdot\frac{1}{(\alpha_i- \alpha_j)^2 \prod_{p\neq i,j}( \alpha_i- \alpha_p)(\alpha_j- \alpha_p)\cdot
	(\alpha_j- \alpha_k)^2 \prod_{p\neq j,k}( \alpha_k- \alpha_p)(\alpha_j- \alpha_p)}\nn\\
	&=&-\frac{\alpha_i^{n} \alpha_j^{n-3} (d\alpha_j)\prod_{p=0}^{d-1}\big(p\alpha_j+(d-p) \alpha_i\big)
	\prod_{p=0}^{d-1}\big(p\alpha_j+(d-p) \alpha_k\big)}
	{\prod_{p\neq i}(\alpha_i- \alpha_p)\cdot
	\prod_{q\neq j}(\alpha_j- \alpha_q)
	\cdot	\prod_{r\neq k}(\alpha_k- \alpha_r)}.
\end{eqnarray}

\begin{equation}\label{eq-graph-localization-2}
 \Gamma_{i,j12,k}=\xy
(0,0); (10,0), **@{-};(20,0), **@{-};
(10,0);(8,-3),**@{.};(8,-4)*+{1};
(10,0);(12,-3),**@{.};(12,-4)*+{2};
(0,0)*+{\bullet};(10,0)*+{\bullet};(20,0)*+{\bullet};(5,2);
(0,3)*+{i};(10,3)*+{j};(20,3)*+{k};
\endxy,\quad \mbox{with } i,k\neq j,
\end{equation}

\begin{eqnarray}\label{eq-localization-Contribution-2}
	&&\int_{\Mbar_{0,4}}\Big(\textcolor{blue}{\psi_3}\alpha_j^{2n-3} (d\alpha_j)\prod_{p=0}^{d-1}\big(p\alpha_j+(d-p) \alpha_i\big)
	\prod_{p=0}^{d-1}\big(p\alpha_j+(d-p) \alpha_k\big)\nn\\
	&&\cdot\frac{\prod_{p\neq i}(\alpha_i- \alpha_p)
	\cdot \big(\prod_{p\neq j}(\alpha_j- \alpha_p)\big)^2 
	\cdot \prod_{p\neq k}(\alpha_k- \alpha_p)}
	{(\alpha_j- \alpha_i-\psi_1)(\alpha_j- \alpha_k-\psi_2)}\nn\\
	&&\cdot	\frac{1}{\prod_{p\neq i,j}(\alpha_i- \alpha_p)\cdot
	\prod_{q\neq j}(\alpha_j- \alpha_q)
	\cdot	\prod_{r\neq k,j}(\alpha_k- \alpha_r)}\nn\\
	&&\cdot\frac{1}{(\alpha_i- \alpha_j)^2 \prod_{p\neq i,j}( \alpha_i- \alpha_p)(\alpha_j- \alpha_p)\cdot
	(\alpha_j- \alpha_k)^2 \prod_{p\neq j,k}( \alpha_k- \alpha_p)(\alpha_j- \alpha_p)}\Big)\nn\\
	&=& \alpha_j^{2n-3} (d\alpha_j)\prod_{p=0}^{d-1}\big(p\alpha_j+(d-p) \alpha_i\big)
	\prod_{p=0}^{d-1}\big(p\alpha_j+(d-p) \alpha_k\big)\nn\\
	&&\cdot\frac{\prod_{p\neq i}(\alpha_i- \alpha_p)
	\cdot \big(\prod_{p\neq j}(\alpha_j- \alpha_p)\big)^2 
	\cdot \prod_{p\neq k}(\alpha_k- \alpha_p)}
	{(\alpha_j- \alpha_i)(\alpha_j- \alpha_k)}\nn\\
	&&\cdot	\frac{1}{\prod_{p\neq i,j}(\alpha_i- \alpha_p)\cdot
	\prod_{q\neq j}(\alpha_j- \alpha_q)
	\cdot	\prod_{r\neq k,j}(\alpha_k- \alpha_r)}\nn\\
	&&\frac{1}{(\alpha_i- \alpha_j)^2 \prod_{p\neq i,j}( \alpha_i- \alpha_p)(\alpha_j- \alpha_p)\cdot
	(\alpha_j- \alpha_k)^2 \prod_{p\neq j,k}( \alpha_k- \alpha_p)(\alpha_j- \alpha_p)}\nn\\
	&=& \frac{\alpha_j^{2n-3} (d\alpha_j)\prod_{p=0}^{d-1}\big(p\alpha_j+(d-p) \alpha_i\big)
	\prod_{p=0}^{d-1}\big(p\alpha_j+(d-p) \alpha_k\big)}{\prod_{p\neq i}(\alpha_i- \alpha_p)\cdot
	\prod_{q\neq j}(\alpha_j- \alpha_q)
	\cdot	\prod_{r\neq k}(\alpha_k- \alpha_r)}.
\end{eqnarray}

$$ \Gamma_{i1,j,k2}=\xy
(0,0); (10,0), **@{-};(20,0), **@{-};
(0,0);(0,-3),**@{.};(0,-4)*+{1};
(20,0);(20,-3),**@{.};(20,-4)*+{2};
(0,0)*+{\bullet};(10,0)*+{\bullet};(20,0)*+{\bullet};(5,2);
(0,3)*+{i};(10,3)*+{j};(20,3)*+{k};
\endxy,\quad \mbox{with } i,k\neq j,
$$

\begin{eqnarray}\label{eq-localization-Contribution-3}
	&&\textcolor{blue}{(\alpha_j- \alpha_i)} \alpha_i^{n} \alpha_k^{n-3} (d\alpha_j)\prod_{p=0}^{d-1}\big(p\alpha_j+(d-p) \alpha_i\big)
	\prod_{p=0}^{d-1}\big(p\alpha_j+(d-p) \alpha_k\big)\nn\\
	&&\cdot\prod_{p\neq i}(\alpha_i- \alpha_p)
	\cdot \big(\prod_{p\neq j}(\alpha_j- \alpha_p)\big)^2 
	\cdot \prod_{p\neq k}(\alpha_k- \alpha_p)\nn\\
	&&\cdot	\frac{1}{\prod_{p\neq i}(\alpha_i- \alpha_p)\cdot
	(2 \alpha_j- \alpha_i- \alpha_k)
	\prod_{q\neq j}(\alpha_j- \alpha_q)
	\cdot	\prod_{r\neq k}(\alpha_k- \alpha_r)}\nn\\
	&&\cdot\frac{1}{(\alpha_i- \alpha_j)^2 \prod_{p\neq i,j}( \alpha_i- \alpha_p)(\alpha_j- \alpha_p)\cdot
	(\alpha_j- \alpha_k)^2 \prod_{p\neq j,k}( \alpha_k- \alpha_p)(\alpha_j- \alpha_p)}\nn\\
	&=& 	\frac{(\alpha_j- \alpha_i) \alpha_i^{n} \alpha_k^{n-3} (d\alpha_j)\prod_{p=0}^{d-1}\big(p\alpha_j+(d-p) \alpha_i\big)
	\prod_{p=0}^{d-1}\big(p\alpha_j+(d-p) \alpha_k\big)}
	{(2 \alpha_j- \alpha_i- \alpha_k)\prod_{p\neq i}(\alpha_i- \alpha_p)\cdot
	\prod_{q\neq j}(\alpha_j- \alpha_q)
	\cdot	\prod_{r\neq k}(\alpha_k- \alpha_r)}.
\end{eqnarray}

$$ \Gamma_{i12,j}=\xy
(0,0); (10,0), **@{-};
(0,0);(0,-3),**@{.};(0,0);(-3,0),**@{.};
(-4,0)*+{1};(0,-4)*+{2};
(0,0)*+{\bullet};(10,0)*+{\bullet};
(0,3)*+{i};(10,3)*+{j};
\endxy,
$$
The contribution is $0$.

\begin{equation}\label{eq-graph-localization-4}
 	\Gamma_{i1,j2}=\xy
(0,0); (10,0), **@{-};
(0,0);(0,-3),**@{.};(0,-4)*+{1};
(10,0);(10,-3),**@{.};(10,-4)*+{2};
(0,0)*+{\bullet};(10,0)*+{\bullet};
(0,3)*+{i};(10,3)*+{j};
\endxy,\quad \mbox{with } i\neq j,
 \end{equation} 

\begin{eqnarray}\label{eq-localization-Contribution-4}
	&& \textcolor{blue}{\frac{\alpha_j -\alpha_i}{2}} \alpha_i^{n} \alpha_j^{n-3} 
	\prod_{p=0}^{2d}\frac{p\alpha_j+(2d-p) \alpha_i}{2}\nn\\
	&&\cdot\prod_{p\neq i}(\alpha_i- \alpha_p)
	\cdot \prod_{p\neq j}(\alpha_j- \alpha_p)\nn\\
	&&\cdot	\frac{1}{\prod_{p\neq i}(\alpha_i- \alpha_p)\cdot
	\prod_{q\neq j}(\alpha_j- \alpha_q)}\nn\\
	&&\cdot\frac{1}{\frac{(\alpha_i- \alpha_j)^4}{4} 
	\prod_{p\neq i,j}( \alpha_i- \alpha_p)(\frac{\alpha_i+\alpha_j}{2}-\alpha_p)(\alpha_j- \alpha_p)}\nn\\
	&=& -\frac{\alpha_i^{n} \alpha_j^{n-3} 
	\prod_{p=0}^{2d}\big(p\alpha_j+(2d-p) \alpha_i\big)}
	{32(\alpha_i- \alpha_j)^3 
	\prod_{p\neq i,j}( \alpha_i- \alpha_p)(\alpha_i+\alpha_j-2\alpha_p)(\alpha_j- \alpha_p)}.
\end{eqnarray}

Now we sum (\ref{eq-localization-Contribution-1}), (\ref{eq-localization-Contribution-2}), (\ref{eq-localization-Contribution-3}) and (\ref{eq-localization-Contribution-4}). Note that the graph (\ref{eq-graph-localization-4})
has an automorphism of order 2, and so is the graph (\ref{eq-graph-localization-2}) when $k=i$. 
We get
\begin{eqnarray}\label{eq-twistedInvariant-descendant-localization-sum-1}
&& \langle \sfh_{n}\psi,\sfh_{n-3}\rangle_{0,2}\nn\\
&=& \sum_{j=0}^{n+1}\sum_{i\neq j}\sum_{k\neq j}\Big(
	-\frac{\alpha_i^{n} \alpha_j^{n-3} (d\alpha_j)\prod_{p=0}^{d-1}\big(p\alpha_j+(d-p) \alpha_i\big)
	\prod_{p=0}^{d-1}\big(p\alpha_j+(d-p) \alpha_k\big)}
	{\prod_{p\neq i}(\alpha_i- \alpha_p)\cdot
	\prod_{q\neq j}(\alpha_j- \alpha_q)
	\cdot	\prod_{r\neq k}(\alpha_k- \alpha_r)}\nn\\
&& + \frac{1}{2}\frac{\alpha_j^{2n-3} (d\alpha_j)\prod_{p=0}^{d-1}\big(p\alpha_j+(d-p) \alpha_i\big)
	\prod_{p=0}^{d-1}\big(p\alpha_j+(d-p) \alpha_k\big)}{\prod_{p\neq i}(\alpha_i- \alpha_p)\cdot
	\prod_{q\neq j}(\alpha_j- \alpha_q)
	\cdot	\prod_{r\neq k}(\alpha_k- \alpha_r)}\nn\\
&&+ 	\frac{(\alpha_j- \alpha_i) \alpha_i^{n} \alpha_k^{n-3} (d\alpha_j)\prod_{p=0}^{d-1}\big(p\alpha_j+(d-p) \alpha_i\big)
	\prod_{p=0}^{d-1}\big(p\alpha_j+(d-p) \alpha_k\big)}
	{(2 \alpha_j- \alpha_i- \alpha_k)\prod_{p\neq i}(\alpha_i- \alpha_p)\cdot
	\prod_{q\neq j}(\alpha_j- \alpha_q)
	\cdot	\prod_{r\neq k}(\alpha_k- \alpha_r)}\Big)\nn\\
&&-\frac{1}{2}\sum_{i=0}^{n+1}\sum_{j\neq i}\frac{\alpha_i^{n} \alpha_j^{n-3} 
	\prod_{p=0}^{2d}\big(p\alpha_j+(2d-p) \alpha_i\big)}
	{32(\alpha_i- \alpha_j)^3 
	\prod_{p\neq i,j}( \alpha_i- \alpha_p)(\alpha_i+\alpha_j-2\alpha_p)(\alpha_j- \alpha_p)}\nn\\
&=& \sum_{j=0}^{n+1}\sum_{i\neq j}\sum_{k\neq j}
	\frac{(d\alpha_j)\prod_{p=0}^{d-1}\big(p\alpha_j+(d-p) \alpha_i\big)
	\prod_{p=0}^{d-1}\big(p\alpha_j+(d-p) \alpha_k\big)}
	{\prod_{p\neq i}(\alpha_i- \alpha_p)\cdot
	\prod_{q\neq j}(\alpha_j- \alpha_q)
	\cdot	\prod_{r\neq k}(\alpha_k- \alpha_r)}\nn\\
&&	\cdot \big(-\alpha_i^{n} \alpha_j^{n-3}+\frac{\alpha_j^{2n-3}}{2}+
	\frac{(\alpha_j- \alpha_i) \alpha_i^{n} \alpha_k^{n-3} }{2 \alpha_j- \alpha_i- \alpha_k}
	\big)	\nn\\
&&-	\frac{1}{2}\sum_{i=0}^{n+1}\sum_{j\neq i}\frac{\alpha_i^{n} \alpha_j^{n-3} 
	\prod_{p=0}^{2d}\big(p\alpha_j+(2d-p) \alpha_i\big)}
	{32(\alpha_i- \alpha_j)^3 
	\prod_{p\neq i,j}( \alpha_i- \alpha_p)(\alpha_i+\alpha_j-2\alpha_p)(\alpha_j- \alpha_p)}.
\end{eqnarray}
We rewrite
\begin{eqnarray}\label{eq-twistedInvariant-descendant-localization-sum-1-1}
&&\sum_{j=0}^{n+1}\sum_{i\neq j}\sum_{k\neq j}
	\frac{(d\alpha_j)\prod_{p=0}^{d-1}\big(p\alpha_j+(d-p) \alpha_i\big)
	\prod_{p=0}^{d-1}\big(p\alpha_j+(d-p) \alpha_k\big)}
	{\prod_{p\neq i}(\alpha_i- \alpha_p)\cdot
	\prod_{q\neq j}(\alpha_j- \alpha_q)
	\cdot	\prod_{r\neq k}(\alpha_k- \alpha_r)}\nn\\
&&	\cdot \big(-\alpha_i^{n} \alpha_j^{n-3}+\frac{\alpha_j^{2n-3}}{2}+
	\frac{(\alpha_j- \alpha_i) \alpha_i^{n} \alpha_k^{n-3} }{2 \alpha_j- \alpha_i- \alpha_k}
	\big)	\nn\\
&=& \sum_{j=0}^{n+1}\sum_{i\neq j}\sum_{k\neq j}
	\frac{(d\alpha_j)\prod_{p=0}^{d-1}\big(p\alpha_j+(d-p) \alpha_i\big)
	\prod_{p=0}^{d-1}\big(p\alpha_j+(d-p) \alpha_k\big)}
	{\prod_{p\neq i}(\alpha_i- \alpha_p)\cdot
	\prod_{q\neq j}(\alpha_j- \alpha_q)
	\cdot	\prod_{r\neq k}(\alpha_k- \alpha_r)}
	\cdot \big(-\alpha_i^{n} \alpha_j^{n-3}+\frac{\alpha_j^{2n-3}}{2}
	\big)	\nn\\
&&+\sum_{i=0}^{n+1}\sum_{k\neq i}\sum_{j\neq i,k}
	\frac{(d\alpha_j)\prod_{p=0}^{d-1}\big(p\alpha_j+(d-p) \alpha_i\big)
	\prod_{p=0}^{d-1}\big(p\alpha_j+(d-p) \alpha_k\big)}
	{\prod_{p\neq i}(\alpha_i- \alpha_p)\cdot
	\prod_{q\neq j}(\alpha_j- \alpha_q)
	\cdot	\prod_{r\neq k}(\alpha_k- \alpha_r)}
	\cdot 
	\frac{(\alpha_j- \alpha_i) \alpha_i^{n} \alpha_k^{n-3} }{2 \alpha_j- \alpha_i- \alpha_k}\nn\\	
&&+ \sum_{i=0}^{n+1}\sum_{j\neq i}
	\frac{(d\alpha_j)\prod_{p=0}^{d-1}\big(p\alpha_j+(d-p) \alpha_i\big)^2}
	{\prod_{p\neq i}(\alpha_i- \alpha_p)^2\cdot
	\prod_{q\neq j}(\alpha_j- \alpha_q)}
	\frac{\alpha_i^{2n-3} }{2}.
\end{eqnarray}

\subsubsection{Residue computations}
In this section we compute some summations using the residue theorem on $\mathbb{P}^1$.
\begin{lemma}\label{lem-residue-localization-5}
Suppose $i\neq k$. 
\begin{eqnarray}\label{eq-residue-localization-5}
&& \sum_{\begin{subarray}{c}0\leq j\leq n+1\\j\neq i,k\end{subarray}}
	\Big(\frac{(d\alpha_j)\prod_{p=0}^{d-1}\big(p\alpha_j+(d-p) \alpha_i\big)
	\prod_{p=0}^{d-1}\big(p\alpha_j+(d-p) \alpha_k\big)}
	{\prod_{q\neq j}(\alpha_j- \alpha_q)}
	\cdot \frac{\alpha_j -\alpha_i }
	{2\alpha_j- \alpha_i- \alpha_k}\Big)\nn\\
&=& -\frac{\prod_{p=0}^{2d}(p \alpha_k+(2d-p)\alpha_i)}{64(\alpha_i- \alpha_k)\prod_{p\neq i,k}(\alpha_i+ \alpha_k-2 \alpha_q)}
	-\frac{(d \alpha_k)^{d+1}\prod_{p=0}^{d-1}\big(p \alpha_k+(d-p) \alpha_i\big)}
	{\prod_{q\neq k}(\alpha_k- \alpha_q)}\nn\\
&&+ \mathrm{Res}_{x=0}\Big(
	\frac{d^3 \alpha_i \alpha_k\prod_{p=1}^{d-1}\big(p +(d-p) \alpha_i x\big)
	\prod_{p=1}^{d-1}\big(p +(d-p) \alpha_k x\big)}
	{x^{4}(2- \alpha_i x- \alpha_k x)(1- \alpha_k x)\prod_{q\neq i,k}(1- \alpha_qx)}
	\Big).		
\end{eqnarray}
\end{lemma}
\begin{proof}
Rewrite the left handside of (\ref{eq-residue-localization-5}) as
\begin{eqnarray*}
&& \sum_{\begin{subarray}{c}0\leq j\leq n+1\\j\neq i,k\end{subarray}}
	\Big(\frac{(d\alpha_j)\prod_{p=0}^{d-1}\big(p\alpha_j+(d-p) \alpha_i\big)
	\prod_{p=0}^{d-1}\big(p\alpha_j+(d-p) \alpha_k\big)}
	{\prod_{q\neq j}(\alpha_j- \alpha_q)}
	\cdot \frac{\alpha_j -\alpha_i }
	{2\alpha_j- \alpha_i- \alpha_k}\Big)\\
&=& \sum_{\begin{subarray}{c}0\leq j\leq n+1\\j\neq i,k\end{subarray}}
	\frac{(d\alpha_j)\prod_{p=0}^{d-1}\big(p\alpha_j+(d-p) \alpha_i\big)
	\prod_{p=0}^{d-1}\big(p\alpha_j+(d-p) \alpha_k\big)}
	{(2\alpha_j- \alpha_i- \alpha_k)(\alpha_j- \alpha_k)\prod_{q\neq i,j,k}(\alpha_j- \alpha_q)}.
\end{eqnarray*}
By the residue theorem,
\begin{eqnarray*}
&& 0=\sum_{x_0=\mathrm{Poles}}\mathrm{Res}_{x=x_0}\Big(
	\frac{(d x)\prod_{p=0}^{d-1}\big(p x+(d-p) \alpha_i\big)
	\prod_{p=0}^{d-1}\big(p x+(d-p) \alpha_k\big)}
	{(2 x- \alpha_i- \alpha_k)(x- \alpha_k)\prod_{q\neq i,k}(x- \alpha_q)}
	\Big)\\
&=& \mathrm{Res}_{x=\frac{\alpha_i+\alpha_k}{2}}\Big(
	\frac{(d x)\prod_{p=0}^{d-1}\big(p x+(d-p) \alpha_i\big)
	\prod_{p=0}^{d-1}\big(p x+(d-p) \alpha_k\big)}
	{(2 x- \alpha_i- \alpha_k)(x- \alpha_k)\prod_{q\neq i,k}(x- \alpha_q)}
	\Big)\\
&& +\mathrm{Res}_{x=\alpha_k}\Big(
	\frac{(d x)\prod_{p=0}^{d-1}\big(p x+(d-p) \alpha_i\big)
	\prod_{p=0}^{d-1}\big(p x+(d-p) \alpha_k\big)}
	{(2 x- \alpha_i- \alpha_k)(x- \alpha_k)\prod_{q\neq i,k}(x- \alpha_q)}
	\Big)\\	
&& +\sum_{\begin{subarray}{c}0\leq j\leq n+1\\j\neq i,k\end{subarray}}
	\mathrm{Res}_{x=\alpha_j}\Big(
	\frac{(d x)\prod_{p=0}^{d-1}\big(p x+(d-p) \alpha_i\big)
	\prod_{p=0}^{d-1}\big(p x+(d-p) \alpha_k\big)}
	{(2 x- \alpha_i- \alpha_k)(x- \alpha_k)\prod_{q\neq i,k}(x- \alpha_q)}
	\Big)\\
&&+ \mathrm{Res}_{x=\infty}\Big(
	\frac{(d x)\prod_{p=0}^{d-1}\big(p x+(d-p) \alpha_i\big)
	\prod_{p=0}^{d-1}\big(p x+(d-p) \alpha_k\big)}
	{(2 x- \alpha_i- \alpha_k)(x- \alpha_k)\prod_{q\neq i,k}(x- \alpha_q)}
	\Big).
\end{eqnarray*}	
Using 
\begin{eqnarray*}
&&\frac{d(\alpha_i+\alpha_k)}{2}\prod_{p=0}^{d-1}\big(\frac{p(\alpha_i+\alpha_k)}{2}+(d-p)\alpha_i\big)\prod_{p=0}^{d-1}\big(\frac{p(\alpha_i+\alpha_k)}{2}+(d-p)\alpha_k\big)\\
&=&2^{-2d-1}\prod_{p=0}^{2d}(p \alpha_k+(2d-p)\alpha_i),
\end{eqnarray*}
we get
\begin{eqnarray*}
&& \sum_{\begin{subarray}{c}0\leq j\leq n+1\\j\neq i,k\end{subarray}}
	\frac{(d\alpha_j)\prod_{p=0}^{d-1}\big(p\alpha_j+(d-p) \alpha_i\big)
	\prod_{p=0}^{d-1}\big(p\alpha_j+(d-p) \alpha_k\big)}
	{(2\alpha_j- \alpha_i- \alpha_k)(\alpha_j- \alpha_k)\prod_{q\neq i,j,k}(\alpha_j- \alpha_q)}\\	
&=& -\frac{\prod_{p=0}^{2d}(p \alpha_k+(2d-p)\alpha_i)}{64(\alpha_i- \alpha_k)\prod_{p\neq i,k}(\alpha_i+ \alpha_k-2 \alpha_q)}
	- \frac{(d \alpha_k)^{d+1}\prod_{p=0}^{d-1}\big(p \alpha_k+(d-p) \alpha_i\big)}
	{(\alpha_k- \alpha_i)\prod_{q\neq i,k}(\alpha_k- \alpha_q)}\\
&&- \mathrm{Res}_{x=0}\Big(
	-\frac{(d x^{-1})\prod_{p=0}^{d-1}\big(p x^{-1}+(d-p) \alpha_i\big)
	\prod_{p=0}^{d-1}\big(p x^{-1}+(d-p) \alpha_k\big)}
	{x^2(2 x^{-1}- \alpha_i- \alpha_k)(x^{-1}- \alpha_k)\prod_{q\neq i,k}(x^{-1}- \alpha_q)}
	\Big).
\end{eqnarray*}
Hence (\ref{eq-residue-localization-5}) follows.
\end{proof}
In the following we omit such details.
\begin{lemma}\label{{lem-residue-localization-2}}
\begin{equation}\label{eq-residue-localization-2}
	\sum_{\begin{subarray}{c}0\leq k\leq n+1\\k\neq i\end{subarray}}
	\frac{\prod_{p=0}^{d-1}\big(p\alpha_i+(d-p) \alpha_k\big)}
	{(\alpha_k- \alpha_i)\prod_{q\neq i,k}(\alpha_k- \alpha_q)}	
	=-\frac{(d \alpha_i)^{d}}
	{\prod_{q\neq i}(\alpha_i- \alpha_q)}.
\end{equation}
\end{lemma}
\begin{proof}
Apply the residue theorem to
\[
\frac{\prod_{p=0}^{d-1}\big(p x+(d-p) \alpha_i\big)}
	{(x- \alpha_i)\prod_{q\neq i}(x- \alpha_q)}.
\]
\end{proof}

\begin{lemma}\label{lem-residue-localization-6}
\begin{eqnarray}\label{eq-residue-localization-6}
&&	\sum_{\begin{subarray}{c}0\leq j\leq n+1\\ j\neq i\end{subarray}}\frac{(d\alpha_j)\prod_{p=0}^{d-1}\big(p\alpha_j+(d-p) \alpha_i\big)^2}
	{\prod_{q\neq j}(\alpha_j- \alpha_q)}\nn\\
&=& -\frac{(d \alpha_i)^{2d+1}}
	{\prod_{q\neq i}(\alpha_i- \alpha_q)}
	+\mathrm{Res}_{x=0}
\frac{d^3 \alpha_i^2\prod_{p=1}^{d-1}\big(p+(d-p) \alpha_i x\big)^2}
	{x^{4}\prod_{q=0}^{n+1}(1- \alpha_q x)}.
\end{eqnarray}
\end{lemma}
\begin{proof}
Apply the residue theorem to
\[
\frac{(d x)\prod_{p=0}^{d-1}\big(p x+(d-p) \alpha_i\big)^2}
	{(x- \alpha_i)\prod_{q\neq i}(x- \alpha_q)}.
\]
\end{proof}

\subsubsection{The summation}
Applying (\ref{eq-residue-localization-2}) to the second sum of (\ref{eq-twistedInvariant-descendant-localization-sum-1-1}), 
 and (\ref{eq-residue-localization-5}) to the third, we get
\begin{eqnarray}\label{eq-twistedInvariant-descendant-localization-sum-2}
&& \langle \sfh_{n}\psi,\sfh_{n-3}\rangle_{0,2}\nn\\
&=& \sum_{j=0}^{n+1}\sum_{i\neq j}
	\frac{(d\alpha_j)\prod_{p=0}^{d-1}\big(p\alpha_j+(d-p) \alpha_i\big)}
	{\prod_{p\neq i}(\alpha_i- \alpha_p)\cdot
	\prod_{q\neq j}(\alpha_j- \alpha_q)}
	\cdot\big(-\frac{(d \alpha_j)^{d}}
	{\prod_{q\neq j}(\alpha_j- \alpha_q)}\big)
	\cdot \big(-\alpha_i^{n} \alpha_j^{n-3}+\frac{\alpha_j^{2n-3}}{2}
	\big)	\nn\\
&&+\sum_{i=0}^{n+1}\sum_{k\neq i}
	\frac{\alpha_i^{n} \alpha_k^{n-3}}
	{\prod_{p\neq i}(\alpha_i- \alpha_p)\cdot
	\prod_{q\neq k}(\alpha_k- \alpha_q)}
	\cdot\bigg( \textcolor{blue}{-\frac{\prod_{p=0}^{2d}(p \alpha_k+(2d-p)\alpha_i)}{64(\alpha_i- \alpha_k)\prod_{p\neq i,k}(\alpha_i+ \alpha_k-2 \alpha_q)}}\nn\\
&&	-\frac{(d \alpha_k)^{d+1}\prod_{p=0}^{d-1}\big(p \alpha_k+(d-p) \alpha_i\big)}
	{\prod_{q\neq k}(\alpha_k- \alpha_q)}
	+ \mathrm{Res}_{x=0}\Big(
	\frac{d^3 \alpha_i \alpha_k\prod_{p=1}^{d-1}\big(p +(d-p) \alpha_i x\big)
	\prod_{p=1}^{d-1}\big(p +(d-p) \alpha_k x\big)}
	{x^{4}(2- \alpha_i x- \alpha_k x)(1- \alpha_k x)\prod_{q\neq i,k}(1- \alpha_qx)}
	\Big)\bigg) \nn\\	
&&+ \sum_{i=0}^{n+1}\sum_{j\neq i}
	\frac{(d\alpha_j)\prod_{p=0}^{d-1}\big(p\alpha_j+(d-p) \alpha_i\big)^2}
	{\prod_{p\neq i}(\alpha_i- \alpha_p)^2\cdot
	\prod_{q\neq j}(\alpha_j- \alpha_q)}
	\frac{\alpha_i^{2n-3} }{2}\nn\\
&&\textcolor{blue}{-	\frac{1}{2}\sum_{i=0}^{n+1}\sum_{j\neq i}\frac{\alpha_i^{n} \alpha_j^{n-3} 
	\prod_{p=0}^{2d}\big(p\alpha_j+(2d-p) \alpha_i\big)}
	{32(\alpha_i- \alpha_j)^3 
	\prod_{p\neq i,j}( \alpha_i- \alpha_p)(\alpha_i+\alpha_j-2\alpha_p)(\alpha_j- \alpha_p)}}\nn\\
&=& \sum_{j=0}^{n+1}\sum_{i\neq j}
	\frac{(d\alpha_j)^{d+1}\prod_{p=0}^{d-1}\big(p\alpha_j+(d-p) \alpha_i\big)}
	{\prod_{p\neq i}(\alpha_i- \alpha_p)\cdot
	\prod_{q\neq j}(\alpha_j- \alpha_q)^2}
	\cdot \big(\textcolor{blue}{\alpha_i^{n} \alpha_j^{n-3}}-\frac{\alpha_j^{2n-3}}{2}
	\big)	\nn\\
&&+\sum_{i=0}^{n+1}\sum_{k\neq i}
	\frac{\alpha_i^{n} \alpha_k^{n-3}}
	{\prod_{p\neq i}(\alpha_i- \alpha_p)\cdot
	\prod_{q\neq k}(\alpha_k- \alpha_q)}
	\cdot\bigg(\textcolor{blue}{-\frac{(d \alpha_k)^{d+1}\prod_{p=0}^{d-1}\big(p \alpha_k+(d-p) \alpha_i\big)}
	{\prod_{q\neq k}(\alpha_k- \alpha_q)}}\nn\\
&&	+ \mathrm{Res}_{x=0}\Big(
	\frac{d^3 \alpha_i \alpha_k\prod_{p=1}^{d-1}\big(p +(d-p) \alpha_i x\big)
	\prod_{p=1}^{d-1}\big(p +(d-p) \alpha_k x\big)}
	{x^{4}(2- \alpha_i x- \alpha_k x)(1- \alpha_k x)\prod_{q\neq i,k}(1- \alpha_qx)}
	\Big)\bigg) \nn\\	
&&+ \sum_{i=0}^{n+1}\sum_{j\neq i}
	\frac{(d\alpha_j)\prod_{p=0}^{d-1}\big(p\alpha_j+(d-p) \alpha_i\big)^2}
	{\prod_{p\neq i}(\alpha_i- \alpha_p)^2\cdot
	\prod_{q\neq j}(\alpha_j- \alpha_q)}
	\frac{\alpha_i^{2n-3} }{2}\nn\\
&=& -\frac{1}{2}\sum_{j=0}^{n+1}\sum_{i\neq j}
	\frac{(d\alpha_j)^{d+1}\alpha_j^{2n-3}\prod_{p=0}^{d-1}\big(p\alpha_j+(d-p) \alpha_i\big)}
	{\prod_{p\neq i}(\alpha_i- \alpha_p)\cdot
	\prod_{q\neq j}(\alpha_j- \alpha_q)^2}	\nn\\
&&+\sum_{i=0}^{n+1}\sum_{k\neq i}
	\bigg(\frac{d^3\alpha_i^{n+1} \alpha_k^{n-2}}
	{\prod_{p\neq i}(\alpha_i- \alpha_p)\cdot
	\prod_{q\neq k}(\alpha_k- \alpha_q)}
	\mathrm{Res}_{x=0}\Big(
	\frac{\prod_{p=1}^{d-1}\big(p +(d-p) \alpha_i x\big)
	\prod_{p=1}^{d-1}\big(p +(d-p) \alpha_k x\big)}
	{x^{4}(2- \alpha_i x- \alpha_k x)(1- \alpha_k x)\prod_{q\neq i,k}(1- \alpha_qx)}
	\Big)\bigg) \nn\\	
&&+\frac{1}{2} \sum_{i=0}^{n+1}
	\Big(\frac{\alpha_i^{2n-3}}
	{\prod_{p\neq i}(\alpha_i- \alpha_p)^2}
	\sum_{j\neq i}\frac{(d\alpha_j)\prod_{p=0}^{d-1}\big(p\alpha_j+(d-p) \alpha_i\big)^2}
	{\prod_{q\neq j}(\alpha_j- \alpha_q)}\Big),
\end{eqnarray}
where the second and the third equality are obtained by cancellation of the terms in \textcolor{blue}{blue}. 
Applying (\ref{eq-residue-localization-2}) to the first sum of (\ref{eq-twistedInvariant-descendant-localization-sum-2})  we have
\begin{eqnarray}\label{eq-twistedInvariant-descendant-localization-subsum-1}
&&-\frac{1}{2}\sum_{j=0}^{n+1}\sum_{i\neq j}
	\frac{(d\alpha_j)^{d+1}\alpha_j^{2n-3}\prod_{p=0}^{d-1}\big(p\alpha_j+(d-p) \alpha_i\big)}
	{\prod_{p\neq i}(\alpha_i- \alpha_p)\cdot
	\prod_{q\neq j}(\alpha_j- \alpha_q)^2}\nn\\
&=& -\frac{1}{2}\sum_{j=0}^{n+1}
	\frac{(d\alpha_j)^{d+1}\alpha_j^{2n-3}}
	{\prod_{q\neq j}(\alpha_j- \alpha_q)^2}	
	\sum_{i\neq j}\frac{\prod_{p=0}^{d-1}\big(p\alpha_j+(d-p) \alpha_i\big)}{\prod_{p\neq i}(\alpha_i- \alpha_p)}\nn\\
&=& -\frac{1}{2}\sum_{j=0}^{n+1}
	\frac{(d\alpha_j)^{d+1}\alpha_j^{2n-3}}
	{\prod_{q\neq j}(\alpha_j- \alpha_q)^2}	\big(-\frac{(d \alpha_j)^{d}}
	{\prod_{q\neq j}(\alpha_j- \alpha_q)}\big)\nn\\
&=& \frac{1}{2}\sum_{j=0}^{n+1}
	\frac{(d\alpha_j)^{2d+1}\alpha_j^{2n-3}}
	{\prod_{q\neq j}(\alpha_j- \alpha_q)^3}.
\end{eqnarray}
To deal with the second sum of (\ref{eq-twistedInvariant-descendant-localization-sum-2})
we rewrite
\begin{eqnarray*}
&&	\mathrm{Res}_{x=0}\Big(
	\frac{\prod_{p=1}^{d-1}\big(p +(d-p) \alpha_i x\big)
	\prod_{p=1}^{d-1}\big(p +(d-p) \alpha_k x\big)}
	{x^{4}(2- \alpha_i x- \alpha_k x)(1- \alpha_k x)\prod_{q\neq i,k}(1- \alpha_qx)}
	\Big)\bigg)\\
&=& \mathrm{Res}_{x=0}\Big(
	\frac{(1- \alpha_i x)\prod_{p=1}^{d-1}\big(p +(d-p) \alpha_i x\big)
	\prod_{p=1}^{d-1}\big(p +(d-p) \alpha_k x\big)}
	{x^{4}(2- \alpha_i x- \alpha_k x)\prod_{q=0}^{n+1}(1- \alpha_qx)}
	\Big)\bigg),
\end{eqnarray*}
and thus get
\begin{eqnarray}\label{eq-twistedInvariant-descendant-localization-subsum-2}
&&\sum_{i=0}^{n+1}\sum_{k\neq i}
	\bigg(\frac{d^3\alpha_i^{n+1} \alpha_k^{n-2}}
	{\prod_{p\neq i}(\alpha_i- \alpha_p)\cdot
	\prod_{q\neq k}(\alpha_k- \alpha_q)}
	\mathrm{Res}_{x=0}\Big(
	\frac{\prod_{p=1}^{d-1}\big(p +(d-p) \alpha_i x\big)
	\prod_{p=1}^{d-1}\big(p +(d-p) \alpha_k x\big)}
	{x^{4}(2- \alpha_i x- \alpha_k x)(1- \alpha_k x)\prod_{q\neq i,k}(1- \alpha_qx)}
	\Big)\bigg)\nn\\
&=&	\sum_{i=0}^{n+1}\sum_{k=0}^{n+1}
	\bigg(\frac{d^3\alpha_i^{n+1} \alpha_k^{n-2}}
	{\prod_{p\neq i}(\alpha_i- \alpha_p)\cdot
	\prod_{q\neq k}(\alpha_k- \alpha_q)}\nn\\
&&	\cdot\mathrm{Res}_{x=0}\Big(
	\frac{(1- \alpha_i x)\prod_{p=1}^{d-1}\big(p +(d-p) \alpha_i x\big)
	\prod_{p=1}^{d-1}\big(p +(d-p) \alpha_k x\big)}
	{x^{4}(2- \alpha_i x- \alpha_k x)\prod_{q=0}^{n+1}(1- \alpha_qx)}
	\Big)\bigg)\nn\\
&& -\sum_{i=0}^{n+1}
	\bigg(\frac{d^3\alpha_i^{2n-1}}
	{\prod_{p\neq i}(\alpha_i- \alpha_p)^2}
	\mathrm{Res}_{x=0}\Big(
	\frac{\prod_{p=1}^{d-1}\big(p +(d-p) \alpha_i x\big)^2}
	{2x^{4}\prod_{q=0}^{n+1}(1- \alpha_qx)}
	\Big)\bigg)	.
\end{eqnarray}
Applying (\ref{eq-residue-localization-6}) to the third sum of (\ref{eq-twistedInvariant-descendant-localization-sum-2}) we get
\begin{eqnarray}\label{eq-twistedInvariant-descendant-localization-subsum-3}
&& \frac{1}{2} \sum_{i=0}^{n+1}
	\Big(\frac{\alpha_i^{2n-3}}
	{\prod_{p\neq i}(\alpha_i- \alpha_p)^2}
	\sum_{j\neq i}\frac{(d\alpha_j)\prod_{p=0}^{d-1}\big(p\alpha_j+(d-p) \alpha_i\big)^2}
	{\prod_{q\neq j}(\alpha_j- \alpha_q)}\Big)\nn\\
&=& \frac{1}{2} \sum_{i=0}^{n+1}
	\bigg(\frac{\alpha_i^{2n-3}}
	{\prod_{p\neq i}(\alpha_i- \alpha_p)^2}\Big(-\frac{(d \alpha_i)^{2d+1}}
	{\prod_{q\neq i}(\alpha_i- \alpha_q)}
	+\mathrm{Res}_{x=0}
\frac{d^3 \alpha_i^2\prod_{p=1}^{d-1}\big(p+(d-p) \alpha_i x\big)^2}
	{x^{4}\prod_{q=0}^{n+1}(1- \alpha_q x)}\Big)\bigg)\nn\\
&=&  -\frac{1}{2} \sum_{i=0}^{n+1}
	\frac{(d \alpha_i)^{2d+1}\alpha_i^{2n-3}}
	{\prod_{p\neq i}(\alpha_i- \alpha_p)^3}	\nn\\
&&	+\frac{1}{2} \sum_{i=0}^{n+1}\Big(\frac{d^3\alpha_i^{2n-1}}
	{\prod_{p\neq i}(\alpha_i- \alpha_p)^2}
	\mathrm{Res}_{x=0}
\frac{\prod_{p=1}^{d-1}\big(p+(d-p) \alpha_i x\big)^2}
	{x^{4}\prod_{q=0}^{n+1}(1- \alpha_q x)}\Big).
\end{eqnarray}	
It follows from (\ref{eq-twistedInvariant-descendant-localization-sum-2}), (\ref{eq-twistedInvariant-descendant-localization-subsum-1}), (\ref{eq-twistedInvariant-descendant-localization-subsum-2}) and (\ref{eq-twistedInvariant-descendant-localization-subsum-3}) that
\begin{eqnarray}\label{eq-twistedInvariant-descendant-localization-sum-3}
&& \langle \sfh_{n}\psi,\sfh_{n-3}\rangle_{0,2}\nn\\
&=& \sum_{i=0}^{n+1}\sum_{k=0}^{n+1}
	\bigg(\frac{d^3\alpha_i^{n+1} \alpha_k^{n-2}}
	{\prod_{p\neq i}(\alpha_i- \alpha_p)\cdot
	\prod_{q\neq k}(\alpha_k- \alpha_q)} \nn\\
&&	\cdot\mathrm{Res}_{x=0}\Big(
	\frac{(1- \alpha_i x)\prod_{p=1}^{d-1}\big(p +(d-p) \alpha_i x\big)
	\prod_{p=1}^{d-1}\big(p +(d-p) \alpha_k x\big)}
	{x^{4}(2- \alpha_i x- \alpha_k x)\prod_{q=0}^{n+1}(1- \alpha_qx)}
	\Big)\bigg).
\end{eqnarray}
We expand 
\[
\mathrm{Res}_{x=0}\Big(
	\frac{(1- \alpha_i x)\prod_{p=1}^{d-1}\big(p +(d-p) \alpha_i x\big)
	\prod_{p=1}^{d-1}\big(p +(d-p) \alpha_k x\big)}
	{x^{4}(2- \alpha_i x- \alpha_k x)\prod_{q=0}^{n+1}(1- \alpha_qx)}
	\Big)
\]
as a polynomial of $\alpha_i$ and $\alpha_k$ of degree $\leq 3$, whose coefficients are constants or symmetric polynomials of $\alpha_0,\dots,\alpha_{n+1}$ that are independent of $i$ and $k$. Since
\[
\sum_{k=0}^{n+1}\frac{\alpha_k^{n-2}\cdot \alpha_k^{b}}{\prod_{q\neq k}(\alpha_k- \alpha_q)}=\delta_{b,3}
\]
for $b\leq 3$,  we obtain
\begin{eqnarray*}
&&\sum_{k=0}^{n+1}
	\bigg(\frac{\alpha_k^{n-2}}
	{\prod_{q\neq k}(\alpha_k- \alpha_q)} \nn\\
&&	\cdot\mathrm{Res}_{x=0}\Big(
	\frac{(1- \alpha_i x)\prod_{p=1}^{d-1}\big(p +(d-p) \alpha_i x\big)
	\prod_{p=1}^{d-1}\big(p +(d-p) \alpha_k x\big)}
	{x^{4}(2- \alpha_i x- \alpha_k x)\prod_{q=0}^{n+1}(1- \alpha_qx)}
	\Big)	\\
&=& \mathrm{Coeff}_{\alpha_k^3}\mathrm{Coeff}_{x^3}	\Big(
	\frac{(1- \alpha_i x)\prod_{p=1}^{d-1}\big(p +(d-p) \alpha_i x\big)
	\prod_{p=1}^{d-1}\big(p +(d-p) \alpha_k x\big)}
	{(2- \alpha_i x- \alpha_k x)\prod_{q=0}^{n+1}(1- \alpha_qx)}
	\Big)	\\	
&=& (d-1)! \mathrm{Coeff}_{\alpha_k^3}\mathrm{Coeff}_{x^3}	\Big(
	\frac{\prod_{p=1}^{d-1}\big(p +(d-p) \alpha_k x\big)}
	{2- \alpha_k x}	\Big)\\
&=& (d-1)!	\mathrm{Coeff}_{x^3}	\Big(
	\frac{\prod_{p=1}^{d-1}\big(p +(d-p) x\big)}{2-  x}	\Big)\\
&=& \frac{\big((d-1)!\big)^2}{2}  \mathrm{Coeff}_{x^3}\Big(\big(1+\sum_{p=1}^{d-1}\frac{(d-p)x}{p}+ \sum_{1\leq p<q\leq d-1}\frac{(d-p)(d-q)x^2}{pq}\\
&&	+ \sum_{1\leq p<q<r\leq d-1}\frac{(d-p)(d-q)(d-r)x^3}{pqr}
	+O(x^4)\big)\big(1+\frac{x}{2}+\frac{x^2}{4}+\frac{x^3}{8}+O(x^4)\big)\Big)\\
&=& \frac{\big((d-1)!\big)^2}{2} \Big(\frac{1}{8}+\frac{1}{4}\sum_{p=1}^{d-1}\frac{d-p}{p}+ \frac{1}{2}\sum_{1\leq p<q\leq d-1}\frac{(d-p)(d-q)}{pq}\\
&&	+ \sum_{1\leq p<q<r\leq d-1}\frac{(d-p)(d-q)(d-r)}{pqr}\Big).
\end{eqnarray*}
The proof of (\ref{eq-twistedInvariant-descendant-localization-result}), and therefore that of Theorem \ref{thm-F20-WM-a(n,d)=(n-1)/2}, is completed.

\end{appendix}

\subsection*{Acknowledgment}
I am grateful to Huai-Liang Chang and Hua-Zhong Ke for enlightening discussions and encouragement to complete this work.  I thank Sergey Galkin and Nicolas Perrin for very helpful discussions on Fano variety of lines and semisimplicity of quantum cohomology. I thank R. Pandharipande for helpful discussions and suggesting the terminology \emph{double root recursion}.   I also thank James Carlson, Huijun Fan, Si-Qi Liu, Christopher Lyons, Giosuè Muratore, Maxim Smirnov, Yang Su,  Gang Tian, Xin Wang, Jinxing Xu, Ze Xu, Lei Zhang and Jian Zhou for  discussions on various related topics. 
 I am grateful to the authors of Macaulay2 \cite{GS-Mac} and its packages for providing  a convenient programming language.
  This work is supported by China postdoctoral science special foundation 2015T80007, NSFC 12371063, and NSFC 11701579.

\textsc{School of Sciences, Great Bay University, Dongguan 523000, P.R. China}

 \emph{E-mail address:}  huxw06@gmail.com

\end{document}